\documentclass[11pt]{amsart}

\usepackage{amsmath,amsthm,amssymb,amsfonts, amscd, bm, graphicx}

\usepackage[margin=1in,centering]{geometry}

\usepackage[cmtip,all]{xy}

\allowdisplaybreaks

\usepackage{mathrsfs}
\usepackage{bbm}
\usepackage{extarrows} 
\usepackage{cancel}

\usepackage[pagebackref=true, colorlinks]{hyperref}

\hypersetup{pdffitwindow=true,linkcolor=blue,citecolor=blue,urlcolor=blue,menucolor=blue}

\usepackage{comment}
\usepackage{enumerate}

\usepackage{booktabs,longtable,array}
\usepackage[nopatch=eqnum]{microtype}

\newtheorem{thm}{Theorem}[section]
\newtheorem{lem}[thm]{Lemma}
\newtheorem{prop}[thm]{Proposition}
\newtheorem{cor}[thm]{Corollary}
\theoremstyle{definition}

\theoremstyle{remark}
\newtheorem{rem}[thm]{Remark}

\theoremstyle{plain}

\numberwithin{equation}{section}
\numberwithin{figure}{section}

\def\beq{\begin{equation}}
\def\eeq{\end{equation}}

\newcommand*{\mmiddle}[1]{\mathrel{}\middle#1\mathrel{}}

\newcommand{\N}{\mathbf N}
\newcommand{\Z}{\mathbf Z}
\newcommand{\R}{\mathbf R}
\newcommand{\C}{\mathbf C}
\newcommand{\Q}{\mathbf Q}
\newcommand{\E}{\mathbf{E}}
\renewcommand{\P}{\mathbf{P}}

\renewcommand{\1}{\mathbf{1}}

\newcommand{\cI}{\mathcal{I}}
\newcommand{\cX}{\mathcal{X}}
\newcommand{\cK}{\mathcal{K}}
\newcommand{\cQ}{\mathcal{Q}}
\newcommand{\cH}{\mathcal{H}}
\newcommand{\bT}{\mathbf{T}}
\newcommand{\cL}{\mathcal{L}}
\newcommand{\cR}{\mathcal{R}}
\newcommand{\cE}{\mathcal{E}}
\newcommand{\cB}{\mathcal{B}}
\newcommand{\cD}{\mathcal{D}}
\newcommand{\cA}{\mathcal{A}}
\newcommand{\cG}{\mathcal{G}}
\newcommand{\cS}{\mathcal{S}}
\newcommand{\cZ}{\mathcal{Z}}
\newcommand{\cC}{\mathcal{C}}
\newcommand{\cW}{\mathcal{W}}
\newcommand{\fC}{\mathfrak{C}}

\newcommand{\fd}{\mathfrak{d}}
\newcommand{\AC}{A^{\CLE}}

\newcommand*{\piic}{\p_{\mathrm{IIC}}}

\newcommand{\cXmu}{\cX^{(m)}_{Q}}
\newcommand{\cXmun}{\cX^{(m)}_{Q,n}}

\DeclareMathOperator{\SLE}{SLE}
\DeclareMathOperator{\CLE}{CLE}

\DeclareMathOperator{\diam}{diam}

\DeclareMathOperator{\GH}{GH}
\DeclareMathOperator{\GHf}{GHf}
\DeclareMathOperator{\GHP}{GHP}
\DeclareMathOperator{\GHPf}{GHPf}
\DeclareMathOperator{\rGHPf}{rGHPf}
\DeclareMathOperator{\Hc}{Hc}
\DeclareMathOperator{\cur}{c}
\DeclareMathOperator{\Pro}{P}
\DeclareMathOperator{\Haus}{H}
\DeclareMathOperator{\out}{out}
\DeclareMathOperator{\inn}{in}

\DeclareMathOperator{\dist}{dist}
\DeclareMathOperator{\intt}{int}
\DeclareMathOperator{\med}{median}

\newcommand{\pth}{{\mathrm {path}}}
\newcommand{\geo}{{\mathrm {geo}}}
\newcommand{\res}{{\mathrm {res}}}
\newcommand{\dpt}{d^{\pth}}
\newcommand{\dg}{d^{\geo}}
\newcommand{\dr}{d^{\res}}

\newcommand{\Dg}{D^{\geo}}
\newcommand{\Dr}{D^{\res}}

\newcommand{\p}{\mathbf{P}}
\newcommand{\median}{{\mathfrak{q}}}
\newcommand{\altq}{\mathbf{q}}
\newcommand{\fp}{\mathfrak{p}}
\newcommand{\dprime}{l}

\begin{document}

\title[Random walk and the intrinsic metric on planar critical percolation]{The scaling limit of random walk and\\ the intrinsic metric on planar critical percolation}

\author{Irina \DJ ankovi\'c, Maarten Markering, Jason Miller, and Yizheng Yuan}
\thanks{Department of Pure Mathematics and Mathematical Statistics, University of Cambridge}
\thanks{Faculty of Mathematics, University of Vienna}

\begin{abstract}
We consider critical site percolation ($p=p_c=1/2$) on the triangular lattice $\bT$ in two dimensions.  We show that the simple random walk on the clusters of open vertices converges in the scaling limit to a continuous diffusion which lives in the gasket of a conformal loop ensemble with parameter $\kappa = 6$ ($\CLE_6$), the so-called \emph{$\CLE_6$ Brownian motion}.  We also show that the intrinsic (i.e., chemical distance) metric converges in the scaling limit to the \emph{geodesic $\CLE_6$ metric}. As a consequence, we deduce the existence of the chemical distance exponent, the resistance exponent, and the spectral dimension of the critical percolation clusters. Moreover, we show that the exponents satisfy the Einstein relations.
\end{abstract}

\date{\today}
\maketitle

\setcounter{tocdepth}{1}

\tableofcontents

\parindent 0 pt
\setlength{\parskip}{0.20cm plus1mm minus1mm}

\addtocontents{toc}{
  \protect\contentsline{part}{\hyperref[sec:introduction]{Part 1}}{}{}
}

\section{Introduction}
\label{sec:introduction}

\subsection{Overview}
\label{sec:overview}

The main focus of this paper is the site percolation model on the triangular lattice~$\bT$ in two dimensions.  Recall that it is defined by fixing a parameter $p \in [0,1]$ and then declaring each vertex to be either open or closed independently with probability $p$.  Given a percolation configuration, a \emph{cluster} of open (resp.\ closed) vertices is a connected component in the subgraph of~$\bT$ where two vertices are connected by an edge if they are adjacent in $\bT$ and they are both open (resp.\ closed).  Recall that the \emph{critical probability} $p_c$ is given by $\inf\{ p \in [0,1] : \theta(p) > 0\}$, where $\theta(p) = \p[ 0 \longleftrightarrow \infty]$ and the notation $0 \longleftrightarrow \infty$ means that the cluster of open vertices containing the origin is infinite, and for $\bT$ we have that $p_c = 1/2$.  In the case that $p > p_c$, there is almost surely a unique infinite cluster and it has been shown that its behavior is strongly Euclidean at large scales in several different senses while for $p=p_c$ the large scale geometry is highly non-trivial and fractal.  In particular, it was shown by Smirnov \cite{smirnov2001critical} and Camia--Newman \cite{camia2006two} that the cluster interfaces converge in the scaling limit to the conformal loop ensemble with $\kappa=6$ ($\CLE_6$) \cite{camia2006two,sheffield2009exploration}, the loop variant of Schramm's $\SLE_6$ curves \cite{schramm2000scaling}.

In this paper, we will be concerned with establishing scaling limit results for:
\begin{enumerate}[(i)]
\item\label{it:srw_problem} the \emph{simple random walk} and
\item\label{it:metric_problem} the \emph{intrinsic} (i.e., chemical distance) metric
\end{enumerate}
on the clusters of open sites for critical percolation on $\bT$.

In the case of~\eqref{it:srw_problem}, we mean that we consider a cluster of open vertices and then study the random walk on that cluster which moves from a site to each of its neighbors at rate $1$ (so the associated discrete time random walk moves from a site to each of its neighbors with equal probability).  The study of this model was popularized by de Gennes in 1976 \cite{degenes1976percolation} in which he referred to it as the ``\emph{Ant in the Labyrinth}'', the ant corresponding to the random walker and the labyrinth corresponding to the random environment described by the percolation configuration that the random walker lives in.  As we will describe just below, this has been the inspiration for a large number of subsequent works.  Our contribution will be to show that the simple random walk on large clusters converges in the scaling limit to the \emph{$\CLE_6$ Brownian motion}, the canonical diffusion on a $\CLE_6$ constructed in \cite{amy2025tightness,my2025diffusion}.  In order to prove this convergence result, we will use the framework of resistance metrics developed by Kigami  \cite{kigami2001analysis,kigami2012resistance} as well as the work of Croydon \cite{croydon2008convergence} which reduces the problem of process convergence to the problem of resistance metric convergence.

In the case of~\eqref{it:metric_problem}, we will consider a cluster of open vertices and then study the graph metric where there is a unit length edge connecting any pair of adjacent open sites.  We will show that the graph metric converges in the scaling limit to the \emph{geodesic $\CLE_6$ metric}, the canonical conformally covariant geodesic metric that one can associate with a $\CLE_6$ constructed in \cite{amy2025tightness,my2025metric}.

The CLE$_6$ in a simply connected domain $D \subsetneq \C$ is a random countable collection of loops contained in $\overline{D}$, introduced in \cite{camia2006two,sheffield2009exploration}. The CLE$_6$ loops intersect themselves and can intersect each other and the domain boundary, but they do not cross upon intersecting. Each of the loop segments looks locally like an SLE$_6$ curve \cite{schramm2000scaling}. The CLE$_6$ loops describe the scaling limit of the collection of cluster boundaries in critical percolation. The intersection points of the loops correspond to the so-called \emph{$4$-arm points} or \emph{pivotal points} in the percolation. We note that in this paper we will consider the \emph{nested CLE$_6$} in which each loop contains a further collection of loops on its inside. The \emph{non-nested CLE$_6$}, in contrast, consists only of the outermost loops. The nested CLE$_6$ is obtained from the non-nested CLE$_6$ by repeatedly adding conditionally independent CLE$_6$'s in each of the connected components inside the loops.

The \emph{CLE$_6$ gaskets} describe the scaling limit of the clusters themselves. The gasket of a non-nested CLE$_6$ is the set of points that are connected to $\partial D$ by a path not crossing any of the loops (i.e.,the set of points that do not lie inside the loops). In a nested CLE$_6$, we can define for each loop $\gamma$ the gasket inside $\gamma$ as the set of points on or inside $\gamma$ that are connected to $\gamma$ by a path not crossing any of the loops inside $\gamma$ (i.e., the gasket of the non-nested CLE$_6$ inside $\gamma$). The CLE$_6$ gaskets are fractal sets and their Hausdorff dimension (with respect to the Euclidean metric) is shown in \cite{ssw2009radii,miller2016dimension} to be $91/48$.

Note that since the CLE$_6$ loops are self-intersecting, there will be distinct points on the percolation clusters (the points on the opposite arms of a pivotal) that converge to the same point on the Euclidean plane in the scaling limit. These points need to be distinguished in the CLE$_6$ metric space. Therefore we will view the CLE$_6$ gasket as an abstract metric space with an embedding in the plane where each point in the plane has one or two preimages in the gasket. See Section~\ref{subsec:statements} below.

The main goal of this paper is to show the convergence of the simple random walk and the intrinsic metric, respectively, to the $\CLE_6$ Brownian motion and the geodesic $\CLE_6$ metric on all macroscopic clusters.

\subsection{Related work}
\label{subsec:prior_work}

Both the simple random walk and the intrinsic metric on percolation clusters have been the subject of a substantial amount of prior work which we now review.  In each case, one can think of putting the existing results into one of three categories: supercritical percolation, critical percolation in two dimensions, and critical percolation in high dimensions.  We will also mention results for related models.

\subsubsection{Ant in the labyrinth}
We direct the reader to \cite{kumagai2014random} for an overview of random walks on critical percolation and related models.  Early investigations of the random walk on percolation include experiments by Last and Thouless \cite{last1971percolation}, who studied the conductivity of a sheet of graphite with randomly placed holes in it, and Brandt \cite{brandt1975solubilities}, who numerically studied the diffusion of gases in glasses using random walk on percolation. Kesten was the first to rigorously study random walk on percolation. In 1986, he introduced the \emph{infinite incipient cluster} (IIC) in two dimensions \cite{kesten1986incipient}. This is  the random configuration which arises by starting off with critical percolation, conditioning on the event that the origin is connected to the ball of radius $n$ by an open path, and then taking a limit as $n \to \infty$.  In the same year, he showed that simple random walk on the IIC is subdiffusive \cite{ kesten1986subdiffusive}, meaning that there exists $\varepsilon > 0$ so that it takes at least $n^{2+\varepsilon}$ steps for it to reach Euclidean distance $n$. This is much longer than the $n^2$ steps it takes for random walk on the entire lattice.  The reason for this phenomenon is that one can think of a large critical percolation cluster as consisting of a long backbone together with many smaller clusters hanging off it. Since the number of vertices in the backbone is negligible compared to the size of the entire cluster, the walk spends most of its time exploring the dead end clusters hanging off the backbone. Thus, the study of random walk on percolation and effective resistance can tell us a great deal about the structure of critical clusters. We refer to \cite{nolin2023backbone} for more background on the percolation backbone and its relation to the random walk. Recently, Ganguly and Lee \cite{ganguly2020chemical} improved Kesten's subdiffusivity result by showing that it holds with the intrinsic metric in place of the Euclidean metric.

A number of recent works have been aimed at the scaling limit of random walk on high-dimensional critical percolation clusters.  For example, Ben Arous, Cabezas and Fribergh \cite{ben2019scaling} obtained the scaling limit of random walk on the range of branching random walk in high dimensions and recently in the case of critical lattice trees \cite{bcf2025highdimensional}.  These are both simpler models which are conjectured to be in the same universality class as random walk on high-dimensional critical percolation clusters.

The \emph{Alexander-Orbach conjecture}  states that the \emph{spectral dimension} of the IIC for critical percolation is $4/3$. Recall that the spectral dimension of an infinite connected graph $G$ is defined as 
\begin{equation}
    d_s = d_s(G)=-2\lim_{n\rightarrow\infty}\frac{p_{2n}(x,x)}{\log n},
\end{equation}
where $p_n$ is the $n$-step transition kernel of simple random walk on $G$.  When $d_s=4/3$, we have that $p_{2n}(x,x) \sim n^{-2/3+o(1)}$ as $n \to \infty$.  This conjecture was confirmed for critical percolation in high dimensions by Kozma and Nachmias \cite{kozma2009alexander}. It is expected to be false in dimensions $d<6$, however, numerical simulations suggest that it is remarkably close to being true \cite{avraham2000diffusion}. In particular, the numerical estimate for planar critical percolation is roughly $1.318\pm0.001$, which is not far from the conjectured value.

In the case of supercritical percolation, the clusters converge to the Euclidean space, and the annealed invariance principle towards an ordinary Brownian motion was proved in \cite{mfgw1985clt}.  This was extended to a quenched invariance principle by Sidoravicius and Sznitman \cite{sidoravicius2004quenched} for $d \geq 4$ and independently by Berger and Biskup \cite{berger2007quenched} and Mathieu and Piatnitski \cite{mathieu2007quenched} for $d \geq 2$.  Gaussian heat kernel estimates in this setting were proved by Barlow \cite{barlow2004random}. In the intervening years, Grimmett, Kesten, and Zhang \cite{grimmett1993random} showed that random walk on supercritical percolation is recurrent in dimension $2$ and transient in dimensions $d\geq3$.

Another line of research from the last $40$ years has focused on constructing random processes on deterministic fractals, with the aim of developing a better understanding of random walk on low-dimensional critical percolation clusters. Some of the early investigations were by physicists in the 1980s, see \cite{alexander1982density, gefen1983fractal}. In 1987, the mathematicians Kusuoka \cite{kusuoka1987diffusions} and Goldstein \cite{goldstein1987random} constructed the scaling limit of the random walk on the \emph{Sierpinski gasket}. Barlow and Perkins \cite{barlow1988brownian} called this process \emph{Brownian motion on the Sierpinski gasket} and further studied its properties. This spurred a great deal of activity by e.g., Lindstrom \cite{lindstrom1990brownian}, Kigami \cite{kigami2001analysis}, and Barlow and Bass \cite{barlow1988brownian}, studying random walks on other deterministic fractals. We refer to Barlow's St.\ Flour lecture notes from 1995 \cite{barlow1995diffusions} for an overview of the work on random walk on (deterministic) fractals.

Kigami developed a useful framework for constructing and studying stochastic processes on low-dimensional fractals, which connects these processes to resistance forms and resistance metrics. See \cite{kigami2001analysis, kigami2012resistance} for an overview. It is precisely this framework that we rely on to obtain the scaling limit of random walk on critical percolation. In particular, we use a result by Croydon \cite{croydon2008convergence} to transfer the scaling limit of the resistance metric to that of the random walk.

\subsubsection{Intrinsic metric}\label{sec:background-intrinsic}
The study of intrinsic (e.g., chemical) distances and scaling limits of large percolation clusters has a rich history. In the case of supercritical percolation, increasingly strong results proved by G\"artner and Molchanov, Antal and Pisztora, and Garet and Marchand \cite{gartner1990parabolic, antal1996chemical, Garet2007large} showed that the intrinsic metric behaves like the Euclidean metric in all dimensions.  This implies that the scaling limit of a large supercritical percolation cluster is simply Euclidean space.  These works complement the scaling limit results for random walk towards Brownian motion in this setting mentioned above \cite{mfgw1985clt,sidoravicius2004quenched,berger2007quenched,mathieu2007quenched}. For high-dimensional ($d>11$ for nearest-neighbor and $d>6$ for spread-out) critical percolation, it was shown using lace expansion by Hara and Slade \cite{hara1990mean} that the intrinsic distance between two points at Euclidean distance $n$ is of order $n^2$, and the limiting distribution of the rescaled distance between two points was obtained recently in \cite{cchs2025chemicalhighd}. However, existence of the large cluster scaling limit is still an open problem. Blanc, Renaudie, Broutin and Nachmias \cite{blanc2024scaling} established the scaling limit of percolation on the hypercube. See also recent work by Chatterjee, Chinmay, Hanson and Sosoe \cite{chatterjee2025robust} on the IIC in high dimensions. In low dimensions, the chemical distance is not well-understood, even by physicists. For planar critical percolation, it is conjectured (e.g., in Schramm's 2006 ICM contribution \cite[Problem 3.3]{schramm2011conformally}) that the intrinsic distance between two points at Euclidean distance $n$ should also be of order $n^{\beta}$ for some exponent $\beta$. In 1999, Aizenman and Burchard \cite{aizenman1999holder} showed that $\beta>1$. About twenty years later, Damron, Hanson and Sosoe \cite{damron2021strict} obtained an upper bound of $\beta < 4/3$. Various simulations are consistent with $\beta\approx1.131$ \cite{grassberger1999pair, herrmann1988fractal, zhou2012shortest}. Further numerical work by Pose, Schrenk, Araujo and Hermann \cite{pose2014shortest} compared the geodesics to $\SLE_\kappa$ curves with $\kappa\approx1.04$.

\subsubsection{Other models}

Barlow, Croydon and Kumagai \cite{barlow2010random} obtained a (subsequential) scaling limit of random walk on the planar uniform spanning tree. In high dimensions, the analogous result was proved by Archer, Nachmias and Shalev \cite{ans2024ust}. On a tree, the resistance metric is the same as the geodesic metric, which makes it more straightforward to obtain the scaling limit of the random walk. See \cite{archer2020infinite, archer2021brownian, archer2023scaling, archer2025quenched, athreya2017invariance, croydon2008convergence, krebs1995brownian} for some examples of random walk scaling limits on other tree-like critical models.

A related model is \emph{long-range} percolation, where every pair of sites has a positive probability of being connected, with the connection probability behaving like a power law in the distance between vertices.  Berger--Tokushige and Crawford--Sly \cite{berger2024scaling, crawford2012simple} established the scaling limit of random walk on long-range supercritical percolation in all dimensions.  It was very recently shown by Ding, Fang and Huang \cite{ding2023uniqueness} that the geodesic metric on critical long-range percolation has a scaling limit, improving on a subsequential convergence result by B\"aumler \cite{baumler2023distances}.  We refer to the discussion in \cite{ding2023uniqueness} for more background on long-range percolation and the intrinsic distance.

\subsection{Main results}
\label{subsec:statements}
We now give the statement of our main results. In order to be concrete, we take $\bT$\phantomsection{}\label{def:bT} to be the triangular lattice centered at $0$ which contains the horizontal axis and for each $n \in \N$ we take $\bT_n=\tfrac{1}{n}\bT$\phantomsection{}\label{def:bT_n} as mentioned above. Throughout the paper, we will work in the \emph{triangular coordinate system} where the vertical axis is the line intersecting the horizontal axis at $0$ at an angle of $\frac{\pi}{3}$. In particular, the point with triangular coordinates $(1,1)$ denotes the Euclidean point $(3/2,\sqrt{3}/2)$, and the notation $[0,1]^2$\phantomsection{}\label{def:unit-square} refers to the parallelogram enclosed by the triangular vertices $(0,0),(1,0),(0,1),(1,1)$, which correspond to the Euclidean vertices $(0,0),(1,0),(1/2,\sqrt{3}/2),(3/2,\sqrt{3}/2)$. Define $\Lambda_n$\phantomsection{}\label{def:Lambda_n} to be the $n\times n$ subgraph of $\bT_n$ restricted to vertices in $[0,1)^2$. By a slight abuse of notation, we also write $\Lambda_n$ to refer to the vertex set of $\Lambda_n$ and, as is standard in the percolation literature, we refer to the vertices of $\Lambda_n$ as \emph{sites}.

Let $\Omega_n:=\{0,1\}^{\Lambda_n}$ be the set of site percolation configurations on $\Lambda_n$ with $\omega(x)=0$ (resp.\ $\omega(x) = 1$) denoting the site $x\in\Lambda_n$ being closed (resp.\ open).  Throughout, we consider \emph{critical percolation} so that each site is open independently with probability $1/2$. We say two open sites $x,y\in\Lambda_n$ are \emph{connected} if there exists a nearest-neighbor path in $\Lambda_n$ from $x$ to $y$ consisting only of open sites.  This induces an equivalence relation on the set of open sites. The equivalence classes are called \emph{open clusters}; \emph{closed clusters} are defined analogously.

Let $\cX^{(1)}_n,\cX^{(2)}_n,\ldots$\phantomsection{}\label{def:cXnm} be the (open and closed) clusters of $\Lambda_n$ induced by $\omega\in\Omega_n$, ordered descending in Euclidean diameter and breaking ties arbitrarily. We define the following metrics on $\Lambda_n$. For $x,y\in\cX^{(m)}_n$, let $P_{n}(x,y)$\phantomsection{}\label{def:Pn} be the set of paths in $\cX^{(m)}_n$ between $x$ and $y$. Let $\diam(\gamma)$\phantomsection{}\label{def:diam} be the Euclidean diameter of a path $\gamma\in P_{n}(x,y)$. Define the \emph{path metric} $\dpt_n$ as 
\begin{equation}\label{eq:dpt_n}
    \dpt_{n}(x,y)=\inf_{\gamma\in P_{n}(x,y)}\diam(\gamma).
\end{equation}
We use the path metric as an auxiliary reference metric to distinguish points on the clusters whose Euclidean distance converges to zero but are separated by a macroscopic closed path.

We recall the definition of the effective resistance metric of a graph $G=(V,E)$ which is assumed to be simple, finite, and connected. For each $A,B\subset V$
define the \emph{effective resistance} between $A$ and $B$ to be
\[
D^{\res}_G(A,B) = \bigl( \min \{ \cE(f,f) \mid f\colon V \to \R,\,f|_A\equiv 1,\,f|_B\equiv 0 \} \bigr)^{-1}
\]
where
\[
\cE(f,f) = \sum_{\{u,v\} \in E} (f(u)-f(v))^2 .
\]
For $x,y\in V$, we write $D_G^{\res}(x,y):=D_G^{\res}(\{x\},\{y\})$. Then $D^{\res}_G$ is a metric on $V$.

We let $\Dg_n$\phantomsection{}\label{def:Dgeon} be the \emph{geodesic} or \emph{graph} distance on all clusters $\cX^{(m)}_n$, where neighboring sites are considered to have distance 1. Finally, we let $\Dr_n$\phantomsection{}\label{def:Dresn} be the \emph{resistance} metric on all $\cX^{(m)}_n$, where the conductances between neighboring sites are also set to 1. Let $(\median_n^{\geo})_{n\in\N}$\phantomsection{}\label{def:mediangeo} and $(\median_n^{\res})_{n\in\N}$\phantomsection{}\label{def:medianres}  
be sequences of deterministic constants, to be specified later. Define the rescaled metrics
\begin{equation}\label{eq:def_dgeon_dresn}
\dg_n:=(\median_n^{\geo})^{-1}\Dg_n \quad\text{and}\quad \dr_n:=(\median_n^{\res})^{-1}\Dr_n.
\end{equation}
Lastly, we let 
\begin{equation}\label{eq:counting_measure_cluster}
    \mu^{(m)}_n:=\frac{1}{n^2A(1,n)}\sum_{x\in\cX^{(m)}_n}\delta_x
\end{equation}
be the normalized counting measure on $\cX^{(m)}_n$, where $A(1,n)$ is the \emph{one-arm} probability as defined in Section~\ref{subsec:arm-exponents}.

Let $\Gamma$\phantomsection{}\label{def:Gamma} be a (nested) $\CLE_6$ in $(0,1)^2$. We define the set $\Lambda=\Lambda(\Gamma)$\phantomsection{}\label{def:Lambda} to be the set of \emph{prime ends}\footnote{Our notion of prime ends differs slightly from the classical notion of prime ends for domains $D \subseteq \C$, however if $D$ is a domain with locally connected boundary $\partial D$, then its prime ends in the classical sense correspond to the prime ends of $\partial D$ in our notion.} of $\Gamma$, i.e., $\Lambda$ is the completion of the set of points in $(0,1)^2$ that are not on a loop in $\Gamma$ with respect to the metric $\dpt$ defined below in~\eqref{eq:dpt-first}. Let $\cX^{(1)},\cX^{(2)},\ldots \subset\Lambda$\phantomsection{}\label{def:cXm} be the interior gaskets of $\Gamma$, ordered descending in Euclidean diameter, i.e.\ the connected components of $\Lambda$.

We also define the path metric on $\Lambda$ as follows. For $x,y\in\Lambda$, let $P(x,y)$\phantomsection{}\label{def:P(x,y)} be the set of admissible paths in $\Lambda$ between $x$ and $y$, where we say a path is \emph{admissible} if it is contained in $(0,1)^2$ and does not cross any loops of $\Gamma$. Note that by our definition, if e.g.\ $\pi$ is the natural map which projects a prime end to its corresponding point in Euclidean space and $\pi(x)=\pi(y)$ is on a loop $\gamma$ but $x$ is a prime end inside of $\gamma$, and $y$ is a prime end outside of $\gamma$, then no path between $x$ and $y$ is considered admissible. The topology of the clusters of $\Lambda$ is encoded by the path metric $\dpt$ defined by
\begin{equation}\label{eq:dpt-first}
    \dpt(x,y)=\inf_{\gamma\in P(x,y)}\diam(\gamma).
\end{equation}

Throughout the paper, we mostly consider $(\cX^{(m)}_n,\dpt_n)$, $(\cX^{(m)}_n,\dg_n)$, and $(\cX^{(m)}_n,\dr_n)$ as abstract metric spaces, equipped with an embedding $\pi_n\colon\Lambda_n\to[0,1]^2$\phantomsection{}\label{def:pin}. Similarly, we consider $(\cX^{(m)},\dpt)$ as an abstract metric space with an embedding $\pi\colon\Lambda\to[0,1]^2$\phantomsection{}\label{def:pi}. In the continuum, it follows from the results of \cite{amy2025tightness,my2025metric,my2025diffusion} that we can define on $\cX^{(m)}$ a conformally covariant geodesic metric $\dg$\phantomsection{}\label{def:dg} (the geodesic $\CLE_6$ metric) and a resistance metric $\dr$\phantomsection{}\label{def:dr} whose associated resistance form is conformally covariant (the $\CLE_6$~resistance metric). Each of the metrics $\dg$ and $\dr$ is defined up to a global deterministic constant. We will review the basic properties of $\dg$ and $\dr$ in Section~\ref{subsec:cle-metrics}. We also let $\mu^{(m)}$\phantomsection{}\label{def:mum} be (a multiple of) the canonical conformally covariant $\CLE_6$ gasket measure on $\cX^{(m)}$ from \cite{miller2024existence, gps}, whose properties we will recall in Section~\ref{subsubsec:measure-conv}.

\begin{thm}\label{thm:main-convergence}
Let $\dg$, $\dr$, and $\mu$ be suitable multiples of the geodesic CLE$_6$ metric, the CLE$_6$ resistance metric, and CLE$_6$ gasket measure, respectively. There exist sequences $(\median_n^{\geo})_{n\in\N},(\median_n^{\res})_{n\in\N}$ such that for all $m$,
\begin{equation}
    \begin{split}
        (\cX^{(m)}_n,\dpt_n,\mu^{(m)}_n,\pi_n)&\rightarrow (\cX^{(m)},\dpt,\mu^{(m)},\pi),\\
        (\cX^{(m)}_n,\dg_n,\mu^{(m)}_n,\pi_n)&\rightarrow (\cX^{(m)},\dg,\mu^{(m)},\pi),\\
        (\cX^{(m)}_n,\dr_n,\mu^{(m)}_n,\pi_n)&\rightarrow (\cX^{(m)},\dr,\mu^{(m)},\pi)\\
    \end{split}
\end{equation}
weakly with respect to the Gromov-Hausdorff-Prokhorov-function topology as $n\rightarrow\infty$.
\end{thm}

The definition of the normalizing constants $\median_n^{\geo}$ and $\median_n^{\res}$ is given in Section~\ref{subsubsec:apriori-sketch}. As a result of this, we obtain existence of the scaling exponent for the intrinsic distance and effective resistance between the left and right sides of the parallelogram.
\begin{cor}\label{cor:scaling-exponent-distance}
    Let $S_n$ be the geodesic distance (resp.\ effective resistance) between the left and the right side of $\Lambda_n$ which we denote by $\partial_L \Lambda_n, \partial_R \Lambda_n$, respectively. Let $\beta \in \{\beta^\geo,\beta^\res\}$ be the scaling exponent of the CLE$_6$ geodesic (resp.\ resistance) metric defined in \cite{my2025metric,my2025diffusion}. For each $q \in (0,1)$, we have
    \[
    \E[S_n \mid \partial_L \Lambda_n \leftrightarrow \partial_R \Lambda_n ] 
    \asymp \operatorname{quant}_q[S_n \mid \partial_L \Lambda_n \leftrightarrow \partial_R \Lambda_n ]
    \asymp \median_n
    = n^{\beta+o(1)} 
    \quad\text{as } n \to \infty ,
    \]
    where $\operatorname{quant}_q$ denotes the $q$-quantile, $\median_n$ denotes $\median_n^{\geo}$ (resp.\ $\median_n^{\res}$), and the multiplicative constants in~$\asymp$ depend on $q$.
\end{cor}

Similarly, we obtain the existence of the exponents for the chemical distance and effective resistance on the IIC.
\begin{cor}\label{co:exponent_iic}
    Let $\piic$ denote the law of the IIC on the entire lattice $\bT$ and let $S_n$ be the geodesic distance (resp.\ effective resistance) between $0$ and $\partial B(0,n)$. Let $\beta \in \{\beta^\geo,\beta^\res\}$ be the scaling exponent of the CLE$_6$ geodesic (resp.\ resistance) metric defined in \cite{my2025metric,my2025diffusion}. Under $\piic$, almost surely,
    \[
    S_n = n^{\beta+o(1)}
    \quad\text{as } n \to \infty .
    \]
    Moreover, for each $\varepsilon > 0$ and $\alpha > 0$ there exists $c(\varepsilon,\alpha)$ such that
    \[
    \piic[ S_n \notin [n^{\beta-\varepsilon},n^{\beta+\varepsilon}] ] \le cn^{-\alpha} .
    \]
\end{cor}
In this paper, $B(0,n)$ refers to the parallelogram of side length $2n$ centered at $0$, but of course the statement is also true if $B(0,n)$ is the Euclidean ball.

From the one-arm exponent for critical percolation established in \cite{lawler2002onearm}, it is known that the volume growth exponent equals $91/48$. By the results from \cite{kumagai2008heat}, this implies the following statements regarding (the existence of) the spectral dimension of the IIC.
\begin{cor}\label{co:spectral-dimension}
    Let $\beta$ be the scaling exponent of the CLE$_6$ resistance metric defined in~\cite{my2025diffusion}. Let $\mathrm{P}_0$ denote the law of the simple random walk on the IIC starting at $0$, let $\mathrm{p}_{2t}(0,0)$ denote its return probability after $2t$ steps, and let $\tau_n$ denote its hitting time of $\partial B(0,n)$. Under $\piic$, almost surely,
    \begin{align*}
    \mathrm{p}_{2t}(0,0) &= t^{-\frac{91/48}{91/48+\beta}+o(1)} \quad\text{as } t\to\infty ,\\
    \quad \mathrm{E}_0 \tau_n &= n^{91/48+\beta+o(1)} \quad\text{as } n\to\infty .
    \end{align*}
\end{cor}
Note that Corollary~\ref{co:spectral-dimension} states that the Einstein relations hold for the IIC in planar critical percolation.

As shown in~\cite{croydon2008convergence}, the Gromov-Hausdorff-Prokhorov convergence of a sequence of metric measure spaces where each metric is a resistance metric implies the convergence of the corresponding stochastic processes.  Thus by combining Theorem~\ref{thm:main-convergence} with \cite{croydon2008convergence}, we obtain the following convergence result for the simple random walk on $\cX_n$.

\begin{thm}\label{thm:main-random-walk}
There exists a sequence $(b_n)_{n\in\N}$ such that the following is true. Let $m \in \N$ and sample $\rho_n \in \cX^{(m)}_n$ and $\rho \in \cX^{(m)}$ such that $(\cX^{(m)}_n,\dpt_n,\pi_n,\rho_n) \to (\cX^{(m)},\dpt,\pi,\rho)$ in law in the rooted GHf topology. Then there exists a coupling such that the following holds almost surely. Let $X^{(n)}$ be the continuous-time simple random walk on $\cX^{(m)}_n$ that jumps to each open neighboring site at rate $b_n$, started from $\rho_n$. Then there exists a diffusion process $X$ on $\cX^{(m)}$ with $X_0=\rho$ such that $\pi_n(X^{(n)})$ converges in law to $\pi(X)$ as $n\rightarrow\infty$ in the space of c\`adl\`ag processes on $[0,1]^2$ equipped with the topology of uniform convergence on compact sets.
\end{thm}

We will see that the jump rate $b_n$ can be taken to be $n^2A(1,n)\median_n^{\res} = n^{91/96+\beta+o(n)}$. The limiting stochastic process $X$ is the \emph{$\CLE_6$ Brownian motion} as constructed in \cite{amy2025tightness,my2025diffusion}. See Section~\ref{subsec:resistance-metrics-walks} for a precise definition.

\subsection{Discussion}\label{sec:discussion}  Let us now make some comments about our main theorem statements.

\subsubsection*{Choice of lattice} The theorems above are stated for critical site percolation on the triangular lattice, since this is the model for which the convergence towards $\SLE_6$/$\CLE_6$ has been shown. However, the proofs of tightness in this paper only use certain arm exponents for percolation, which are also shown for the square lattice $\Z^2$.  The $\SLE$/$\CLE$ machinery is only used to prove the uniqueness of the subsequential limits. Therefore, our scaling limit result would apply to other lattices if the convergence to $\SLE_6$/$\CLE_6$ were shown.

\subsubsection*{Random walk speed} In this paper, we show convergence of the random walk that jumps to each of its open neighbors at the same rate $b_n$. This is also known as the \emph{variable-speed} random walk (VSRW), since the total jump rate at each site is proportional to the number of open neighbors. Another natural model is the \emph{constant-speed} random walk (CSRW), where the total jump rate is constant at each site. The stationary measure of the VSRW is the uniform measure on $\cX_n$, which is proportional to $\mu_n$. Or, in the language of resistance metrics: the VSRW is the unique \emph{Hunt process} associated with $(\cX_n,\dr_n,\mu_n)$ (see Section~\ref{subsec:resistance-metrics-walks} for more background).  This allows us to use the results of \cite{gps}, which give the convergence of the uniform measure on critical percolation clusters. The stationary measure of the CSRW is the \emph{degree-weighted} probability measure on $\cX_n$, so one would have to adapt the results from \cite{gps} to this setting. This should not be too challenging, since both measures only differ by at most a factor of $6$ (the maximum number of neighbors at a site) on each set.

\subsubsection*{Alexander--Orbach} The resistance estimates we show in this paper can be used to obtain bounds on the spectral dimension of large critical percolation clusters in two dimensions. Corollary~\ref{co:spectral-dimension} states that the spectral dimension of the IIC equals $\frac{2D}{D+\beta}$ where $D=91/48$ is the volume growth exponent and $\beta$ is the resistance growth exponent.

We do not obtain a precise value for $\beta$, but it is known that $\beta\in [3/4,4/3]$ \cite{my2025diffusion}, implying an upper bound on the spectral dimension of approximately $1.433$. So although this does not disprove the Alexander-Orbach conjecture, we do obtain an upper bound that is not far from the conjectured value of $4/3$. To disprove the conjecture, one would have to improve the lower bound of the resistance growth $\beta$ to $91/96 \approx0.948$.

\subsubsection*{Explicit lengths} In this paper (see Corollary~\ref{cor:scaling-exponent-distance}), we show that distances in critical planar percolation scale like $n^{\beta+o(1)}$ for some exponent $\beta$, thus partially answering \cite[Problem 3.3]{schramm2011conformally}. We do not obtain an explicit expression for $\median_n$ in terms of $n$ and thus do not obtain any precise estimate for the exponent $\beta$.

\subsection*{Acknowledgements}
J.M.\ and Y.Y.\ were supported by ERC starting grant SPRS (804116) and from ERC consolidator grant ARPF (Horizon Europe UKRI G120614). Y.Y.\ in addition received support from the Royal Society. M.M.\ was supported by the University of Cambridge Harding Distinguished Postgraduate Scholarship Programme. I.Đ.\ was supported by a Cambridge International Scholarship from the Cambridge Trust.

\subsection*{Outline}\label{sec:outline}
The rest of the paper is organized as follows. We first state the preliminary results needed for the rest of the paper in Section~\ref{sec:preliminaries}. In Section~\ref{sec:sketch}, we give a detailed overview of the proofs of Theorems~\ref{thm:main-convergence} and~\ref{thm:main-random-walk} and Corollaries~\ref{cor:scaling-exponent-distance}--\ref{co:spectral-dimension}. The rest of the paper is devoted to the proofs of these theorems. In Section~\ref{sec:convergence-path-metric-proof}, we show the convergence of the percolation clusters as metric spaces defined by $\dpt_n$. As a first step towards tightness, we prove in Section~\ref{sec:polynomial_tails} polynomial tails for the geodesic distances and effective resistances in various setups. We then improve these results to superpolynomial tails in Section~\ref{sec:superpolynomial-concentration} using spatial independence arguments. We then conclude the proof of tightness in Section~\ref{sec:tightness}. In Section~\ref{sec:uniqueness-proof}, we identify the subsequential limits with the geodesic CLE$_6$ metric and the CLE$_6$ resistance metric, respectively, yielding the full convergence of the metrics in Theorem~\ref{thm:main-convergence}. In Section~\ref{sec:ghpf-proof}, we explain how the results of \cite{gps} and \cite{my2026upcoming} imply the convergence of $\mu_n^{(m)}$ to $\mu^{(m)}$, completing the proof of Theorems~\ref{thm:main-convergence} and~\ref{thm:main-random-walk}. Finally, in Section~\ref{sec:scaling-exponent-proof}, we deduce the existence of the exponents in the Corollaries~\ref{cor:scaling-exponent-distance}--\ref{co:spectral-dimension}. We refer to Section~\ref{sec:sketch} for a more detailed outline of the arguments from Sections~\ref{sec:convergence-path-metric-proof}--\ref{sec:scaling-exponent-proof}.

\section{Preliminaries}
\label{sec:preliminaries}

The purpose of this section is to collect a number of preliminaries.  First, we will describe the notation that we will use in the rest of this article in Section~\ref{subsec:notation}. In Section~\ref{subsec:topologies} we will review the variants of the Gromov-Hausdorff topology which we will consider. In Section~\ref{subsec:arm-exponents}, we will recall the basic results about arm exponents for critical percolation. We will then give the statements of the results regarding the convergence of percolation in Section~\ref{subsec:percolation_convergence_results}. The purpose of Section~\ref{subsec:cle-metrics} is to review the basic properties of $\CLE$ metrics (both geodesic and resistance).  We finally recall the relationship between resistance metrics and stochastic processes in Section~\ref{subsec:resistance-metrics-walks}.

\subsection{Notation}
\label{subsec:notation}

Throughout the paper, we refer to points in $\R^2$ by their triangular coordinates as explained at the beginning of Section~\ref{subsec:statements}, in particular $[a,b) \times [c,d)$ denotes a parallelogram with interior angles $\pi/3,2\pi/3$.

Let $R\in\N$ and let $\cD^{R}_k=\{(iR^{-k},jR^{-k}):0\leq i,j< R^k\}$\phantomsection{}\label{def:D_kR} be the set of $R$-adic points in $[0,1)^2$ at level $k$. Define $\cB_k^R:=\{[iR^{-k},(i+1)R^{-k})\times[jR^{-k},(j+1)R^{-k}):0\leq i,j< R^k\}$\phantomsection{}\label{def:B_kR}, i.e., the set of $R$-adic boxes at level $k$. Let $\cQ_k^R:=\{\text{int}(B_1\cup B_2\cup\ldots)\text{ simply connected}: B_1,B_2,\ldots\in\cB_k^R\}$\phantomsection{}\label{def:cQ_kR}. Further define $\cB^R:=\bigcup_k\cB_k^R$\phantomsection{}\label{def:cBR}, $\cQ^R:=\bigcup_{k}\cQ_k^R$\phantomsection{}\label{def:cQR} and $\cD^R:=\bigcup_k\cD_k^R$\phantomsection{}\label{def:cDR}. We also let $\cB^R_{k,n}:=\{[\frac{i}{n}R^{\kappa-k},\frac{i+1}{n}R^{\kappa-k})\times[\frac{j}{n}R^{\kappa-k},\frac{j+1}{n}R^{\kappa-k}):0\leq i,j< R^{k-1}\}$\phantomsection{}\label{def:cBR_k,n} and $\cD_{k,n}^R:=\{(\frac{i}{n}R^{\kappa-k},\frac{j}{n}R^{\kappa-k})\}_{0\leq i,j<R^{k-1}}$\phantomsection{}\label{def:cDR_k,n}, where $\kappa\in\N$\phantomsection{}\label{def:kappa} is such that $R^{\kappa-1}<n\leq R^{\kappa}$. The boxes in $\cB^{R}_{k,n}$ are also called $R$-adic boxes at level $k$. The distinction with $\cB^R_k$ should be clear from the context. When $R=2$, we suppress $R$ from the notation and simply write $\cD_k,\cB_k,\cQ\ldots$\phantomsection{}\label{def:suppress_R}.  

For $x\in[0,1]^2$ and $r_2>r_1>0$, we write $B(x,r_1)=x+[-r_1,r_1)^2$\phantomsection{}\label{def:box} and $A(x;r_1,r_2)=B(x,r_2)\setminus B(x,r_1)$\phantomsection{}\label{def:annulus} for the box centered at $x$ of side length $2r_1$ and for the parallelogram-shaped annulus of side lengths $2r_1$ and $2r_2$ centered at $x$ respectively. For some box $B=B(x,r)$ and $r'>0$, we write $r'B:=B(x,r'\cdot r)$. Throughout the paper, we will sometimes let a box $B$ denote either a subset of $[0,1]^2$ as above or the subgraph of $\Lambda_n$ induced by the vertices inside $B$. The exact meaning should be clear from the context.

For $Q\in\cQ$\phantomsection{}\label{def:Q-in-cQ}, we let $\Lambda_{Q,n}$\phantomsection{}\label{def:Lambda_Q,n} be the set of sites in $\bT_n$ that are contained in $Q$ and $\Omega_{Q,n}:=\{0,1\}^{\Lambda_{Q,n}}$ be the set of percolation configurations on $\Lambda_{Q,n}$. If $\Gamma_Q$ is a collection of $\CLE_6$ loops in $Q$, we let $\Lambda_Q=\Lambda_Q(\Gamma_Q)$\phantomsection{}\label{def:Lambda_Q} be the set of prime ends associated with $\Gamma_Q$. In the paper, we might sometimes write $Q$ to denote $\Lambda_Q$ by an abuse of notation when it is clear from the context what is meant. We refer to Section~\ref{subsubsec:convergence-boundary-loops} for more background on CLE and the relationship between $\Gamma$ and $\Gamma_Q$.

We say two points $x,y\in\Lambda_Q$ are in the same cluster if there exists a path from $x$ to $y$ in the interior of $Q$ that does not cross a loop in $\Gamma_Q$, or in other words, if there exists a path from $x$ to $y$ that is continuous with respect to $\dpt_Q$. We label the clusters of $\Gamma_Q$ in descending order of Euclidean diameter as $\cX_{Q}^{(1)},\cX_Q^{(2)},\ldots$\phantomsection{}\label{def:cXQm}. The (open and closed) percolation clusters of $\Gamma_{Q,n}$ are labeled similarly as $\cX_{Q,n}^{(1)},\cX_{Q,n}^{(2)},\ldots$\phantomsection{}\label{def:cXQmn}, breaking ties arbitrarily. We further denote by $\cX$\phantomsection{}\label{def:cX-first} (resp.\ $\cX_n$\phantomsection{}\label{def:cXn-first}) the outermost open $\CLE_6$ (resp.\ percolation) cluster surrounding the point $(\frac{1}{2},\frac{1}{2})$.\footnote{We say that a point is surrounded by a loop $\gamma \in \Gamma$ (or is inside $\gamma$) if the winding number of $\gamma$ about $x$ is $\pm 1$ (depending on the choice of orientation).}

Let $d_{Q,n}$ be any (path, geodesic, or resistance) metric on $\Lambda_{Q,n}$. We say, e.g., the (collection of) metric space(s) $(\Lambda_{Q,n},d_{Q,n})$ converges to $(\Lambda_Q,d_Q)$ in some metric $\Delta$ if $(\cX^{(m)}_{Q,n},d_{Q,n})$ converge to $(\cX^{(m)}_{Q},d_{Q})$ in $\Delta$ for each $m$. This convergence is metrized by $\Delta^\infty((\Lambda_{Q,n},d_{Q,n}),(\Lambda_Q,d_Q))=\sum_{m=1}^\infty2^{-m}\!\left(\Delta((\cX^{(m)}_{Q,n},d_{Q,n}),(\cX^{(m)}_Q,d_Q))\wedge 1\right)$. In other words, we say $(\Lambda_{Q,n},d_{Q,n})$ converges to $(\Lambda_Q,d_Q)$ if, for every $m$, the $m$th largest (in Euclidean diameter) cluster of $\Lambda_{Q,n}$ converges to the $m$th largest cluster of $\Lambda_Q$.

We let $\dist$\phantomsection{}\label{def:dist}, $\diam$\phantomsection{}\label{def:diam-second} denote the Euclidean distance and diameter on $[0,1]^2$.

\subsection{Variants of the Gromov-Hausdorff metric}
\label{subsec:topologies}

We will now describe some variants of the Gromov-Hausdorff (GH) metric that we will consider in this paper. To begin with, we need to recall the definition of the ordinary GH metric (see \cite[Chapter~27]{villani2009optimal} for a more in depth review), which we recall defines a distance between compact metric spaces.  The starting point is the definition of the ordinary Hausdorff metric.  Suppose that $(W,d_W)$ is a metric space and $K_1,K_2 \subseteq W$ are non-empty and compact.  Then the Hausdorff distance between $K_1,K_2$ is given by
\begin{equation}\label{eq:hausdorff}
    \Delta_{\Haus}(K_1,K_2) = \inf\{ \epsilon > 0 : \forall x_i \in K_i, \exists x_{3-i} \in K_{3-i} \text{ s.t.\ } d_W(x_i,x_{3-i}) < \epsilon \text{ for } i=1,2 \}.
\end{equation}
Now suppose that $(X_i,d_i)$ for $i=1,2$ are non-empty compact metric spaces.  Then the GH distance between $(X_1,d_1)$ and $(X_2,d_2)$ is defined as
\begin{equation}\label{eq:def-gh-embeddings}
    \Delta_{\GH}((X_1,d_1),(X_2,d_2)):=\inf_{\psi_1,\psi_2} \Delta_{\Haus}(\psi_1(X_1),\psi_2(X_2)),
\end{equation}
where $\Delta_{\Haus}$ is the Hausdorff distance and the infimum is taken over all metric spaces $(W,d_W)$ and isometric embeddings $\psi_i \colon X_i\to W$ for $i=1,2$. It is sometimes more convenient to express the GH distance in terms of \emph{correspondences}. A correspondence is a subset $\cR\subseteq X_1\times X_2$ such that for every $x\in X_1$, there exists $y\in X_{2}$ such that $(x,y)\in\cR$ and for every $y\in X_2$, there exists $x\in X_1$ such that $(x,y)\in\cR$. Then
\begin{equation}\label{eq:def-gh-correspondence}
    \Delta_{\GH}((X_1,d_1),(X_2,d_2))=\frac{1}{2}\inf_{\cR}\sup_{(x,y),(x',y')\in \cR}|d_1(x,x')-d_2(y,y')|,
\end{equation}
where the infimum is over all correspondences $\cR$. Furthermore, since the $(X_i,d_i)$ are compact, there exist isometric embeddings $\psi_1,\psi_2$ and a correspondence $\cR$ such that $\psi_1,\psi_2$ achieve the infimum in~\eqref{eq:def-gh-embeddings} and $\cR$ achieves the infimum in~\eqref{eq:def-gh-correspondence}. 
In addition, if $(X_n,d_n)_{n\in\N}$ converges to some $(X,d)$ in the GH distance, then there exists some compact metric space $(W,d_W)$ and isometries $\psi_n\colon X_n\to W$, $\psi\colon X\to W$ such that $\psi_n(X_n)$ converges to $\psi(X)$ in Hausdorff distance.

Lastly, we can still define the GH distance for pseudo-metric spaces. Then if $\Delta_{\GH}((X_n,d_n),(X,d))\rightarrow0$ and $(X,d)$ is a pseudo-metric space, then $(X_n,d_n)$ converges in GH distance to the quotient metric space $(X/{\sim},d)$, where $x\sim y$ in $X$ if and only if $d(x,y)=0$.

Next, we define the \emph{Gromov-Hausdorff-function} (GHf) distance (see \cite[Appendix~A]{amy2025tightness} for more details). Let $(Y,d_Y)$ be a proper metric space (in the sense that every closed ball is compact). Suppose that $X_i\subseteq W$ are non-empty compact subsets for $i=1,2$ and let $f_i\colon X_i\to Y$ be bounded continuous functions. Then we define $\Delta_f(f_1,f_2)$ to be the infimum over $\delta > 0$ such that for each $x_i \in X_i$ there is $x_{3-i} \in X_{3-i}$ with $d_W(x_i,x_{3-i})\leq\delta$ and $d_Y(f_i(x_i),f_{3-i}(x_{3-i}))\leq\delta$. Now let $(X_i,d_i)$ for $i=1,2$ be general compact spaces (so not necessarily a subspace of the same space $W$). Then the GHf distance between $(X_1,d_1,f_1)$ and $(X_2,d_2,f_2)$ is defined as 
\begin{equation}\label{eq:ghf}
    \Delta_{\GHf}((X_1,d_1,f_1),(X_2,d_2,f_2)):=\inf_{\psi_1,\psi_2}\!\left\{\Delta_{\Haus}(\psi_1(X_1),\psi_2(X_2))+\Delta_f(f_1\circ\psi_1^{-1},f_2\circ\psi_2^{-1}) \right\}.
\end{equation}

We now turn our attention to metric spaces equipped with a measure. Suppose that we have the setting described just above and for each $i=1,2$, let $\mu_i$ be a finite measure on $X_i$. Then the \emph{Gromov-Hausdorff-Prokhorov} (GHP) distance between $(X_1,d_1,\mu_1)$ and $(X_2,d_2,\mu_2)$ is given by
\begin{equation}\label{eq:ghp}
    \Delta_{\GHP}((X_1,d_1,\mu_1),(X_2,d_2,\mu_2)) :=\inf_{\psi_1,\psi_2}\!\left\{\Delta_{\Haus}(\psi_1(X_1),\psi_2(X_2))+\Delta_{\Pro}(\mu_1\circ\psi_1^{-1},\mu_2 \circ \psi_2^{-1})\right\},
\end{equation}
where $\Delta_{\Pro}$ is the Prokhorov distance. We also define the \emph{Gromov-Hausdorff-Prokhorov-function} (GHPf) distance as
\begin{equation}\label{eq:ghpf}
    \begin{split}
        &\Delta_{\GHPf}((X_1,d_1,\mu_1,f_1),(X_2,d_2,\mu_2,f_2))\\
        &:=\inf_{\psi_1,\psi_2}\!\left\{\Delta_{\Haus}(\psi_1(X_1),\psi_2(X_2))+\Delta_{\Pro}(\mu_1\circ\psi_1^{-1},\mu_2\circ\psi_2^{-1})+\Delta_{f}(f_1\circ\psi_1^{-1},f_2\circ\psi_2^{-1}\right\}.
    \end{split}
\end{equation}
Lastly, we define the rooted Gromov-Hausdorff-Prokhorov-function distance (rGHPf) between function-equipped metric measure spaces which are marked by a distinguished point (i.e., the root). Let $\rho_i\in X_i$ for $i=1,2$. Then
\begin{equation}\label{eq:rghpf}
    \begin{split}
        & \Delta_{\rGHPf}((X_1,d_1,\mu_1,f_1,\rho_1),(X_2,d_2,\mu_2,f_2,\rho_2))
        :=\inf_{\psi_1,\psi_2}\big\{ \Delta_{\Haus}(\psi_1(X_1),\psi_2(X_2))+\\
         &\Delta_{\Pro}(\mu_1\circ\psi_1^{-1},\mu_2\circ \psi_2^{-1})+
         \Delta_{f}(f_1\circ\psi_1^{-1},f_2\circ\psi_2^{-1})+d_W(\rho_1,\rho_2)\big\}.
    \end{split}
\end{equation}

\subsection{Arm exponents}\label{subsec:arm-exponents}
We gather all the results on the arm exponents in annuli that we will later use. That is, consider the annulus determined by two boxes $b= B(0,r)$ and $B= B(0, R)$, where $r< R$ are positive integers. Recall that the boxes and annuli in this paper are always considered to be parallelograms. Let $j\in \mathbf{N}$ and let $\sigma$ be a sequence of $j$ colors (open/closed). We consider the probability $A_{j,\sigma}(b,B)=A_{j,\sigma}(r, R)$\phantomsection{}\label{def:arm-percolation} of the event that when considering percolation on $\bT$, there exist $j$ disjoint crossings of the annulus $B\setminus b$ that, in clockwise order, have the colors prescribed by $\sigma$. If $j=1$, we simply write $A_1(r,R)$ and omit the $\sigma$. We refer to these crossings as \emph{arms}. 

We can similarly define the arm probabilities $A^{\CLE}_{j,\sigma}(r,R)$\phantomsection{}\label{def:arm-cle} for $\CLE_6$ and $0<r<R$. It is shown in \cite[Section~2.4]{gps} that $A^{\CLE}_{j,\sigma}(r,R) = \lim_{n\to\infty}A_{j,\sigma}(nr,nR)$.

We will repeatedly use the following quasi-multiplicativity result shown in \cite[Proposition 16]{nolin2008near}:
\begin{prop}\label{lem:quasi-multiplicativity-arms} Fix some $j\in \mathbf{N}$ and let $\sigma$ be any sequence of $j$ colors. Then there exist $n_0(j)\in \mathbf{N}$ and an absolute constant $C\in (0,1)$ such that $$ C\cdot A_{j,\sigma}(n_1,n_2)\cdot A_{j,\sigma}(n_2,n_3) \leq A_{j,\sigma}(n_1,n_3) \leq  A_{j,\sigma}(n_1,n_2)\cdot A_{j,\sigma}(n_2,n_3) $$ for all $n_0(j)\leq n_1 < n_2 <n_3.$ \label{prop:quasi_multiplicativity}
\end{prop}
\noindent The $n_0(j)$ is there for technical reasons; the $j$ arms need to have enough space to approach the box of radius $n_1$.

We are also interested in the probability $A^h_{j,\sigma}(b,B)=A^h_{j,\sigma}(r, R)$\phantomsection{}\label{def:arm-percolation-half} of the event that there are $j$ arms of colors prescribed by $\sigma$ that cross the annulus $A(0;r,R)$ and stay in the upper half-plane $\{(x,y) \mid y>0\}$. Quasi-multiplicativity holds for half-plane events as well -- the proof of Proposition~\ref{prop:quasi_multiplicativity} given in \cite{nolin2008near} also applies to this case.

 For non-constant sequences $\sigma$, we know by \cite[Theorem 4]{smirnov2001criticalexponents} that for $j>1$ $$A_{j,\sigma}(r,R)=R^{-(j^2-1)/12+o(1)}, \textrm{ as }R\rightarrow \infty.$$ Further, \cite[Theorem 3]{smirnov2001criticalexponents} says that for \textit{any} $\sigma$ we know that $$A_{j,\sigma}^h(r,R)=R^{-j(j+1)/6+o(1)}, \textrm{ as }R\rightarrow \infty.$$ 
 
 Combining these results and Proposition~\ref{prop:quasi_multiplicativity} we get the following lemma.

\begin{lem}
\label{lem:arm_exponents}
For every $\zeta>0$, $j\geq 1$, there exists $\rho_{\zeta,j}>0$ big enough such that
\begin{enumerate}
\item[(a)] If $j > 1$, for every non-constant sequence of colors $\sigma$ and
for all $R\geq r \geq \rho_{\zeta,j}$ 
$$\!\left(\frac{R}{r}\right)^{-\alpha_j-\zeta}\leq\mathbf{P} (j \textrm{ arms of colors given by }\sigma \textrm{ in } A(r,R) )\leq \!\left(\frac{R}{r}\right)^{-\alpha_j+\zeta},$$
where $\alpha_j=\frac{j^2-1}{12}$ are the full-plane polychromatic exponents.
\item[(b)] For every sequence of colors $\sigma$ and for all $R\geq r\geq \rho_{\zeta,j}$ 
\begin{equation*}
    \begin{split}
        \!\left(\frac{R}{r}\right)^{-\alpha^h_j-\zeta}\leq\mathbf{P} (j \textrm{ arms of colors given by }\sigma \textrm{ in } A(r,R) {\textrm{ in the half-plane}})
        \leq \!\left(\frac{R}{r}\right)^{-\alpha^h_j+\zeta},
    \end{split}
\end{equation*}
where $\alpha^h_j=\frac{j(j+1)}{6}$ are the half-plane exponents. 
\end{enumerate}
\end{lem}

In the half-plane, there is no difference between polychromatic and monochromatic events -- the reason is that one can explore the crossings from ``below'' one-by-one and flip all colors above some of them without changing the probabilities of configurations. However, we cannot do this in the full-plane case, as the lowest crossing is no longer well defined. The following lemma, given in \cite[Theorem 4 and Theorem 5]{beffara2011monochromatic}, tells us that monochromatic exponents also exist and bounds them in terms of the polychromatic exponents.

\begin{lem}
For every $\zeta>0$, $j\geq 2$, and a constant sequence of colors $\sigma$, there exists $\rho_{\zeta,j}'>0$ big enough such that for all $R\geq r \geq \rho_{\zeta,j}'$ 
$$\!\left(\frac{R}{r}\right)^{-\alpha_j'-\zeta}\leq\mathbf{P} (j \textrm{ arms of colors prescribed by }\sigma \textrm{ in } A(r,R) )\leq \!\left(\frac{R}{r}\right)^{-\alpha_j'+\zeta},$$
where $\alpha_j'$ are the monochromatic exponents that satisfy $\alpha_j<\alpha_j'<\alpha_{j+1}$.\label{lem:arm_exponents_monochromatic}
\end{lem}

\subsection{Percolation convergence results}
\label{subsec:percolation_convergence_results}

We are now going to review prior convergence results for percolation that will be useful in this work.

\subsubsection{Convergence of interface loops and the quad-crossing space}
\label{subsubsec:convergence-boundary-loops}
Let $Q\in\cQ$ be a union of dyadic boxes (see Section~\ref{subsec:notation}). Recall that $\Lambda_{Q,n}$ is the set of sites in $\bT_n$ that are contained in $Q$, and let $\partial\Lambda_{Q,n}$ be the set of sites in $\bT_n$ that are not contained in $Q$, but have a neighbor in $\Lambda_{Q,n}$. If we consider all the sites in $\partial\Lambda_{Q,n}$ to be closed, we obtain a collection $\Gamma_{Q,n}$ of loops which are contained in $Q$, where each loop $\gamma\in\Gamma_{Q,n}$\phantomsection{}\label{def:Gamma_Q,n} is the interface loop of a percolation cluster (i.e., the interface between open and closed). If we let $W$ be the (not necessarily simply connected) union of finitely many $Q\in\cQ$, we define the collections $\Gamma_{W,n}$ and $\Gamma_{W^c,n}$ analogously. If $Q=(0,1)^2$, we just write $\Gamma_n:=\Gamma_{(0,1)^2,n}$\phantomsection{}\label{def:Gamma_n}.

For two curves $\eta,\eta'\colon[0,1]\to[0,1]^2$, we define the uniform metric modulo parameterization as
\begin{equation}\label{eq:delta-cur-paths}
    \Delta_{\cur}(\eta,\eta')=\inf_{\psi} \sup_{t\in[0,1]}\|\eta(\psi(t))-\eta'(t)\|_2,
\end{equation}
where the infimum is over all homeomorphisms $\psi \colon [0,1] \to [0,1]$ that fix 0 and 1 and $\|\cdot\|_2$ refers to the Euclidean norm. For two loops $\gamma,\gamma'\colon S^1\to[0,1]^2$, we also define the uniform metric modulo parameterization as
\begin{equation}\label{eq:delta-cur-loops}
    \Delta_{\cur}(\gamma,\gamma')=\inf_{\psi} \sup_{t\in[0,1]}\|\gamma(\psi(t))-\gamma'(t)\|_2,
\end{equation}
where the infimum is now over all orientation preserving homeomorphisms $\psi \colon S^1 \to S^1$ from the unit circle to itself and $\|\cdot\|_2$ refers to the Euclidean norm. For two collections of curves $\Gamma,\Gamma'$, we let $\Delta_{\Hc}(\Gamma,\Gamma')$\phantomsection{}\label{def:delta-Hc} be the Hausdorff distance between $\Gamma$ and $\Gamma'$ with respect to $\Delta_{\cur}$.

Camia and Newman \cite{camia2006two} showed that the percolation interface loops converge weakly to $\CLE_6$.
\begin{thm}
\label{th:camia-newman}
The collection $(\Gamma_{n})_{n\in\N}$ converges weakly to the collection of $\CLE_6$ loops $\Gamma$ with respect to $\Delta_{\Hc}$ as $n\rightarrow\infty$.
\end{thm}

For this paper, we will need to consider not only the percolation on the entire box $\Lambda_n$, but also on subboxes $\Lambda_{Q,n}$. A different way to encode a percolation configuration other than by its interface loops, is by considering the collection of \emph{quads} it crosses, which are homeomorphisms $R\colon[0,1]^2\to\C$. A quad is crossed if there exists an open path within $R([0,1]^2)$ that intersects both the left-hand side (i.e., $R(\{0\}\times[0,1])$) and right-hand side (i.e., $R(\{1\}\times[0,1])$) of the quad. Similarly, each collection of $\CLE_6$-loops $\Gamma$ also induces a collection of crossed quads. It was shown by Garban, Pete and Schramm in \cite[Section~2.3]{gps} that the collection of percolation-quads also converges to the collection of $\CLE_6$-quads in an appropriate topology on the space of collections of quads. Furthermore, the quad-limit is almost surely measurable with respect to the $\CLE_6$ loops $\Gamma$. Holden and Sun showed in \cite[Theorem~6.10]{holden2023convergence} that the converse is also true. 

For a collection $\cS$ of quads and a simply connected open subdomain $U$, we can simply define the induced collection $\cS_{U}$\phantomsection{}\label{def:Gamma_Q} as all the quads in $\cS$ that lie inside $U$. This defines a $\CLE_6$ in $U$ in the quad-crossing sense. By the measurability results stated above, we can then identify this with a collection of $\CLE_6$ loops $\Gamma_U$. 
For general subdomains~$U$, 
one can still obtain a collection of loops by applying the measurability result to each simply connected subdomain of $U$. The following result proved in~\cite[Theorem~1.19]{schramm2011scaling} shows that $\CLE_6$ is measurable with respect to any nice subdomain and its complement. For notational simplicity, we write $\Gamma_{U^c}=\Gamma_{\overline{U}^c}$. In \cite[Corollary 1.20]{schramm2011scaling}, percolation is also characterized as a \emph{black noise} as defined by Tsirelson in \cite{tsirelson2004nonclassical, tsirelson2004scaling, tsirelson2014noise}. This implies that percolation in disjoint domains is independent.

\begin{lem}\label{lem:measurability-cle}
Let $U\in[0,1]^2$ be an open and simply connected domain such that its boundary $\partial U$ is a rectifiable loop with finitely many double points. Then the $\CLE_6$ $\Gamma$ is almost surely measurable with respect to the $\sigma$-algebra generated by $\Gamma_U$ and $\Gamma_{U^c}$. Furthermore, $\Gamma_U$ and $\Gamma_{U^c}$ are independent.
\end{lem}
The following convergence result is shown in~\cite[Section~2.3]{gps}.
\begin{thm}\label{thm:convergence-cle-percolation-full}
There exists a coupling of percolation on $\Lambda_n$ and the collection of $\CLE_6$ loops $\Gamma$ such that almost surely for all $Q\in\cQ$, the sequences $(\Gamma_{Q,n})_{n\in\N}$ and $(\Gamma_{Q^c,n})_{n\in\N}$ converge to $\Gamma_Q$ and $\Gamma_{Q^c}$ with respect to the quad-crossing topology introduced in \cite{schramm2011scaling}.
\end{thm}
In Proposition~\ref{prop:convergence-ghf-path}, we show that we can further couple percolation on $\Lambda_n$ and $\Gamma$ such that almost surely for all $Q\in\cQ$, the sequence $(\Lambda_{Q,n},\dpt_{Q,n},\pi_n)_{n\in\N}$ converges to $(\Lambda_Q,\dpt_Q,\pi)$ in the $\GHf$ distance.

In this paper, we will not need the precise quad crossing machinery. Rather, we only need to know that we can define $\CLE_6$ on subdomains in such a way that Lemma~\ref{lem:measurability-cle} and Theorem~\ref{thm:convergence-cle-percolation-full} hold. Therefore, we do not give the precise definitions regarding the space of quads and its topology. Instead, we refer to the cited works above for more background.

We finish this part with the following lemma, which will be used throughout the paper. The lemma states that loops that touch in the $\CLE_6$ will also eventually touch in the percolation configuration under a suitable coupling.

\begin{lem}\label{lem:touching-loops}
There exists a coupling such that almost surely $(\Lambda_n,\dpt_n,\pi_n) \to (\Lambda,\dpt,\pi)$ in GHf and the following holds. Suppose that $(\Lambda_n,\dpt_n,\pi_n)$ and $(\Lambda,\dpt,\pi)$ are isometrically embedded into a common metric space $(W,d_W)$ such that the sequence converges under the $\Delta_{\Haus}+\Delta_f$ metric in~\eqref{eq:ghf}. Let $\gamma,\gamma'$ be two different $\CLE_6$ loops that touch at a prime end $z\in\Lambda$ that is outside the loops $\gamma,\gamma'$. Let $\gamma_n,\gamma'_n$ be percolation interface loops converging to $\gamma$ and $\gamma'$ respectively. Then there exists $n_0$ such that for all $n\geq n_0$, the loops $\gamma_n$ and $\gamma'_n$ touch (i.e., come within distance $1$ of each other) at a point $z_n$ such that $z_n\rightarrow z$ as $n\rightarrow\infty$ in $d_W$ (and hence also $\pi_n(z_n) \to \pi(z)$ in $\R^2$). The analogous statement holds for a loop $\gamma$ that touches itself at a point $z$, and for a loop $\gamma$ that touches the boundary $\partial [0,1]^2$.
\end{lem}
\begin{proof}
Let $\cA_6(x,\varepsilon,\rho)$ denote a polychromatic $6$-arm event in the annulus $A(x;\varepsilon,\rho)$. By \cite[Lemma~2.9]{gps}, we have that $(\Gamma_n,\1_{\cA_6(x,\varepsilon,\rho)}(\Gamma_n))$ converges in law to $(\Gamma,\1_{\cA_6(x,\varepsilon,\rho)}(\Gamma))$. Combining this with Proposition~\ref{prop:convergence-ghf-path}, we can couple $(\Gamma_n,\Gamma)$ so that $((\Lambda_n,\dpt_n,\pi_n),(\1_{\cA_6(x_i,\varepsilon_i,\rho_i)}(\Gamma_n))_{i\in\N})$ converges almost surely to $((\Lambda,\dpt,\pi),(\1_{\cA_6(x_i,\varepsilon_i,\rho_i)}(\Gamma))_{i\in\N})$ for a countable dense set of annuli.

Suppose there exist two $\CLE_6$ loops $\gamma,\gamma'$ that touch at some point $z$, but such that their approximants $\gamma_n,\gamma'_n$ do not touch for all $n$. Fix any $0 < \rho < \min\{\diam(\gamma),\diam(\gamma')\}$. Let $\varepsilon>0$. Then there exists $n_0$ such that for all $n\geq n_0$, $\gamma_n$ and $\gamma'_n$ come within distance $\varepsilon$ of $z$. Then the box $B(z,\varepsilon)$ has $6$ arms to distance $\rho/2$. Since the $6$-arm exponent is larger than $2$, the probability that there is such a box goes to $0$ as we let $\varepsilon \searrow 0$ and $n \to \infty$. In the coupling described above, almost surely, for $\varepsilon$ sufficiently small and $n$ sufficiently large, there is no such box in the percolation configuration. This is a contradiction, so $\gamma_n$ and $\gamma'_n$ must touch within distance $\rho$ of $z$ for all $n$ sufficiently large. By considering smaller and smaller constants $\rho > 0$, we can find touching points $z_n$ that also satisfy $\|z-z_n\|_2\rightarrow0$.

It remains to show that $z_n\rightarrow z$ in $d_W$. Note that since $z$ is a touching point of two loops, there also exist $z'\in\Lambda$ inside the loop $\gamma'$ and $z''\in\Lambda$ inside the loop $\gamma$ with $\pi(z)=\pi(z')=\pi(z'')$, i.e., $z$ is one of three \emph{prime ends}. By definition of Gromov-Hausdorff convergence, each subsequence $I$ has a further subsequence $J$ such that $(z_n)_{n\in J}$ converges in $d_W$. This limit must be one of $z,z',z''$. Since none of the $z_n$ are inside $\gamma_n$ and $\gamma_n''$, we must have that the limit is $z$. We conclude that $z_n\rightarrow z$ in~$d_W$.

For a loop $\gamma$ touching itself or the boundary, the proof is similar. In the boundary case, we use instead that the half-plane $3$-arm exponent is larger than $1$.
\end{proof}

\subsubsection{Convergence of the cluster measure}
\label{subsubsec:measure-conv}
In order to prove the GHP convergence from Theorem~\ref{thm:main-convergence}, we first prove GH convergence of the metric spaces and then use prior results on convergence of the one-arm measure on percolation clusters to strengthen this to GHP convergence. In this section, we state the required cluster measure results proved in \cite{gps} and \cite{my2026upcoming}.

Let us start with the annulus counting measure. Let $A$ be a topological annulus in $[0,1]^2$ with piecewise smooth boundary. Let $A^{\inn}$ (resp.\ $A^{\out}$) denote its inner (resp.\ outer) face, i.e.\ the bounded (resp.\ unbounded) connected component of $\R^2 \setminus A$. Define
\begin{equation}\label{eq:annulus_measure}
    \mu_n^A:= \frac{1}{n^2 A_1(1,n)}\sum_{x\in\Lambda_n\cap A^{\inn}}\1_{\{x\leftrightarrow A^{\out}\}}\delta_x ,
\end{equation}
i.e.\ the normalized counting measure of the set of points in $\Lambda_n \cap A^{\inn}$ that are connected to the outer boundary of $A$.

\begin{thm}[{\cite[Theorem 5.1]{gps}}]\label{th:gps}
There exists a coupling of $\Gamma_n$ and $\Gamma$ such that the sequence $(\mu^A_n)_{n\in\N}$ converges almost surely in Prokhorov distance to some measure $\mu^A$\phantomsection{}\label{def:mu^A}. The measure $\mu^A$ is a function of the collection $\Gamma$ of $\CLE_6$ loops. 
\end{thm}
Let $B\subset[0,1]^2$ be a domain with rectifiable boundary and let $Y^A_k(B)$\phantomsection{}\label{def:YA_k} and $Y^A_{k,n}(B)$\phantomsection{}\label{def:YA_k,n} be the number of boxes $Q$ in $\cB_k$ that are contained in $B\cap A^{\inn}$ and such that $2Q\leftrightarrow A^{\out}$ in respectively the $\CLE_6$ and the percolation configuration. 
The following description of $\mu^A$ is essentially proved in \cite[Section 4.6]{gps}. Note that this result was stated for the pivotal measure rather than the cluster measure, but it was remarked in Section 5.1 of that paper that the analogous results also hold for the cluster measure.
\begin{prop}\label{prop:char-gps}
There exists a constant $c>0$ such that for each (deterministic) domain $B\subset[0,1]^2$ with rectifiable boundary,
\begin{equation}\label{eq:muA(B)}
    \mu^A(B)=\lim_{k\rightarrow\infty}\frac{cY^A_k(B)}{2^{2k} A^{\CLE}_1(2^{-k},1)}\qquad\text{in $L^2$},
\end{equation}
where $A^{\CLE}_1(2^{-k},1)$ is the $\CLE_6$ one-arm probability as defined in Section~\ref{subsec:arm-exponents}.
\end{prop}
\begin{proof}
It is shown in \cite[Section 4.6]{gps} that there is a coupling under which the following hold. Firstly, there exists $\zeta(n,k)$ that tends to 0 as $k\rightarrow\infty$ and $\frac{n}{2^k}\rightarrow\infty$ and such that
\begin{equation}\label{eq:gps_discrete_measapprox}
   \left\|\mu^A_n(B)-\frac{cY^A_{k,n}(B)}{2^{2k}A_1(n2^{-k},n)}\right\|_2\leq\zeta(n,k).
\end{equation}
Furthermore, for all $k$,
\begin{equation}
    \frac{cY^A_{k,n}(B)}{2^{2k}A_1(n2^{-k},n)} \to \frac{cY^A_k(B)}{2^{2k}\AC_1(2^{-k},1)} \quad\text{in probability,}\qquad n\rightarrow\infty.
\end{equation}
Lastly, $(\frac{cY^A_k(B)}{2^{2k}\AC_1(2^{-k},1)})_{k\in\N}$ is a Cauchy sequence in $L^2$, and so 
\begin{equation}\label{eq:gps_lim_measapprox}
    \widetilde{\mu}^A(B):=\lim_{k\rightarrow\infty}\frac{cY^A_k(B)}{2^{2k}\AC_1(2^{-k},1)}\qquad\text{in $L^2$}
\end{equation}
exists. Therefore we just need to argue that $\mu^A(B) = \widetilde{\mu}^A(B)$ almost surely. Combining the above, taking $n\rightarrow\infty$ in~\eqref{eq:gps_discrete_measapprox} and then $k\rightarrow\infty$ in~\eqref{eq:gps_lim_measapprox}, it follows that for all $B$,
\begin{equation}
    \widetilde{\mu}^A(B)=\lim_{n\rightarrow\infty}\mu_n^A(B)\qquad\text{in probability}.
\end{equation}
Now recall that $\mu^A$ is the a.s.\ limit of $(\mu^{A}_n)_{n\in\N}$ with respect to the L\'evy-Prokhorov distance. It is shown in \cite[Corollary~4.16]{gps} that under a suitable coupling, for each $B$ with rectifiable boundary, we have the $L^2$-convergence $\mu_n^A(B) \to \mu^A(B)$ and that $\mu^A(\partial B) = 0$ almost surely. This implies that $\mu^A(B) = \widetilde{\mu}^A(B)$ almost surely and completes the proof.
\end{proof}

\begin{lem}\label{lem:measure_2arm}
    Fix an annulus $A$ and $r>0$. For each $\delta > 0$, let $X_\delta$ be the set of points $x \in A^{\inn}$ such that there are two disjoint (open or closed) arms crossing $B(x,r)\setminus B(x,\delta)$. Then
    \[
    \lim_{\delta\searrow 0} \E[\mu^A(X_\delta)] = 0 .
    \]
    In particular, the set of points with two disjoint arms to distance $r$ has measure $0$ almost surely.
\end{lem}
\begin{proof}
    Let $\widetilde{X}_\delta$ be the union of boxes $B \in \cB_\delta$ such that there are two disjoint arms from $2B$ to distance $r$. Then $X_\delta \subseteq \widetilde{X}_\delta$. By applying Proposition~\ref{prop:char-gps} to each of the finite number of boxes in $\cB_\delta$ we see for each fixed $\delta$ that
    \[
    \E[\mu^A(\widetilde{X}_\delta)] = \lim_{k\rightarrow\infty}\frac{c\E[Y^A_{k}(\widetilde{X}_\delta)]}{2^{2k} \AC_1(2^{-k},1)}
    \]
    By the BK inequality, we have $\AC_2(\delta,r) \le \AC_1(\delta,r)^2$, and by Lemma~\ref{lem:quasi-multiplicativity-arms} we have $\AC_1(2^{-k},1) \asymp \AC_1(2^{-k},\delta)\AC_1(\delta,r)$. Combining everything gives
    \[\begin{split}
        \E[Y^A_{k}(\widetilde{X}_\delta)]
        &= \sum_{B \in \cB_\delta} \sum_{Q \in \cB_k, Q \subseteq B} \P[B \subseteq \widetilde{X}_\delta ,\, Q \leftrightarrow A^{\out}] \\
        &\asymp 2^{2k} \AC_1(2^{-k},\delta)\AC_2(\delta,r) \\
        &\le 2^{2k} \AC_1(2^{-k},\delta)\AC_1(\delta,r)^2 \\
        &\lesssim 2^{2k} \AC_1(2^{-k},r)\AC_1(\delta,r) .
    \end{split}\]
    Therefore we have
    \[
    \E[\mu^A(X_\delta)] \le \E[\mu^A(\widetilde{X}_\delta)] \lesssim \AC_1(\delta,r),
    \]
    which tends to $0$ as $\delta\rightarrow0$.
\end{proof}

Let us turn to the $\CLE_6$ gasket measure. The conformally covariant volume measure on the $\CLE_6$ gasket was first constructed in \cite{miller2024existence}. It is uniquely characterized up to a deterministic constant which we fix arbitrarily. For each cluster~$\cX^{(m)}$ as defined in Section~\ref{subsec:statements}, we let $\mu^{(m)}$\phantomsection{}\label{def:mum-second} denote the gasket measure on $\cX^{(m)}$. The following characterization is shown in \cite{my2026upcoming}. Let $B\subseteq [0,1]^2$ be a box or ball and let $Y^{(m)}_k(B)$\phantomsection{}\label{def:Ymk} be the number of boxes $Q \in \cB_k$ such that $Q\subset B$ and such that $2Q$ has non-empty intersection with $\pi(\cX^{(m)})$.
\begin{prop}\label{prop:char-measure-my}
Let $\mu^{(m)}$ be the $\CLE_6$ gasket measure on $\cX^{(m)}$. Then there exists a constant $c>0$ such that for each box or ball $B\subseteq [0,1]^2$,
\begin{equation}\label{eq:expression-mu}
    \mu^{(m)}(B)=\lim_{k\rightarrow\infty}\frac{cY^{(m)}_k(B)}{2^{(91/48)k}}\qquad\text{almost surely}.
\end{equation}
\end{prop}
Since the exact renormalization constant is shown in \cite{my2025diffusion}, we obtain the following asymptotic expression for the CLE$_6$ one-arm probability by a comparison.
\begin{cor}\label{co:cle_one_arm}
There is a constant $c>0$ such that
\[
\AC(2^{-k},1) = c2^{-(5/48)k}+o(2^{-(5/48)k}) \quad\text{as } k\to\infty .
\]
\end{cor}
The proof of this statement is given in Section~\ref{sec:ghpf-proof}.

\subsection{CLE metrics}\label{subsec:cle-metrics}
The limiting measures $\dpt$ and $\dg$ are examples of a class of (pseudo\nobreakdash-)metrics introduced in \cite{amy2025tightness} called \emph{$\CLE$ metrics}, which we define in this section. Note that the definitions in this paper are in fact slightly stronger than those in \cite{amy2025tightness}. This is because \cite{amy2025tightness} work in a very general framework, whereas we are only interested in continuous geodesic and resistance $\CLE_6$ metrics.

We first generalize the metric $\dpt$ from~\eqref{eq:dpt-first} to encode also the topologies of regions whose interior consist of collections of domains. Let $U\subset[0,1]^2$ and let $\cG$ be any collection of loops in $\R^2$. Define $P(x,y;U;\cG)$ to be the set of paths from $x\in U$ to $y\in U$ that stay inside $U$ that do not cross any loop in $\cG$. Define 
\begin{equation}\label{eq:def-dptU}
    \dpt_{U,\cG}(x,y)=\inf\{\diam(\gamma): \gamma\in P(x,y;U;\cG)\},
\end{equation}
where $\diam$ denotes diameter with respect to the Euclidean metric. If $\cG$ is the entire collection of $\CLE_6$ loops $\Gamma$, then we simply write $\dpt_{U,\Gamma}=\dpt_{U}$.

Recall that the convergence results in Section~\ref{subsec:statements} were stated simultaneously for all clusters. In order to identify the limit with the $\CLE_6$ metric, it suffices to consider a specific cluster which we denote by $\cX$. The choice of $\cX$ is arbitrary, but for the sake of consistency with \cite{amy2025tightness, my2025diffusion, my2025metric}, we choose $\cX$\phantomsection{}\label{def:cX-second} to be the outermost open cluster in $\Gamma$ surrounding the center $(\frac{1}{2},\frac{1}{2})$. By the local absolute continuity between the different interior clusters, the $\CLE_6$ metric on each cluster $\cX^{(m)}$ can be recovered from the $\CLE_6$ metric on the cluster $\cX$, see \cite{my2025metric}. We define $\cX_n$\phantomsection{}\label{def:cXn-second} similarly for the percolation configuration. Therefore it suffices to prove the convergence of the geodesic and resistance metrics on $\cX_n$ to those on $\cX$.

To be precise, let $\cL$\phantomsection{}\label{def:cL} be the outermost loop in $\Gamma$ surrounding the point~$(\frac{1}{2},\frac{1}{2})$ (in the sense that its winding number about~$(1/2,1/2)$ is $\pm 1$) and let $\cX$ be the gasket of $\Gamma$ whose exterior boundary is $\cL$ as defined in Section~\ref{subsec:statements}. Write $\Gamma_{\cL}$\phantomsection{}\label{def:Gamma_cL} to be the collection of loops in $\Gamma$ that are in the next nesting level inside $\cL$. Let $D$ be the regions surrounded by $\cL$ (i.e., the union of connected components of $(0,1)^2\setminus\cL$ with winding number 1). For each open, simply connected $U\subset[0,1]^2$\phantomsection{}\label{def:U}, let $\Gamma_{U^*}$\phantomsection{}\label{def:GammaU*} be the collection of loops in $\Gamma$ that are inside $\cL$ and entirely contained in $U$, and let $U^*$\phantomsection{}\label{def:U*} be the set of points in $[0,1]^2$ that are inside $\cL$ but not on or inside any other loop of $\Gamma_{\cL}\setminus\Gamma_{U^*}$. We view $U^*$ as the metric space equipped with the metric $\dpt_{\overline{U^*},\Gamma\setminus\Gamma_{U^*}}$. Note that $\Gamma_{U^*}$ is indeed a $\CLE_6$ inside $U^*$ as defined in Section~\ref{subsubsec:convergence-boundary-loops}.

Suppose $\cG=\{\gamma_1,\ldots,\gamma_m\}\subset\Gamma_{\cL}$ and $V\subset D$ is a union of Euclidean connected components of $D\setminus\cG$ that are not inside the loops $\gamma_1,\ldots,\gamma_m$. We let $\fC$\phantomsection{}\label{def:fC} be the collection of regions $V\subset D$ as described above (for some choice of $\cG$) such that $\overline{V}$ is simply connected with respect to $\dpt_{\overline{V},\{\cL\}\cup\cG}$ (where $\overline{V}$ denotes the completion with respect to $\dpt_{\overline{V},\{\cL\}\cup\cG}$). For $U\subset[0,1]^2$ open, we let $\fC_U=\{V\in\fC: \overline{V}\subset U\}$.

Let $(\cX\cap\overline{V},\fd_V)_{V\in\fC}$\phantomsection{}\label{def:fd-first} be a collection of resistance or geodesic (pseudo-)metric spaces coupled with~$\Gamma$. Here the notation $\cX\cap U$ is shorthand for the set of sites $x\in\cX\subset\Lambda$ such that $\pi(x)\in U$. We say it is a $\CLE_6$ metric if the following hold almost surely.

\textbf{Continuity:} For each $V\in\fC$, we have that $\fd_V$ is continuous with respect to $\dpt_{\overline V}$.

\textbf{Markovian property:} Let $U\subset[0,1]^2$ be open and simply connected. The conditional law of the collection $(\cX\cap \overline{V},\fd_{V})_{V\in\fC_U}$ given $\{\cL\}\cup\Gamma_{\cL}\setminus\Gamma_{U^*}$ and $(\cX\cap \overline{V'},\fd_{V'})_{V'\in\fC_{[0,1]^2\setminus\overline{U}}}$ is almost surely measurable with respect to $U^*$.

\textbf{Translation invariance:} 
Let $D_1$ be the support of $(\cX\cap \overline{V},\fd_V)_{V\in\fC_U}$ and $D_2$ be the support of $U^*$. Let $\mu_U(\cdot\mid\cdot)$ be a probability kernel from $D_2$ to $D_1$ such that $\mu_U(\cdot\mid U^*)$ is a.s.\ equal to the conditional law of $(\cX\cap \overline{V},\fd_V)_{V\in\fC_U}$ given $U^*$. Then for each $z \in \C$ such that $U+z\subset[0,1]^2$ one can choose versions of $\mu_U$ and $\mu_{U+z}$ such that for all events $A$ in $D_1$ and $u \in D_2$ we have
\begin{equation}
    \mu_U(A\mid u)=\mu_{U+z}(A+z\mid u+z).
\end{equation}
where $A+z$ denotes the event that the translated metrics $(\cX\cap \overline{V}-z,\fd_{V}(\cdot+z,\cdot+z))_{V\in\fC_{U+z}}$ are in the event $A$.

\textbf{Compatibility:} Let $V,V'\in\fC$, $V\subset V'$, and $x,y\in\cX\cap\overline{V}$ such that for every $u\in V'\setminus V$ there is a point $z\in\overline{V}$ that separates $u$ from $x,y$ in $\cX\cap\overline{V'}$. Then
\begin{equation}
    \fd_{V'}(x,y)=\fd_V(x,y).
\end{equation}

\textbf{Monotonicity:} Let $V,V'\in\fC$, $V\subset V'$. Then, for each $x,y \in \cX \cap \overline{V}$,
\begin{equation}
    \fd_{V'}(x,y) \le \fd_V(x,y) .
\end{equation}

\textbf{Series law:} Let $V\in\fC$ and let $x,y,z\in V$ such that $z$ separates $x$ from $y$ in $\cX\cap\overline{V}$. Then
\begin{equation}
    \fd_V(x,y)=\fd_V(x,z)+\fd_V(z,y).
\end{equation}

\textbf{Generalized parallel law:} Let $V\in\fC$ and let $x,y,z_1,\ldots,z_M\in V$ such that $x$ and $y$ are separated in $\cX\cap\overline{V}\setminus\{z_1,\ldots,z_M\}$. Let $K_x$ be the connected component of $\cX\cap\overline{V}\setminus\{z_1,\ldots,z_M\}$ containing~$x$, and let $V_x \in \fC$ be the smallest region such that $\overline{V_x} \supseteq K_x$. Then 
\begin{equation}
    \fd_V(x,y)\geq\frac{1}{M}\min_i\fd_{V_x}(x,z_i).
\end{equation}

The following is shown in \cite[Proposition 6.14]{amy2025tightness}.
\begin{prop}\label{prop:positivity}
Let $(\fd_V)_{V\in\fC}$ be a $\CLE$ (pseudo-)metric. Then either $\fd_V=0$ for each $V\in\fC$ or it is in fact a true metric, i.e., $\fd_V(x,y)>0$ for each $V\in\fC$ and $x,y\in\cX\cap\overline{V}$ with $x\neq y$.
\end{prop}

We now turn our attention to specifically geodesic and resistance $\CLE$ metrics. We say a collection $(\fd_V)_{V\in\fC}$ is a \emph{geodesic $\CLE$ metric}\footnote{In \cite{my2025metric}, this definition is called a \emph{weak} geodesic CLE metric but it is shown that each weak geodesic CLE metric is a (strong) geodesic CLE metric, see Theorem~\ref{th:uniqueness-geodesic-cle}.} if it is a CLE$_6$ metric in the sense defined above and additionally the following holds. For every $V\in\fC$ and $x,y\in\cX\cap\overline{V}$, we have
\begin{equation}
    \fd_V(x,y)=\inf_{\gamma\in P(x,y;\overline{V})}L_{\fd}(\gamma),
\end{equation}
where 
\begin{equation}
    L_{\fd}(\gamma)=\sup_{0\leq t_0\leq\ldots\leq t_M\leq 1}\sum_{i=1}^M \fd_D(\gamma(t_{i-1}),\gamma(t_{i}))
\end{equation}
is the \emph{length} of $\gamma$ with respect to $\fd_D$.

The following is shown in \cite{my2025metric}. Recall the conformally covariant geodesic metric $\dg$ from \cite{my2025metric} introduced in Section~\ref{subsec:statements}.
\begin{thm}\label{th:uniqueness-geodesic-cle}
    Suppose that we have the setup described above and $(\fd_V)_{V\in\fC}$ is a geodesic CLE$_6$ metric. Then there exists a deterministic constant $c$ such that $\fd_D = c\,\dg|_{\cX}$ almost surely.
\end{thm}

We now define the \emph{$\CLE_6$ resistance metric}. First, we recall the definition of a resistance metric \cite{kigami2001analysis}. A metric $R$ on a non-empty set $X$ is called a \emph{resistance metric} if for each finite subset $A \subseteq X$ there is a (unique) symmetric weight function $w\colon A \times A \to [0,\infty)$ such that the restriction of $R$ to $A$ agrees with the effective resistance metric associated with the weighted graph $G = (A,w)$.

We say $(\fd_V)_{V\in\fC}$ is a \emph{$\CLE_6$ resistance metric}\footnote{In \cite{my2025diffusion}, this definition is called the resistance metric of a \emph{weak} CLE resistance form but it is shown that each weak CLE resistance form is a (strong) CLE resistance form, see Theorem~\ref{th:uniqueness-resistance-cle}.} if it is a CLE$_6$ metric in the sense defined above and additionally the following properties are satisfied.
\begin{enumerate}[(a)]
    \item For each $V \in \fC$, we have that $(\cX\cap\overline{V},\fd_V)$ is a resistance metric space.
    \item Suppose that $V,V_1,V_2\in\fC$ are such that $\cX\cap\overline{V} = (\cX\cap\overline{V_1}) \cup (\cX\cap\overline{V_2})$ and $(\cX\cap\overline{V_1}) \cap (\cX\cap\overline{V_2}) = \{a,b\}$. There exists a countable set $\cZ_V=\{z_1,z_2,\ldots\} \subseteq \cX\cap\overline{V}$ that contains $a,b$, is dense with respect to $\dpt_{\overline{V}}$, and such that the following holds for sufficiently large $m$. Let $A = \{z_1,\ldots,z_m\}$, $A_1 = A \cap \overline{V_1}$, $A_2 = A \cap \overline{V_2}$, and let $w$ (resp.\ $w_1$, $w_2$) be the weight function associated with $\fd_V(\cdot,\cdot)|_{A}$ (resp.\ $\fd_{V_1}(\cdot,\cdot)|_{A_1}$, $\fd_{V_2}(\cdot,\cdot)|_{A_2}$). Then
\begin{itemize}
    \item $w(x,y)=w_1(x,y)$ for each $x,y\in A_1$;
    \item $w(x,y)=w_2(x,y)$ for each $x,y\in A_2$;
    \item $w(x,y)=0$ for each $x\in A_1 \setminus \{a,b\}$, $y\in A_2 \setminus \{a,b\}$;
    \item $w(a,b)=0$.
\end{itemize}
\end{enumerate}

The following is shown in \cite{my2025diffusion}. Recall the conformally covariant resistance metric $\dr$ from \cite{my2025diffusion} introduced in Section~\ref{subsec:statements}.
\begin{thm}\label{th:uniqueness-resistance-cle}
    Suppose that we have the setup described above and $(\fd_V)_{V\in\fC}$ is a CLE$_6$ resistance metric. Then there exists a deterministic constant $c$ such that $\fd_D = c\,\dr|_{\cX}$ almost surely.
\end{thm}

\subsection{Resistance metrics and stochastic processes}\label{subsec:resistance-metrics-walks}
An important tool in the proof of the scaling limit of random walk on $\cX_n$ is the connection between resistance metrics and stochastic processes. We refer to \cite{croydon2018scaling} for a definition of resistance metrics on general spaces. See \cite[Chapter 9]{levin2017markov} for a more accessible perspective from the viewpoint of random walks. The following result was proved by Croydon in \cite{croydon2018scaling}. Let $(K,d,\mu,f,\rho)$ be a rooted metric-measure-function space. We assume that $(K,d)$ is non-empty and compact, $d$ is a resistance metric,\footnote{Since $(K,d)$ is compact, $d$ is a \emph{regular} resistance metric, see \cite[Corollary~6.4]{kigami2012resistance}.} $\mu$ is a finite Borel measure with full support on $(K,d)$, and $\rho \in K$. See \cite{croydon2018scaling} for the definition of a regular resistance metric. Then there exists a unique \emph{Hunt process} $X$ associated with $(K,d,\mu,\rho)$ with $X_0=\rho$. A Hunt process is a strong Markov process that is c\`adl\`ag and \emph{quasi-left continuous}, i.e., if $(\sigma_n)_{n\in\N}$ is an increasing sequence of stopping times with limit $\sigma$, then 
\begin{equation}
    \P(\lim_{n\rightarrow\infty}X_{\sigma_n}=X_\sigma\mid \sigma<\infty)=1.
\end{equation}
If $K$ is the vertex set of a finite graph and $d$ the associated resistance metric, then the associated Hunt process is the continuous-time random walk with jump rates such that $\mu$ is its stationary measure.
\begin{thm}[{\cite[Theorem 7.1]{croydon2018scaling}}]\label{th:resistance-random-walk-convergence}
Let $(K_n,d_n,\mu_n,f_n,\rho_n)_{n\in\N}$ be a sequence of rooted metric-measure-function spaces where $f_n\colon(K_n,d_n)\to(Y,d_Y)$ and $f\colon (K,d)\to(Y,d_Y)$ are continuous and $(Y,d_Y)$ is some compact metric space. Assume that
\[
\Delta_{\rGHPf}((K_n,d_n,\mu_n,f_n,\rho_n), (K,d,\mu,f,\rho)) \to 0 .
\]
Let $X_n$ and $X$ be the Hunt processes associated with $(K_n,d_n,\mu_n,\rho_n)$ and $(K,d,\mu,\rho)$, respectively. Then $f_n(X_n)$ converge weakly to $f(X)$ as $n\rightarrow\infty$ in the space of cadlag processes on $Y$ equipped with the Skorokhod $J_1$-topology.
\end{thm}
Note that Croydon does not use the term $\GHPf$ convergence, but rather calls it \emph{spatial $\GHP$ convergence.} The random walk convergence from Theorem~\ref{thm:main-random-walk} is thus a direct consequence of Theorem~\ref{thm:main-convergence} and Theorem~\ref{th:resistance-random-walk-convergence}. (Note that since the limiting process $X$ is continuous, convergence in the Skorokhod $J_1$-topology is equivalent to the uniform convergence on compact time intervals.)

\section{Outline of the proof}\label{sec:sketch}
The proof of Theorem~\ref{thm:main-convergence} has three main components. First, we prove that the sequence of path metric spaces $(\cX^{(m)}_n,\dpt_n,\pi_n)_{n\in\N}$ converges weakly with respect to the GHf topology. We give a sketch of this part in Section~\ref{sec:convergence-path-sketch}. The majority of this paper is devoted to proving tightness of the normalized geodesic and resistance metrics with respect to the GHf topology, which is the second part. A sketch of this step is given in Section~\ref{sec:tightness-sketch}. In that section, we also define the normalizing constants $(\median_n)_{n\in\N}$. For the third part, we show that every subsequential weak limit is a $\CLE_6$ metric by checking it has the properties stated in Section~\ref{subsec:cle-metrics}. The uniqueness results from Theorems~\ref{th:uniqueness-geodesic-cle} and~\ref{th:uniqueness-resistance-cle} then imply that all subsequential limits are in fact the same, which completes the proof of GHf convergence of the geodesic and resistance metrics. See Section~\ref{sec:uniqueness-sketch} for an overview. To complete the proof of Theorem~\ref{thm:main-convergence}, it only remains to strengthen the GHf convergence to GHPf convergence, using the measure convergence from Section~\ref{subsubsec:measure-conv}. Convergence of the random walk from Theorem~\ref{thm:main-random-walk} then follows from the GHPf convergence of the resistance metric combined with Theorem~\ref{th:resistance-random-walk-convergence}. These last two steps are sketched in Sections~\ref{sec:sketch-measure} and~\ref{sec:sketch-random-walk}.

We give a unified proof of Theorem~\ref{thm:main-convergence} that works for both the geodesic and resistance metric simultaneously. The proof relies only on the following properties of the metrics. Let $D_n$\phantomsection{}\label{def:D_n} denote either $D_n^{\geo}$ or $D_n^{\res}$ and let $d_n$\phantomsection{}\label{def:d_n} denote the rescaled metrics $\dg_n$ or $\dr_n$ respectively, as defined in Section~\ref{sec:introduction}. Let $x,y\in\Lambda_n$ be sites that belong to the same open cluster. 
\begin{enumerate}[(A)]
\item\label{property:serial}  If there exists a site $a\in\Lambda_n$ such that every open path from $x$ to $y$ has to pass through $a$, then $D_n(x,y)=D_n(x,a)+D_n(a,y)$. 
\item\label{property:parallel} If there exists $M\in \N$ and a set of sites $a_1,\ldots, a_M$ such that every open path from $x$ to $y$ has to pass through at least one $a_i$, then we have that $D_n(x,y) \geq \frac{1}{M} \min_{1\leq i\leq M} D_n(x,a_i)$.
\item\label{property:monotonicity} If $G$ is a subgraph of $\Lambda_n$ containing both sites $x,y$, then $D_{G}(x,y)\geq D_n(x,y)$, where $D_{G}$ denotes the induced geodesic/resistance metric on $G$.
\end{enumerate}
The same statements hold for $d_n$. The three properties are clearly true for the geodesic metric. For the resistance metric, \eqref{property:serial} is the series law, \eqref{property:parallel} follows from the parallel law, and~\eqref{property:monotonicity} is Rayleigh's monotonicity principle. For the majority of the paper, we will be using the notation $D_n$ or $d_n$ and not refer to the geodesic or resistance metrics specifically. Sometimes, we will drop the $n$ and just write $D$ when the underlying graph is clear from the context.

\subsection{Convergence of the path metric}\label{sec:convergence-path-sketch}

Path metric convergence is the most straightforward and relies only on Camia and Newman's percolation cluster interface convergence as stated in Theorem~\ref{th:camia-newman}.

Let $Q\in\cQ$. Define 
\begin{equation}\phantomsection{}\label{def:dpt_{Q,n}}
    \dpt_{Q,n}:=\inf_{\gamma\in P_{Q,n}(x,y)}\diam(\gamma),
\end{equation}
where $P_{Q,n}(x,y)$ is the set of all admissible (i.e., containing only all open or all closed vertices) paths from $x$ to $y$ inside $\Lambda_{Q,n}$. Note that if $P_{Q,n}(x,y)$ is empty, then $\dpt_{Q,n}(x,y)=\infty$. For a collection of $\CLE_6$-loops $\Gamma_Q$ in $Q$, define $\dpt_Q$ analogously (recall~\eqref{eq:def-dptU}). We prove GHf convergence of the path metric (and the projection map). We note that although the definition of the path metric seems a bit artificial, it is useful as a reference metric in the proof of the tightness of the geodesic and the resistance metrics. In particular, given the Proposition~\ref{prop:convergence-ghf-path} below, it will suffice to bound $\dg_n$ and $\dr_n$ in terms of~$\dpt_n$.

\begin{prop}\label{prop:convergence-ghf-path}
There exists a coupling between $\Gamma_n$ and $\Gamma$ such that the following holds almost surely. For each $Q\in\cQ$, the sequence $(\Lambda_{Q,n},\dpt_{Q,n},\pi_n)_{n
\in\N}$ converges to $(\Lambda_Q,\dpt_Q,\pi)$ in the $\GHf$ distance.
\end{prop}

Recall from Section~\ref{subsec:notation} that convergence of $(\Lambda_{Q,n},\dpt_{Q,n},\pi_n)_{n\in\N}$ in $\GHf$ means that the individual clusters converge in $\GHf$.

In Section~\ref{subsubsec:convergence-boundary-loops}, we described two different encodings of percolation: by interface loops and by quad crossings. A different way to encode a percolation configuration is by considering the collection of open paths in $[0,1]^2$, equipped with the topology $\Delta_{\cur}$. It was shown by Aizenman and Burchard in \cite{aizenman1999holder} that the collection of open percolation paths is tight with respect to $\Delta_{\cur}$. Proposition~\ref{prop:convergence-ghf-path} then follows from uniqueness of the subsequential limits. This is widely believed to follow from convergence of the interface loops, but the authors are not aware of any existing proof in the literature. Therefore, we provide a short proof of Proposition~\ref{prop:convergence-ghf-path}. We first show that the sequence above is uniformly totally bounded with arbitrarily large probability, which implies that it is tight. This follows from the fact that each Euclidean ball in $[0,1]^2$ cannot have too many disjoint open arms coming out of it. Thus, every Euclidean ball has a bounded number of points in it that are connected but only by a path of large diameter. This implies that we can cover each open cluster by a bounded number of $\dpt_n$-balls. To prove convergence of the sequence, we then show that the subsequential limiting path distance between two points $x$ and $y$ equals the $\CLE_6$ path distance $\dpt(x,y)$. This is a consequence of the interface loop convergence proved by Camia and Newman as stated in Section~\ref{subsubsec:convergence-boundary-loops}, but requires some precision. The full proof is given in Section~\ref{sec:convergence-path-metric-proof}.

\subsection{Tightness of the geodesic and resistance metrics}\label{sec:tightness-sketch}
We now focus our attention on proving convergence of the geodesic and resistance metrics. We first show that they are tight in an appropriate topology.

We will define a sequence of positive real numbers $\median^{\geo}_n$ resp.\ $\median^{\geo}_n$. The precise definition is given in Sections~\ref{subsubsec:apriori-sketch}  and~\ref{subsubsec:holder-sketch} below; it is the $(1-p)$th quantile of an appropriate random variable defined in terms of the metric $D^{\geo}_n$ resp.\ $D^{\res}_n$ for a suitable $p \in (0,1)$.

Furthermore, for each $Q\in\cQ$, we let $D_{Q,n}$\phantomsection{}\label{def:DQn} denote either the geodesic or resistance metric on percolation on $\Lambda_{Q,n}$, where we set $D_{Q,n}(x,y)=\infty$ if $x$ and $y$ are not in the same cluster. We define
\begin{equation}\label{eq:def_dQn}
    d_{Q,n}:=\median_n^{-1}D_{Q,n}
\end{equation}
where $\median_n$\phantomsection{}\label{def:median_n} denotes either $\median^{\geo}_n$ or $\median^{\res}_n$.

We now specify the topology with respect to which we prove tightness. Let $Y=[0,1]^2\times [0,1]^2\times[0,\infty]$, let~$d_Y$ be the product metric on $Y$, and consider the functions $f_{Q,n}\colon\Lambda^2_{Q,n}\to Y$ given by $f_{Q,n}(x,y)=(\pi_n(x),\pi_n(y),d_{Q,n}(x,y))$. To ease the notation, will view $(\Lambda_{Q,n},\dpt_{Q,n},d_{Q,n},\pi_n)$ as the tuple $(\Lambda^2_{Q,n},(\dpt_{Q,n})^2,f_{Q,n})$ which lives in the space of compact metric spaces equipped with a continuous function. The majority of the paper will be devoted to proving tightness of the sequence $(\Lambda_{Q,n},\dpt_{Q,n},d_{Q,n},\pi_n)_{n\in\N}$ in the GHf topology. As in Section~\ref{subsec:notation}, we mean by this the convergence of $(\cX^{(m)}_{Q,n},\dpt_{Q,n},d_{Q,n},\pi_n)_{m \in \N}$ in the product of the GHf topologies as $n \to \infty$.
\begin{prop}\label{prop:GHf-tightness-dU}

For each $Q\in\cQ$, the sequence $(\Lambda_{Q,n},\dpt_{Q,n},d_{Q,n},\pi_n)_{n\in\N}$ is tight with respect to the GHf distance. Each subsequential limit is of the form $(\Lambda_Q,\dpt_Q,d_Q,\pi)$, where $d_Q$ is an almost surely $\eta$-H\"older continuous (pseudo-)metric with respect to $\dpt_Q$ for some $\eta > 0$. 
\end{prop}

\begin{rem}
    It is shown in \cite{my2025metric,my2025diffusion} that the metrics $\dg$ resp.\ $\dr$ are locally H\"older continuous for any exponent $\eta < \beta^\geo$ resp.\ $\eta < \beta^\res$ where the exponents satisfy $\beta^\geo \ge 1$, $\beta^\res \ge 3/4$. We will prove that it is in fact globally $\eta$-H\"older continuous (also near the boundary of $Q$) for any $\eta < 5/8$.
\end{rem}

From this, we obtain subsequential GHf convergence of $d_{Q,n}$.
\begin{prop}\label{prop:main-conv-dU-GH}
For each $Q\in\cQ$ and each subsequence of $(\Lambda_{Q,n},d_{Q,n},\pi_n)_{n\in\N}$, there exists a further subsequence that converges weakly with respect to the GHf distance.
\end{prop}

Proposition~\ref{prop:main-conv-dU-GH} follows almost immediately from Proposition~\ref{prop:GHf-tightness-dU}. In the proof, we show in particular that the GHf convergence of $(\Lambda_{Q,n},\dpt_{Q,n},d_{Q,n})$ implies GH convergence of $(\Lambda_{Q,n},d_{Q,n})$.

The proof of Proposition~\ref{prop:GHf-tightness-dU} has four main steps. To prove Proposition~\ref{prop:GHf-tightness-dU}, we want to show that $d_n$ is H\"older continuous with respect to $\dpt_n$ with a constant and exponent that is uniform in~$n$ (stated below in Proposition~\ref{prop:holderct}). Then, since $\dpt_n$ converges, we can conclude tightness of $d_n$.

The proof of the H\"older estimate follows a chaining argument. For each annulus $A(z,r,2r) \subseteq [0,1]^2$ we will prove an upper bound of the type
\begin{equation}\label{eq:chaining_bound}
    \p(d_n(v_i,v_j) > r^\eta) = O(r^\alpha)
\end{equation}
where $v_i,v_j$ are some pivotal points on open annulus crossings (see Section~\ref{subsubsec:annulus-sketch}). From this, we deduce a bound for $d_n(x,y)$ in terms of $\dpt_n(x,y)$ as follows. Let $\gamma$ be a path between $x$ and $y$ with smallest diameter. We cover the path $\gamma$ with annuli of sizes approximately proportional to their respective distances to $x$ and $y$ (with the largest one having size approximately comparable to $\dpt_n(x,y)$). The path $\gamma$ then forms a series of crossings of these annuli. By applying the bound~\eqref{eq:chaining_bound} and then `stringing together' the distances across all annuli we obtain a bound of the type $d_n(x,y) \lesssim \dpt_n(x,y)^\eta$ on an event with large probability. This requires taking a union bound of~\eqref{eq:chaining_bound} for all annuli covering the space $[0,1]^2$, and therefore we need to prove~\eqref{eq:chaining_bound} for a sufficiently large exponent $\alpha$. In fact, we will prove that it has superpolynomial tails.

The first step in the proof of~\eqref{eq:chaining_bound} is to prove an \emph{a priori estimate}. Instead of bounding the distance between connected annulus crossings, we upper bound the distance between two points in a toy version of the problem. For the a priori bound, we obtain polynomial concentration on this distance. See Section~\ref{subsubsec:apriori-sketch} for a sketch of this step. In that section, we also explain how the normalizing constants~$\median_n$ are defined. 
Informally, $\median_n$ should be thought of as the median distance between two points in the same cluster at Euclidean distance of order $n$. The next step consists of turning this a priori bound into a polynomial right tail bound for the distance between two crossings of an annulus. See Section~\ref{subsubsec:annulus-sketch} for a sketch. We then show that any polynomial bound with some exponent $\alpha$, say, can be bootstrapped into a bound with exponent $\alpha'>\alpha$. Thus, we obtain a bound that is stronger than any polynomial. We then also turn the a priori bound in the toy model into a bound on the distance between connected annulus crossings. This step is sketched in Section~\ref{subsubsec:subpolynomial-sketch}. The final step is given by the path covering procedure to obtain our H\"older bound, see Section~\ref{subsubsec:holder-sketch}. The corresponding proofs are contained in Sections~\ref{sec:apriori}--\ref{sec:holder}. The final proofs of Propositions~\ref{prop:GHf-tightness-dU} and~\ref{prop:main-conv-dU-GH} are given in Section~\ref{sec:conv-CLE}.

\subsubsection{A priori estimate}\label{subsubsec:apriori-sketch}
In this section, we define the auxiliary sequence $(\altq_n)_{n\in\N}$ that depends on a parameter $p\in(0,1)$. The normalizing constant $\median_n$ will be defined as $\altq_n(p)$ for a particular choice $p\in(0,1)$, which we will indicate as $\fp$. We also give a sketch of the first step of the proof of Proposition~\ref{prop:GHf-tightness-dU}, which we call the \emph{a priori estimate}.

\begin{figure}[ht]
    \centering
    \includegraphics[scale=0.8]{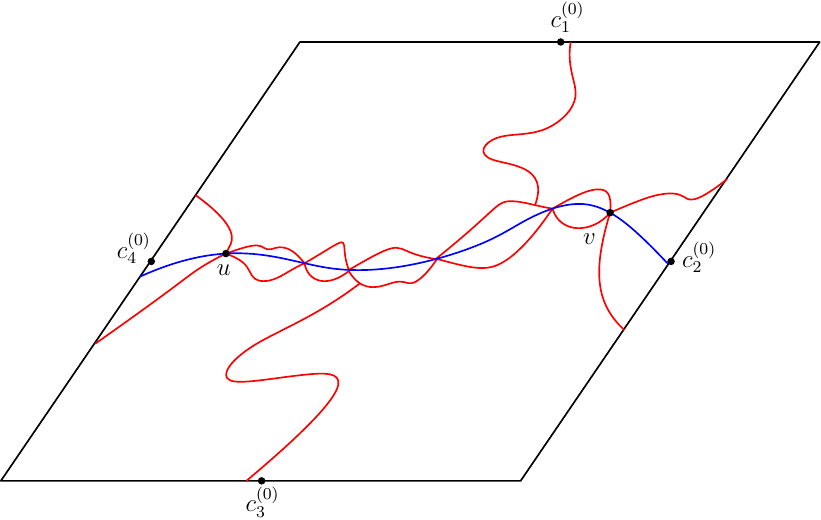}
    \caption{An illustration of the event $\cE$. The blue curve represents an open crossing from left to right. The red curves represent two closed left-to-right crossings that come within distance $2$ but do not intersect.}
    \label{fig:global_event_B(0)}
\end{figure}

We first define the auxiliary normalizing constants $\altq_n$ as follows. Let $c_1^{(0)},\ldots,c_4^{(0)}$ denote the midpoints of the sides of the box $[0,1]^2$, numbered clockwise starting from the top side. Let $\delta=1/10000$\phantomsection{}\label{def:delta} and define $\cE=\cE(\Lambda_n)$\phantomsection{}\label{def:cE} to be the event that the following events occur for percolation on $\Lambda_n$:
\begin{enumerate}[(i)]
\item\label{property:E1} there exist 3 alternating left-to-right crossings such that the top and bottom crossings are closed and the middle crossing is open;
\item\label{property:E2} the closed clusters of the top and bottom crossings come within (graph) distance 2 of each other;
\item\label{property:E3} the closed cluster of the top crossing hits the top side within distance $\delta $ of $c_1^{(0)}$ and does not come within distance $2\delta $ of $c_2^{(0)},c_3^{(0)},c_4^{(0)}$;
\item\label{property:E4} the closed cluster of the bottom crossing hits the bottom side within distance $\delta$ of $c_3^{(0)}$ and does not come within distance $2\delta$ of $c_1^{(0)},c_2^{(0)},c_4^{(0)}$;
\item\label{property:E5} the open cluster of the middle crossing hits the left and right sides within distance $\delta$ of $c_2^{(0)}$ and $c_4^{(0)}$ respectively and does not come within distance $2\delta$ of $c_1^{(0)}$ and $c_3^{(0)}$
\end{enumerate}
We refer to Figure~\ref{fig:global_event_B(0)} for an illustration of this event. Note that the probability of $\mathcal{E}(\Lambda_n)$ is uniformly bounded away from $0$ and $1$ due to the RSW theorem (see Lemma~\ref{lem:E(Lambda_n)_nontrivial}).

Let $u=u\!\left(\Lambda_n\right)$ and $v=v\!\left(\Lambda_n\right)$ be the left-most and right-most points that neighbor both the top and bottom closed clusters of $\Lambda_n$. Define $X_{\Lambda_n}=D(u,v)$\phantomsection{}\label{def:Xn} to be the distance between these points if they exist. Note that the items~\eqref{property:E3}--\eqref{property:E5} in the definition of the event $\cE(\Lambda_n)$ imply that the distance $D(u,v)$ does not depend on the percolation configuration within the $2\delta$-neighborhoods of $c_1^{(0)},c_2^{(0)},c_3^{(0)},c_4^{(0)}$. This will be useful later in Lemma~\ref{lem:perfectbounds}. 
More generally, for any box $B\in\cB^R_{k,n}$ (recall Section~\ref{subsec:notation}), we define the event $\cE(B)$ and random variable $X_B$ analogously.

We let $X_n$ be a random variable whose law is given by $X_{\Lambda_n}$ conditional on the event $\mathcal{E}(\Lambda_n)$.  Let $p\in(0,1/2]$ and let 
\begin{equation}\label{eq:def-median_n}
    \altq_n=\altq_n(p)=\inf\{x:\P(X_n\leq x)\geq 1-p\}
\end{equation}
be the $(1-p)$-th quantile of $X_n$. We have $\P(X_n> \altq_n)\leq p$ by construction. 

We define normalizing constant $\median_n$ as
\begin{equation}\label{eq:def-median_n-second}
    \median_n:=\altq_n(\fp),
\end{equation}
where $\fp\in(0,1)$ is some small number specified in Section~\ref{subsubsec:holder-sketch}.

So a priori, we prove Theorem~\ref{thm:main-convergence} when $\median_n$ is chosen to be a sufficiently large quantile of $X_n$. 
A posteriori in the proof of Corollary~\ref{cor:scaling-exponent-distance}, we will see that in fact the scaling limit exists for all $p$, and choosing a different value for $p$ only results in a deterministic factor in the limit. 

Throughout the rest of the paper, we will often implicitly take $p$ small enough that the desired results are true, and omit the dependence of $\altq_n$ on~$p$. Recall that $\kappa$ is defined to be the unique integer such that $R^{\kappa-1}<n\leq R^\kappa$. The immediate goal is to show the following.

\begin{prop}\label{prop:apriori}
Let $\varepsilon>0$. There exist $R_{0}=R_{0}(\varepsilon)$ and $C=C(\varepsilon)$ such that for each $R \ge R_0$ there exists $p_0=p_0(\varepsilon,R)$ such that for all $p<p_0$, $k \geq 1$, and $n> R^k$,
\begin{equation}
\P\!\left(X_{R^{\kappa-k}}> \altq_n\right)\leq CR^{-\!\left(\frac{3}{4}-\varepsilon\right)k}.
\end{equation}
\end{prop}
For the proof, we divide $[0,1]^2$ into smaller $R$-adic boxes in $\cB^R_{k,n}$ with $R^{\kappa-k}\times R^{\kappa-k}$ vertices for some $k\geq1$. We roughly say a box $B$ is \emph{important} if there are closed arms from the top and bottom sides of $B$ to the centers of the top and bottom sides of $[0,1]^2$ and open arms from the left and right sides of $B$ to the centers of left and right sides of $[0,1]^2$. It is \emph{bad} if inside $B$, the event $\cE(B)$ happens, the open and closed arms inside $B$ defined by the event $\cE(B)$ link up with the open and closed arms outside that go to the sides of $[0,1]^2$, and furthermore that $X_B> \altq_n$. Note that $X_{R^{\kappa-k}}$ is distributed as $X_B$ conditional on $\cE(B)$. Points~\eqref{property:E3}--\eqref{property:E5} in the definition of the event $\cE$ guarantee that the clusters inside $B$ link up with the arms outside of $B$ with probability bounded away from 0. They also guarantee that $u$ and $v$ always have 4 long arms, which will be useful in the annulus covering parts of the paper. Combining all this, we show 
\begin{equation}
    \P(\text{$B$ is bad})\asymp\P(\text{$B$ is important})\P(X_{R^{\kappa-k}}> \altq_n).
\end{equation}
Furthermore, if a box $B$ is bad, then $\cE(\Lambda_n)$ occurs and $X_{\Lambda_n}> \altq_n$. Thus, if the probability that there is at least one bad box is large, then the probability that $X_n> \altq_n$ must be large, but this probability is at most $p$ by definition of $\altq_n$. So in order to bound $\P(X_{R^{\kappa-k}}> \altq_n)$, we want to lower bound the number of important boxes.

Since the 4-arm exponent for planar critical percolation is $\frac{5}{4}$ (see also Section~\ref{subsec:arm-exponents}), the probability that a box at level $k$ is important is roughly $R^{-\frac{5}{4}k}$. So the expected number of important boxes is $R^{2k}R^{-\frac{5}{4}k}=R^{\frac{3}{4}k}$. Since the probability that a box is bad is bounded from above by $p$, it follows that $\P(X_{R^{\kappa-k}}> \altq_n)\lesssim R^{-\frac{3}{4}k}$. To make this argument precise, we use a second moment method to estimate the number of important/bad boxes. We bound the second moment by a multiscale argument decorrelating the events that two boxes that are far away are both important.

\subsubsection{Annulus crossing bounds}\label{subsubsec:annulus-sketch}
We now want to transform the right tail bound on $X_n$ into a more useful bound on the probability that the distance between two crossings of the annulus that are connected exceeds $\altq_n$. Let $\varepsilon>0$ and set $\alpha=\frac{3}{4}-\varepsilon$\phantomsection{}\label{def:alpha-first}. Let $R=R(\varepsilon)>8$ and $p=p(\varepsilon,R)$ be as in Proposition~\ref{prop:apriori} . The goal is to prove the following result.
\begin{prop}\label{prop:annulus-crossing}
Let $f\in\N$ and consider the annulus $A(0;2n^{-1}R^{\kappa-f},4n^{-1}R^{\kappa-f})$. Let $s<1$ and define $G=G(s,f)$ to be the event that for each pair of open annulus crossings that are that are connected inside $B(0,2n^{-1}R^{\kappa-f})$, the geodesic resp.\ resistance distance between them within $B(0,4n^{-1}R^{\kappa-f})$ is at most $s\altq_n$. Then for every $\zeta,p'>0$, there exists $C=C(\zeta,p')>0$ independent of $f,s$ and some event $F$ depending only on the sites in $A(0;n^{-1}R^{\kappa-f},n^{-1}R^{\kappa-f+1})$ such that $\P(F)>1-p'$ uniformly in $n,f$ and for all $R$ sufficiently large depending on $\zeta$, and $p$ small enough depending on $R$,
\begin{equation}\label{eq:distance-between-crossings-proposition}
    \P(G^{c}\mid F)\leq Cs^{-(\frac{3}{4}+\zeta)\alpha}R^{-\alpha f}.
\end{equation}
\end{prop}
In fact, we prove something slightly stronger than the above. See Proposition~\ref{lem:distancebetweencrossings} for a precise statement of the desired result. The idea behind the proof is as follows. Assume first that there are only two open crossings $\eta_1,\eta_2$ of the annulus $A(0;2n^{-1}R^{\kappa-f},4n^{-1}R^{\kappa-f})$. Then there also exist two closed paths $\phi_1,\phi_2$ that touch at $v_1\in\eta_1$ and $v_2\in\eta_2$. So $\{v_1\}\cup\phi_1\cup\{v_2\}\cup\phi_2$ form a loop inside $A(0;2n^{-1}R^{\kappa-f},4n^{-1}R^{\kappa-f})$ surrounding $B(0;2n^{-1}R^{\kappa-f})$. Consider the event $E$ that $\eta_1,\eta_2$ are connected with open paths to the centers of the left- and right-hand sides of $B(0,\frac{1}{2}n^{-1}R^{\kappa-f+1})$ respectively, and $\phi_1,\phi_2$ are connected with closed paths to the top and bottom sides of $B(0,\frac{1}{2}n^{-1}R^{\kappa-f+1})$ respectively. Then $E$ implies that the event $\cE(B(0,\frac{1}{2}n^{-1}R^{\kappa-f+1}))$ occurs. Furthermore, $v_1$ and $v_2$ are pivotal sites for global left-to-right crossing of $B(0,\frac{1}{2}n^{-1}R^{\kappa-f+1})$. So $X_{B(0,\frac{1}{2}n^{-1}R^{\kappa-f+1})}\geq D(v_i,v_j)$ on $E$. Here $D(v_i,v_j)$ refers to the (geodesic or resistance) distance inside the region enclosed by the loop $\{v_1\}\cup\phi_1\cup\{v_2\}\cup\phi_2$.

What if $E$ does not occur? Then with high probability, $\eta_1,\eta_2,\phi_1,\phi_2$ are still almost connected to the sides of $B(0,\frac{1}{2}n^{-1}R^{\kappa-f+1})$ by bounded sequences of touching, but not intersecting, clusters. If we resample the percolation configuration in small boxes around these touch-points, then with positive probability all touching loops become intersecting and the event $E$ occurs. Furthermore, we can show that if the configuration is nice (this is the event $F$), then the successful resampling probability is uniformly bounded away from $0$, and further the distance $D(v_i,v_j)$ remains unchanged. Combining these two facts, we obtain
\begin{equation}
    \P(D(v_1,v_2)> \altq_n,F)\lesssim\P(D(v_1,v_2)>\altq_n,F,E)\lesssim\P(X_{R^{\kappa-f+1}}>\altq_n)\lesssim R^{-\alpha f},
\end{equation}
where the last inequality is precisely the a priori bound from Proposition~\ref{prop:apriori}.

What if we have $m>2$ open crossings? Then there are points $v_1,\ldots,v_m$ and closed paths $\phi_1,\phi_m$ such that $\{v_1\}\cup\phi_1\cup\ldots\cup\{v_m\}\cup\phi_m$ form a loop in the annulus. If we resample in small boxes around the points $v_3,\ldots,v_m$, then the paths $\phi_3,\ldots,\phi_m$ link up with probability bounded away from 0 like before and we end up with only two open crossings. Again, we can define $F$ such that this resampling procedure does not affect $D(v_i,v_j)$ on $F$. We can then proceeds as before.

The addition of the extra factor $s<1$ will be important in later parts of the paper where we will sum the distances across multiple annuli. We first prove the upper bound $\P(X_{R^{\kappa-k}}>s\altq_n)\lesssim s^{-(\frac{3}{4}+\zeta)\alpha}R^{-\alpha k}$ by showing that $\altq_{R^{-k}n}\lesssim R^{-(3/4-\varepsilon)k}\altq_n$. The rest of the proof then follows as before.

\subsubsection{Superpolynomial concentration}\label{subsubsec:subpolynomial-sketch}
Having established the polynomial upper bound on the right tails of $X_n$ with some exponent $\alpha$ that is almost $\frac{3}{4}$, we then show that we can bootstrap this estimate and obtain concentration with some exponent $\alpha+\beta$. From this, we obtain a superpolynomial upper bound. More precisely, we prove the following. Recall that $\altq_n = \altq_n(p)$ depends implicitly on $p$.
\begin{prop}\label{prop:subpolynomialconc}
There exist $\varepsilon_0>0$, $\beta>0$, $R_0 > 0$, and a function $p_0(R)$ such that the following holds for each $R \ge R_0$ and $p < p_0(R)$. Assume that 
for some $\alpha>\frac{3}{4}-\varepsilon_0$, there exists $C_0=C_0(\alpha,R,p)$ such that for all $k\in\N$ and $n> R^k$, we have
\begin{equation}\label{eq:assumptionsubpoly}
    \P(X_{R^{\kappa-k}}> \altq_n)\leq C_0R^{-\alpha k}.
\end{equation}
Then 
there exists $C_1=C_1(\alpha,R,p)$ such that for all $k\in\N$ and $n> R^k$, we have
\begin{equation}
    \P(X_{R^{\kappa-k}}> \altq_n)\leq C_1R^{-(\alpha+\beta) k}.
\end{equation}
\end{prop}

Note that $\beta$ does not depend on $\alpha$. Thus, we immediately obtain a superpolynomial right tail bound on the annulus crossing distance.
\begin{cor}\label{cor:subpolynomialconc-general}
The statement in equation~\eqref{eq:assumptionsubpoly} and hence Proposition~\ref{prop:annulus-crossing} (and Proposition~\ref{lem:distancebetweencrossings}) hold for any exponent $\alpha>0$, upon choosing some $R$ large enough and $p=p(R)$ small enough, independently of $\alpha$.
\end{cor}
Throughout the proof of Proposition~\ref{prop:subpolynomialconc}, we fix some $\alpha>\frac{3}{4}-\varepsilon_0$ and take $\altq_n=\altq_n(p)$, where we suppose that $R$ is sufficiently large and $p$ sufficiently small such that Proposition~\ref{prop:apriori} holds for $\varepsilon_0$. Furthermore, at certain steps of the proof, we will need to further increase $R$ (and thus also further decrease $p$ depending on $R$) such that certain estimates hold. Crucially, the necessary value of $R$ does not depend on $\alpha$, so the bootstrapping procedure is valid for a fixed large $R$ and small $p$. 

We now give a sketch of the proof. Let $k \in \N$ and consider percolation on $\Lambda_n$ restricted to some box $B^{(k)}\in\cB^R_{k,n}$\phantomsection\label{def:Bk-second}. The goal is to bound $\P(X_{R^{\kappa-k}}>\altq_n)\lesssim\P(X_{B^{(k)}}>\altq_n,\cE(B^{(k)}))$. We say a site $x$ is \emph{pivotal} if every open path from the left side of $B^{(k)}$ to the right side passes through $x$. We will also call them \emph{pinch points}. Assume the event $\mathcal{E}(B^{(k)})$ holds, so that there exist pivotal sites. Let $u,v$ be the left- and rightmost pivotal points and let $(x_0=u,x_1,x_2,\ldots,v)$ be the sequence of pivotal points from left to right. Write $X_{x_i,x_j}=D(x_i,x_j)$ to be the (geodesic or effective resistance) distance between $x_i$ and $x_j$. In particular, on the event $\mathcal{E}(B^{(k)})$, we have $X_{B^{(k)}}=X_{u,v}$. Define $\tau_0=0$ and for $\ell\geq 1$ define
\begin{equation} \label{eq:tau_definition}
    \tau_{\ell}=\min\{i\geq\tau_{\ell-1}\colon X_{x_{\tau_{\ell-1},x_i}}> \altq_n\}.
\end{equation}
If the set above is empty, set $\tau_{\ell}=\infty$. Let
\begin{equation}
    N=\min\{\ell\geq 1\colon \tau_{\ell}=\infty\}.
\end{equation}
Then, because the $x_i$ are pivotal sites, we know that $D(u,v)=D(x_{\tau_0},x_{\tau_1})+D(x_{\tau_1},x_{\tau_2})+\cdots + D(x_{\tau_N-1},v)$ by the series law (property~\eqref{property:serial} of the metric). In particular,
\begin{equation}\label{eq:upperboundXdecomp}
    X_{u,v}\leq N\altq_n+\sum_{\ell=1}^{N-1}X_{{\tau_{\ell}-1},\tau_{\ell}}.
\end{equation}
We bound both terms separately. We first show the following bound on the second term.
\begin{prop}
There exist $\beta>0$ and a constant $C(R,\alpha)>0$ such that
\begin{equation}
    \begin{split}
        \P\!\left(\exists i\colon X_{x_{i},x_{i+1}}> \altq_n\right)\leq C(R,\alpha)\cdot R^{-(\alpha+\beta)k},
    \end{split}
\end{equation}
for all $n> R^k$.
\label{prop:no_individual_bubble_long}
\end{prop}

The intuition behind the proof is as follows. In between two consecutive pivotals $x_i$ and $x_{i+1}$, there is a \emph{bubble} bounded by two disjoint open interface paths $\gamma_1$ and $\gamma_2$. Assume the bubble is roughly of size $R^{-d}$ for $d\geq k$. If we cover the bubble with annuli of size a bit smaller than $R^{-d}$, then $\gamma_1$ and $\gamma_2$ form crossings of a string of annuli. For an annulus of size $R^{-f}$, the probability of the (geodesic resp.\ resistance) distance between two crossings being large (a bit less than $\altq_n$) is just a bit more than~$R^{-\alpha f}$. We say such an annulus is bad. So the probability of there being two disjoint bad annuli is just a bit more than $R^{-2\alpha f}$. This already gives us an exponent strictly larger than~$\alpha$, so we may assume without loss of generality that there is at most one bad annulus. If, say, $\gamma_1$ does not pass through this annulus, then $\gamma_1$ can be covered by annulus crossings that all have small distance to each other (see Figure~\ref{fig:linking_boxes}), and we can `string together' short distances between these crossings to conclude that $D(x_i,x_{i+1})\leq \altq_n$ by choosing the correct exponents. We refer to the proof of Lemma \ref{lem:crossingbubble} for an explanation of this stringing procedure. If both $\gamma_1$ and $\gamma_2$ pass through the bad annulus, then we have a 5-arm event occurring for that annulus. Multiplying the probability that an annulus is bad and the probability of the 5-arm event, we obtain a bound of $R^{-(\alpha+\beta')f}$ for some $\beta'>0$. Summing over all annuli of level $f$ slightly larger than $d$, we obtain that the probability that $X_{x_i,x_{i+1}}>\altq_n$ for a bubble of size $d$ is at most $R^{-(\alpha+\beta)d}$ for some $\beta>0$. By a 4-arm exponent calculation, it can be shown that the number of bubbles of size $R^{-d}$ is of order $R^{\frac{3}{4}d}$. Summing over all such bubbles and using that $\alpha+\beta>\frac{3}{4}$ when $\alpha$ is chosen close enough to $\frac{3}{4}$, we obtain the desired bound. Here we need that $\alpha$ can be chosen arbitrarily close to $\frac{3}{4}$. 

Having shown the above, we obtain a bound of $X_{u,v}\le N\altq_n+N\altq_n=2N\altq_n$ in~\eqref{eq:upperboundXdecomp} with probability at least $1-O(R^{-(\alpha+\beta)k})$. The following proposition bounds $N$ by $2$ with probability at least $1-O(R^{-(\alpha+\beta)k})$. 
\begin{prop} There exists $\beta>0$ and $C=C(\alpha,R)$ such that for all $k$ and $n> R^k$,
\[\mathbf{P}(\tau_2 <\infty)\leq C R^{-(\alpha+\beta)k}.\] \label{prop:bound_on_stopping_time}
\end{prop}
Note that if the $\tau_i-\tau_{i-1}$ were i.i.d., then we would have $\P(\tau_2<\infty)=\P(\tau_1<\infty)^2\lesssim R^{-2\alpha k}$. We show that this is almost true, at least with some exponent that is strictly larger than $\alpha$. At time~$\tau_1$, we can explore the domain that is enclosed by the top and bottom touching clusters up until~$x:=x_{\tau_1}$. We call the complement of that domain $C$. We consider the \emph{Whitney decomposition} of $C$, which is a covering of $C$ with boxes such that the size of a Whitney box is proportional to its distance to the explored domain. We will first show that there exists a global event $H$ that happens with very high probability, and on the event $H$ it will be possible to find an open path $\gamma$ from~$x$ to~$v$ that can be covered by a small number of $R$-adic Whitney boxes of each level $i\geq k+1$. The potential obstruction is that if $\gamma$ comes close to the boundary of $C$, then the number of Whitney boxes that cover $\gamma$ is large, since the Whitney boxes close to the boundary of $C$ are small. However, the boxes close to the boundary of $C$ are also close to the explored bubbles, which leads to these Whitney boxes having 6 so-called \emph{semi-arms}. The event $H$ bounds the number of Whitney boxes with 6 semi-arms. Following this, we will show that with high probability all of these boxes are good, i.e., open crossings of the boxes have short geodesic resp.\ resistance distance. This is possible since the number of boxes in the covering of $\gamma$ is not too big. Finally, we will string together short crossings across each one of these boxes to obtain an upper bound on the remaining distance $X_{x,v}$. There is a technical issue in showing that the covering boxes are good. To bound the distance between annulus crossings for annuli that cover the path from $x$ to $v$, we would like to use the annulus crossing bound from Proposition~\ref{prop:annulus-crossing} \emph{conditionally on $\{\tau_1 < \infty\}$} (since this is what improves the exponent from $\alpha$ to $\alpha+\beta$). However, as we explain in Section~\ref{subsec:bootstrapping}, the conditioning on $\cE(B^{(k)}) \cap \{\tau_1 < \infty\}$ induces a bias, and the percolation in the domain $C$ is no longer i.i.d.\ Bernoulli. Thus, we also need to prove some results to deal with decorrelating the events $F,G$ from Proposition \ref{prop:annulus-crossing} for an annulus in $C$ and the event $\cE(B^{(k)}) \cap \{\tau_1 < \infty\}$.

\subsubsection{H\"older continuity}\label{subsubsec:holder-sketch}
The last technical result we need for proving tightness is H\"older continuity of the geodesic/resistance metrics on percolation with respect to the path metric with a constant that is uniform over $n$.
\begin{prop}\label{prop:holderct}
For all $\eta<\frac{5}{8}$ and $\varepsilon>0$, there exist $p=p(\eta)$ and $C=C(\eta,\varepsilon)>0$ such that for all $n$,
\begin{equation}
    \P\!\left(\forall x,y\in \Lambda_n\colon\,\altq_n^{-1}D_n(x,y)\leq C\dpt_n(x,y)^\eta\right)>1-\varepsilon.
\end{equation}
\end{prop}
We are now ready to define $\median_n$. To prove Theorem~\ref{thm:main-convergence}, we will need that $d_n=\median_n^{-1}D_n$ is uniformly H\"older continuous with some H\"older exponent. Since the precise H\"older exponent is not important, we simply choose $\eta=1/4$ and set $\fp:=p(1/4)$ and $\median_n=\altq_n(\fp)$\phantomsection{}\label{def:median_n-hoelder}.

The proof of Proposition~\ref{prop:holderct} uses many ingredients from previous proofs in this paper. To bound $D_{n}(x,y)$, we consider again a Whitney decomposition, this time of the set $\Lambda_n\setminus\{x,y\}$. Let $\gamma$ be a path from $x$ to $y$ with minimal diameter. Again, we would like to show that we can cover $\gamma$ in not too many Whitney boxes and then use the same stringing argument to conclude that $D_{n}(x,y)$ is small. However, there is a technical issue: we can only apply Proposition~\ref{prop:annulus-crossing} to annuli that are fully contained inside~$\Lambda_n$. So the stringing argument will not work if our Whitney boxes covering $\gamma$ are too close to the boundary of $\Lambda_n$. To deal with this, we take the size of the Whitney boxes to be proportional to the distance to $\partial\Lambda_n\cup\{x,y\}$. 
This is expensive: if $\gamma$ stays too close to the boundary then we need too many boxes, and this will result in a bound that is too weak. However, unless there is a closed blocking path coming close to the boundary, we can reroute $\gamma$ to be away from the boundary. If there is such a blocking path, then the Whitney box must have 3 long arms. We can bound the number of such boxes by a 3-arm half-plane exponent calculation.

Let us now explain how to get the exponent $\frac{5}{8}$. Let $k_0=\lfloor-\log_R\dpt_n(x,y)\rfloor$. If the path $\gamma$ does not come near the boundary, the number of Whitney boxes covering $\gamma$ at each level $k$ would be bounded. From~\eqref{eq:distance-between-crossings-proposition}, we get that with high enough probability, the distance between two crossings of an annulus at $R$-adic level $k$ is at most $R^{-(\frac{3}{4}-o(1))k}\altq_n$. Since the number of such good annuli at level~$k$ is bounded for each level, we would obtain
\begin{equation}\label{eq:chaining_sum}
    d_n(x,y)\lesssim\sum_{k\geq k_0}[\textrm{$\#$ of Whitney boxes at level $k$}]\cdot R^{-(\frac{3}{4}-o(1))k}\lesssim R^{-(\frac{3}{4}-o(1))k_0}\approx\dpt_n(x,y)^{\frac{3}{4}-o(1)}.
\end{equation}
In the case $\gamma$ comes close to the boundary of $\Lambda_n$, we cannot guarantee that the number of Whitney boxes of each level $k$ is bounded, but since they are all contained in a box of Euclidean size $R^{-k_0}$, their number (i.e.\ boxes of size $R^{-k}$ and Euclidean distance to the boundary proportional to $R^{-k}$) is trivially bounded by $R^{k-k_0}$. This is not enough to sum up the distance in~\eqref{eq:chaining_sum}, but as we mentioned before, if $\gamma$ is forced to go through a small Whitney box near the boundary, then that box has $3$ long arms. By an arm exponent calculation, we show that with high probability, there are at most $R^{k/3}$ such Whitney boxes. Therefore, if we estimate the number of Whitney boxes by $R^{k-k_0}$ for $k \le \frac{3}{2}k_0$, say, and by $R^{k/3}$ for $k > \frac{3}{2}k_0$, then we obtain in~\eqref{eq:chaining_bound} the estimate $d_n(x,y)\lesssim  R^{-\frac{5}{8}k_0}\approx \dpt_n(x,y)^{\frac{5}{8}}$.

This estimate is not optimal, but for the results in this paper any exponent is good enough, so we made no effort to push the exponent any further.

\subsection{Identification of the limit}\label{sec:uniqueness-sketch}
Let $(\median_n)_{n\in\N}$ be as defined in Sections~\ref{subsubsec:apriori-sketch} and~\ref{subsubsec:holder-sketch}. By Propositions~\ref{prop:GHf-tightness-dU} and~\ref{prop:main-conv-dU-GH}, for every subsequence $I\subset\N$, we can choose a further subsequence $J\subset I\subset\N$ such that $(\Lambda_{Q,n},\dpt_{Q,n},d_{Q,n},\pi_n)$ converge weakly in the GHf topology along $J$ for all $Q\in\cQ$. 

We denote the limiting metric of $\dg_{Q,n}$ resp.\ $\dr_{Q,n}$ by $\dg_Q$\phantomsection{}\label{def:dg_Q} resp.\ $\dr_Q$\phantomsection{}\label{def:dr_Q}. If $Q=(0,1)^2$, we simply write $\dg=\dg_Q$\phantomsection{}\label{def:dg_scaling_limit} resp.\ $\dr=\dr_Q$\phantomsection{}\label{def:dr_scaling_limit}. The goal of this section is to show that $\dg$ and $\dr$ are indeed the (unique up to multiplicative constants) conformally covariant geodesic and resistance metrics respectively on the clusters of $\Lambda$, as previously introduced in Section~\ref{subsec:statements}.

In the following, we let $d=\dg$\phantomsection{}\label{def:def-d} or $d=\dr$. Note that for both the geodesic and resistance metrics, we have $d_Q\leq d_{Q'}$ on $\Lambda_{Q'}$ if $Q'\subset Q$ (since we have $d_{Q,n}\leq d_{Q',n}$ by construction, and the inequality is preserved in the GHf limit by \cite[Lemma~A.10]{amy2025tightness}). To identify the limit with the geodesic CLE$_6$ metric resp.\ the CLE$_6$ resistance metric, we verify the conditions in Section~\ref{subsec:cle-metrics}. For this, we consider the restriction of $d$ to the cluster $\cX$ (recall the definition of $\cX$ and $\fC$ from Section~\ref{subsec:cle-metrics}) and define the collection of metrics $(\fd_V)_{V \in \fC}$ as follows:
\begin{equation}\label{eq:def-dU}
    \fd_V(x,y)=\sup_{\substack{Q\Supset V\\Q\in\cQ}}d_Q(x,y), \quad x,y \in \cX \cap \overline{V}.
\end{equation}
Upon verifying that the collection of metrics $(\fd_V)_{V \in \fC}$ satisfy the definition in Section~\ref{subsec:cle-metrics}, we obtain from Theorem~\ref{th:uniqueness-geodesic-cle} resp.~\ref{th:uniqueness-resistance-cle} that the restriction of $d$ to the cluster $\cX$ is (a multiple of) the geodesic CLE$_6$ metric resp.\ the CLE$_6$ resistance metric, and that the weak convergence of $(\cX_n,\dpt_{Q,n},d_{Q,n},\pi_n)$ holds for the full sequence. It follows by the local absolute continuity (see \cite[Section~7.2]{my2025metric}) that the global metric $d$ (defined on all clusters) agrees with (a global multiple of) the geodesic CLE$_6$ metric resp.\ the CLE$_6$ resistance metric.

\begin{thm}\label{th:dU-is-CLE-metric}
 The collection $(\cX\cap\overline{V},\fd_V)_{V\in\fC}$ as defined in~\eqref{eq:def-dU} defines a $\CLE_6$ metric. If $d=\dg$, then $d$ is indeed a geodesic $\CLE_6$-metric. Furthermore, if $d=\dr$, then $d$ is a resistance $\CLE_6$-metric.
\end{thm}

For the proof of Theorem~\ref{th:dU-is-CLE-metric}, we show that $(\fd_V)_{V\in\fC}$ satisfies all of the axioms of a $\CLE_6$ metric. This is done in Section~\ref{sec:uniqueness-proof}. Analogues of the axioms are clearly true for the finite metrics $d_{Q,n}$, so it just remains to show that they hold in the GHf limit.
The continuity of $\fd_V$ with respect to $\dpt_{\overline{V}}$ follows immediately from Proposition~\ref{prop:main-conv-dU-GH}. The monotonicity follows immediately from~\eqref{eq:def-dU}. Compatibility is stated in Proposition~\ref{prop:compatibility}. The series and generalized parallel laws are stated in Proposition~\ref{prop:series} and~\ref{prop:parallel}. The Markovian Property is stated in Proposition~\ref{prop:markov}. Translation invariance is stated in Proposition~\ref{prop:translation-invariance}. This shows that $(\fd_V)_{V\in\fC}$ is a $\CLE_6$-metric. If $d=\dg$, it is shown in Section~\ref{sec:geodesic-cle} that $d$ is a geodesic $\CLE_6$-metric. If $d=\dr$, it is shown in Section~\ref{sec:edge-weight-compatibility} that $d$ is a resistance $\CLE_6$-metric.

We obtain convergence of the geodesic and resistance metrics to $\dg$ and $\dr$ as a corollary. 
\begin{cor}\label{cor:final-convergence-metrics}
There exist true metrics $\dg$ and $\dr$ on all clusters $\cX^{(m)}$ such that
\begin{equation}
    \begin{split}
        (\cX^{(m)}_n,\dg_n,\pi_n)&\rightarrow(\cX^{(m)},\dg,\pi),\\
        (\cX^{(m)}_n,\dr_n,\pi_n)&\rightarrow(\cX^{(m)},\dr,\pi)\\
    \end{split}
\end{equation}
weakly with respect to $\Delta_{\GHf}$ as $n\rightarrow\infty$.
\end{cor}
For the cluster $\cX$, this corollary follows almost immediately from Theorems~\ref{th:uniqueness-geodesic-cle} and~\ref{th:uniqueness-resistance-cle}, which imply that all subsequential limits are multiples of each other. It remains to show that the subsequential limiting metric $d$ is non-trivial and that all subsequential limits are actually equal and not just multiples. We can then carry over convergence for $\cX$ over to all clusters $\cX^{(m)}$ by a local absolute continuity argument. See Section~\ref{sec:final-convergence-metrics} for a proof of Corollary~\ref{cor:final-convergence-metrics}.

\subsection{Measure and projection convergence}\label{sec:sketch-measure}

Having proved Corollary~\ref{cor:final-convergence-metrics}, the final ingredient to obtain Theorem~\ref{thm:main-convergence} is the convergence of the counting measure on percolation clusters. This was essentially proved by Garban, Pete, and Schramm in \cite{gps} as stated in Section~\ref{subsubsec:measure-conv}. However, their result concerns the counting measure for an annulus defined in~\eqref{eq:annulus_measure}. Thus, we first need to express the counting measure on $\cX^{(m)}_n$ in terms of counting measures for annuli. We are then able to transfer GHf convergence and measure convergence to the desired GHPf convergence. In order to show that the limiting measure equals the CLE$_6$ gasket measure $\mu$, we show that it also satisfies the characterization from Proposition~\ref{prop:char-measure-my} by making use of the analogous property of the annulus measure in Proposition~\ref{prop:char-gps}. This is stated as Proposition~\ref{prop:char-mu-tilde}. See Section~\ref{sec:ghpf-proof} for the details.

\subsection{Random walk convergence}\label{sec:sketch-random-walk}
The proof of Theorem~\ref{thm:main-random-walk} is now straightforward, and so we give it here.
\begin{proof}[Proof of Theorem~\ref{thm:main-random-walk}]
Let $X$ be the Hunt process associated with $(\cX^{(m)},\dr,\mu^{(m)},\rho)$, as described in Section~\ref{subsec:resistance-metrics-walks}. Since $\dr,\mu^{(m)}$ are measurable with respect to the $\Gamma$, we have that $(\cX^{(m)}_n,\dr_n,\mu^{(m)}_n,\pi_n,\rho_n)$ converges in law to $(\cX^{(m)},\dr,\mu^{(m)},\pi,\rho)$, so that we can couple them such that the sequence converges in $\rGHPf$ almost surely. The conclusion of Theorem~\ref{thm:main-random-walk} then follows from Theorem~\ref{th:resistance-random-walk-convergence}.
\end{proof}

\subsection{Existence of distance, resistance, and displacement exponents}\label{subsec:scaling-exponent-proof-sketch} 
The only results left to prove are Corollaries~\ref{cor:scaling-exponent-distance}--\ref{co:spectral-dimension}. The fact that $\median_n=n^{\beta+o(1)}$ follows mostly from Theorem~\ref{thm:main-convergence} and the scale-covariance $R^\beta d_Q(\cdot,\cdot;\Gamma)=d_{RQ}(R\cdot,R\cdot;R\Gamma)$ of the CLE$_6$ metrics proved in~\cite{my2025metric,my2025diffusion}. In order to show that $\E[S_n]\asymp\operatorname{quant}_q[S_n]\asymp \median_n$, we need a more quantitative version of the H\"older bound in Proposition~\ref{prop:holderct}, which we show in Section~\ref{sec:holder}. See Section~\ref{sec:scaling-exponent-proof} for a full proof of Corollary~\ref{cor:scaling-exponent-distance}. The existence of the spectral dimension then follows almost immediately from Kumagai and Misumi's \cite{kumagai2008heat} results relating the spectral dimension of a random graph to its resistance and volume growth.

\section{Convergence of the path metric}\label{sec:convergence-path-metric-proof}
In this section, we prove Proposition~\ref{prop:convergence-ghf-path}.
\subsection{Tightness}
We first show tightness of the path metric over all clusters and subdomains of $[0,1]^2$. Recall the definition of $\cX_{Q,n}^{(m)}$ from Section~\ref{subsec:notation}, which is the $m$-th largest cluster by diameter of $\Lambda_{Q,n}$.
\begin{prop}\label{prop:finaltightness}
The collection $\{(\cX^{(m)}_{Q,n},\dpt_{Q,n},\pi_n)\colon n\in\N,\,m\in\N,\,Q\in\cQ\}$ is tight with respect to the GHf topology.
\end{prop}

\begin{lem}\label{lem:uniftotallybounded}
Let $\delta>0$. Let $N^p_n(\delta,m,Q)$ be the minimum number of $\delta$-balls in the metric $\dpt_{Q,n}$ needed to cover $\cX^{(m)}_{Q,n}$. Let $N'(\delta)$ denote some integer-valued function of $\delta$. Let $E_n(N')$ be the event that for all $\delta>0$, $m\in\N$, and $Q\in\cQ$, we have $N^p_n(\delta,m,Q)\leq N'(\delta)$. Then for all $\varepsilon>0$, there exists a function $N'_{\varepsilon}$ such that for all $n$,
\begin{equation}
\P(E_n(N'_{\varepsilon}))>1-\varepsilon.
\end{equation}
\end{lem}

\begin{figure}[ht]
    \centering
    \includegraphics[scale=0.75]{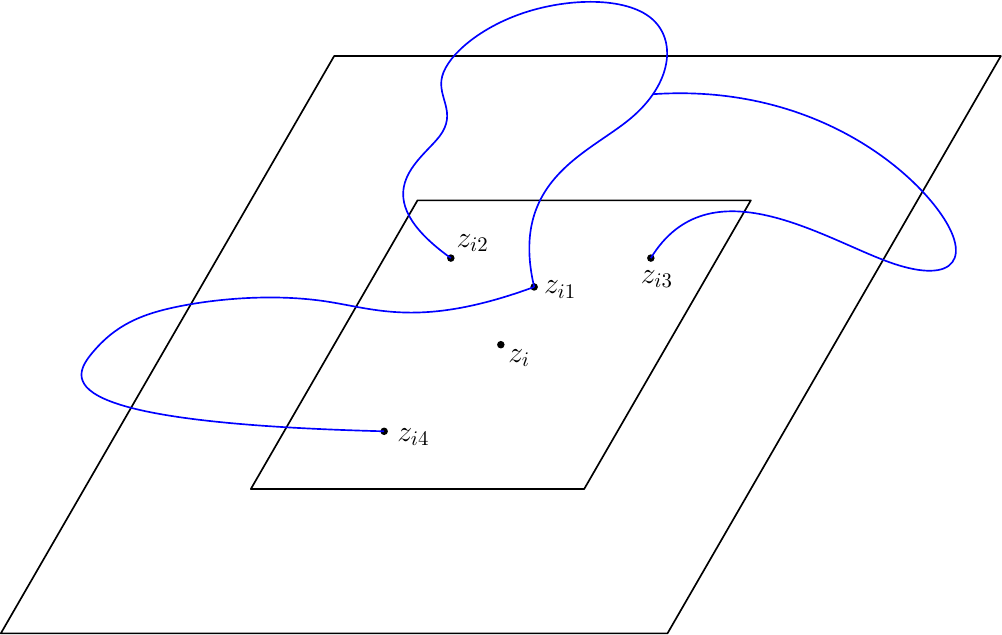}
    \caption{Four points $z_{i1},\ldots,z_{i4}$ together with an annulus such that the pairwise path diameter distances are at least the size of the annulus. Because $\dpt_n(z_{i1},z_{i2})$ is greater than the diameter of the annulus, it must be that any open path connecting these two points leaves the annulus. Moreover, since these points are in the same cluster, such an open path exists. In particular, we can find $2$ disjoint open crossings of our annulus -- the section of the path going from $z_{i1}$ to the exterior of the annulus, and the section of the path going from the exterior to $z_{i2}$. These must be disjoint as otherwise we would be able to find an open path connecting $z_{i1}$ and $z_{i2}$ within $B(z_i,\ell^{-1})$. Every subsequent point $(z_{i3}, z_{i4},\ldots)$ adds one more open crossing. Indeed, there is an open path between $z_{ij}$ to $z_{i1}$ that leaves the annulus. If it doesn't create any additional crossings, then $z_{ij}$ must be connected to one of $z_{i1},z_{i2},\ldots,z_{i,j-1}$ within $B(z_i,\ell^{-1})$, contradicting our original assumption.}
    \label{fig:uniformly_totally_bounded}
\end{figure}

\begin{proof}
Note that it suffices to prove the statement for all $\delta$ of the form $\delta=\frac{1}{\ell}$ for $\ell\in\N$. Consider some $n\in \N$ and let $z_1,\ldots,z_{\ell^2}\in\Lambda_n$ such that the boxes $B(z_i,\ell^{-1}/2)$ cover $[0,1]^2$. Let $K>0$, whose exact value we will determine later. Let $F_{\ell}=F_{\ell}(K)$ be the event that for all $1\leq i\leq\ell^2$, there are at most $K$ disjoint open crossings of $A(z_i;\ell^{-1}/2,\ell^{-1})$ for percolation in $\Lambda_n$. Let $Q\in\cQ$ and $m\in\N$. Suppose that there exist $z_{i1},z_{i2},\ldots,z_{i(K+1)}\in\cX_{Q,n}^{(m)}\cap B(z_i,\ell^{-1}/2)$ such that $\dpt(z_{ij},z_{ij'})\geq\frac{2}{\ell}$ for all $1\leq j,j'\leq K+1$. Then the annulus $A(z_i;\ell^{-1}/2,\ell^{-1})$ must have $K+1$ disjoint open crossings (see Figure~\ref{fig:uniformly_totally_bounded}). So on the event $F_{\ell}(K)$, for each box $B(z_i,\ell^{-1}/2)$, we can choose at most $K$ sites $z_{i1},\ldots,z_{iK}$ such that 
\begin{equation}
	\max_{v\in B(z_i,\ell^2)\cap\cX^{(m)}_{Q,n}}\min_{1\leq j\leq K}\dpt_{Q,n}(v,z_{ij})\leq2\ell^{-1}.
\end{equation}
So on $F_{\ell}(K)$, the space $\cX^{(m)}_{Q,n}$ can be covered by at most $K\ell^2$ many $\dpt_{Q,n}$-balls of radius $2\ell^{-1}$, i.e., $N_n^p(2\ell^{-1},m,Q)\leq K\ell^2$. Let $0<p_0<1$ be a uniform (in n) upper bound on the probability that there exists an open crossing of $A(z;\ell^{-1}/2,\ell^{-1})\cap[0,1)^2$ in $\Lambda_n$. Then by the BK inequality,
\begin{equation}
	\P(F_{\ell}(K)^c)\leq\ell^2p_0^K<\varepsilon2^{-\ell}
\end{equation}
for $K=K(\ell,\varepsilon)$ sufficiently large. Let $N'_{\varepsilon}(2\ell^{-1})=K\ell^2$. Then
\begin{equation}
	\P(E_n(N'_\varepsilon)^c)\leq\sum_{\ell\geq1}\P(F_{\ell}(K)^c)<\varepsilon,
\end{equation}
which concludes the proof.
\end{proof}

\begin{proof}[Proof of Proposition~\ref{prop:finaltightness}] 
Note that on the event $E_n(N'_{\varepsilon})$ defined in Lemma~\ref{lem:uniftotallybounded}, each $(\cX^{(m)}_{Q,n},\dpt_{Q,n})$ where $m\in\N$, $Q\in\cQ$ is totally bounded as quantified by the function $N'_\varepsilon$, and hence is contained in a compact set $\cK(N'_\varepsilon)$ in the Gromov-Hausdorff topology. Furthermore, for any two points $x_n,y_n\in \cX^{(m)}_{Q,n}$ we have $\|\pi_n(x_n)-\pi_n(y_n)\|_2\leq\dpt_{Q,n}(x_n,y_n)$. So by \cite[Proposition A.4]{amy2025tightness}, the collection $\{(\cX^{(m)}_{Q,n},\dpt_{Q,n},\pi_n)\colon (\cX^{(m)}_{Q,n},\dpt_{Q,n}) \in \cK(N'_\varepsilon)\}$ is relatively compact with respect to the GHf topology. Tightness thus follows immediately from Lemma~\ref{lem:uniftotallybounded}. 
\end{proof}

\subsection{Uniqueness} 
The previous proposition implies that $(\Lambda_{Q,n},\dpt_{Q,n},\pi_n)_{n\in \N}$ has subsequential limits. We now show that every limit in fact equals $(\Lambda_Q,\dpt_Q,\pi)$.
\begin{proof}[Proof of Proposition~\ref{prop:convergence-ghf-path}]

By Proposition~\ref{prop:finaltightness}, Theorem~\ref{th:camia-newman}, and the Skorokhod representation theorem, for each subsequence $I\subset\N$, there exists a further subsequence $J\subset I$ and a coupling such that for all $Q\in\cQ$, $(\Gamma_{Q,n})_{n\in\N}$ converges to $\Gamma_Q$ and such that $(\Lambda_{Q,n},\dpt_{Q,n},\pi_n)_{n\in J}$ has a GHf limit, which we call $(\tilde{\Lambda}_Q,\tilde{d}^{\pth}_Q,\tilde{\pi})$. To complete the proof, we construct an isometry $\phi\colon\tilde{\Lambda}_Q\to\Lambda_Q$ for each $Q\in\cQ$. Throughout the proof, we assume that $(\Lambda_{Q,n},\dpt_{Q,n})_{n\in J}$ and $(\tilde{\Lambda}_Q,\tilde{d}^{\pth}_Q)$ are isometrically embedded in a common metric space in accordance with the GHf convergence. Further, by \cite[Theorem A.1]{aizenman1999holder}, switching to a subsequence if necessary, we can assume that the family of open paths in $\Lambda_{Q,n}$ converges in $\Delta_{\Hc}$.

Let $\tilde{x},\tilde{y}\in\tilde{\Lambda}_Q$ be distinct points in the same component of $\tilde{\Lambda}_Q$. There exist sequences $(x_n)_{n\in J}$ and $(y_n)_{n\in J}$ where $x_n$ and $y_n$ are in the same component of $\Lambda_{Q,n}$ that converge to $\tilde x$ and $\tilde y$, respectively, and such that their projections converge. 
Let $\eta_n$ be a path between $x_n$ and $y_n$ with minimal diameter that lies entirely inside $Q$, i.e.\ $\dpt_{Q,n}(x_n,y_n) = \diam(\eta_n)$. 
Let $\eta$ be any subsequential limit. By the convergence with respect to $\Delta_{\cur}$, we have that $\eta(0)=\tilde\pi(\tilde x)$ and $\eta(1)=\tilde\pi(\tilde y)$. Furthermore, the diameter of $\eta$ equals $\tilde{d}^{\pth}_{Q}(\tilde x,\tilde y)$. By Lemma~\ref{lem:touching-loops}, $\eta$ can be chosen to lie entirely inside $Q$. 
We show that $\eta$ is an open path in the $\CLE_6$-gasket. Let $\gamma\in\Gamma$ be a $\CLE_6$ loop. 
Let $\gamma_n\in\Gamma_n$ converge to $\gamma$. Then for all $n$ large enough, $x_n$ is either on the inside of $\gamma_n$ or on the outside. Assume the former. Since $x_n$ and $y_n$ are connected, $y_n$ must also be on the inside of $\gamma_n$. So $\eta_n$ lies entirely inside $\gamma_n$. Then $\eta$ must also lie entirely inside $\gamma$, so $\gamma$ does not cross $\eta$. Now assume $x_n$ and $y_n$ lie outside $\gamma_n$. Then by the same reasoning, $\eta$ must lie outside $\gamma$ and so $\gamma$ does not cross $\eta$. It follows that $\eta$ is an open $\CLE_6$-path. Furthermore, the endpoints $\eta(0),\eta(1)$ define prime ends $x,y\in\Lambda$ such that $\eta$ is an open path from $x$ to $y$. Let $x=\phi(\tilde x)$ and $y=\phi(\tilde y)$. Note that the map $\phi$ is indeed well-defined, i.e., the choice of $x$ does not depend on $\tilde y$. Indeed, let $\tilde{y}'\in\tilde\Lambda_Q$ and let $(y'_n)_{n\in J}$ converge to $\tilde y'$ in $\tilde d^{\pth}_Q$. Let $\eta'_n$ be an open path between $x_n$ and $y'_n$ with minimal diameter inside $Q$ and let $\eta'$ be a subsequential limit. Then by the same argument as above, $\eta'$ is an open $\CLE_6$ path. Moreover, the concatenation of $\eta_n,\eta'_n$ does not cross any loop of $\Gamma_{Q,n}$, hence the concatenation of $\eta,\eta'$ does not cross any loop of $\Gamma_Q$. Therefore $\eta(0),\eta'(0)$ define the same prime end. Similarly, our definition of $\phi(\tilde x)$ does not depend on the choice of the sequence $(x_n)$ converging to $\tilde x$, since for a different sequence $\tilde x_n \to \tilde x$ we can apply the same argument with a path $\tilde\eta_n$ with $\diam(\tilde\eta_n) = \dpt_{Q,n}(x_n,\tilde x_n) \to 0$. We conclude that $\dpt_Q(\phi(\tilde x),\phi(\tilde y))\leq\tilde d^{\pth}_Q(\tilde x,\tilde y)$.

We now show the opposite inequality. Let $x,y\in\Lambda$ and let $\eta$ be a simple $\CLE_6$-open path from $x$ to $y$ with smallest diameter. Let $\varepsilon>0$ and consider the $\varepsilon$-fattening $B(\eta,\varepsilon):=\{z\in[0,1]^2:\dist(z,\eta)<\varepsilon\}$. We claim that for all $n$ large enough, there exist $x'_n(\varepsilon),y'_n(\varepsilon)\in\Lambda_{Q,n}$ such that $\dist(\pi_n(x'_n(\varepsilon)),\pi(x))$, $\dist(\pi_n(y'_n(\varepsilon)),\pi(y))<\varepsilon$ and such that $x'_n$ and $y'_n$ are connected inside $B(\eta,\varepsilon)$. Indeed, assume that is not the case, i.e.\ for infinitely many $n$, the $\varepsilon$-neighborhoods of $\pi(x)$ and of $\pi(y)$ are disconnected within $B(\eta,\varepsilon)$ by some percolation loop $\gamma_n$. In particular, there are $0<t_1<t_2<1$ such that $\gamma_n$ crosses $\eta[t_1,t_2]$ and exits $B(\eta,\varepsilon)$. 

Without loss of generality (considering a subsequence if necessary), we may assume that $\gamma_n$ converges to some non-empty loop $\gamma$. The loop $\gamma$ must then cross $\eta$, which is a contradiction. Now, consider a sequence $\varepsilon_k \searrow 0$ and let $n_k$ be large enough so that there is a path $\eta_{n_k}$ from $x'_{n_k}(\varepsilon_k)$ to $y'_{n_k}(\varepsilon_k)$ within $B(\eta,\varepsilon_k)$. Considering a further subsequence if necessary, we may assume that $x'_{n_k},y'_{n_k}$ converge in $\tilde d^{\pth}_{Q}$ to some $x',y' \in \tilde\Lambda_Q$, respectively, and $\eta_{n_k}$ converges in $\Delta_{\cur}$ to $\tilde\eta$. We have $\pi(\tilde\eta) = \pi(\eta)$, and since $\eta$ is a simple path, it follows that $\tilde\eta=\eta$. 
By the first part of the proof, we see that $x = \phi(x')$, $y=\phi(y')$, and we have $\tilde{d}^{\pth}_{Q}(x',y') = \lim_{k\to\infty}\dpt_{Q,n_k}(x'_{n_k}(\varepsilon_k), y'_{n_k}(\varepsilon_k)) \leq\diam(\eta) = \dpt_Q(x,y)$.

It follows from the above that $\phi$ is surjective. The above also shows that $\phi$ is injective. Indeed, let $\tilde x,\tilde y\in\tilde\Lambda_Q$ and assume that $\phi(\tilde x)=\phi(\tilde y)$. Then $\tilde d^{\pth}_Q(\tilde x,\tilde y)\leq \dpt_Q(\phi(\tilde x),\phi(\tilde y))=0$, so $\tilde x=\tilde y$. Since we also showed that $\dpt(\phi(\tilde x),\phi(\tilde y))=\tilde d^{\pth}(\tilde x,\tilde y)$, the proof is complete.
\end{proof}

\addtocontents{toc}{
  \protect\contentsline{part}{\hyperref[sec:polynomial_tails]{Part 2}}{}{}
}

\section{Polynomial tail estimates}\label{sec:polynomial_tails}

\subsection{A priori estimate}\label{sec:apriori}
In this section, we give a proof of Proposition~\ref{prop:apriori}.

Recall the definition of the event $\mathcal{E}(\Lambda_n)$ from Section~\ref{subsubsec:apriori-sketch}. Before starting the proof of Proposition~\ref{prop:apriori}, let us show that the event $\mathcal{E}(\Lambda_n)$ occurs with probability bounded away from 0.
\begin{lem}
Fix $\delta\in(0,1/4)$. There exists a constant $c=c(\delta)>0$ such that for all $n\in \mathbf{N}$ with $\delta n\geq 2$, we have $\P(\mathcal{E}(\Lambda_n))\geq c$.\footnote{The restriction on $n$ is simply there to ensure that the set of sites of the lattice $\Lambda_n$ that are at most $\delta$ away from the $c_i^{(0)}$'s is non-empty.} \label{lem:E(Lambda_n)_nontrivial}
\end{lem}
\begin{proof} Fix $\delta$ and let $n\geq 2/\delta$. Define the event $\mathcal{E}_1(\Lambda_n)$ on $\Lambda_n$ as follows:
\begin{enumerate}
\item[$\bullet$] There is a closed left-to-right crossing inside $[0,1]\times\!\left[0,\frac{1}{6}\right]$;
\item[$\bullet$] There is a closed top-to-bottom crossing inside $\!\left[\frac{1}{2}-\delta,\frac{1}{2}+\delta\right]\times \!\left[0,\frac{1}{6}\right]$;
\item[$\bullet$] There is a closed circuit in the semiannulus $A(c_3^{(0)};\delta,2\delta)\cap \Lambda_n$;
\item[$\bullet$] There is an open left-to-right crossing inside $[0,1]\times (1/6,1/3]$.
\end{enumerate}
See Figure~\ref{fig:E(Lambda_n)_is_nontrivial} for illustration. By the standard RSW estimates and Harris' inequality, we see that there exists a constant $c_1(\delta)>0$ such that $\P(\mathcal{E}_1(\Lambda_n))\geq c_1$ for all $n\geq 2\delta$. On this event, one can explore the closed cluster adjacent to the bottom side of $\Lambda_n$; let $\mathcal{S}$ denote the set of explored sites - this includes the closed cluster adjacent to the bottom and all the neighboring sites which must be open. Condition on $\mathcal{S}=S$. The remaining percolation on $\Lambda_n \setminus S$ is still i.i.d.\ critical percolation. Define the event $\mathcal{E}_2(\Lambda_n,S)$ on $\Lambda_n\setminus S$ as follows:
\begin{enumerate}
\item[$\bullet$] There exist open circuits in semiannuli $A(c_2^{(0)};\delta,2\delta)\cap \Lambda_n$ and $A(c_2^{(0)};\delta,2\delta)\cap \Lambda_n$;
\item[$\bullet$] There exist open left-to-right paths in $[0,2\delta]\times\!\left[\frac{1}{2}-\delta,\frac{1}{2}+\delta\right]$ and $[1-2\delta,1]\times\!\left[\frac{1}{2}-\delta,\frac{1}{2}+\delta\right]$;
\item[$\bullet$] There exist open top-to-bottom paths in $[0,2\delta]\times\!\left[0,\frac{1}{2}+\delta\right] \setminus S$ and $[1-2\delta,1]\times\!\left[0,\frac{1}{2}+\delta\right]\setminus S$;
\item[$\bullet$] There is a closed top-to bottom path in $[1/2-\delta,1/2+\delta]\times[0,1] \setminus S$;
\item[$\bullet$] There is a closed left-to-right path in $[0,1]\times[5/6,1]$;
\item[$\bullet$] There exists an open circuit inside the semiannulus $A(c_1^{(0)}; \delta,2\delta)$.
\end{enumerate}
\begin{figure}[ht]
    \centering
    \includegraphics[scale=0.7]{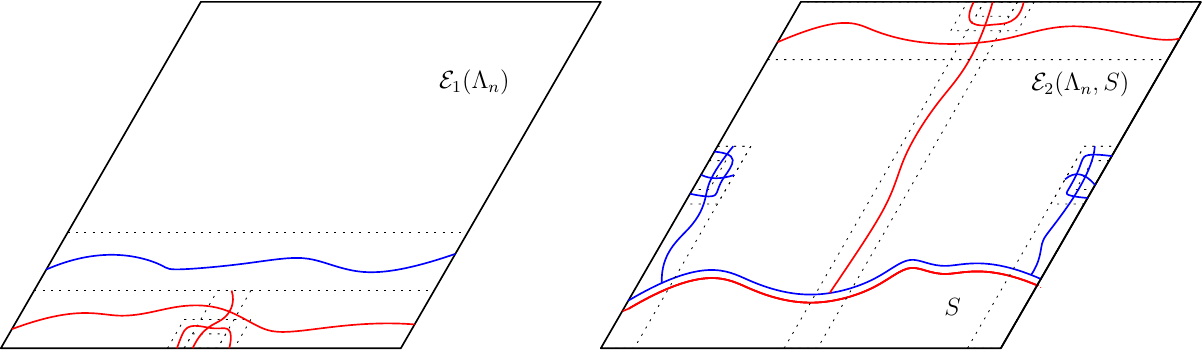}
    \caption{An illustration of events $\mathcal{E}_1(\Lambda_n)$ and $\mathcal{E}_2(\Lambda_n,S)$. The red and blue interface on the bottom of the second picture belong to the set $S$ - this is the boundary of the closed cluster adjacent to the bottom. }
    \label{fig:E(Lambda_n)_is_nontrivial}
\end{figure}
Note that for any possible $S$, $\P(\mathcal{E}_2(\Lambda_n,S))\geq \P(\mathcal{E}_2(\Lambda_n,\emptyset))$, and the latter is bounded away from zero by a constant $c_2=c_2(\delta)>0$, again by standard RSW estimates. Moreover, the event $\mathcal{E}_1(\Lambda_n)\cap\{\mathcal{S}=S\}\cap \mathcal{E}_2(\Lambda_n,S)$ implies $\mathcal{E}(\Lambda_n)$. Thus, because $\mathcal{E}_1(\Lambda_n)\cap\{\mathcal{S}=S\}$ can be verified on $S$, we have that 
\begin{equation}
\begin{split}
\P(\mathcal{E}(\Lambda_n)) &\geq \sum_{S} \P(\mathcal{E}_1(\Lambda_n)\cap \{\mathcal{S}=S\}\cap \mathcal{E}_2(\Lambda_n,S))\\  
&= \sum_S \P(\mathcal{E}_2(\Lambda_n,S)\,|\,\mathcal{E}_1(\Lambda_n)\cap\{\mathcal{S}=S\}) \cdot \P(\mathcal{E}_1(\Lambda_n)\cap \{\mathcal{S}=S\})\\
&= \sum_S \P(\mathcal{E}_2(\Lambda_n,S)) \cdot \P(\mathcal{E}_1(\Lambda_n)\cap\{\mathcal{S}=S\}) \geq c_2 \cdot \P(\mathcal{E}_1(\Lambda_n)) = c_2 \cdot c_1 >0,
\end{split}
\end{equation} as wanted.
\end{proof}

Due to the nature of the event $\mathcal{E}$, it will be useful to have a result concerning the probability that there exist arms in an annulus that hit inner and outer boundaries of this annulus at prescribed boundary segments. The following lemma, which follows immediately from \cite[Theorem 10]{nolin2008near}, addresses this situation. More precisely, it tells us that the probability that there are $4$ arms in an annulus such that each arms has a small extension and that their endpoints are at prescribed intervals at the inner and outer boundaries of our annulus is comparable to the probability that we simply have $4$ arms in this annulus when we consider critical percolation. Namely, let $b$ and $B$ be two boxes in the triangular grid of radii $r$ and $R$ respectively, such that $b \subseteq B$ and $\textrm{dist}(b,\partial B) \geq \frac{R}{3}$. Consider the (non-centered) annulus determined by $b$ and $B$. 

For $\delta$ as in the definition of the event $\mathcal{E}$, let $i_1,\ldots, i_4$ denote subintervals of the sides of box $b$, each of length $2\delta r$, centered at the midpoints of the respective sides and listed in clockwise order starting from the top side. Define $I_1,\ldots, I_4$ similarly for the box $B$. Suppose there exists a closed arm $\eta_1$ going from $i_1$ to $I_1$. Say $z_1$ and $Z_1$ are its endpoints on $b$ and $B$ respectively. We say that $\eta_1$ has a small extension on the outside, if:
\begin{enumerate}
\item[$\bullet$] There exists a closed horizontal crossing $C_1$ of the parallelogram $Z_1+[-\delta R,\delta R]\times [0,\sqrt{\delta}R]$
\item[$\bullet$] $C_1$ is connected to $\eta_1$ via a closed path fully contained within the box $B(Z_1,\sqrt{\delta}R)$.
\end{enumerate}
We can analogously define what it means for $\eta_2,\eta_3$ and $\eta_4$ to have small extensions on the outside, as well as small extensions on the inside. See Figure~\ref{fig:arms_with_extensions_in_the_annulus} below. 
\begin{figure}[ht]
    \centering
    \includegraphics[scale=0.7]{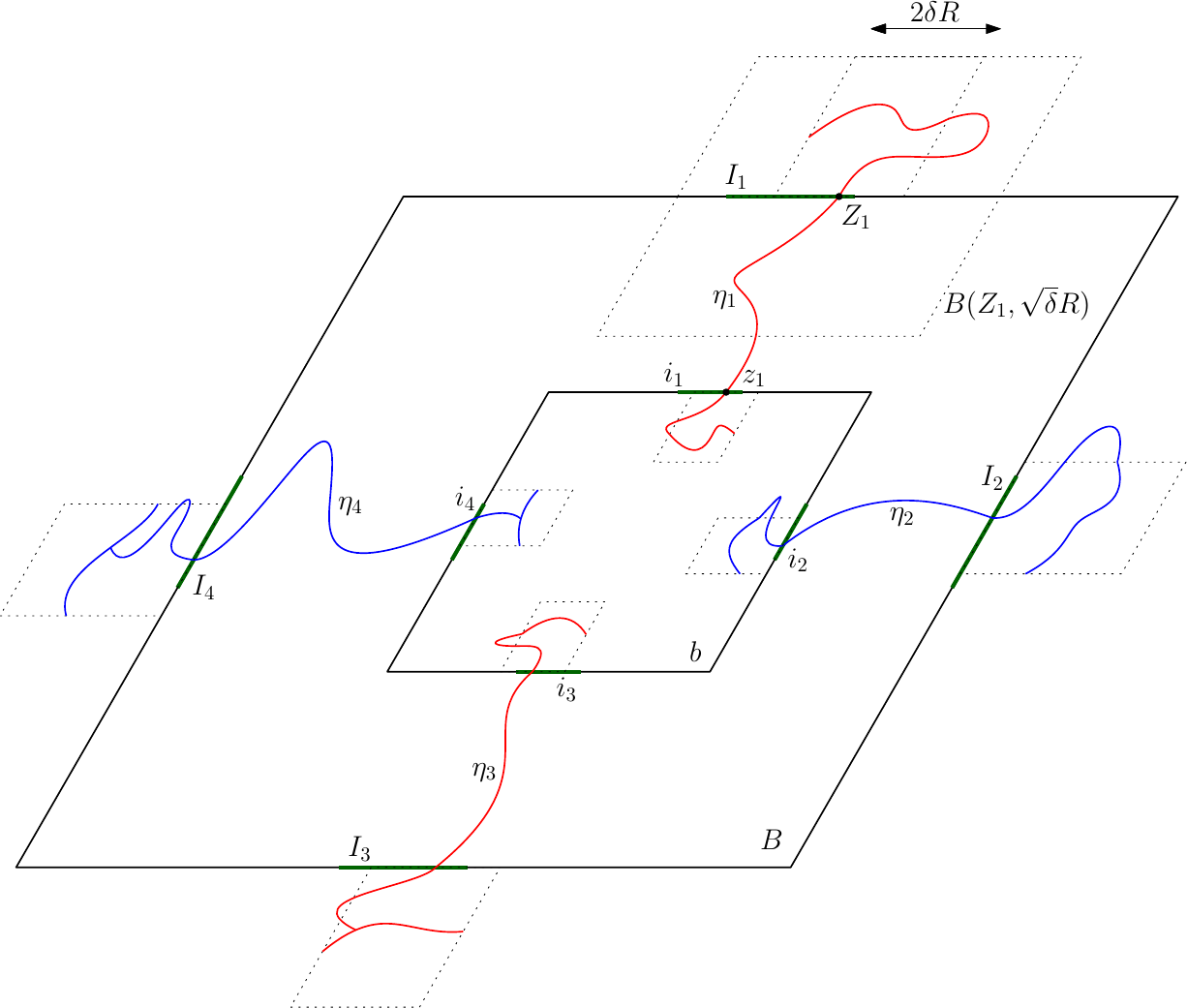}
    \caption{The figure above illustrates the event $\widehat{A_4}(b,B)$ for a (not necessarily centered) annulus $B\setminus b$. Green segments on the outer (resp.\ inner) boundaries represent the segments of length $2\delta R$ (resp.\ $2\delta r$) centered at the midpoints of the $4$ sides of this boundary. Red and blue curves represent the closed (resp.\ open) arms crossing this annulus. Note that for all $1\leq j\leq 4$, $\eta_j$ starts at $i_j$ and ends at $I_j$. Moreover, each $\eta_j$ has a small extension outside of $B$ and within $b$. For the extension of $\eta_1$ on the outside, we can see that there exists a closed left-to-right crossing of $Z_1+ [-\delta R,\delta R]\times [0,\sqrt{\delta}R]$, and this crossing is connected to $\eta_1$ within the larger box $B(Z_1,\sqrt{\delta}R)$ pictured above.}
    \label{fig:arms_with_extensions_in_the_annulus}
\end{figure}

Having this definition in mind, we will consider the following two events:
\begin{equation}
\begin{split}
A_4(b,B) &=\{\textrm{there are 4 arms of alternating colors between } \partial b \leftrightarrow\partial B\} \\
\widehat{A_4}(b,B) &= \!\left\{\begin{array}{cc}
\textrm{there are 4 arms } \eta_1,\ldots \eta_4 \textrm{ of alternating colors between }\partial b \leftrightarrow\partial B \\ \textrm{ s.t.\ } \eta_1 \textrm{ is closed, the endpoints of } \eta_a \textrm{ belong to } i_a,I_a \textrm{ for } a=1,\ldots,4, \\ \textrm{and all } \eta_a \textrm{ have small extensions both on the inside and on the outside } 
\end{array}\right\}.
\label{eq:defn_of_event_with_arm_extensions}
\end{split}
\end{equation}
In the remainder of the paper, we fix
\[
\delta= \frac{1}{10000} .
\]
For such $\delta$, \cite[Theorem 10]{nolin2008near} tells us that: 
\begin{lem}There exists a constant $C \in (0,1)$ such that for all $\delta^{-2}\leq r <\frac{R}{3}$ and boxes $b,B$ as above we have $$C\cdot \P(A_4(b,B)) \leq \P(\widehat{A_4}(b,B)) \leq \P(A_4(b,B)).$$ \label{lem:arms_with_prescribed_endpoints_comparable}
\end{lem}

Throughout the remaining of this section and Section~\ref{sec:superpolynomial-concentration}, we will fix some $R=R(\varepsilon)\in \mathbf{N}$ sufficiently large, so that the assumptions of all the necessary lemmas are satisfied.

\subsubsection{Proof of Proposition~\ref{prop:apriori}}
\medskip \noindent 
Recall the definition of $\cB^R_{k,n}$ from Section~\ref{subsec:notation}. Note that the set $\cB^R_{k,n}$ partitions the $(n+1)\times(n+1)$ graph $\Lambda_n$ into $R^{2(k-1)}$ identical $R^{\kappa-k}\times R^{\kappa-k}$ subgraphs, together with the strips $\Lambda_n\cap\!\left([n^{-1}R^{\kappa-1},1)\times[0,1)\right)$ and $\Lambda_n\cap\!\left([0,1)\times[n^{-1}R^{\kappa-1},1)\right)$. Consider a fixed box $B\in \cB^R_{k,n}$ and let $B^{(k-1)},\ldots,B^{(1)},B^{(0)}$\phantomsection\label{def:Bk-first} be its unique $R$-adic superboxes at levels $k-1,\ldots,1,0$.

In particular, $B^{(1)}=[0,n^{-1}R^{\kappa-1})^2$ and we define $B^{(0)}:=[0,1)^2$. 
We say that $B$ is \emph{away from the boundary} if each of its superboxes $B^{(i)}$ is at distance at least $\frac{1}{4n}R^{\kappa-i+1}$ away from the sides of the superbox $B^{(i-1)}$ at the next level for $i\geq 2$. In the following, we will only consider boxes $B\in \cB^R_{k,n}$ that are away from the boundary.

We say $B$ is $B^{(l)}$-\textit{important} with $l\geq1$ if there exist 4 arms $\eta_1^{(k)},\ldots,\eta_4^{(k)}$ of alternating colors from $\partial B$ to $\partial B^{(l)}$ such that
\begin{enumerate}
\item[$\bullet$] $\eta_1^{(k)}$ is a closed/red crossing;
\item[$\bullet$] for all $l\leq i\leq k$ and $0\leq j\leq 4$, $\eta_j^{(k)}$ passes through superbox $B^{(i)}$ within distance $\delta n^{-1}R^{\kappa-i}$ of $c_j^{(i)}$;
\item[$\bullet$] for all $l\leq i\leq k$ and $0\leq j\leq 4$, the (open/closed) cluster of $\eta_j^{(k)}$ does not come within distance $2\delta n^{-1}R^{\kappa-i}$ of $c_a^{(i)}$ for $a\neq j$. 
\end{enumerate}

\begin{figure}[ht]
    \centering
    \includegraphics[scale=0.7]{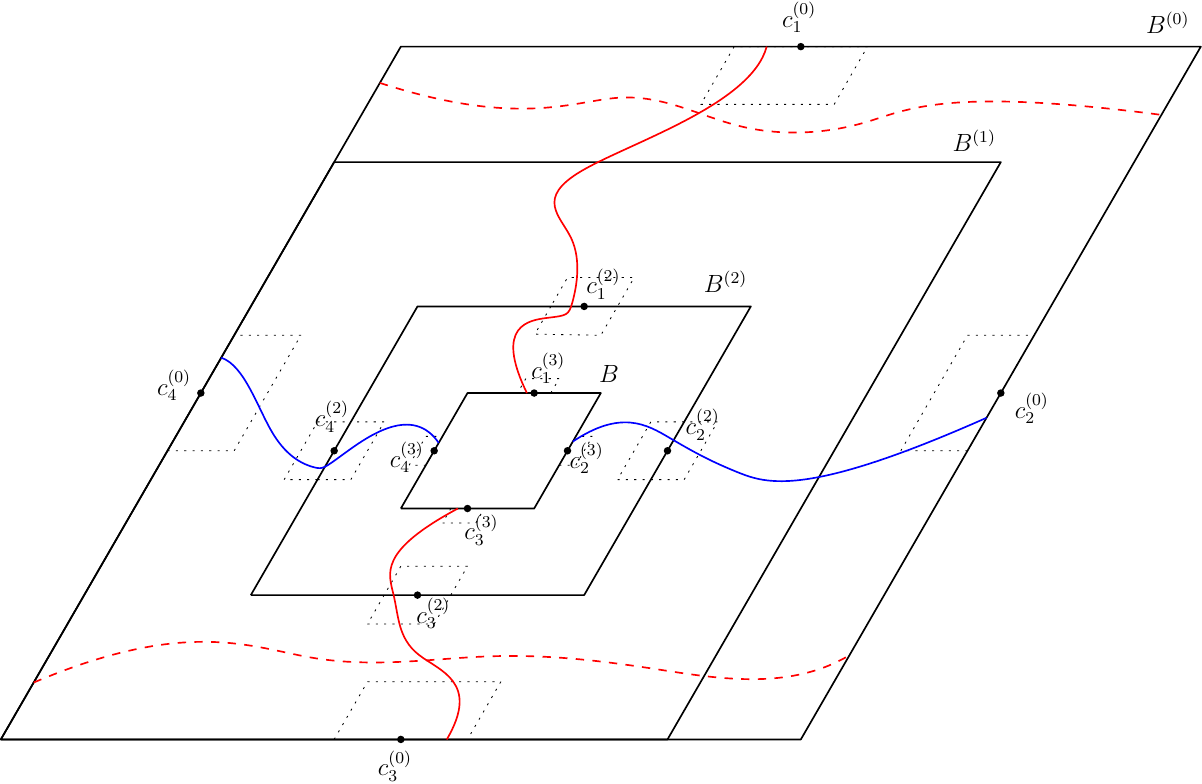}
    \caption{An illustration showing the event that a box $B$ of level $3$ is important. The blue and red curves represent the open (resp.\ closed) crossings of our annuli. Note that these curves come close to centers of all superboxes of $B$, except for $B^{(1)}$. The two dashed red curves represent the closed horizontal crossings in the outermost annulus. The two dashed red lines represent the closed horizontal crossings in the outermost annulus. Without them, the picture would represent the event that $B$ is $B^{(0)}$-important.}
    \label{fig:good_event}
\end{figure}

A box of level at least $2$ is $B^{(0)}$-important if the above hold for all $2\leq i\leq k$, and, furthermore, {every crossing} $\eta_j^{(k)}$ hits the boundary of $B^{(0)}=[0,1]^2$ within distance $\delta$ of $c_j^{(0)}$ and the cluster of $\eta_j^{(k)}$ does not come within distance $2\delta$ of $c_a^{(0)}$ for $a\neq j$. 
We define $B$ to simply be \emph{important} if it is $B^{(0)}$-important and if there exist two closed left-to-right crossings that stay within the strip of width $\frac{1}{4}$ from the top and bottom sides respectively. See Figure~\ref{fig:good_event}.

\medskip \noindent We say that $B$ is $B^{(l)}$-\textit{bad} if 
\begin{enumerate}
\item[$\bullet$] $B$ is $B^{(l)}$-important;
\item[$\bullet$] inside of $B$, the event $\mathcal{E}(B)$ occurs;
\item[$\bullet$] the top closed crossing inside $B$ is connected to $\eta_1^{(k)}$;
\item[$\bullet$] the bottom closed crossing inside $B$ is connected to $\eta_3^{(k)}$;
\item[$\bullet$] the open crossing in $B$ is connected to $\eta_2^{(k)}$ and $\eta_4^{(k)}$;
\item[$\bullet$] $ X_B\geq \altq_n$.
\end{enumerate}
We say the ``arms link up'' if the third to fifth events occur. Finally, say $B$ is \textit{bad} if in the above definition we substitute the condition that $B$ is $B^{(l)}$-important with the condition that $B$ is important.
\begin{figure}[ht]
    \centering
    \includegraphics[scale=0.7]{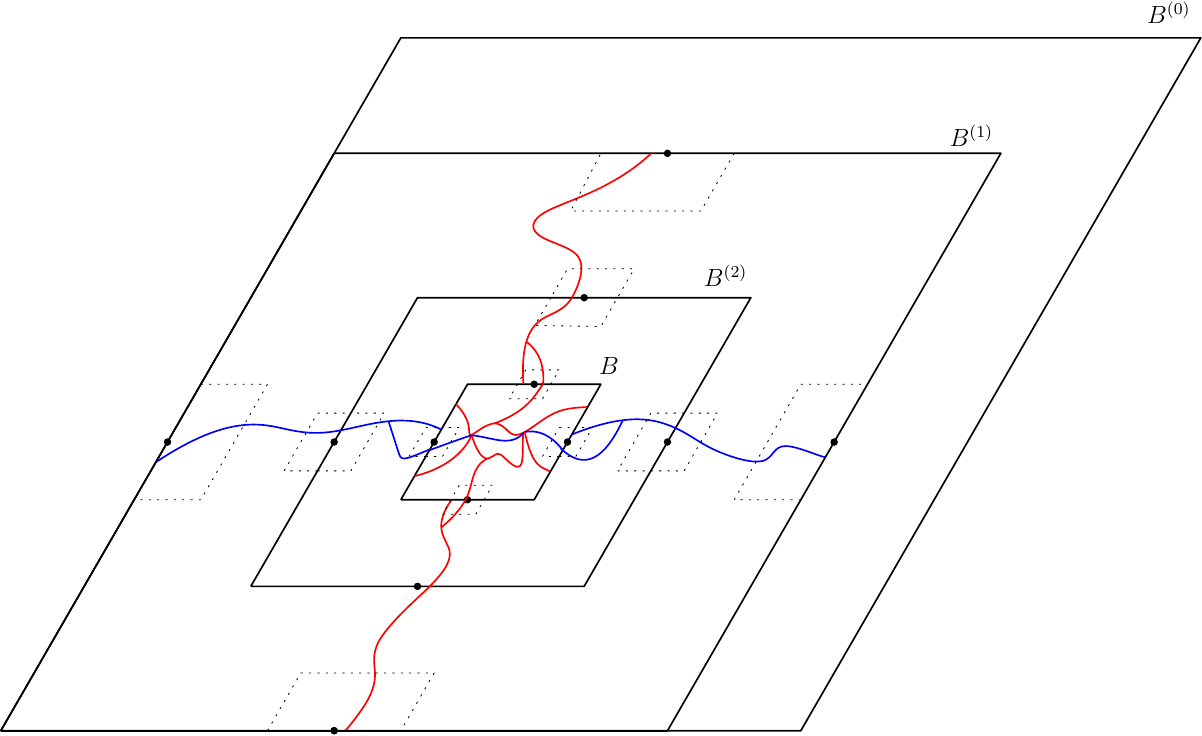}
    \caption{An illustration of the event that a box $B$ of level $3$ is $B^{(1)}$-bad. If we wanted $B$ to be $B^{(0)}$-bad, the crossings in the final annulus $B^{(1)}\setminus B^{(2)}$ should be replaced by crossings in $B^{(0)}\setminus B^{(2)}$, which now need not come close to midpoints of the sides of $B^{(1)}$, but, they must come close to midpoints of the sides of $B^{(0)}$. If in addition we wanted $B$ to be bad, we would need to include two red left-to-right crossings inside of $B^{(0)}$, as on Figure~\ref{fig:good_event}.}
    \label{fig:perfect_event}
\end{figure}

Let $T_k$ be the number of bad boxes at level $k\geq 2$. Note that if there exists a bad box $B$ at some level $k\geq 2$, then the event $\mathcal{E}\!\left(B^{(0)}\right)$ occurs and $X_{B^{(0)}}\geq X_B$. Indeed, the two extra closed crossings guarantee that the top and bottom closed arms are part of closed clusters that cross the box. The arms linking up implies that these clusters touch. Furthermore, every open path from $u\!\left(B^{(0)}\right)$ to $v\!\left(B^{(0)}\right)$ must pass through $u(B)$ and $v(B)$. Thus, by property~\eqref{property:serial} of our metric $D$, we know that $D(u,v)\geq D(u(B),v(B))$, as wanted. Thus, for every $k\geq1$,
\begin{equation}
\begin{split}
p\geq\P(X_n> \altq_n)\geq\P\!\left(X_{B^{(0)}}> \altq_n,\mathcal{E}\!\left(B^{(0)}\right)\right) &\geq \P(T_k>0) \geq \frac{\E[T_k]^2}{\E[T_k^2]}. \label{eqn:second_moment_bound}
\end{split}
\end{equation}
The goal is to obtain a lower bound for $\E[T_k]$ and an upper bound for $\E[T_k^2]$. We first prove some lemmas.

\begin{lem}\label{lem:4arms_polynomial_concentration}
Let $\zeta>0$. Then there exists $R(\zeta)$ such that for all $R\geq R(\zeta)$, there exists $c=c(R)$ such that for all $k\geq2$, $n> R^k$ and $B\in\cB^R_{k,n}$,
\begin{equation}
\begin{split}
\P(B\textrm{ is important})\geq cR^{-\!\left(\frac{5}{4}+\zeta\right)k}.
\end{split}
\end{equation}
\end{lem}

\begin{proof} We will subdivide the proof into four steps.

\textit{Step 1:} There exists $\tilde{R}(\zeta)$ such that for all $R\geq \tilde{R}(\zeta)$, for $B\in\cB^R_{2,n}$ and every $n> R^2$,  we have $\P(B \textrm{ is } B^{(1)}\textrm{-important})\geq R^{-\frac{5}{4}-\zeta}$.

\smallskip The main ingredient in this step is Lemma~\ref{lem:arms_with_prescribed_endpoints_comparable}. Let $B^{(1)}_\delta\subseteq B^{(1)}$ be the box obtained by removing 4 strips of width $\sqrt{\delta}n^{-1}R^{\kappa-1}$ from each side of $B^{(1)}$. Similarly, let $B^{(2)}_\delta\supseteq B^{(2)}$ be the box obtained by adding 4 strips of width $\sqrt\delta n^{-1}R^{\kappa-2}$ to each side of $B^{(2)}=B$. Let $c_{\delta,j}^{(1)},c_{\delta,j}^{(2)}$ for $j=1,\ldots,4$ be the midpoints of the sides of $B_\delta^{(1)}$ and $B_\delta^{(2)}$. Recall the definition~\eqref{eq:defn_of_event_with_arm_extensions} of the event $\widehat{A_4}(B_\delta^{(2)},B_\delta^{(1)})$ and assume that it holds. Suppose that $\eta_1$ is the closed arm that connects the top sides of the boxes and let $x$ be its endpoint on $B^{(1)}_{\delta,1}$. Then there exists a closed horizontal crossing $\gamma$ of $x+\!\left[-\frac{\delta}{2}n^{-1}R^{\kappa-1}, \frac{\delta}{2}n^{-1}R^{\kappa-1}\right] \times [0,\sqrt \delta n^{-1}R^{\kappa-1})$ that is connected to $\eta_1$ inside of $B\!\left(c_{\delta,1}^{(1)},2\sqrt \delta n^{-1}R^{\kappa-1}\right)$.  Consider the parallelogram $c_{\delta,1}^{(1)}+[-2\sqrt \delta n^{-1}R^{\kappa-1}, 2\sqrt \delta n^{-1}R^{\kappa-1}]\times [0,\sqrt \delta n^{-1}R^{\kappa-1})$. Suppose inside of it the following event occurs:
\begin{enumerate}
\item[$\bullet$] Subdivide $c_{\delta,1}^{(1)}+[-\delta n^{-1}R^{\kappa-1},\delta n^{-1}R^{\kappa-1}]\times [0,\sqrt \delta n^{-1}R^{\kappa-1})$ into 3 disjoint vertical parallelograms of dimensions $\frac{2}{3}\delta n^{-1}R^{\kappa-1}\times \sqrt \delta n^{-1}R^{\kappa-1}$. Assume each of them has a vertical closed crossing.
\item[$\bullet$] Assume that there is a closed circuit inside the semiannulus $A\!\left(c_{1}^{(1)};2\delta n^{-1}R^{\kappa-1},4\delta n^{-1}R^{\kappa-1}\right)\allowbreak \cap B^{(1)}$.
\end{enumerate}
Call this event $H\!\left(c_1^{(1)}\right)$. Note that by RSW and positive correlation the probability of this event is bounded away from zero by an absolute constant $c(\delta)$, depending only on $\delta$. Define $H\!\left(c_j^{(1)}\right),H\!\left(c_j^{(2)}\right)$ analogously. Since these $8$ events depend on disjoint sets of sites, hence are independent, their intersection, called $H$, has probability bounded away from zero by $c^8$. See Figure~\ref{fig:picture_for_lemma_5_1}.

\begin{figure}[ht]
    \centering
    \includegraphics[scale=0.77]{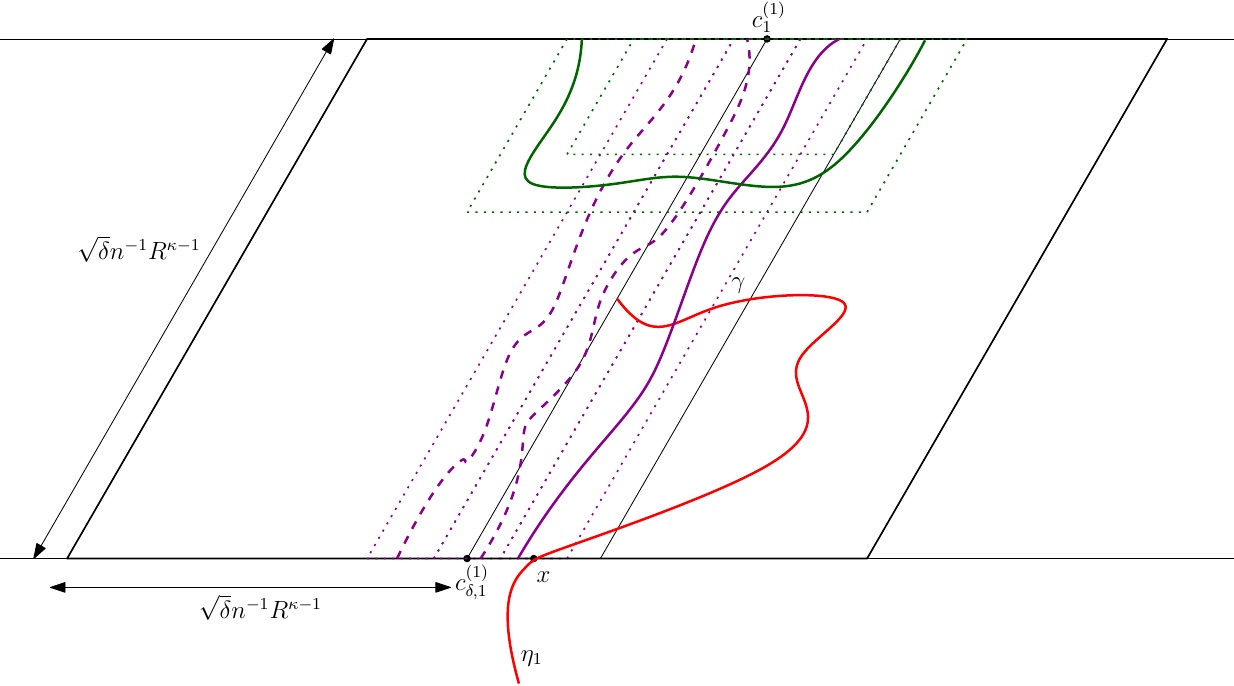}
    \caption{An illustration of events $H$ and $\widehat{A_{4}}\!\left(B_{\delta}^{(2)}, B_{\delta}^{(1)}\right)$ near the point $c_1^{(1)}$. The red curve represents the closed arm $\eta_1$ with an extension $\gamma$ coming from the event $\widehat{A_{4}}\!\left(B_{\delta}^{(2)}, B_{\delta}^{(1)}\right)$. The green curve is the closed circuit inside of the annulus $A(c_1^{(1)};\delta n^{-1}R^{\kappa-1},2\delta n^{-1}R^{\kappa-1})$, whose existence is guaranteed by $H$. $H$ also gives us the existence of the three purple curves, which represent the vertical closed crossings of the three parallelograms partitioning $c_{\delta,1}^{(1)}+[-\delta n^{-1}R^{\kappa-1},\delta n^{-1}R^{\kappa-1}]\times [0,\sqrt \delta n^{-1}R^{\kappa-1})$. The dashed ones do not intersect $\gamma$, whereas the full one does. Choosing $R$ big enough so that $R-\sqrt{\delta}>\frac{2}{3}R$ guarantees that at least one of these three crossings intersect $\gamma$. Moreover the closed circuit ensures that the open clusters corresponding to $\eta_2$ and $\eta_4$ do not come $2\delta n^{-1}R^{\kappa-1}$-close to $c^{(1)}_1$. So indeed, $H\cap \widehat{A_{4}}\!\left(B_{\delta}^{(2)}, B_{\delta}^{(1)}\right)$ implies that $B$ is $B^{(1)}$-important.}
    \label{fig:picture_for_lemma_5_1}
\end{figure}

Note that on the event $\widehat{A_4}\!\left(B_\delta^{(2)},B_\delta^{(1)}\right)\cap H$, $B^{(2)}$ is $B^{(1)}$-important. Moreover, conditional on the percolation configuration outside these $8$ parallelograms, the events $H$ and $\widehat{A_4}\!\left(B_\delta^{(2)},B_\delta^{(1)}\right)$ are positively correlated. Denote by $\mathcal{A}$ the set of sites in $B^{(1)}\setminus B^{(2)}$ that do not belong to any of the $8$ parallelograms constructed above, and let $\omega_{\mathcal{A}}$ denote the percolation configuration on $\mathcal{A}$. Then 
\begin{equation}
\begin{split}
 \P(B^{(2)} \textrm{ is } B^{(1)}\textrm{-important }| \, \omega_{\mathcal{A}})&\geq\P\!\left(\widehat{A_4}(B_\delta^{(2)},B_\delta^{(1)})\cap H \, | \, \omega_{\mathcal{A}}\right) \geq \P\!\left(\widehat{A_4}(B_\delta^{(2)},B_\delta^{(1)}) \, | \, \omega_{\mathcal{A}}\right) \cdot \P\!\left(H \, | \, \omega_{\mathcal{A}}\right) \\ &= \P\!\left(\widehat{A_4}(B_\delta^{(2)},B_\delta^{(1)}) \, | \, \omega_{\mathcal{A}}\right) \cdot \P\!\left(H\right) \geq c^8 \cdot \P\!\left(\widehat{A_4}(B_\delta^{(2)},B_\delta^{(1)}) \, | \, \omega_{\mathcal{A}}\right). \nonumber
\end{split}
\end{equation}
Summing over all possible configurations $\omega_\mathcal{A}$ and applying Lemmas~\ref{lem:arms_with_prescribed_endpoints_comparable} and~\ref{lem:arm_exponents}, we get that there exists $\tilde{R}=\tilde{R(\zeta)}$ such that for all $R>\tilde{R}$ we have \begin{equation}
\begin{split}
\P(B^{(2)} \textrm{ is } B^{(1)}\textrm{-important}) &\geq c^{8}\cdot \P\!\left(\widehat{A_4}(B_\delta^{(2)},B_\delta^{(1)}) \right) \geq c^8 C \cdot \P\!\left(A_4(B_\delta^{(2)},B_\delta^{(1)})\right)  \\ &\geq c^8 C \cdot  R^{-5/4-\zeta/2} \geq R^{-5/4-\zeta}. \nonumber
\end{split}
\end{equation}

\textit{Step 2:} There exists $\tilde{R}(\zeta)$ such that for all $R\geq \tilde{R}(\zeta)$, for all $k\geq 2$, $n> R^k$ and $B\in\cB^R_{k,n}$,  we have $\P\!\left(B \textrm{ is } B^{(1)}\textrm{-important}\right)\geq R^{-(\frac{5}{4}+\zeta)(k-1)}$.

\smallskip Conditional on the percolation configuration outside the $4(k+1)$ boxes $B\!\left(c_i^{(j)},2\delta n^{-1} R^{\kappa-j}\right) \cap (\Lambda_n\setminus B)$ for $j=1,\ldots, k$, $i=1,\ldots,4$, the event $\bigcap_{j=2}^{k} \!\left\{B^{(j)} \textrm{ is }B^{(j-1)}\textrm{-important}\right\}$ and the event that there exist $4(k+1)$ appropriately colored circuits in each $A\!\left(c_i^{(j)};\delta n^{-1}R^{\kappa-j},2\delta n^{-1}R^{\kappa-j}\right)$ are positively correlated. Moreover, these two events imply that $B$ is $B^{(1)}$-important. Hence, exactly as in Step 1 we get that 
\begin{equation}
\begin{split}
\P(B \textrm{ is } B^{(1)}\textrm{-important}) &\geq \P\!\left(\bigcap_{j=2}^{k} \!\left\{B^{(j)} \textrm{ is }B^{(j-1)}\textrm{-important}\right\} \bigg| \, \omega_\mathcal{A}\right)  \cdot \P(4(k+1) \textrm{ circuits}). \nonumber
\end{split}
\end{equation}
Note again by RSW estimates that the probability of each circuit existing is bounded away from zero by an absolute constant $c$. Hence, using independence and applying Step 1 $k-1$ times, we get that there exists $R(\zeta)$ big enough such that for all $R>R(\zeta)$,
\begin{equation*}
    \begin{split}
        \P(B \textrm{ is } B^{(1)}\textrm{-important})\geq{}& c^{4k}\cdot \prod_{j=2}^{k}\P\!\left( B^{(j)} \textrm{ is }B^{(j-1)}\textrm{-important}\right)  \\
        \geq{}& c^{4k}\cdot \!\left(R^{-5/4-\zeta/2}\right)^{k-1} \geq R^{-(5/4+\zeta)(k-1)}.
    \end{split}
\end{equation*}

\textit{Step 3:} There exists $\tilde{R}(\zeta)$ such that for all $R\geq \tilde{R}(\zeta)$, there exists $c=c(R)$ such that for all $k\geq 2$, $n> R^k$ and $B\in\cB^R_{k,n}$,  we have $\P\!\left(B \textrm{ is } B^{(0)}\textrm{-important}\right)\geq cR^{-(\frac{5}{4}+\zeta)k}$.

\smallskip An analogous proof to the one we have used in Step 1 shows that $\P(B^{(2)} \textrm{ is } B^{(0)}\textrm{-important})\geq R^{-5/4-\zeta}$. Moreover, an analogous argument to Step $2$ then yields that

\begin{equation}
    \P(B\textrm{ is }B^{(0)}\textrm{-important})\geq c_1R^{-(5/4+\zeta))(k-1)},
\end{equation} for some constant $c_1=c_1(R)$ depending only on $R$.

\textit{Step 4:} The same formulas hold for $B$ being important rather than $B^{(0)}$-important.

\smallskip To ensure that $B$ is important, we only need to add two horizontal closed crossings inside of $B^{(0)}$. We will add them in the strips of width $\frac{1}{4R}\leq \frac{1}{4}\sqrt{\delta}$ from the top and bottom sides of $B^{(0)}$ (instead of $1/4$ as in the definition) to ensure that this event has positive correlation with events $H\!\left(c_1^{(0)}\right),H\!\left(c_3^{(0)}\right)$ and is independent of $H\!\left(c_2^{(0)}\right),H\!\left(c_4^{(0)}\right)$ and $H\!\left(c_j^{(k)}\right)$ for every $k\geq 2, j=1,\ldots,4$. Since the probability of these two closed crossings existing is bounded away from 0 by a constant that depends only on $R$, we have $\P(B\textrm{ is important})\gtrsim_R\P(B\textrm{ is }B^{(0)}\textrm{-important})$, which completes the proof.
\end{proof}

\begin{lem}\label{lem:perfectbounds}
Let $B\in\cB^R_{k,n}$ and let $B^{(i)}$ be one of its superboxes. There exists a constant $c=c(\delta)>0$ such that for all $R$ sufficiently large, $k\geq 2$, $0\leq i<k$ and $n> R^k$,
\begin{equation}
    \begin{split}
        c\P(X_{R^{\kappa-k}} > \altq_n) \P(B\textrm{ is $B^{(i)}$-important})\leq{}& \P(B\textrm{ is $B^{(i)}$-bad})\\
        \leq{}&\P(X_{R^{\kappa-k}} > \altq_n) \P(B\textrm{ is $B^{(i)}$-important}) 
    \end{split}
\end{equation}
\end{lem}
\begin{proof}
By the definition of a bad box and independence of the percolation configuration inside and outside of $B$, we have
\begin{equation}\label{eq:perfectdecomp}
\begin{split}
 \P\!\left(B\textrm{ is $B^{(i)}$-bad}\right)&= \P(B\textrm{ is $B^{(i)}$-important})\P(\mathcal{E}(B))\P(X_B> \altq_n\mid \mathcal{E}(B))\\
 &\cdot\P(\textrm{arms link up}\mid \mathcal{E}(B),\,B\textrm{ is $B^{(i)}$-important},\,X_B> \altq_n).
\end{split}
\end{equation}
Lemma~\ref{lem:E(Lambda_n)_nontrivial} tells us that $\P(\mathcal{E}(B))$ is bounded uniformly from below over all $n$ and $B$ by a constant depending only on $\delta$ (take $R\geq 2/\delta$ so that the lemma applies). Furthermore, $\P(X_B> \altq_n\mid \mathcal{E}(B))=\P(X_{R^{\kappa-k}}> \altq_n)$. We now show that the last term in~\eqref{eq:perfectdecomp} is bounded from below by a universal constant.

Consider the box $B(c_j^{(k)},2\delta n^{-1}R^{\kappa-k})$ of radius $2\delta n^{-1}R^{\kappa-k}$ around $c_j^{(k)}$ and let $A(c_j^{(k)};\delta n^{-1}R^{\kappa-k},\allowbreak 2\delta n^{-1}R^{\kappa-k})$ be the annulus of radii $\delta n^{-1} R^{\kappa-k}$ and $2\delta n^{-1} R^{\kappa-k}$. As in the proof of the previous lemma, we can show that conditional on the percolation configuration outside the 4 boxes $B(c_j^{(k)};2\delta n^{-1}R^{\kappa-k})$, the event that the arms link up and the event $\mathcal{E}(B)\cap\{B\textrm{ is }B^{(i)}\textrm{-important}\}\cap\{X_B> \altq_n\}$ are positively correlated. Here we use that the cluster of the top closed crossing does not come within distance $2\delta n^{-1}R^{\kappa-k}$ of $c_2^{(k)},c_3^{(k)},c_4^{(k)}$ and the similar holds for the other two crossings. Since $X_B$ is the distance between two meeting points of the closed top and bottom clusters, this also implies that $X_B$ is independent of the percolation configurations inside the balls $B(c_j^{(k)},2\delta n^{-1}R^{\kappa-k})$. We conclude the proof by noting that the arms link up if there exist loops of the correct color inside each annulus $A(c_j^{(k)};\delta n^{-1}R^{\kappa-k},2\delta n^{-1}R^{\kappa-k})$. The probability of this event is bounded from below by a universal constant by RSW estimates.
\end{proof}

We now bound the second moment of $T_k$. Let $B_1,B_2\in\cB^R_{k,n}$ be $R$-adic boxes at level $k$, and let $B_{1,2}^{(d)}$ be their smallest common $R$-adic superbox, which is at some level $d\geq1$. Let $B_1^{(d+1)},B_2^{(d+1)}$ be the $(d+1)$-level subboxes of $B_{1,2}^{(d)}$ containing $B_1$ and $B_2$ respectively. Then
\begin{equation}\label{eq:uppertwoperfect}
\begin{split}
&\mathbf{P}(B_1,B_2 \textrm{ are bad})\\
&{=} 
\mathbf{P}(B_{1,2}^{(d)} \textrm{ is important}) \mathbf{P}(B_1,B_2 \textrm{ are bad } | B_{1,2}^{(d)} \textrm{ is important} )  \\
&\leq  \mathbf{P}(B_{1,2}^{(d)} \textrm{ is important}) \mathbf{P}(B_1 \textrm{ is }B_1^{(d+1)}\textrm{-bad, }  B_2 \textrm{ is }B_2^{(d+1)}\textrm{-bad} | B_{1,2}^{(d)} \textrm{ is important} )  \\
&= \mathbf{P}(B_{1,2}^{(d)} \textrm{ is important})\mathbf{P}(B_1 \textrm{ is }B_1^{(d+1)}\textrm{-bad })\mathbf{P}(  B_2 \textrm{ is }B_2^{(d+1)}\textrm{-bad} ). 
\end{split}
\end{equation}

Furthermore,
\begin{equation}
\begin{split}
\mathbf{P}(B_1 \textrm{ is bad} )={}& \mathbf{P}(B_1 \textrm{ is bad } | B_1 \textrm{ is } B_1^{(d+1)}\textrm{-bad, } B_1^{(d+1)} \textrm{ is important}) \\ &\cdot \mathbf{P}(B_1 \textrm{ is } B_1^{(d+1)}\textrm{-bad})\mathbf{P}(B_1^{(d+1)} \textrm{ is important})  \\
={}& \mathbf{P}(\textrm{arms link up} | B_1 \textrm{ is } B_1^{(d+1)}\textrm{-bad, } B_1^{(d+1)} \textrm{ is important}) \\ &\cdot \mathbf{P}(B_1 \textrm{ is } B_1^{(d+1)}\textrm{-bad})\mathbf{P}(B_1^{(d+1)} \textrm{ is important}).
\end{split}
\end{equation}
For the first equality, we use that the percolation configurations inside and outside $B_1^{(d+1)}$ are independent. Now by the exact same argument as in the proof of Lemma~\ref{lem:perfectbounds}, we obtain
\begin{equation*}\label{eq:lowerperfect}
    \begin{split}
        \P(B_1\text{ is bad})\geq&\,C_0 \mathbf{P}(B_1 \textrm{ is } B_1^{(d+1)}\textrm{-bad})\mathbf{P}(B_1^{(d+1)} \textrm{ is important})\\
        \geq& \,C_0 \cdot C(R) \cdot 
\mathbf{P}(B_1 \textrm{ is } B_1^{(d+1)}\textrm{-bad})\mathbf{P}(B_1^{(d)} \textrm{ is important})
    \end{split}
\end{equation*}
for some absolute constant $C_0$ and a constant $C(R)$ depending only on $R$. The last inequality follows from Lemma~\ref{lem:4arms_polynomial_concentration} combined with another application of the argument in the proof of Lemma~\ref{lem:perfectbounds}.

It follows from equations~\eqref{eq:uppertwoperfect} and~\eqref{eq:lowerperfect} and Lemma~\ref{lem:perfectbounds} that there exist $C=C(R)>0$ and $C'=C'(R)>0$ such that for all $k,n$,
\begin{equation}
\begin{split}
\mathbf{E}[T_k^2] ={}& \mathbf{E} \!\left[ \!\left(\sum_{i} \mathbf{1}_{B_i \textrm{ is bad}}\right)^2 \right] \\
={}&\sum_{i\neq j} \mathbf{P}(B_i,B_j \textrm{ are bad}) + \sum_{i} \mathbf{P}(B_i \textrm{ is bad}) \\
\leq{}& \sum_{i\neq j} \mathbf{P}(B_{i,j}^{(d)} \textrm{ is important}) \mathbf{P}\!\left(B_i \textrm{ is } B_i^{(d+1)}\textrm{-bad}\right)\mathbf{P}\!\left(B_j \textrm{ is } B_j^{(d+1)}\textrm{-bad}\right)\\
&+\sum_i\P\!\left(B_i \textrm{ is important}\right)\P\!\left(X_{R^{-k}n}> \altq_n\right) \\
\leq{}& \sum_{i\neq j} \mathbf{P}(B_{i,j}^{(d)} \textrm{ is important}) \mathbf{P}\!\left(B_i \textrm{ is } B_i^{(d+1)}\textrm{-important}\right)\mathbf{P}\!\left(B_j \textrm{ is } B_j^{(d+1)}\textrm{-important}\right)\\
&\cdot \P(X_{R^{-k}n}> \altq_n)^2 +\sum_i\P\!\left(B_i \textrm{ is important}\right)\P\!\left(X_{R^{-k}n}> \altq_n\right) \label{eqn:upper_bound_on_E[T_k^2]}
\end{split}
\end{equation}
and
\begin{equation}
\begin{split}
\mathbf{E}[T_k]^2 ={}& \!\left(\sum_i \mathbf{P}(B_i \textrm{ is bad})\right)^2 \\
={}&\sum_{i\neq j} \mathbf{P}(B_i \textrm{ is bad})\mathbf{P}(B_j \textrm{ is bad}) + \sum_{i} \mathbf{P}(B_i \textrm{ is bad})^2  \\
\geq{}& C\sum_{i \neq j} \mathbf{P}(B_i \textrm{ is } B_i^{(d+1)}\textrm{-bad})\mathbf{P}(B_i^{(d)} \textrm{ is important})\mathbf{P}(B_j \textrm{ is } B_j^{(d+1)}\textrm{-bad})\mathbf{P}(B_j^{(d)} \textrm{ is important})\\
\geq{}& C'\sum_{i\neq j} \mathbf{P}(B_{i,j}^{(d)} \textrm{ is important})^2 \mathbf{P}\!\left(B_i \textrm{ is } B_i^{(d+1)}\textrm{-important}\right)\mathbf{P}\!\left(B_j \textrm{ is } B_j^{(d+1)}\textrm{-important}\right)\\
&\cdot\P\!\left(X_{R^{-k}n}> \altq_n\right)^2. \label{eqn:lower_bound_on_E[T_k]^2}
\end{split}
\end{equation}
In the above, we sum over all $R$-adic subboxes $B_i$ and $B_j$ at level $k$ that are away from the boundary. These equations give us an upper bound on $\E[T_k^2]$ and a lower bound on $\E[T_k]^2$. The goal is now to substitute these bounds into~\eqref{eqn:second_moment_bound} and obtain a bound on $\P(X_{R^{\kappa-k}}> \altq_n)$. In order to do so, we will need the following lemma telling us that two sums in~\eqref{eqn:upper_bound_on_E[T_k^2]} and~\eqref{eqn:lower_bound_on_E[T_k]^2} are comparable.

\begin{lem}\label{lem:comparison-sums}
There exists a constant $C=C(R)>0$ such that for all $k\geq 2$ and $n>R^k$,
\begin{equation}
\begin{split}
&\sum_{i\neq j} \mathbf{P}(B_{i,j}^{(d)} \textrm{ is important})^2 \mathbf{P}\!\left(B_i \textrm{ is } B_i^{(d+1)}\textrm{-important}\right)\mathbf{P}\!\left(B_j \textrm{ is } B_j^{(d+1)}\textrm{-important}\right)\\
&\geq C\sum_{i\neq j} \mathbf{P}(B_{i,j}^{(d)} \textrm{ is important}) \mathbf{P}\!\left(B_i \textrm{ is } B_i^{(d+1)}\textrm{-important}\right)\mathbf{P}\!\left(B_j \textrm{ is } B_j^{(d+1)}\textrm{-important}\right).
\end{split}
\end{equation}
\end{lem}
\begin{proof} Note first that if $d\geq 2$, then
\begin{equation}
\begin{split}
&\mathbf{P}\!\left( B_i \textrm{ is important}\right)\\
={}& \mathbf{P}\!\left( B_i \textrm{ is } B_i^{(d+1)} \textrm{-important}\right) \cdot \mathbf{P}\!\left( B_{i,j}^{(d)} \textrm{ is }\textrm{important}\right) \cdot  \mathbf{P}\!\left( B_i^{(d+1)} \textrm{ is } B_i^{(d)} \textrm{-important}\right)\\&\cdot \mathbf{P}\!\left( \textrm{arms link up } \big| B_i \textrm{ is } B_i^{(d+1)} \textrm{-important, }B_{i,j}^{(d)} \textrm{ is }\textrm{important, } B_i^{(d+1)} \textrm{ is } B_i^{(d)} \textrm{-important}  \right). \label{eqn:subpolynomial_concentration_equivalence}
\end{split}
\end{equation}

By Lemma~\ref{lem:4arms_polynomial_concentration} and Lemma~\ref{lem:arm_exponents}, we know that
\begin{equation}
    R^{-5/4-\zeta} \leq \mathbf{P}\!\left( B_i^{(d+1)} \textrm{ is } B_i^{(d)} \textrm{-important}\right) \leq R^{-5/4+\zeta}.
\end{equation}
Moreover, an argument similar to the one in Lemma~\ref{lem:perfectbounds} gives us that the linking probability is bounded away from zero by an absolute constant. Thus, for $d\geq 2$ we have 
\begin{equation}
    \P(B_i \textrm{ is important}) \asymp_R \P(B_i \textrm{ is } B_{i}^{(d+1)}\textrm{-important})\cdot\P(B_{i,j}^{(d)} \textrm{ is important}).
\end{equation}
Using Lemma~\ref{lem:4arms_polynomial_concentration}, we can extend this formula to the case $d=1$. 

With the above formula in hand, we can see that in order to show the desired inequality, it suffices to show that there exists a constant $C'=C'(R)$ such that for all $n,k$
\begin{equation}
\begin{split}
&\sum_{i\neq j} \mathbf{P}\!\left( B_i \textrm{ is important}\right)\mathbf{P}\!\left( B_j \textrm{ is important}\right)\\
\geq{}& C' \cdot \sum_{i\neq j} \mathbf{P}\!\left( B_{i,j}^{(d)} \textrm{ is important} \right)^{-1} \mathbf{P}\!\left( B_i \textrm{ is important}\right)\mathbf{P}\!\left( B_j \textrm{ is important}\right). 
\end{split}
\end{equation}
Denote by $P(i)=\mathbf{P}(B_i \textrm{ is important})$ and let $l(i,j)=d$ if the smallest common superbox of $B_i$ and $B_j$ is of $R$-adic level $d$. Again, using Lemma~\ref{lem:4arms_polynomial_concentration} to estimate $\mathbf{P}\!\left(B_{i,j}^{(d)} \textrm{ is important}\right) \geq R^{-(5/4+\zeta)d}$ for all $0\leq d< k$, we note that it suffices to show that we can find a constant $C'$ such that
\begin{equation}
\sum_{i\neq j} P(i) P(j) \geq C' \cdot \sum_{0\leq d <k}\sum_{\substack{i\neq j\\l(i,j)=d}} R^{(5/4+\zeta)d} P(i) P(j).
\end{equation}
We will now show that this is true. For each box $B_j$ at level $k$, let $\overrightarrow{n}(j)=(n_1,\ldots, n_k)$ denote the vector of relative positions of the box $B_{j}^{(d+1)}$ within $B_{j}^{(d)}$, for $1\leq d<k$. Fix $d\in \{1,\ldots k\}$. Note that for every $B_j$, there exist at most $(R/2)^{2d}$ boxes $B_{j'}$ that are away from the boundary and such that $\overrightarrow{n}(j)$ and $\overrightarrow{n}(j')$ differ only in the first $d$ coordinates. 

Moreover, we claim that if $j$ and $j'$ do share the same last $k-d$ relative coordinates, then \begin{equation}P(j)\leq R^{2\zeta d} c^{-d} \label{eq:boundonlastk-dcoordssame}P(j'),\end{equation} where $c$ is the linking constant coming from Lemma~\ref{lem:perfectbounds}. To see this, note that the event that $B_j$ is important is simply the product of probabilities that each of the $k-1$ (not necessarily centered) annuli surrounding $B_j$ has the important 4-arm event and that we have $k-2$ linkings at each border of consecutive annuli given the 4-arm events. Because the last $k-d$ coordinates of $B_j$ and $B_{j'}$ match, we know that the crossing events of those annuli have exactly the same probability, as well as linking probabilities at those levels. Thus the only difference comes from crossing events for annuli at levels $0,2,3,\ldots, d$ and the corresponding linking events. By Lemma~\ref{lem:4arms_polynomial_concentration} and Lemma~\ref{lem:arm_exponents}, we know that crossing probabilities at each level can differ by at most a factor of $R^{2\zeta}$, hence the factor $R^{-2\zeta d}$ in the bound above. Finally, by Lemma~\ref{lem:perfectbounds} , we know that the linking probabilities for these annuli are bounded away from zero by an absolute constant $c>0$, hence the factor $c^{-d}$. This shows our claim.

Summing~\eqref{eq:boundonlastk-dcoordssame} over all suitable $j'$ we find that \begin{equation}
\!\left(\frac{R}{2}\right)^{2d} P(j) \leq R^{2\zeta d} c^{-d} \sum_{j'}P(j').
\end{equation}
Thus
\begin{equation}
\begin{split}
\sum_{\substack{i\neq j\\l(i,j)=d}} R^{(5/4+\zeta)d} P(i) P(j) &=  R^{(5/4+\zeta)d} \sum_i P(i) \sum_{\substack{j\neq i \\ l(i,j)=d }} P(j) \\
&\leq R^{(5/4+\zeta)d} \sum_i \!\left[ P(i) \sum_{\substack{j\neq i \\ l(i,j)=d }}  \!\left(\frac{R}{2}\right)^{-2d}\cdot R^{2\zeta d} c^{-d} \sum_{j'}P(j') \right] \\
&\leq R^{-(3/4-3\zeta)d} \cdot 2^{2d} \cdot c^{-d} \sum_{i\neq j} P(i) P(j),
\end{split}
\end{equation}
where the sum over $j'$ is taken over all $j'$ whose coordinates differ from coordinates of $j$ in only possibly the first $d$ places. The last inequality then follows since $j'$ (possibly) differs from $j$ in the first $d$ coordinates, while $j$ and $i$ share these $d$ coordinates and differ in at least one remaining coordinate. Thus this sum can be upper bounded by the sum over all positions differing from $i$. Summing the above inequality over $1\leq d <k$ we get that
\begin{equation}
\begin{split}
\sum_{1\leq d <k}\sum_{\substack{i\neq j\\l(i,j)=d}} R^{(5/4+\zeta)d} P(i) P(j) &\leq \sum_{1\leq d<k} R^{-(3/4-3\zeta)d} \cdot 2^{2d} \cdot c^{-d} \sum_{i\neq j} P(i) P(j) \\
&= \sum_{i\neq j} P(i) P(j) \cdot \sum_{1\leq d<k} \!\left(\frac{4}{c}\cdot R^{-(3/4-3\zeta)}\right)^d \\
&\leq \sum_{i\neq j} P(i) P(j) \cdot \sum_{d=1}^{\infty} \!\left(\frac{4}{c}\cdot R^{-(3/4-3\zeta)}\right)^d.
\end{split}
\end{equation}
For $R$ big enough (e.g.\ bigger than the absolute constant $(8/c)^2$) and for $\zeta<1/12$, the series on the RHS is summable (and in particular less than 2), hence the claim follows.
\end{proof}

We are now ready to complete the proof of Proposition~\ref{prop:apriori}.
\begin{proof}[Proof of Proposition~\ref{prop:apriori}]
Let
\begin{equation}
    f_R(k):=\sum_{i\neq j} \mathbf{P}(B^{(d)} \textrm{ is important}) \mathbf{P}\!\left(B_i \textrm{ is } B_i^{(d+1)}\textrm{-important}\right)\mathbf{P}\!\left(B_j \textrm{ is } B_j^{(d+1)}\textrm{-important}\right) \label{eqn:defn_of_f_R}
\end{equation}
and
\begin{equation}
    g_R(k):=\sum_i\P\!\left(B_i \textrm{ is important}\right), \label{eqn:defn_of_g_R}
\end{equation}
again summing over all boxes $i$ and $j$ at level $k$ that are away from the boundary. Then substituting~\eqref{eqn:upper_bound_on_E[T_k^2]} and~\eqref{eqn:lower_bound_on_E[T_k]^2} into~\eqref{eqn:second_moment_bound} and using Lemma~\ref{lem:comparison-sums}, we obtain that there exists $c=c(R)$ such that for all $k\geq2$ and $n>R^k$,
\begin{equation}
\begin{split}
p\geq\frac{\E[T_k]^2}{\E[T_k^2]}\geq\frac{cf_R(k)\P\!\left(X_{R^{\kappa-k}}> \altq_n\right)^2}{f_R(k)\P\!\left(X_{R^{\kappa-k}}> \altq_n\right)^2+g_R(k)\P\!\left(X_{R^{\kappa-k}}> \altq_n\right)}. \label{eqn:lower_bound_for_second_moment_bound}
\end{split}
\end{equation}
Rearranging this equation gives
\begin{equation}
\begin{split}
\P\!\left(X_{R^{\kappa-k}}> \altq_n\right)\leq\frac{p}{c-p}g_R(k)f_R(k)^{-1},
\end{split}
\end{equation}
provided $p<c$. Note that $f_R(k)\geq g_R(k)^2-g_R(k)$. Indeed, if $\tilde{T}_k$ denotes the number of \textit{important} boxes of level $k$, then $g_R(k)=\mathbf{E}[\tilde{T}_k]$ and by the same argument as in~\eqref{eqn:upper_bound_on_E[T_k^2]} we get that $g_R(k)^2=\E[\tilde{T}_k]^2\leq\E[\tilde{T}_k^2]\leq f_R(k)+g_R(k).$

\noindent Finally, the number of $R$-adic boxes at level $k$ that are away from the boundary is $\!\left(\frac{R}{2}\right)^{2(k-1)}$. By Lemma~\ref{lem:4arms_polynomial_concentration}, we know that $g_R(k)\gtrsim\!\left(\frac{R}{2}\right)^{2(k-1)}R^{\!\left(-\frac{5}{4}-\zeta\right)k}\geq 2$ for all $\zeta<1/4$ and $R\geq \max\{R(\zeta),16\}$. Therefore, for every $\varepsilon>0$, there exists $\zeta>0$, $R=R(\zeta,\varepsilon)$ and $c=c(R)$ such that for all $k\geq 2$, $n>R^k$ and $p<c$, 
\begin{equation}
\begin{split}
\P\!\left(X_{R^{\kappa-k}}> \altq_n\right)\leq \frac{p}{c-p}\cdot \frac{1}{g_R(k)-1}\leq\frac{2p}{c-p}\!\left(\frac{R}{2}\right)^{-2k}R^{\!\left(\frac{5}{4}+\zeta\right)k}\leq\frac{2p}{c-p}R^{-\!\left(\frac{3}{4}-\varepsilon\right)k}, \label{eqn:final_bound_on_g_R/f_R}
\end{split}
\end{equation}
which completes the proof. Note that the value of $c$ may change from line to line but only depends on $R$. For the case when $k=1$, we may simply assume that the constant $C$ in the statement of the proposition is such that $C R^{-\!\left(\frac{3}{4}-\varepsilon\right)}\geq 1$.  
\end{proof}

\begin{rem}
In the last step of the proof above, any $p<c$ suffices. Note that $c$ depends on $R$. Later on, we choose $p$ to be some arbitrarily small number. This additional restriction is needed later in the proof of Lemma~\ref{lem:compXscales}. Furthermore, since $R$ will be chosen to be some sufficiently large constant, we will have to choose $p$ sufficiently small accordingly.
\end{rem}

\subsection{Annulus crossing bounds}\label{subsec:distance-crossings}
The goal of this section is to prove Proposition~\ref{lem:distancebetweencrossings}, which is a more precise version of Proposition~\ref{prop:annulus-crossing}.

Throughout the section, we fix some $\alpha>\frac{3}{4}-{\varepsilon_0}$ and assume that the bound~\eqref{eq:assumptionsubpoly} holds. {To be precise, we assume that for all sufficiently large $R\in\N$, sufficiently small $p$ (depending on $R$), there exists $C_0=C_0(\alpha,R)$ such that for all $k\in\N$ and $n> R^k$, we have
\begin{equation}\label{eq:assumptionsubpolyrepeated}
    \P(X_{R^{\kappa-k}}> \altq_n)\leq C_0R^{-\alpha k}.
\end{equation}
}
At the moment, we know that this assumption holds for any $\alpha<\frac{3}{4}$ by Proposition~\ref{prop:apriori}. In Section~\ref{sec:superpolynomial-concentration} below we will bootstrap the a priori estimate to obtain superpolynomial concentration. In turn, we will be able to use the results of this section for general $\alpha$. 

First, we need the following lemma that is a consequence of the assumption~\eqref{eq:assumptionsubpoly}.
\begin{lem}\label{lem:compXscales}
For all $\zeta>0$, there exist $R=R(\zeta)$ big enough, $p$ small enough (depending on $R$) and $C=C(R,\zeta,\alpha)>0$ such that for all $0< s\leq 1$, $k\geq1$ and $n> R^k$, 
we have
\begin{equation}
    \P(X_{R^{\kappa-k}}> {s^{\frac{3}{4}-\zeta}}\altq_n)\leq C {s^{-\alpha}}R^{-\alpha k}.
\end{equation}
\end{lem}

\begin{proof}
Fix $\zeta>0$. We first show that there exist $R=R(\zeta)$ large enough, $k_0=k_0(R)\in \mathbf{N}$ and {$\tilde{c}=\tilde{c}(\zeta,R)$} 
such that for all $k\geq k_0$ and $n> R^{k}$, we have
\begin{equation}\label{eq:comp_mn_scales}
    \frac{\altq_{n}}{\altq_{R^{\kappa-k}}}\geq {\tilde{c}R^{(\frac{3}{4}-\zeta)k}}.
\end{equation}

To that end, it suffices to show that
\begin{equation}
    \P\!\left(X_n\leq {\tilde{c}R^{(\frac{3}{4}-\zeta)k}}\altq_{R^{\kappa-k}}\right)\leq 1-p\leq\P(X_n\leq \altq_n)
\end{equation}
for {$R,k$ sufficiently large and $c$ sufficiently small}.

Fix some ${R}={R}(\zeta)$ so that all the lemmas in Section~\ref{sec:apriori} hold. Let ${\kappa}$ be such that ${R}^{{\kappa}-1}<n\leq {R}^{{\kappa}}$. Recall the definition of a bad (${R}$-adic) box given at the beginning of Section~\ref{sec:apriori}. Say that an ${R}$-adic box $B\subseteq B^{(1)}$ at level ${k}$ is \textit{poor} if it satisfies all the conditions for a box to be bad, except that instead of requiring $X_B> \altq_n$, we require $X_B> \altq_{{R}^{{\kappa}-{k}}}-\frac{1}{12}$. 
Let $T_{{k}}'$ be the number of poor $R$-adic boxes at level ${k}$ in $B^{(1)}$. 

Note first that for both the geodesic and resistance metric on any subgraph of $\bT$, the distance between any two distinct vertices is at least $\frac{1}{6}$ by the parallel law. So either $X_{n}\geq\frac{1}{6}$ or $X_n=0$. The latter is true if there is only one pivotal site for the event that there exists an open left-to-right crossing. This pivotal site then has 6 arms (4 open and 2 closed) to the boundary of $\Lambda_n$. Moreover, all the arms have length at least $\delta$ due to the definition of $\cE$. Since the 6-arm exponent is larger than 2, the probability of such a point existing tends to 0 as $n\rightarrow\infty$, and so it eventually must be smaller than $1/2$. Thus, there exists\footnote{We also choose this $n_0$ large enough in terms of $\delta$ so that $\altq_n(p)$ is well-defined for all $n\geq n_0$ -- it suffices to assume $n_0\geq 10/\delta$ to ensure that $\P(\cE(\Lambda_n))>0$.} $n_0=n_0(\delta)$ such that for any $n\geq n_0$ and $0<p<1/2$ we have that $\altq_n(p)\geq\frac{1}{6}$.

It follows that $\altq_{{R}^{{\kappa}-{k}}}-{\frac{1}{12}}\geq \frac{1}{2}\altq_{{R}^{{\kappa}-{k}}}$ for  all $n>{R}^{{k}}$, since then ${R}^{{\kappa}-{k}}\geq {R}$. Thus, for $\epsilon>0$ to be chosen later, we have by Paley-Zygmund that 
\begin{equation}
\begin{split}
    \P\!\left(X_n\leq \epsilon {R}^{\!\left(\frac{3}{4}-\zeta\right){k}}\altq_{{R}^{{\kappa}-{k}}}\right)\leq&\, \P \!\left( T_{{k}}'\leq 2\epsilon {R}^{\!\left(\frac{3}{4}-\zeta\right){k}}\,\big|\,\cE\right)\\
    =&\, \P\!\left(T_{{k}}'\leq \theta_{{k}}\E [T_{{k}}']\mid\cE\right)\\
    \leq&\, 1-\P\!\left(T_{{k}}'>\theta_{{k}}\E[T_{{k}}'],\,\cE\right)\\
    \leq&\, 1-(1-\theta_{{k}})^2 \frac{\E[T_{{k}}']^2}{\E[T_{{k}}'^{2}]},
\end{split}
\end{equation}
where $\theta_{{k}}=\frac{2\epsilon {R}^{\!\left(\frac{3}{4}-\zeta\right){k}}}{\E[T_{{k}}']}$ and in the last line we use that $T'_{{k}}>0$ implies the event $\cE$. Recall the definitions of $f_{{R}}$ and $g_{{R}}$ in~\eqref{eqn:defn_of_f_R} and~\eqref{eqn:defn_of_g_R}. Exactly as in the previous section one can obtain an analogue of~\eqref{eqn:lower_bound_for_second_moment_bound} for $T_{{k}}'$:
\begin{equation}
\begin{split}
\frac{\E[T_{{k}}']^2}{\E[T_{{k}}'^{2}]}&\geq \frac{cf_{{R}}({k}) \,\P\!\left(X_{{R}^{{\kappa}-{k}}}> \altq_{{R}^{{\kappa}-{k}}}-{\frac{1}{12}}\right)^2}{f_{{R}}({k}) \P\!\left(X_{{R}^{{\kappa}-{k}}}> \altq_{{R}^{{\kappa}-{k}}}-\frac{1}{12}\right)^2 + g_{{R}}({k})\P\!\left(X_{{R}^{{\kappa}-{k}}}> \altq_{{R}^{{\kappa}-{k}}}-{\frac{1}{12}}\right)} \\
&= \frac{c}{1+\frac{g_{{R}}({k})}{f_{{R}}({k})}\frac{1}{\P\!\left(X_{{R}^{{\kappa}-{k}}}> \altq_{{R}^{{\kappa}-{k}}}-{\frac{1}{12}}\right)}} \geq \frac{c}{1+ 2 {R}^{-{k}/2}/p},
\end{split}
\end{equation}
where $c=c({R})>0$ is a constant whose value may change from line to line. The last inequality follows by bounding $\P\!\left(X_{{R}^{{\kappa}-{k}}}> \altq_{{R}^{{\kappa}-{k}}}-{\frac{1}{12}}\right)\geq p$, and using the same bounds used in~\eqref{eqn:final_bound_on_g_R/f_R} to bound $\frac{g_{{R}}}{f_{{R}}}$ from above (fix $\varepsilon=1/4$ and take ${R}$ large enough). Finally, note that by arguments similar to the ones in the derivation of~\eqref{eqn:lower_bound_on_E[T_k]^2}, there exists a constant $C=C(R)>0$ such that 
\begin{equation}
\begin{split}
\E[T_{{k}}']&=\sum_{i} \P(B_i \textrm{ is poor}) \geq C \cdot \sum_i \P(B_i \textrm{ is important})\cdot \P\!\left(X_{{R}^{{\kappa}-{k}}}> \altq_{{R}^{{\kappa}-{k}}}-{\frac{1}{12}}\right) \\ &\geq C p \cdot \!\left(\frac{{R}}{2}\right)^{2{k}}\cdot {R}^{-(5/4+\zeta/3){k}}\geq Cp \cdot {R}^{(3/4-\zeta/2){k}},
\end{split}
\end{equation}
where the last line follows from Lemma~\ref{lem:4arms_polynomial_concentration} for ${R}\geq {R}(\zeta/3)$ big enough. Choose $\epsilon=Cp/2$. Then $\theta_{{k}}\leq {R}^{-{k}\zeta/2}\leq {R}^{-\zeta/2}\leq \frac{1}{2}$ for ${R}={R}(\zeta)$ large enough.

Combining these results, we get that $$\P\!\left(X_n\leq \epsilon {R}^{\!\left(\frac{3}{4}-\zeta\right){k}}\altq_{{R}^{{\kappa}-{k}}}\right) \leq 1- \frac{1}{4} \frac{c({R})}{1+2{R}^{-{k}/2}/p}.$$ Let $k_0=k_0({R})$ be the smallest positive integer such that ${R}^{-k_0/2}\leq \frac{c({R})}{16}$. Then for all ${k}\geq k_0$, the RHS of the above equation is smaller than $1-p$, provided that $p<\frac{c({R})}{8}$.  Thus, {for $\tilde{c}=\tilde{c}(R,\zeta)$ sufficiently small (i.e.\ $\tilde{c}\leq \epsilon=C(R)p(R)/2$), we have} 
\begin{equation}
    \P\!\left(X_n\leq {\tilde{c}R^{(\frac{3}{4}-\zeta)k}}\altq_{{R}^{{\kappa}-{k}}}\right)\leq 1-p\leq\P(X_n\leq \altq_n), \label{eq:bound_on_tilde_R}
\end{equation} for all ${k}\geq k_0$ and $n$ such that $n>{R}^{k}$. 

Now fix any {$\tilde{c}\leq 1$ small enough and $k_0$ large enough such that~\eqref{eq:comp_mn_scales} holds.} 
Let $0<s\leq 1$ and let $\lambda\in\N$ be the smallest integer such that $\tilde{c}^{-1}R^{-\!\left(\frac{3}{4}-\zeta\right)\lambda}\leq s^{\frac{3}{4}-\zeta}$ and $\lambda\geq k_0$. 
In particular, $s^{\frac{3}{4}-\zeta}\leq \tilde{c}^{-1}R^{-\!\left(\frac{3}{4}-\zeta\right)(\lambda-k_0)}$ i.e.\ $R^{\lambda}\leq C{s^{-1}}$ for some constant $C=C(R,\zeta)$. It follows that for all $n>R^{\lambda}$
\begin{equation}
    \begin{split}
        \P(X_{R^{\kappa-k}}>{s^{\frac{3}{4}-\zeta}}\altq_n)\leq&\,\P(X_{R^{\kappa-k}}>{c^{-1}R^{-(\frac{3}{4}-\zeta)\lambda}}\altq_n)\\
        \leq&\,\P(X_{R^{\kappa-k}}>\altq_{R^{\kappa-\lambda}})\\
        \leq&\,CR^{-\alpha(k-\lambda)}\\
        \leq& C{s^{\alpha}}R^{-\alpha k},
    \end{split}
\end{equation}
where the value of $C=C(R,\zeta,\alpha)$ may change from line to line. In the second line we have used~\eqref{eq:comp_mn_scales}, while in the third we are using the standing assumption~\eqref{eq:assumptionsubpoly}. On the other hand, for $n\leq R^{\lambda}$ we have that $R^k<n\leq R^{\lambda}$, making the above inequality trivially true for any $C\geq 1$. This concludes the proof.
\end{proof}

\begin{cor}\label{cor:lower_bound_on_m_n}
There exists $n_0$ such that the following holds. For all $\zeta>0$, there exist $\tilde{C}(\zeta)>0$ and $p_0>0$ such that for all $p\leq p_0$ we have
\[ \altq_n\geq \tilde{C}(\zeta) n^{3/4-\zeta} \]
for all $n\geq n_0$ sufficiently large.\footnote{Like before, we need $n\geq n_0$ so that $\altq_n$ is defined and positive.} 
\end{cor}
\begin{proof} 
Suppose that $n_0=n_0(\delta)$ is chosen as in the proof of Lemma \ref{lem:compXscales} above, so that $\altq_n$ is well defined for all $n\geq n_0$ and so that $\altq_n(p)\geq 1/6$ for all $n\geq n_0$ and $p\in (0,1/2)$. By equation~\eqref{eq:comp_mn_scales}, we have that for all $\zeta>0$, there exist $R_0\geq n_0$ sufficiently large depending on $\zeta$ and $p_0<1/2$ sufficiently small depending on $R_0,\zeta$ such that for all $n> R_0^{k_0}$
\begin{equation}
    \altq_n(p_0)\geq {\tilde{c}R_0^{(\frac{3}{4}-\zeta)(\kappa-1)}}\altq_{R_0}(p_0)\geq\tilde{C}(\zeta,R_0)n^{(\frac{3}{4}-\zeta)}\cdot{\frac{1}{6}}\geq\tilde{C}(\zeta,R_0)n^{(\frac{3}{4}-\zeta)}.
\end{equation}
In the above, we have used that $\altq_{R_0}(p_0)\geq {\frac{1}{6}}$ (we chose $R_0\geq n_0$) and that $\tilde{c}$ depends only on $R_0$ and $\zeta$. Since $k_0=k_0(R_0)$, it follows that the above inequality holds for all $n\geq n_0$ upon possibly decreasing the constant $\tilde{C}(\zeta,R_0)$ -- it suffices to choose $\tilde{C}(\zeta,R_0)\leq R_0^{-\left(3/4-\zeta\right)k_0}/6$, as $\altq_n(p_0)\geq 1/6$ for $n\geq n_0$.
Finally, note that $\altq_n$ is decreasing in $p$, so the statement also holds for all $p\leq p_0$.
\end{proof}

From now on, given any $R\in \N$, we fix $p$ and, in turn, $\altq_n$ so that the statements of the above lemmas hold and they depend only on the choice of $R$.

We now turn our attention to open crossings of an annulus. The following lemmas control the number of open crossings of an annulus, as well as the distance between the said open crossings. In all the lemmas, we consider percolation on $\Lambda_n$ restricted to the annulus $A(x;r,2r)\subset[0,1]^2$ for some $x\in[0,1]^2$ and $r>0$.

\begin{figure}[ht]
    \centering
    \includegraphics[width=1\linewidth]{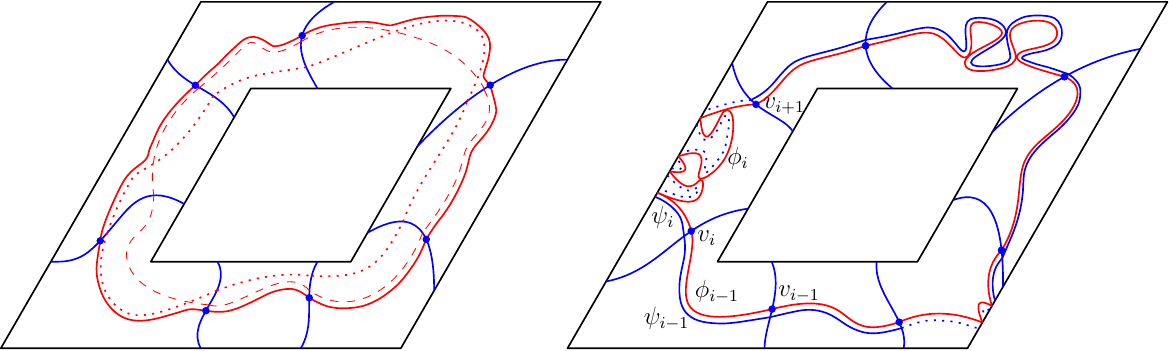} 
    \caption{On the left-hand side are pictured two closed loops with the least amount of defects -- the dotted and the dashed red lines. The full red line represents the concatenation of the outermost parts of these loops.  On the right-hand side picture, the red loop represents the outermost loop $\phi$. Note that this loop can be self-intersecting. The blue sites represent the defects $v_1,\ldots, v_m$. One can also see the blue layer of open sites neighboring the sites of $\phi$ represented by both full and dotted blue curves above. It exists precisely next to the segments of $\phi$ that do not belong to the outer boundary $\partial B(x,2r)$. The full blue segments represent precisely the paths $\psi_i$.}

    \label{fig:outermost_loop_phi}
\end{figure}

The first of these lemmas is purely combinatorial, and is a direct consequence of Menger's theorem (see e.g.\ \cite[Theorem 3.3.1]{diestel2012graph}).
\begin{lem}\label{lem:decompositioncrossings} 
Fix some $r>0$. Assume that there are exactly $m$ disjoint open crossings $\eta_1,\ldots,\eta_m$ of the annulus $A(x;r,2r)$. Then there exists a sequence of sites $v_1,\ldots,v_m$ and closed paths $\phi_1,\ldots,\phi_m$ such that $\phi:=\{v_1\}\cup\phi_1\cup\ldots\cup\{v_m\}\cup\phi_m$ forms a loop.
\end{lem}
The outermost such loop with the least number of `defects' is well defined. Indeed, if we consider two such closed loops with $m$ defects, we can obtain another one by concatenating the outermost closed segments of both loops. The new loop still has $m$ defects, and we still have disjoint open crossings of the annulus through each of the defects. See Figure~\ref{fig:outermost_loop_phi}. We will choose $\phi$ to be this outermost loop and let $\{v_1\},\phi_1,\ldots,\{v_m\},\phi_m$ be the associated open defects and closed segments in clockwise order. Here $\phi$ is allowed and will in general be a non-simple loop. Finally, let $\operatorname{int} \phi$ denote the set of sites in the bounded component of $\phi^c$. 

Note that because $\phi$ is the outermost loop, it must consist of closed sites belonging to cluster boundaries -- the closed sites of $\phi$ are either at the outer boundary of the annulus, or there is another layer of open sites neighboring the sites of $\phi$ -- see Figure~\ref{fig:outermost_loop_phi}. Moreover, the open paths neighboring $\phi_i$ and $\phi_{j}$ for $i\neq j$ are disjoint. Indeed, if they intersected, they would have to intersect in the exterior of $\phi$. This would create a loop with the same number of defects as $\phi$ that is outside of $\phi$, contradicting that $\phi$ is the outermost such loop. Denote by $\psi_i$ the open path neighboring $\phi_i$ that starts at $v_i$ and ends the first time it hits either $v_{i+1}$ or the boundary of $B(x,2r)$, depending on whether the closed cluster of $\phi_i$ touches the outer boundary of the annulus or not. 

We will ultimately want to show that if any two $v_i,v_j$ are connected inside $\operatorname{int}\phi$, then their (geodesic or resistance) distance \textit{inside of} $\operatorname{int}\phi$ is not large. Using a suitable resampling procedure, we will be able to compare this distance to a distance from the a priori estimate. 

The following three lemmas explain the typical behavior of $\phi$. The first one tells us that with high probability $m$ is finite and the defect points $v_1,\ldots,v_m$ are macroscopically far away from each other and the annulus boundary. The second one tells us that the exterior of $\phi$ is sufficiently nice. The last lemma tells us that every defect point is `shielded away' from $\operatorname{int}\phi$ by closed clusters, which ensures that we can resample our percolation configuration near these points without affecting the distances in $\operatorname{int}\phi$.

\begin{lem}\label{lem:number_crossings} 
The following holds for $R \in \N$ large enough. 
Let $F_1=F_1(M)$ be the event that there exist at most $M$ disjoint open crossings of $A(x;r,2r)$. On this event, let $\phi$ and $v_1,\ldots, v_m$ be the outermost almost closed loop and the associated $m\leq M$ defects. Let $F_2=F_2(a)$ be the event that the defects $v_1,\ldots, v_m$ are all at $L^1$-distance at least $ar$ from the annulus boundary. Let $F_3=F_3(a)$ be the event that $v_1,\ldots,v_m$ are all at $L^1$-distance at least $ar$ apart. Let $0<p'<1$. Then there exists $M_0=M_0(p')$ such that for all $M\geq M_0$ there exists $a_0=a_0(p',M)$ such that for all $a<a_0$ and all radii $r\geq \frac{R}{an}$ we have $\P(F_1\cap F_2\cap F_3)\geq 1-p'$.
\end{lem}

\begin{rem} 
In order for the statements of the lemmas in this section to hold, we need to choose~$R$ large enough so that the arm estimates apply for annuli of radii larger than $R$. In the setting of the above lemma, we will apply annulus estimates for annuli of (graph) radius $anr$ -- thus, it is sufficient to assume that $anr\geq R$ to make the proof work.\footnote{The same result in fact holds for all annuli of radius $anr\geq 10$, say. Indeed, the arm exponent estimates hold modulo an additional constant factor for all annuli that have `enough space' to have $6$ arms. 
}
\end{rem}

\begin{proof}[Proof of Lemma~\ref{lem:number_crossings}]
Note that the probability that there exists an open crossing of the annulus is bounded from above by some $p_0<1$ uniformly in $n,r$. Then by the BK inequality, $\P(F_1^c)\leq p_0^M<\frac{p'}{3}$ for $M$ sufficiently large.

We now turn our attention to bounding $\P(F_2)$. Consider $12a^{-1}$ evenly spaced points $z$ around the inner and outer boundaries of the annulus. On the event $F_2^c$, there exists such a $z$ for which $v_i\in B(z,ar)$ for some $1\leq i\leq M$. There are $3$ `arms' emanating from $v_i$ - one open arm going to the opposite boundary of the annulus (depending on whether $z$ is on the inner/outer boundary), and two closed `arms' corresponding to the circuit $\phi$. These are not fully closed as $\phi$ contains open sites (points $v_j$), so our closed arms have 'defects'. On the event $F_1$, there are at most $M$ such defects. Moreover, these three arms go at distance of at least $\frac{1}{2}r$, and, due to the lozenge shape of our annulus, are fully contained in a sector of angle $\frac{5\pi}{3},\frac{4\pi}{3}$ (if we are on the inner boundary) or $\pi$ (if we are on the outer boundary). See Figure~\ref{fig:number_of_crossings} for an illustration of this event. The half-plane three arm exponent is known to be $2$ (see e.g.\ \cite[Theorem 3]{smirnov2001criticalexponents}). Using the conformal invariance between the half-plane and the wedge domains, one can see that the three arm exponent of a $\frac{5\pi}{3}$-wedge is $2\cdot \frac{3}{5}$.\footnote{One can do this proof without conformal invariance, using only arm estimates. Namely, for an $a$-box centered on the outer boundary of our annulus, the arms always belong to a half-plane, thus we can just directly use half-plane exponents. For an $a$-box centered at a point $z$ at the inner boundary of the annulus, let $xar$ with $x\in\N$ be the distance of $z$ to the nearest corner of the inner boundary of the annulus. Then there are $3$ arms in the half-plane inside $A(z;ar,xar)$. There are also $3$ arms in the plane inside $A(z;xar,r/2)$. Using Lemma~\ref{lem:arm_exponents} we find that the probability of this event is $\lesssim\sum_{x=1}^{\frac{1}{2}a^{-1}}\!\left(\frac{xar}{ar}\right)^{-2+\zeta}\cdot \!\left(\frac{r/2}{xar}\right)^{-2/3+\zeta}\lesssim a^{2/3-\zeta}\cdot \sum_{x=1}^{\frac{1}{2}a^{-1}}x^{-4/3}\lesssim a^{2/3-\zeta}$. Here we have used that the series $\sum_{x\geq 1}x^{-4/3}$ is summable. This event will keep showing up throughout the paper. For the sake of brevity, we will deal with it using conformal invariance every subsequent time, but one should keep in mind that this alternative proof exists.} 
Moreover, Nolin proved in \cite[Proposition 18]{nolin2008near} that defects on the arms only add a logarithmic correction to the probability that we have regular arms without defects.

Thus, if we choose $\zeta =1/6$ and $R\geq \rho_{\zeta,3}$ from Lemma~\ref{lem:arm_exponents}, we know that for all $r$ such that $r\cdot an\geq R$ we have that $$\P(F_1\cap F_2^c) \leq 12a^{-1}\cdot \!\left(1+\log \!\left(\frac{rn}{arn}\right)\right)^M \!\left( \frac{\frac{1}{2}rn}{arn}\right)^{-6/5+\zeta} \lesssim \!\left( 1-\log a\right)^M \cdot a^{1/5-\zeta}\leq \frac{p'}{3},$$ for $a=a(p',M)$ sufficiently small.

\begin{figure}[ht]
    \centering
    \includegraphics[width=0.9\linewidth]{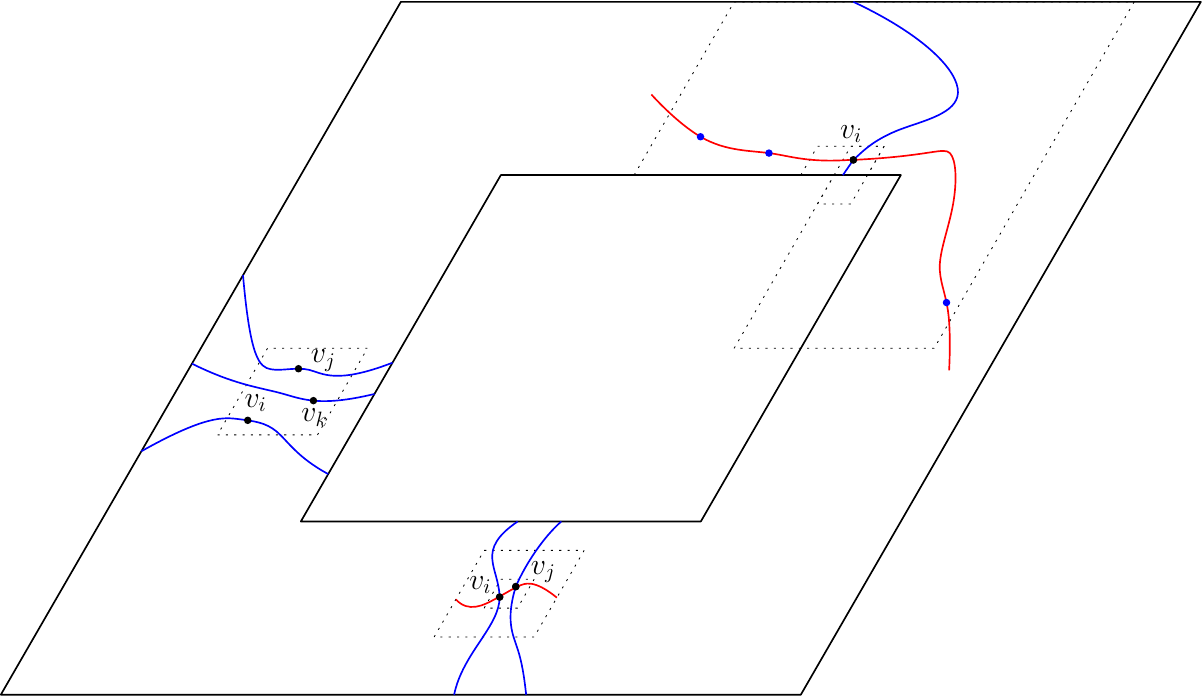}
    \caption{A sketch of the events $F_{2}^c$ in the top-right corner, $F_{31}$ on the left and $F_{32}$ on the bottom. On $F_2^c$, there are 2 closed arms with defects (blue points on the picture above) and 1 open arm from the sides of the box going all the way to the boundary. These arms are contained in a wedge domain of angle at most $5\pi/3$ as pictured above. On $F_{31}$, there are 6 open blue arms from the sides of the box to the boundary. On $F_{32}$, there are 4 open blue arms from the sides of the box to the boundary and 2 closed red arms going far. The size of the boxes is not to scale.}
    \label{fig:number_of_crossings}
\end{figure}

We finally bound $\P(F_3)$. Note that there are disjoint open annulus crossings through each $v_i$ and closed paths between each $v_i$ and $v_{i+1}$. Assume that we are working on $F_3$ -- there exists a pair $v_i,v_{j}$, with $i\neq j$ such that $\|v_i-v_{j}\|_1\leq ar$. Then either there exists $v_{k}$, $k\neq i,j$ such that at least one of $\|v_i-v_{k}\|_1,\|v_{j}-v_{k}\|_1\leq\sqrt{a}r$ holds or not. We call the first event $F_{31}$ and the second event $F_{32}$. 
See Figure~\ref{fig:number_of_crossings} for an illustration. On $F_{31}$, there are 6 disjoint open arms from a box of radius $2\sqrt{a}r$ to the boundary of the annulus. On $F_{32}$, there are 4 disjoint open arms from a box of radius $ar$ to the boundary of the annulus and 2 disjoint closed arms to radius $\sqrt{a}r$. 

We bound the probability of $F_{31}$. Consider a grid of evenly spaced sites $z$ in the annulus at distance $\sqrt{a}r$. Then under $F_{31}\cap F_2\!\left(a^{1/100}\right)$, there exists a $z$ such that $B(z,4\sqrt{a}r)$ has 6 disjoint open arms to the boundary of the annulus, which is at least $\frac{1}{2}a^{1/100}r$ away from $z$. Hence, if we choose $\zeta=1/100$ and $R>\rho'_{\zeta,6}$ as in Lemma~\ref{lem:arm_exponents_monochromatic}, then we know that for all $r$ such that $\sqrt{a}rn\geq R$ we have $$ \P\!\left(F_{31}\cap F_2\!\left(a^{1/100}\right)\right) \lesssim a^{-1} \cdot \!\left( \frac{a^{1/100}}{a^{1/2}}\right)^{-35/12+\zeta} \leq a^{5/12}\leq \frac{p'}{6}$$ for $a=a(p')$ small enough. We can similarly bound $\P(F_{32})$. For $\zeta=1/100$ and $R>R(1/100)$, Lemma~\ref{lem:arm_exponents} tells us that for all $r$ such that $arn\geq R$ we have $$\P\!\left(F_{32}\cap F_2\!\left(a^{1/100}\right)\right)\lesssim a^{-2} \cdot \!\left(\frac{a^{1/2}}{a}\right)^{-35/12+\zeta}\cdot \!\left(\frac{a^{1/100}}{a^{1/2}}\right)^{-5/4+\zeta}\leq a^{1/24} \leq \frac{p'}{6}$$ for $a=a(p')$ small enough.

Thus, if we choose $M=M(p')$ big enough so that $\P(F_1^c)<\frac{p'}{3}$ and $a=a(M,p')$ small enough so that $\P\!\left(F_1\cap F_2^c\!\left(a^{1/100}\right)\right)<\frac{p'}{3}$ and $\P\!\left(F_3^c\cap F_2\!\left(a^{1/100}\right)\right)<\frac{p'}{3}$, then we get that for all $r$ with $arn\geq R$ we have $\P\!\left(\!\left(F_1(M)\cap F_2(a)\cap F_3(a) \right)^c\right)\leq \P\!\left(\!\left(F_1(M)\cap F_2(a^{1/100})\cap F_3(a) \right)^c\right)\leq p'$, as wanted.
\end{proof}

\begin{lem} The following holds for $R \in \N$ large enough. Fix $0<p'<1$. Suppose that $M,a$ are chosen so that $\P(F_1\cap F_2\cap F_3)\geq 1-p'/2$. For $0<b<a$, consider a grid of evenly spaced sites at distance $br$ of each other, such that boxes of side length $br$ centered at those sites cover $B(x,2r)$. Let $F_4(b)$ be the event that for any $\phi_i$, there exists a sequence of boxes $B_1,\ldots, B_t$, such that:
\begin{enumerate}
\item[$\bullet$] the boxes are consecutive, i.e., $B_{\tilde{t}}$ and $B_{\tilde{t}+1}$ share a side for all $1\leq \tilde t <t$;
\item[$\bullet$] $B_t$ intersects $\partial B(x,2r)$, $B_1\cap \phi_i \neq \emptyset$ and $\bar{B}_1\cap \psi_i \neq \emptyset$ where $\bar{B}_1$ denotes the set of all sites of $B_1$ as well as all the sites neighboring those;
\item[$\bullet$] $4B_{\tilde{t}}$ does not intersect $\phi \setminus \phi_i$ for any $1\leq \tilde{t}\leq t$. \label{lem:corridor_of_boxes}
\end{enumerate}

Then there exists $b_0=b_0(p',M,a)$ such that for all $b<b_0$ we have that $\P(F_1\cap \cdots \cap F_4)\geq 1-p'$ holds for all annuli of radius $r$ with $r\geq \frac{R}{bn}$. 
\label{lem:exterior_is_nice}
\end{lem}
\begin{proof} We use similar ideas to those in the proof of the previous lemma. Assume that $F_1\cap F_2\cap F_3$ holds with the appropriate constants $M$ and $a$. Consider a grid of evenly spaced sites $z$ in the annulus at distance $br$, and consider the boxes centered at these points of side length $br$.

Fix any $1\leq i\leq m$. We split the proof into two cases, depending on whether the closed cluster of $\phi_i$ hits the outer boundary of $A(x;r,2r)$ or not. If it does, consider the first point of $\phi_i$ that belongs to $\partial B(x,2r)$ (it neighbors the endpoint of $\psi_i$ on $\partial B(x,2r)$) and let $B_1$ be the box in our $br$-lattice that contains this point. Then, $B_1$ on its own satisfies the conditions of the lemma, unless $4B_1$ intersects $\phi\setminus \phi_i$. Call this event $F_{40}$. On $F_{40}$ there exists $j\neq i$ such that $\phi_i$ and $\phi_j$ intersect $4B_1$. Moreover, at least one open crossing of the annulus must pass through this box, since it has to pass in between the $\phi_i$ and $\phi_j$. See Figure~\ref{fig:exterior_is_not_ugly}. Thus, $4B_1$ has three arms that go distance $a$ in the half-plane. Lemma~\ref{lem:arm_exponents} with $\zeta=1/2$ and $R\geq R(\zeta)$ suitably large then tells us that
$$\P(F_1\cap F_2\cap F_3\cap F_{40})\leq \frac{8rn}{brn} \cdot \!\left(\frac{arn}{4brn}\right)^{-2+\zeta}\leq 128 a^{-3/2}b^{1/2}$$ for all $r$ such that $brn\geq R$. The right-hand side can be made smaller than $p'/6$ by choosing $b=b(p',M,a)$ suitably small. 

Now, assume that the cluster of $\phi_i$ does not intersect the boundary of $B(x,2r)$. In this case, $\psi_i$ is an open path connecting $v_i$ and $v_{i+1}$ and so we automatically have $\bar{B}_1\cap \psi_i\neq \emptyset$ if $B_1\cap \phi_i\neq \emptyset$. Consider the lattice of boxes of side length $br$ needed to cover $B(x,2r)$. In analogy to percolation on the square lattice, we say that any such box $B$ such that $4B$ intersects $\phi \setminus \phi_i$ is closed. All other boxes are open. We say that two open boxes are neighboring if they share a side, while two closed boxes are neighboring if they share a vertex. 

We know that there exist disjoint open crossings $\eta_i$ and $\eta_{i+1}$ of the annulus that do not pass through any of the other $v_j$'s for $j\neq i,i+1$. These crossings enclose $\phi_i$ in a sector-like domain of the annulus. In particular, for any closed box $B$, we must have that $4B$ intersects the outside of this sector. On $F_4^c$, there does not exist a path of open boxes connecting $\phi_i$ to the outside boundary of the annulus, so there must exist a blocking path of closed boxes closing off $\phi_i$. In particular, there are two adjacent closed boxes $B$ and $B'$, with $4B$ intersecting one side of the sector and $4B'$ intersecting the other one. Thus, there exists a site $z$ in our grid such that both open crossings pass through the box $B(z,6br)$. Moreover, this box intersects $\phi\setminus \phi_i$ on both sides of the sector.

We will differentiate between two cases: either $z$ is at distance at least $ar$ from the boundary of the annulus or not. Call these events $F_{41}$ and $F_{42}$. On $F_{41}$, there are $6$ arms with defects crossing the annulus $A\!\left(z;6br,\frac{1}{4}a r\right)$. Indeed, there are $4$ open arms without defects coming from the two disjoint open crossings $\eta_{i}$ and $\eta_{i+1}$ of the annulus. Additionally, we have two closed arms that can have at most $1$ defect (only one of the defect points can be in $B\!\left(z,\frac{1}{4}ar\right)$) and they come from the loop $\phi$ -- see Figure~\ref{fig:exterior_is_not_ugly} for an illustration of this event. 
As in the previous lemma, using \cite[Proposition 18]{nolin2008near} and Lemma~\ref{lem:arm_exponents} we get that for all radii $r$ such that $brn\geq R$ we have $$\P(F_1\cap F_2\cap F_3\cap F_{41})\leq 4b^{-2} \cdot \!\left(1+\log \!\left( \frac{arn}{24brn}\right)\right)\cdot \!\left(\frac{arn}{24brn}\right)^{-\frac{35}{12}+\frac{5}{12}}\lesssim a^{-5/2}\cdot \!\left(2-\log b\right) \cdot b^{1/2}.$$ We can make the right-hand side smaller than $p'/6$ by taking $b=b(p',M,a)$ sufficiently small.

\begin{figure}[ht]
    \centering
    \includegraphics[width=1\linewidth]{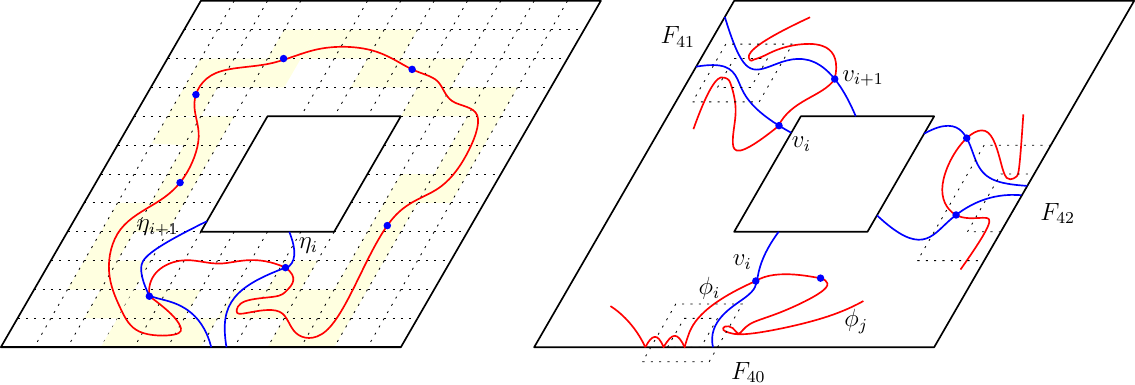}
    \caption{Figure on the left illustrates the event $F_4^c$ in the case that $\phi_i$ does not intersect the boundary of $B(x,2r)$. Yellow boxes above represent the closed boxes. Blue curves represent the open crossings $\eta_i$ and $\eta_{i+1}$ and the area containing $\phi_i$ that is enclosed by these crossings is the relevant sector-like domain. The second figure depicts the events $F_{40}$ on the bottom, $F_{41}$ on the top and $F_{42}$ on the right.}
    
    \label{fig:exterior_is_not_ugly}
\end{figure}

Similarly, on $F_{42}$ there exists a point $w$ in the grid of $12a^{-1}$
evenly spaced points around the inner and outer boundaries of the annulus such that $\eta_i,\eta_{i+1}$ and $\phi$ intersect $B(w,2ar)$. Then we know that there are 4 arms with defects going from this box to distance $\frac{1}{2}r$
in the half-plane corresponding to the two open crossings going through this bottleneck and the loop $\phi$. 
We have dealt with this exact event in the analysis of $\P(F_1\cap F_2^c)$ in the previous lemma. Upon possibly making $a_0=a_0(a,M)$ slightly smaller, we can guarantee that $\P(F_1\cap F_2\cap F_3 \cap F_{42}) \leq \frac{p'}{6}$. In particular, $$\P\!\left((F_1\cap \cdots \cap F_4)^c\right)\leq \P((F_1\cap F_2 \cap F_3)^c)+\P(F_1\cap F_2 \cap F_3\cap (F_{40}\cup F_{41}\cup F_{42}))\leq \frac{p'}{2}+\frac{p'}{6}+\frac{p'}{6}+\frac{p'}{6}=p'.$$ 
\end{proof}

\begin{lem} The following holds for $R \in \N$ large enough. Fix $0<p'<1$. Suppose again that $M,a$ are chosen so that $\P(F_1\cap F_2\cap F_3)\geq 1-p'/2$. For $0<b<a$, let $F_5(b)$ be the event that for every $i\leq m$, there exists $\hat{v}_i\in \operatorname{int}\phi$ such that $\hat{v}_i$ is another open touch-point of the closed clusters of the segments $\phi_{i-1}$ and $\phi_i$. Moreover, we ask that $\hat{v}_i$ is at distance at least $br$ away from $v_i$. For $0<c<b$, let $F_6(c)$ be the event that there exist closed paths from $\phi_{i-1}$ and $\phi_{i}$ to $\hat{v}_i$ that are contained fully outside of the box $B(v_i,c^{1/3}r)$. Also, let $F_7(c)$ be the event that for all $i$, once $\phi_{i-1}$ or $\phi_i$ leave the box $B(v_i,c^{1/3}r)$, they don't enter $B(v_i,cr)$ again. There exists $b_0=b_0(p',M,a)>0$ such that for every $b<b_0$ there exists $c_0=c_0(p',M,a,b)>0$ such that for every $c<c_0$ we have that $\P(F_1\cap F_2\cap F_3 \cap F_5 \cap F_6 \cap F_7) \geq 1-p'$ for all radii $r$ such that $crn\geq R$. 
\label{lem:points_are_shielded_away}
\end{lem}

\begin{proof} Assume that $M,a$ are chosen appropriately and that we are working on the event $F_1\cap F_2\cap F_3$. On $F_5^c$, there are $6$ arms crossing the annulus $A(v_i; br,ar)$: two closed ones coming from the closed clusters of $\phi_{i-1}$ and $\phi_i$, and $4$ open arms corresponding to the open sites neighboring these two clusters. These $4$ arms are non-intersecting. The two open arms in $\operatorname{int} \phi$ are disjoint as we are assuming that the closed clusters of $\phi_{i-1}$ and $\phi_i$ do not have any more touch points. The arms on the exterior of $\phi$ are disjoint since we have chosen $\phi$ to be the outermost loop. See Figure~\ref{fig:there_are_shields}. Arguing as before we find that for all levels $f$ such that $brn\geq R$ we have $$\P(F_1\cap F_2\cap F_3 \cap F_6^c)\lesssim b^{-2}\cdot \!\left(\frac{arn}{brn}\right)^{-35/12+5/12} \leq a^{-5/2}\cdot b^{1/2}.$$ For $b=b(p',M,a)$ sufficiently small this quantity is smaller than $p'/12$.

\begin{figure}[ht]
    \centering
    \includegraphics[width=1\linewidth]{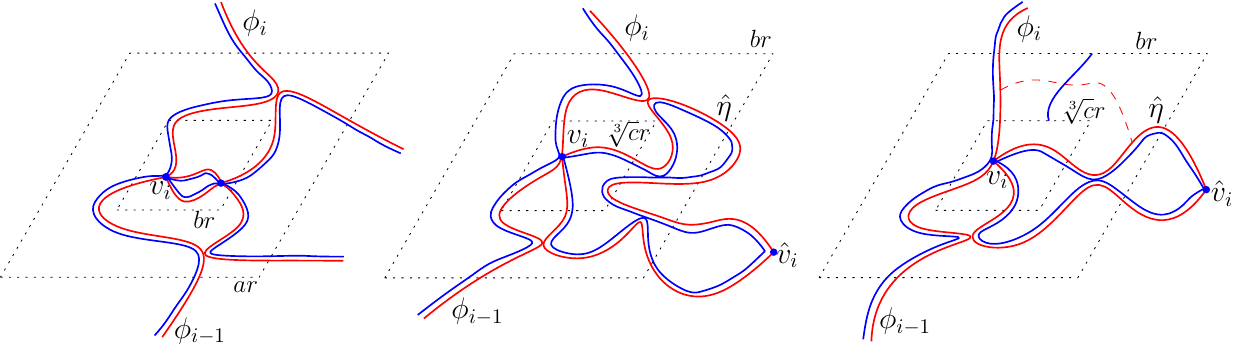}
    \caption{The first of the three figures illustrates the event $F_5(b)^c$ - if clusters of $\phi_{i-1}$ and $\phi_i$ do not touch outside of the box $B(v_i,br)$, then the blue arms pictured above must be disjoint. 
    The second figure pictures the event $F_6(c)^c$ in the case that continuation of the path $\eta$ revisits the box $B(v_i,c^{1/3}r)$. Last figure illustrates the same event, but in the case that $\eta$ does not come back inside $B(v_i,c^{1/3}r)$. If there is no closed path connecting $\phi_i$ and $\eta$ inside of the annulus $A(v_i;c^{1/3}r,br)$ (the dashed red line on the picture), then we must have a blocking path - it is represented by the blue line crossing the dashed red line.} 
    \label{fig:there_are_shields}
\end{figure}

Assume now that $b$ is as above and that we are additionally working on the event $F_5$. On $F_6^c$, there either does not exist a closed path from $\phi_{i}$ to $\hat{v}_i$ or from $\phi_{i-1}$ to $\hat{v}_i$ that is contained fully outside of the box $B(v_i,c^{1/3}r)$. Assume the former holds. Denote by $\eta$ the boundary path from $v_i$ to $\hat{v}_i$ associated to the closed cluster of $\phi_i$. Denote by $\hat{\eta}$ the first segment of this path that crosses the annulus $A(v_i;c^{1/3}r,br)$. Either the continuation of $\eta$ crosses the annulus again and ends up inside $B(v_i,c^{1/3}r)$ again, or it stays fully outside this box. If the former holds, there are $7$ arms crossing the annulus -- see Figure~\ref{fig:there_are_shields}. If the latter holds, there cannot be a closed path inside the annulus connecting $\phi_i$ and $\hat{\eta}$, otherwise there is a closed path connecting $\phi_i$ and $\hat{v}_i$ outside of $B(v_i,c^{1/3}r)$. This means that there must exist an open crossing of the annulus disconnecting these two closed crossings. Once again this yields $7$ arms -- see Figure~\ref{fig:there_are_shields} again. Hence, as before using Lemma~\ref{lem:arm_exponents} we find that for, say, $\zeta=1$ and $R>\rho_{1,7}$, we have $$\P(F_1\cap F_2\cap F_3 \cap F_5 \cap F_6^c) \leq \!\left(c^{1/3}\right)^{-2} \cdot \!\left( \frac{\frac{1}{4}brn}{4c^{1/3}rn}\right)^{-4+\zeta}\lesssim b^{-3} \cdot c^{1/3},$$ for all annuli of radius $r$ such that $c^{1/3}rn\geq R$. By choosing $c=c(p',M,a,b)$ sufficiently small, we can make the above expression less than $p'/12$. 

Finally, on the event $F_7(c)^c$ we see exactly as before that there are $7$ arms in the annulus $\allowbreak A(v_i; cr,\sqrt[3]{c}r)$. Thus, choosing $\zeta = \frac{1}{2}$ and $R\geq \rho_{1/2,7}$, Lemma~\ref{lem:arm_exponents} tells us that  $$\P(F_7(c)^c) \leq c^{-2} \cdot \!\left( \frac{\frac{1}{4}c^{1/3}rn}{4crn}\right)^{-4+\zeta} \lesssim c^{1/3},$$ for all annuli of radius $r$ such that $crn\geq R$. Again, we can make the above expression smaller than $p'/12$ by choosing $c$ sufficiently small.

Combining this with the above results, we get $\P(F_1\cap F_2\cap F_3\cap F_5 \cap F_6\cap F_7) \geq 1-p',$ as required.
\end{proof}

Define $\widehat{F}(A(x;r,2r))=F_1\cap \cdots \cap F_7$. The above lemmas tell us that for every $0<p'<1$, there exist $M=M(p'),a=a(p',M),b=b(p',M,a)$ and $c=c(p',M,a,b)$ such that $\P(\widehat{F})\geq 1-p'$, for all annuli $A(x;r,2r)$ of radius $r$ such that $crn\geq R$. 
Finally, if we denote by $\widehat{F}'(A(x;r,2r))$ the same event as $\widehat{F}$ but with the roles of open and closed clusters reversed, the exact same lemmas tell us that $\widehat{F}'(A(x;r,2r))\geq 1-p'$ too.

We would like to use a resampling procedure to link clusters with positive probability. To control the probability of the resampling being successful, we will define the so-called \emph{well-separation} of annulus crossings.

We now describe the setup for the resampling procedure. Given an annulus $A(w;\tilde{r},2\tilde{r})\subseteq \Lambda_n$, define $\omega_{\textrm{out}}({A(w;\tilde{r},2\tilde{r})})$ as the percolation configuration outside of $B(w,2\tilde{r})$ together with all the cluster interface explorations started from $\partial B(w,2\tilde{r})$ and stopped upon hitting $\partial B(w,\tilde{r})$ or returning to $\partial B(w,2\tilde{r})$. Let $\Lambda_{\textrm{out}}({A(w;\tilde{r},2\tilde{r})})$ denote the set of sites in $B(w,2\tilde{r})^c$ together with the explored sites in $A(w;\tilde{r},2\tilde{r})$. Note that the conditional law of the percolation in $\Lambda_{\textrm{out}}({A(w;\tilde{r},2\tilde{r})})^c$ is still i.i.d.\ Bernoulli percolation.

For $\xi>0$, we say that $\Lambda_{\textrm{out}}({A(w;\tilde{r},2\tilde{r})})$ is \emph{$\xi$-well separated} if for every site $z \in \partial B(w,\tilde{r})$ that is the endpoint of some cluster interface exploration $\gamma$ started from $\partial B(w,2\tilde{r})$ we have that $B(z,\xi \tilde{r}) \cap \Lambda_{\textrm{out}}({A(w;\tilde{r},2\tilde{r})})\subseteq \gamma$. In words, this says that endpoints of cluster interface explorations that reach $\partial B(w,\tilde{r})$ are $\xi \tilde{r}$-far away from the rest of the explorations.
\begin{lem} The following holds for $R \in \N$ large enough. Let $\tilde{p}>0$. Suppose we are given a box $B=B(w,\tilde{r})$. Then, there exists $\xi=\xi(\tilde{p})>0$ such that for all radii $\tilde{r}$ such that $\xi \tilde{r}n\geq R$ we have that $$\P(\Lambda_{\textrm{out}}(2B\setminus B) \textrm{ is }\xi\textrm{-well separated})\geq 1-\tilde{p}.$$ \label{lem:omega_is_well_separated}
\end{lem}
\begin{proof} The proof again follows from a straightforward arm exponent calculation. Consider $4 \xi^{-1}$ evenly spaced points around the inner boundary of the annulus. If $\Lambda_{\textrm{out}}(2B\setminus B)$ is not $\xi$-well separated, then there exists a point $z$ on $\partial B$ such that it is the endpoint of an exploration $\gamma$ starting from $\partial 2B$ such that another such exploration $\gamma'$ enters $B(z,\xi \tilde{r})$. Consider the point $w$ on our grid such that $z\in {B(w,\xi \tilde{r})}$. Since both $\gamma$ and $\gamma'$ consist of open and closed paths, this means that there are at least $3$ arms going from $B(w,2\xi \tilde{r})$ up to distance of at least $\frac{1}{2}\tilde{r}$. Moreover, because these arms belong to the interior of $2B\setminus B$, it means that they are fully contained in a sector of angle at most $\frac{5\pi}{3}$ centered at $w$. Hence, as in the proof of Lemma~\ref{lem:number_crossings}, we find that for $\zeta=1/10$ we have that $$\P(\Lambda_{\textrm{out}}(2B\setminus B) \textrm{ is not } \xi\textrm{-well separated}) \leq 4\xi^{-1} \cdot \!\left(\frac{\frac{1}{2}\tilde{r}n}{\xi \tilde{r}n}\right)^{-6/5+\zeta}\lesssim \xi^{1/10}$$ for any box $B$ of radius satisfying $\xi\tilde{r}n\geq R$. Thus, by choosing $\xi=\xi(\tilde{p})$ small enough we can guarantee that $\P(\Lambda_{\textrm{out}}(2B\setminus B) \textrm{ is } \xi\textrm{-well separated})\geq 1-\tilde{p}$. 
\end{proof}

\begin{lem}\label{le:linking_probability}
For each $\xi>0$ there exists $p_1=p_1(\xi) > 0$ such that the following is true. Suppose that $\xi\tilde{r}n \ge 10$. Let $\Lambda_{\textrm{out}}({A(w;\tilde{r},2\tilde{r})})$ be as defined above, and sample the percolation in the components bounded by $\Lambda_{\textrm{out}}({A(w;\tilde{r},2\tilde{r})})$. Let $J_1,J_2,\ldots$ be the closed boundary components of the inner face of $\Lambda_{\textrm{out}}({A(w;\tilde{r},2\tilde{r})})$, and let $E_{ij}$ be the event that $J_i$ is connected to $J_j$ by a closed path in the inner face. Then
\[ \min_{i,j}\p[ E_{ij} \mid \omega_{\textrm{out}}({A(w;\tilde{r},2\tilde{r})})] \ge p_1 \1_{\{\text{$\Lambda_{\textrm{out}}({A(w;\tilde{r},2\tilde{r})})$ is $\xi$-well separated}\}} . \]
The analogous statement holds for the open boundary components.
\end{lem}
\begin{proof}
This follows from RSW estimates, which apply in our case since the resampled configuration is still i.i.d.\ percolation. By symmetry, it suffices to assume that the target clusters are closed. Let $z$ (resp.\ $z'$) be one of the endpoints of $J_i$ (resp.\ $J_j$). By the definition of $\Lambda_{\textrm{out}}({A(w;\tilde{r},2\tilde{r})})$, they must be the endpoints on $\partial B(w,\tilde{r})$ of interface curves $\gamma,\gamma'$, respectively. The resampling in the inner face will be successful if 
\begin{enumerate}
\item[$\bullet$] In the resampled configuration, there are closed circuits both in $A\!\left(z;\frac{1}{2}\xi\tilde{r},\xi\tilde{r}\right) \setminus \gamma$ and in $A\!\left(z';\frac{1}{2}\xi\tilde{r},\xi\tilde{r}\right) \setminus \gamma'$; 
\item[$\bullet$] There are closed crossing paths in the $4$ strips of width $\frac{1}{2}\xi\tilde{r}$ and length $\tilde{r}$ adjacent to the four sides of box $B(w,\tilde{r})$ in its interior. See Figure~\ref{fig:resampling_successful_positive}.
\end{enumerate}

\begin{figure}[ht]
    \centering
    \includegraphics[width=0.8\linewidth]{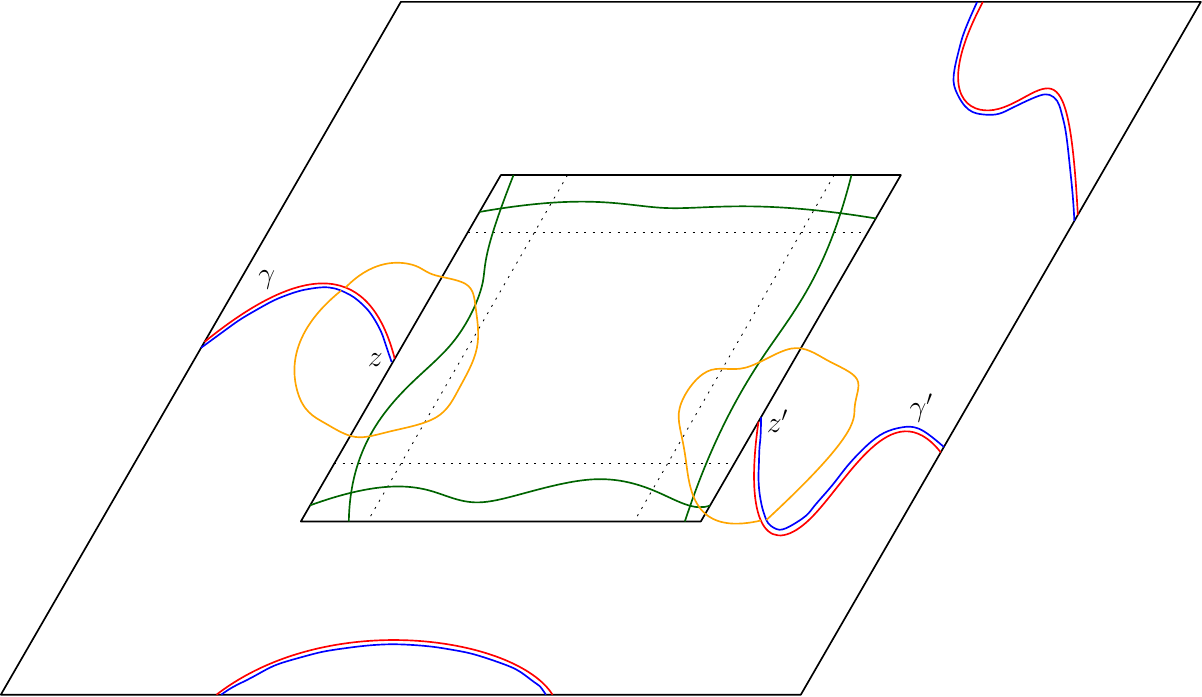}
    \caption{The blue-red curves on the picture represent the explored cluster boundaries starting from $\partial B(w,2\tilde{r})$. The two curves reaching the inner boundary are $\gamma$ and $\gamma'$. The yellow curves represent the closed circuits in $A\!\left(z;\frac{1}{2}\xi\tilde{r},\xi\tilde{r}\right)\setminus \gamma$ and $A\!\left(z';\frac{1}{2}\xi\tilde{r},\xi\tilde{r}\right) \setminus \gamma'$. The green curves in the strips of width $\frac{1}{2}\xi\tilde{r}$ represent the four closed crossings. On this event, the closed parts of $\gamma$ and $\gamma'$ are in the same closed cluster.}
    \label{fig:resampling_successful_positive}
\end{figure}

Given that $\xi\tilde{r}n\geq 10$, we know that the probability of the first item is bounded away from zero by an absolute constant using Russo-Seymour-Welsh. Similarly, each strip has a closed crossing with probability bounded away from zero by a constant depending only on $\xi$. Finally, positive correlation tells us that our event, and hence also the event that the resampling is successful both have probability bounded away from zero by a constant depending only on $\xi$.
\end{proof}

This finishes the string of general results for annuli $A(x;r,2r)\subseteq \Lambda_n$. As advertised before, we will use these results to show that in any annulus there is a finite number of points disconnecting inside from the outside and that the distance between these points is not big. We do so using a suitable resampling procedure and the a priori estimate.

\medskip Consider an $R$-adic annulus $A(x;n^{-1}R^{\kappa-f},n^{-1}R^{\kappa-f+1})\subseteq \Lambda_n$ of level $f$. Fix any $0<p'<1$. We define the event $F=F(A(x;n^{-1}R^{\kappa-f},n^{-1}R^{\kappa-f+1}))$ as the intersection of the following events: \phantomsection\label{defn:definition_of_event_F}
\begin{enumerate}
\item[$\bullet$] $\widehat{F}(A(x;2n^{-1}R^{\kappa-f},4n^{-1}R^{\kappa-f}))$, where $M_1,a_1,b_1,c_1$ are chosen so that
\begin{equation}
    \P(\widehat{F}(A(x;2n^{-1}R^{\kappa-f},4 n^{-1}R^{\kappa-f})))\geq 1-p'/4.
\end{equation}
Recall that we must have $2c_1R^{\kappa-f}\geq R$ for this to be possible. As before, denote by $\phi$ the outermost closed loop and let $v_1,\ldots,v_m$ be the associated defects.
\item[$\bullet$] We cover the annulus $A\!\left(x;2n^{-1}R^{\kappa-f},\frac{1}{2}n^{-1}R^{\kappa-f+1}\right)$ with at most $2R^{4}c_1^{-2}$ boxes of side length $c_1R^{-f}$. Let $M_2$ be large enough and $a_2,b_2,c_2$ small enough depending on $p',M_1,a_1,b_1,\allowbreak c_1$ and $R$ such that for each box $B$ in the covering , we have $\P(\widehat{F}(2B\setminus B))\geq 1-\frac{p'c_1^2}{16R^{4}}$ if $f$ is such that $c_2\cdot c_1R^{-f}n\geq R$. We assume that $\widehat{F}(2B\setminus B)$ and $\widehat{F}'(2B\setminus B)$ hold for all boxes $B$ in our covering. Let $\mathscr{C}^{\bullet}(B)$ and $\mathscr{C}^{\circ}(B)$ denote the collections of outermost loops of closed/open clusters in $2B\setminus B$ whose existence is guaranteed by $\widehat{F}$ and $\widehat{F}'$ respectively. 
\item[$\bullet$] There does not exist a point $w$ in the grid of side length $c_1c_2R^{-f}$ inside the annulus $A(x;\frac{3}{2}n^{-1}R^{\kappa-f},\frac{3}{4}n^{-1}R^{\kappa-f+1})$ such that there are $6$ arms of alternating colors crossing $A(w;2c_2c_1R^{-f},\frac{1}{2}a_2c_1R^{-f})$. As before, using Lemma~\ref{lem:arm_exponents} we see that if $c_2$ is chosen small enough (depending on $p',R$ and $a_2$), then the probability of this event is smaller than $p'/8$ for all annuli for which $c_1c_2R^{-f}n\geq R$. 
\item[$\bullet$] There does not exist a point $w$ in the grid of equally spaced points at distance $c_1R^{-f}$  around $\partial B(x,4n^{-1}R^{\kappa-f})$, such that there are $3$ arms crossing $A\!\left(w;4c_1R^{-f},\frac{b_1}{400}R^{-f}\right)\cap B(x,4n^{-1}R^{\kappa-f})$. 
Again, using half-plane exponents from Lemma~\ref{lem:arm_exponents} we see that by choosing $c_1$ small enough, the probability of this event can be made smaller than $p'/8$.
\item[$\bullet$]Consider also the grid of equally spaced points at distance $\frac{c_1c_2}{8}R^{-f}$ away from each other. Consider all the boxes of side length $\frac{c_1c_2}{8}R^{-f}$ centered at these points that are contained inside of $A(x;\frac{3}{2}n^{-1}R^{\kappa-f},\frac{3}{4}n^{-1}R^{\kappa-f+1})$. The number of such boxes is at most $36(c_1c_2)^{-2}R^{4}$. Let $\tilde{p}=\frac{(c_1c_2)^{2}p'}{144 \cdot R^{4}}$ and let $\xi=\xi(\tilde{p})$ be the value of $\xi$ from Lemma~\ref{lem:omega_is_well_separated}. We require that for all boxes $B$ as above we have that $\Lambda_{\textrm{out}}(2B\setminus B)$ is $\xi$-well separated. We will call this event $E_\textrm{sep}$. By the union bound and Lemma~\ref{lem:omega_is_well_separated}, we know that $\P(E_\textrm{sep})\geq 1-p'/4$ whenever $\xi \frac{c_1c_2}{8}R^{-f}n\geq R$.  
\end{enumerate}
Note that the event $F(A(x;n^{-1}R^{\kappa-f},n^{-1}R^{\kappa-f+1})$ only depends on the percolation configuration inside the annulus $A(x;n^{-1}R^{\kappa-f},n^{-1}R^{\kappa-f+1})$. Thus, for $f,g$ such that $f\neq g$, the events $F(A(x;n^{-1}R^{\kappa-f},n^{-1}R^{\kappa-f+1}))$ and $F(A(x;n^{-1}R^{\kappa-g},n^{-1}R^{\kappa-g+1}))$ are independent.

\begin{lem} The following holds for $R \in \N$ large enough. Fix $p'>0$. Consider any annulus $A(x;n^{-1}R^{\kappa-f},n^{-1}R^{\kappa-f+1}) \subseteq \Lambda_n$, where $x$ is a site in $\Lambda_n$ and $f\in \mathbf{N}$. There exists a choice of parameters $M_1,a_1,b_1,c_1,M_2,a_2,b_2,c_2,\xi$ depending only on $p'$ and $R$ such that $$\P\!\left(F(A(x;n^{-1}R^{\kappa-f},n^{-1}R^{\kappa-f+1}))\right)\geq 1-p'$$ for all annuli such that $\xi\frac{c_1c_2}{8}\cdot R^{\kappa-f-1}\geq R$.
\label{lem:bound_on_F} 
\end{lem}
\begin{proof}
Lemmas~\ref{lem:number_crossings}--\ref{lem:points_are_shielded_away} guarantee that there is a choice of $M_1,a_1,b_1,c_1,M_2,a_2,b_2,c_2$ so that the probabilities of the first two events in the definition of $F$ satisfy the above bounds. Moreover, upon possibly making $c_1$ and $c_2$ smaller, the events in the third and fourth items above also have probability smaller than $p'/8$ by an arm exponent calculation. The discussion in the last item in the definition of $\hat{F}$ shows that there is a choice of $\xi=\xi(p',R,c_1,c_2)$ such that $\P(E_{\textrm{sep}})\geq 1-\frac{p'}{4}$.

Thus, using the union bound, we find that 
\begin{equation}\P(F^c)\leq \frac{p'}{4}+4R^{4}c_1^{-2} \cdot \frac{p'c_1^2}{16R^{4}}+\frac{p'}{8}+\frac{p'}{8}+\frac{p'}{4}=p' \label{eq:bound_on_F}
\end{equation}for all annuli $A(x;n^{-1}R^{\kappa-f},n^{-1}R^{\kappa-f+1})$ such that $c_1c_2R^{-f}n\geq R$ and $\xi \frac{c_1c_2}{8} R^{-f}n\geq R$. Both of these conditions are satisfied if $\xi \frac{c_1c_2}{8}R^{\kappa-f-1}\geq R$, which concludes the proof.  
\end{proof}
We now know that there exists an event $F$ that has probability arbitrarily close to $1$, and on which percolation in our box is well-behaved. As hinted earlier, this will allow us to prove the following lemma, which gives us control on the distance between pivotals in annuli.

\begin{prop}\label{lem:distancebetweencrossings}
Consider percolation on some annulus $A(x;n^{-1}R^{\kappa-f},n^{-1}R^{\kappa-f+1})\subseteq \Lambda_n$. Let $\phi,\{v_i\}_i$ be the loop with defects corresponding to the annulus $A(x;2n^{-1}R^{\kappa-f},4n^{-1}R^{\kappa-f}) \subseteq \Lambda_n$, as defined in Lemma~\ref{lem:decompositioncrossings}. Furthermore, consider the event $F=F(A(x;n^{-1}R^{\kappa-f},n^{-1}R^{\kappa-f+1}))$ defined as above with parameters $M_1,a_1,b_1,c_1,M_2,a_2,b_2,c_2,\xi$. Let $0\leq s\leq 1$. Define $G=G(x,s,f)$ to be the event that $F$ happens and that for every $v_i,v_j$ that are connected inside $\intt\phi$, the (graph resp.\ effective resistance) distance between them in $\operatorname{int} \phi$ is at most $s\altq_n$.

Let $\alpha$ be the exponent of the a priori estimate. Then for all $\zeta>0$, there exists a constant $C=C(\alpha,\zeta,R,M_1,a_1,b_1,c_1,M_2,a_2,b_2,c_2,\xi)$ such that for all $R$ sufficiently large depending on $\zeta$, and $p$ small enough depending on $R$,
\begin{equation}\label{eq:distance-crossings-bound}
    \P(G^{c}\mid F)\leq Cs^{-(\frac{4}{3}+\zeta)\alpha} R^{-\alpha f},
\end{equation} for all levels $f$ such that $\xi\frac{c_1c_2}{8}\cdot R^{\kappa-f-1}\geq R$. 
\end{prop}

Proposition \ref{lem:distancebetweencrossings} will play an important role in Sections \ref{sec:superpolynomial-concentration} and \ref{sec:tightness}. It is the main tool to control distances between two points $x,y\in\Lambda_n$ that are in the same cluster. It will be used as follows. If $\gamma$ is a path between such $x$ and $y$, we will cover $\gamma$ with annuli of various sizes. By a union bound over Lemma~\ref{lem:bound_on_F} and Proposition~\ref{lem:distancebetweencrossings}, on a global event with high probability, we can cover space with annuli for which $F\cap G$ holds. On this event, each open path $\gamma$ crosses each such annulus a bounded number of times, and the distance between each pair of consecutive crossings is bounded by the event $G$. Thus, we can `string together' the distances across all the annuli that $\gamma$ crosses and obtain a bound on $d_n(x,y)$. This stringing argument is first given in the proof Lemma~\ref{lem:crossingbubble}, and we refer to the proof sketch there for an explanation of the stringing procedure.

We note that the value of $R$ (and therefore also $p$) for which~\eqref{eq:distance-crossings-bound} holds may depend on $\zeta$ but not on $\alpha$. Later on, we will apply Proposition~\ref{lem:distancebetweencrossings} with some universal $\zeta$ chosen sufficiently small, independently of $\alpha$.

The main idea in the proof of Proposition~\ref{lem:distancebetweencrossings} is to compare the distance between $v_i$ and $v_j$ inside $\phi$ to the `crossing distance' from the a priori estimate corresponding to the box
\[ B^+=B\!\left(x,\frac{1}{2}n^{-1}R^{\kappa-f+1}\right) . \]

Recall the definition of the event $\mathcal{E}(B^+)$. To make this comparison, we need to:
\begin{enumerate}\phantomsection\label{list:requirements_to_compare_to_apriori}
\item[$\bullet$] `close' all the other pivotals $v_k$ for $k\neq i,j$;
\item[$\bullet$] connect $v_i$ and $v_j$ via open paths to segments of length $2\delta n^{-1} R^{\kappa-f+1}$ centered at midpoints of left and right sides of $B^+$ respectively;
\item[$\bullet$] connect $\phi_i$ and $\phi_j$ via closed paths to segments of length $2\delta n^{-1} R^{\kappa-f+1}$ centered at midpoints of top and bottom sides of $B^+$ respectively;
\item[$\bullet$] create closed left-to-right crossings of $B^+$ above and below the open crossing going through $v_i$ and $v_j$;
\item[$\bullet$] create blocking closed clusters in the two annuli $A(c_1^+;\delta n^{-1} R^{\kappa-f+1}, 2\delta n^{-1} R^{\kappa-f+1})$ and $A(c_3^+;\delta n^{-1} R^{\kappa-f+1}, 2\delta n^{-1} R^{\kappa-f+1})$, where, as before, $c_1^+$ and $c_3^+$ denote the centers of the top and bottom sides of the box $B^+$ respectively;
\item[$\bullet$] create blocking open clusters in the two annuli $A(c_2^+;\delta n^{-1} R^{\kappa-f+1}, 2\delta n^{-1} R^{\kappa-f+1})$ and in $A(c_4^+;\delta n^{-1} R^{\kappa-f+1}, 2\delta n^{-1} R^{\kappa-f+1})$, where $c_2^+$ and $c_4^+$ denote the centers of the right and left sides of the box $B^+$ respectively.  
\end{enumerate}
Fix some $v_i$ and $v_j$ that are connected in $\operatorname{int}\phi$ and assume that $i<j$. Recall that Lemma~\ref{lem:exterior_is_nice} tells us that on the event $F$, for each $k$ we can find a corridor consisting of consecutive boxes of side length $b_1R^{-f}$ connecting $\phi_k$ and $\psi_k$ to the outside boundary of $B(x,4n^{-1}R^{\kappa-f})$, that is at distance at least $b_1R^{-f}$ away from $\phi\setminus \phi_k$. Consider the two corridors $\mathcal{C}_i$ and $\mathcal{C}_j$ corresponding to $\phi_i,\psi_i$ and $\phi_j,\psi_j$ respectively.

We divide the proof into three steps. Our first goal is to use these two corridors to construct four disjoint thinner corridors that, in clockwise order, connect $\psi_i,\phi_i,\psi_j$ and $\phi_j$ to $\partial B^+$ `without interfering with $\intt \phi$'. In the next step of the proof, we will use that on the event $F$, each of the boxes $B$ in these corridors is disconnected from $(2B)^c$ by a finite string of closed (and open) clusters in $2B\setminus B$. This will allow us to find \textit{finite} sequences of neighboring closed/open clusters connecting $\phi_i,\phi_j / \psi_i,\psi_j$ to $\partial B^+$. The final step of the proof is to resample connections between these clusters to form the long closed/open paths that we require to compare the desired distance to the distance in the toy model.

\textbf{Step 1:} \textit{Construction of \textit{disjoint} corridors.}

As mentioned above, we will use the corridors $\mathcal{C}_i$ and $\mathcal{C}_j$ to construct the desired \textit{disjoint} corridors connecting $\phi_i,\phi_j,\psi_i,\psi_j$ to $\partial B^+$. We would like to construct these corridors within $\mathcal{C}_i$ and $\mathcal{C}_j$, respectively. However, it may be the case that $\phi_i\cap \mathcal{C}^i$ consists of too few sites, making it impossible to find two spaced-apart sites of $\phi_i$ and $\psi_i$, respectively, inside of $\mathcal{C}_i$. Therefore we will modify them slightly.

First, we can assume that the boxes in $\mathcal{C}_i$ and $\mathcal{C}_j$ do not intersect $\phi$ except for the first one (otherwise just discard the initial ones). Moreover, we can assume that $\mathcal{C}_i$ and $\mathcal{C}_j$ do not cross upon intersecting (i.e.\ `locally make an X-shape'). Otherwise we let them merge instead, i.e.\ when $\mathcal{C}_j$ intersects $\mathcal{C}_i$ (this can only happen after the first box due to the requirement in Lemma~\ref{lem:exterior_is_nice}), we replace its remaining part with the remaining part of $\mathcal{C}_i$. This modified corridor is still disjoint from $\phi$ except for its first box.

Let $\mathcal{C}_i^{1/2}$ be the corridor in the grid of side length $\frac{1}{2}b_1R^{-f}$ consisting of boxes $\{\frac{1}{2}B:B\in \mathcal{C}_i\}$ together with the intermediate boxes joining $\frac{1}{2}B$ and $\frac{1}{2}B'$ for any two $B,B'$ consecutive in $\mathcal{C}_i$. In particular $\mathcal{C}_i^{1/2}\subseteq \mathcal{C}_i$. If $B^{1/2,i}_1,\ldots, B^{1/2,i}_t$ denote the boxes in the corridor $\mathcal{C}_i^{1/2}$, let $\mathcal{C}_i':=\!\left(3B_1^{1/2,i},B_3^{1/2,i},B_4^{1/2,i},\ldots,B_t^{1/2,i}\right)$. Define $\mathcal{C}_j'$ analogously.
By construction we know that:
\begin{enumerate}
\item[$(1)$] consecutive boxes in $\mathcal{C}_i'$ share a side and distinct boxes in the corridor are disjoint ({except for} boundary intersections);
\item[$(2)$] $\phi_i\cap2B_1^{1/2,i}\neq \emptyset$, $\psi_i\cap\overline{2B_1^{1/2,i}}\neq \emptyset$ and $B^{1/2,i}_t\cap \partial B(x,4n^{-1}R^{\kappa-f})\neq \emptyset$;
\item[$(3)$] No box of $\mathcal{C}_i'$ intersects $\phi\setminus \phi_i$, and no box other than $3B_1^{1/2,i}$ intersects $\phi_i$;
\item[$(4)$] $3B_1^{1/2,i}$ and $\mathcal{C}_j'$ are disjoint.  
\end{enumerate}

The same holds with $i$ and $j$ swapped. Further, $\mathcal{C}_i'$ and $\mathcal{C}_j'$ do not cross in case they intersect.

We claim the following. Inside of the set $B(x,4n^{-1}R^{\kappa-f})^c\cup \!\left(\bigcup \!\left\{B : B\in \mathcal{C}_i'\cup \mathcal{C}_j'\right\}\right)$, we can find $4$ disjoint thinner corridors consisting of adjacent boxes in the grid of side length $c_1R^{-f}$ that are in clockwise order as follows:
\begin{enumerate}
\item[$\bullet$] a corridor connecting $\psi_i$ to the center of the left-hand side of $B^{+}$ that does not come within distance $10 c_1R^{-f}$ from the top and bottom sides of $B^+$;
\item[$\bullet$] a corridor connecting $\phi_i$ to the center of the top side of $B^{+}$;
\item[$\bullet$] a corridor connecting $\psi_j$ to the center of the right-hand side of $B^{+}$ that does not come within distance $10 c_1R^{-f}$ from the top and bottom sides of $B^+$;
\item[$\bullet$] a corridor connecting $\phi_j$ to the center of the bottom side of $B^+$.
\end{enumerate}

\begin{figure}[ht]
    \centering
    \includegraphics[width=1\linewidth]{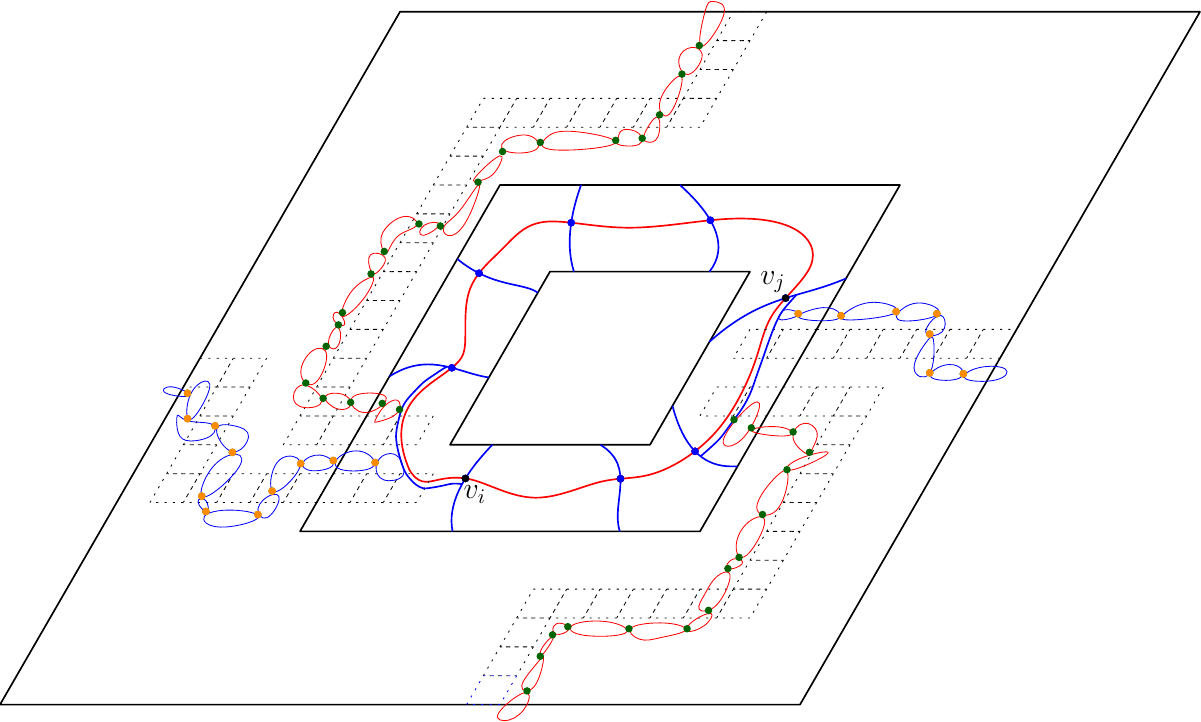}
    \caption{This figure illustrates the 4 corridors connecting different segments of $\phi$ to $\partial B_+$ - these are represented by the dotted boxes. The many thinner blue and red curves represent the sequences of open and closed clusters in the annuli around these boxes that we will use in our resampling procedure. In particular, orange points are the closed sites where two consecutive open clusters touch. We will resample our configuration around these points with the goal of joining these open clusters together. These are the points of $\mathscr{S}^{\circ}$ corresponding to the corridors of $\phi_i$ and $\phi_j$. Similarly, green points represent the points in $\mathscr{S}^{\bullet}$. Finally, blue points are the points of $\{v_k: k\neq i,j\}$. }
    \label{fig:resampling_connections}
\end{figure}

We can additionally guarantee that the first two corridors consist of boxes $B$ such that $2B$ does not intersect $\phi\setminus \phi_i$, and the second two corridors consist of boxes $B$ such that $2B$ does not intersect $\phi\setminus \phi_j$. Moreover, we can choose these corridors so that they are at least at distance $\frac{b_1}{100}R^{-f}$ 
away from each other, since we can choose $c_1\ll b_1$. Finally, we can also choose these so that for all $l=1,\ldots,4$ the corridor going to $c_l^+$ does not come $2\delta n^{-1}R^{\kappa-f+1}$ close to any other $c_{\tilde{l}}^+$. See Figure~\ref{fig:resampling_connections} for an illustration.

\begin{figure}[ht]
    \centering
    \includegraphics[width=1\linewidth]{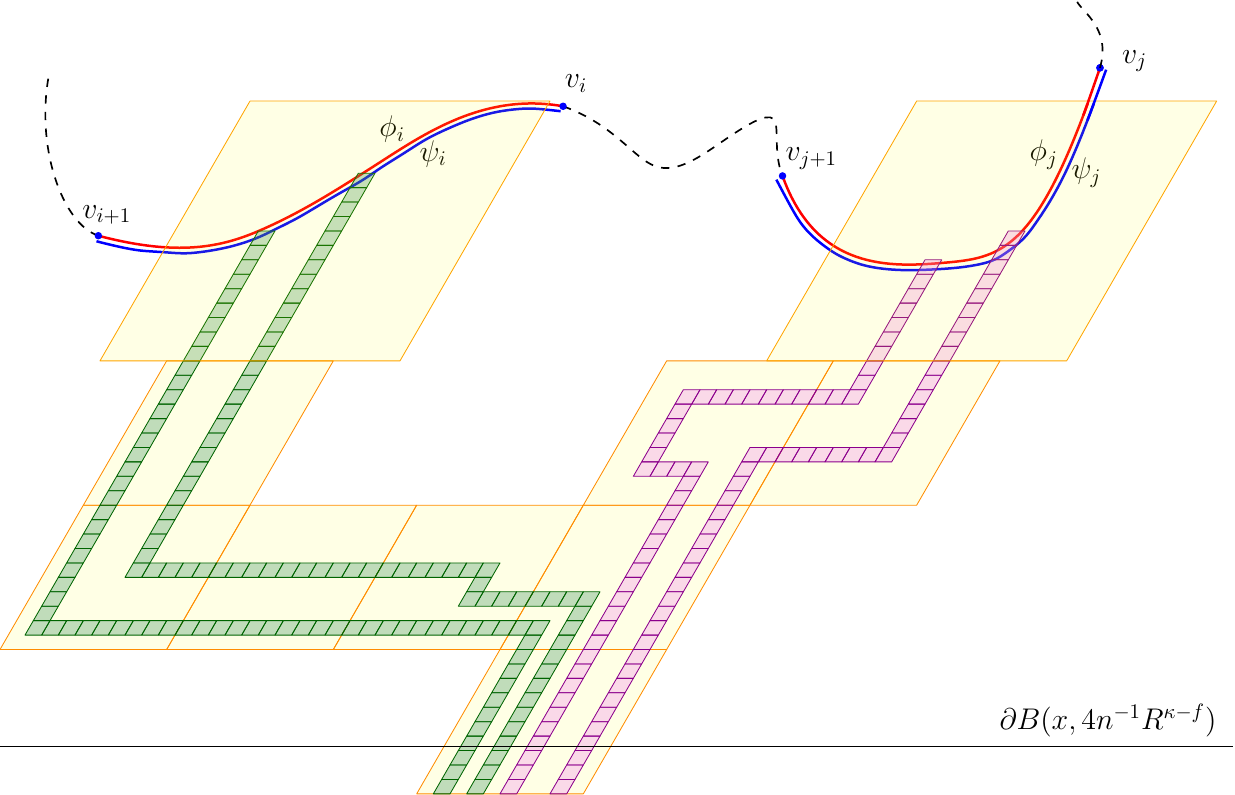}
    \caption{The figure above illustrates how one constructs the 4 thinner corridors in the case that neither $\phi_i$ nor $\phi_j$ intersect $\partial B(x,4n^{-1}R^{\kappa-f})$. The red/blue interfaces represent the curves $\phi_i,\psi_i$ and $\phi_j,\psi_j$, while the dashed line represents the remaining portion of $\phi$. The yellow boxes illustrate the modified corridors $\mathcal{C}_i'$ and $\mathcal{C}_j'$ that join at some point (in the case above $B_6^i=B_4^j$). Note that the first box in each of the yellow corridors is larger - its diameter is three times bigger than the diameter of other boxes (the picture is not to scale). The two green thinner corridors are the ones corresponding to $\phi_i$ and $\psi_i$, while the two purple ones correspond to $\phi_j$ and $\psi_j$. Note that for the most part, all 4 corridors keep their distance to the boundary of corridors $\mathcal{C}_i'$ and $\mathcal{C}_j'$ fixed. There are two exceptions: in the boxes $B_5^i$ and $B_3^j$ we need a slight deviation to make sure that the continuations of all four corridors stay disjoint and far apart; and, we must make a choice of boxes inside of $B_1^i$ until we enter the rest of $\mathcal{C}_i'$ and can follow along the boundary.}
    \label{fig:we_can_choose_box_sequence}
\end{figure}

To convince ourselves that this choice is possible, we need to consider separately the cases when $\phi_i$ intersects $\partial B(x,4n^{-1}R^{\kappa-f})$ and when it does not. Assume the latter holds. Let $B_1^i$ be the first box in the corridor $\mathcal{C}_i'$. Note that there exists at least one box of side length $c_1R^{-f}$ inside of $\frac{2}{3}B_1^i$ that contains a point of $\phi_i$. Since $\phi_i$ covers macroscopic distance $a_1R^{-f} \gg b_1 R^{-f}$, $c_1\ll b_1$ and $\psi_i$ follows along the entire length of $\phi_i$, it means that there exists at least one more box of side length $c_1 R^{-f}$ inside of $B_1^i$ that intersects $\psi_i$. Moreover, we can choose this box so that it is located at distance at least $\frac{b_1}{100}R^{-f}$ 
away from the first box we chose. Then simply continue adding boxes to the first one, keeping the distance from the smaller corridor to the boundary of $\mathcal{C}_i'$ fixed. We only deviate from this pattern possibly right at the start (we choose the first $\sim \frac{b_1 c_1^{-1}}{100}$
boxes of both corridors so that their continuations stay suitably far away) and so that they actually enter $B_2^i$, and at the box leading to the `joining box' of $\mathcal{C}_i'$ and $\mathcal{C}_j'$ to make sure that there is enough space for corridors corresponding to $\phi_j$ and $\psi_j$. See Figure~\ref{fig:we_can_choose_box_sequence} for an illustration. As before, we can edit the $4$ corridors that we have obtained so that no box but the first one in each corridor intersects $\phi$. Then, of the two corridors corresponding to $\phi_i$, choose them to correspond to $\phi_i$ and $\psi_i$ in such a way that $v_i$, $\psi_i$-corridor and $\phi_i$ corridor and $v_{i+1}$ are located in clockwise order along $\phi_i$. Do the same for $\phi_j$ and $\psi_j$. This creates  well-separated corridors connecting $\phi_i,\psi_i,\phi_j$ and $\psi_j$ to $\partial B(x,4n^{-1}R^{\kappa-f})$. Extend these corridors arbitrarily within $A(x;4n^{-1R^{\kappa-f}},\frac{1}{2}n^{-1}R^{-\kappa-k+1})$ so that they come close to the respective midpoints of the target sides of $B^+$ without coming close to each other or other sides of $B^+$. Taking $c_1\ll 1$ ensures that this choice is possible. See Figure~\ref{fig:resampling_connections} for an illustration.

\begin{figure}[ht]
    \centering
    \includegraphics[width=0.8\linewidth]{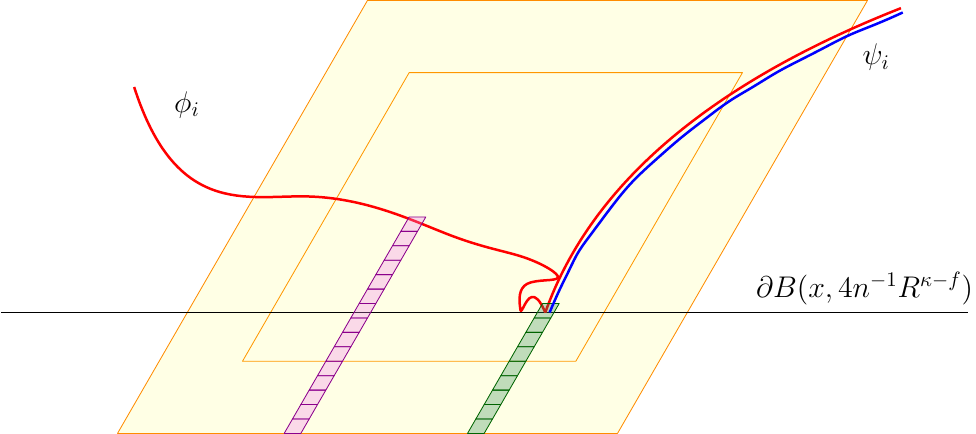}
    \caption{This figure illustrates how we find the 2 corridors corresponding to $\phi_i$ and $\psi_i$ in the case that $\phi_i$ intersects $\partial B(x,4n^{-1}R^{\kappa-f})$. The two yellow boxes represent the box $B_1^i$ that is the single box of corridor $\mathcal{C}_i'$ and $\frac{2}{3}B_1^i$ that is the single box in the original corridor $\mathcal{C}_i$ (the picture is not to scale). The green corridor is the one containing the endpoint of $\psi_i$ belonging to $\partial B(x,4n^{-1}R^{\kappa-f})$. On the picture above, this will be the corridor corresponding to $\psi_i$, while the purple corridor is the one that will correspond to $\phi_i$. Note that we have constructed these corridors by moving straight down along the axis of our coordinate system that is 'perpendicular' to the axis determined by the relevant side of $\partial B(x,4n^{-1}R^{\kappa-f})$. If it were the case that the first purple box we have chosen was directly above the first green one (or if their horizontal coordinates didn't differ by at least $\frac{b_1}{100}R^{-f}$ 
    ), we would first move $\sim \frac{b_1c_{1}^{-1}}{100}$
    steps along the horizontal axis (in our case to the left) and then continue moving down along the 'vertical' axis as determined before. On the picture above, such a corridor would intersect $\phi_i$ one more time at the segment of $\phi_i$ not covered by $\psi_i$. Since we will ultimately truncate this corridor to start from the last box intersecting $\phi_i$, the final picture will look similar to the case already depicted above.}
    \label{fig:we_can_find_paths_if_touches}
\end{figure}

In the case that $\phi_i$ intersects $\partial B(x,4 n^{-1}R^{\kappa-f})$, $\psi_i$ ends at this boundary and we have chosen $\mathcal{C}_i'$ to consist of a single box $B_1^i$ such that the box $\frac{2}{3}B_1^i$ in the grid of side-length $b_1R^{-f}$ contains this point. Consider the box of side length $c_1R^{-f}$ that contains the endpoint of $\psi_i$ belonging to $\partial B(x,4 n^{-1}R^{\kappa-f})$. This box on its own forms a corridor connecting $\psi_i$ to $\partial B(x,4n^{-1}R^{\kappa-f})$. To find a corridor corresponding to $\phi_i$, note again that the section of $\phi_i$ not covered by $\psi_i$ covers a distance of at least $a_1R^{-f}\gg b_1R^{-f}$. Thus, there exists at least one more box of side length $c_1R^{-f}$ inside of $B_1$ that intersects this section of $\phi_i$. Moreover, we can choose this box so that it is at distance at least $\frac{b_1}{100}R^{-f}$
from the first one. Consider the corridor starting at this box moving along one of our two coordinate axes 'straight towards' the target side of $\partial{B(x,4n^{-1}R^{\kappa-f})}$. This, together with the original box containing the endpoint of $\psi_i$ will produce two spaced-apart corridors that join $\psi_i$ and $\phi_i$ to $\partial B(x,4n^{-1}R^{\kappa-f})$ unless the horizontal coordinates of the initial boxes do not differ by more than $\frac{b_1}{100}R^{-f}$. 
In this case, first move horizontally $\sim\frac{b_1c_1^{-1}}{100}$
places (either left or right depending on where the segment of $\phi_i$ not covered by $\psi_i$ is relative to $\psi_i$) and then continue to go downwards. See Figure~\ref{fig:we_can_find_paths_if_touches} for an illustration. As before, we truncate the resulting corridor to start from the last box intersecting $\phi_i$. Finally, we extend these corridors within $A(x;4n^{-1}R^{\kappa-k},n^{-1}R^{\kappa-k+1})$ to reach the respective midpoints of the sides of $B^+$. We do so arbitrarily, with the only requirement being that we 'get away from' $\partial B(x,4n^{-1}R^{\kappa-f})$ immediately upon leaving the box. More precisely, if we consider the box containing the endpoint of $\psi_i$, we choose the next box in the $\psi_i$-corridor to be a neighbor of this first box that does not intersect $B(x,4n^{-1}R^{\kappa-f})$. We choose the remaining boxes $B$ in the corridor to be such that $2B\cap B(x,4n^{-1}R^{\kappa-f})=\emptyset$. In Figure~\ref{fig:we_can_find_paths_if_touches}, this is illustrated by taking a few steps along the axis in the 'outwards' direction as soon as we hit $\partial B(x,4n^{-1}R^{\kappa-f})$ with our corridor. As mentioned above, the rest of the corridor is arbitrary as long as it is away from $\partial B(x,4n^{-1}R^{\kappa-f})$.

\medskip \textbf{Step 2:} \textit{Construction of cluster sequences.}

\medskip We can now use these corridors to find finite strings of consecutive closed/open clusters connecting $\phi_i/\phi_j$ and $\psi_i/\psi_j$ to $\partial B^+$. Namely, on $F$, we know that for every box $B$ in the $c_1R^{-f}$-covering of $B^+$ there exist at most $M_2$ macroscopic closed clusters in $2B\setminus B$ that touch each other and, together with their touch points, disconnect $\partial B$ from $\partial 2B$. The same holds for open clusters. In particular, these clusters can be used to string a sequence of macroscopic closed/open clusters connecting $\phi_i/\psi_i$ or $\phi_j/\psi_j$ to the boundary of $B^+$.

Let us focus on the corridor connecting $\phi_i$ to $c_1^+$ -- the center of the top side of $B^+$. Denote the boxes in this corridor by $B_1,\ldots, B_t$ where $B_1\cap \phi_i$ is non-empty. We will inductively build a string $\mathscr{S}$ of closed clusters connecting $\phi_i$ to the top side of $B^+$ such that consecutive closed clusters in $\mathscr{S}$ touch each other. Start with $\mathscr{S}=\emptyset$. Note that the word `cluster' here refers to the \emph{local} clusters for percolation restricted to a box $B_s$. Two different such clusters might be part of the same \emph{global} cluster when considering percolation on the whole of $B^+$, but we still consider them to be different clusters in the following section.

Because $\phi_i$ has length at least $a_1R^{-f}\gg 2c_1R^{-f}$ and it intersects the interior of $B_1$, we know that $\phi_i$ crosses $2B_1\setminus B_1$. Thus, one of the closed clusters in $\mathscr{C}^{\bullet}(B_1)$ must be a closed local cluster of $\phi_i$. Moreover, box $B_2$ crosses the annulus $2B_1\setminus B_1$ and so some clusters of $\mathscr{C}^{\bullet}(B_1)$ must intersect $B_2$. In fact, this string of clusters must cross the annulus $2B_2\setminus B_2$. Append to $\mathscr{S}$ the consecutive clusters of $\mathscr{C}^{\bullet}(B_1)$ starting from $\phi_i$ going through the interior of $B_2$ and ultimately crossing $2B_2\setminus B_2$. In turn, the clusters of $\mathscr{C}^{\bullet}(B_2)$ must cross the clusters of $\mathscr{S}$. Continue inductively to append clusters of $\mathscr{C}^{\bullet}(B_i)$'s until we get to the final box $B_t$. Assume that we have chosen $c_1<\delta$ (recall $\delta$ from the definition of the event $\mathcal{E}(B^+)$). Then the clusters of $\mathscr{C}^{\bullet}(B_t)$ must intersect the top side of $B^+$ within $\delta n^{-1}R^{\kappa-f+1}$ from the center of the side, since $B_t$ contains the center. We stop appending clusters to $\mathscr{S}$ once we have reached the top side of $B^+$. 

Consider the last cluster in $\mathscr{S}$ that intersects $\phi_i$. We claim that none of the subsequent clusters can intersect $\operatorname{int}\phi$. To see this, note that $\mathscr{S}$ is a string of clusters such that the first cluster intersects $\phi_i$ and the last cluster intersects the outer boundary of the box $B^+$. Importantly, $\partial B^+$ is fully contained in the exterior of $\phi$. Because our corridor is far away from $\phi\setminus \phi_i$, it means that the only way to `escape' $\operatorname{int}\phi$ is through $\phi_i$. In particular, if a cluster in $\mathscr{S}$ intersects $\operatorname{int}\phi$, it means that there must be a subsequent cluster that intersects $\phi_i$, proving our claim.

We now discard all the clusters preceding this one from $\mathscr{S}$. This results in a string of consecutive closed clusters, the first of which is connected to $\phi_i$, the last of which hits the upper side of $B^+$ within $\delta n^{-1}R^{\kappa-f+1}$ from its center and such that none of the clusters, apart from the first one, intersect $\operatorname{int}\phi$.

\begin{figure}[ht]
    \centering
    \includegraphics[scale=0.8]{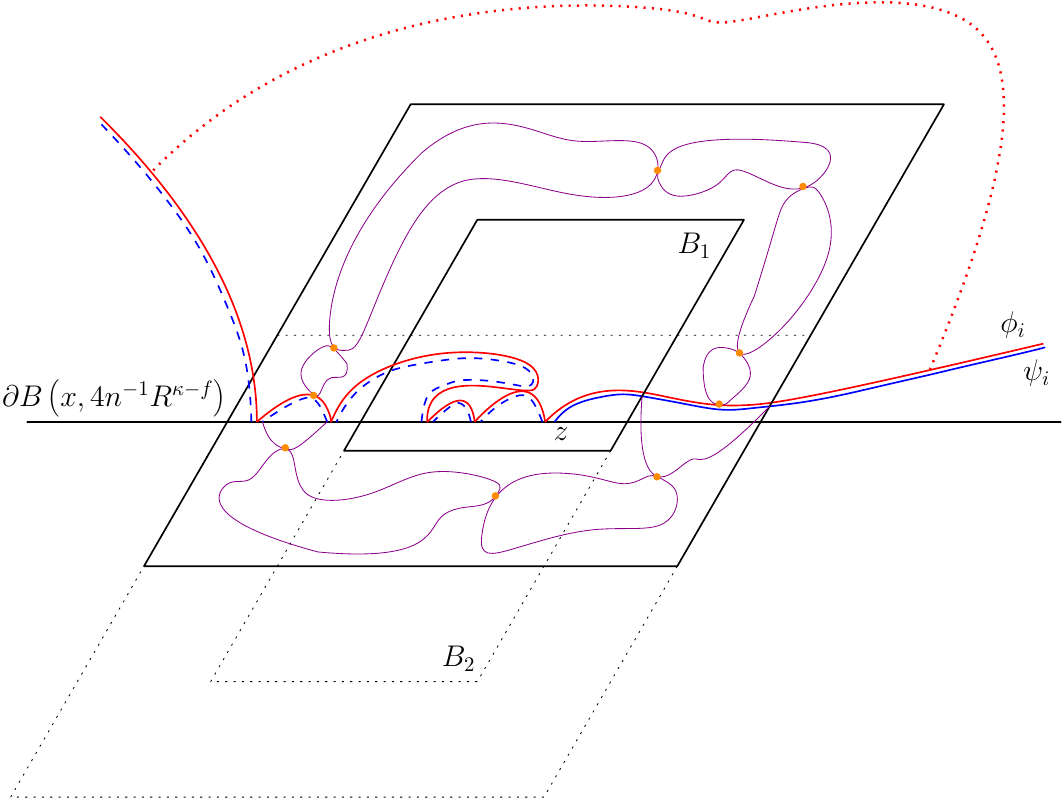}
    \caption{The figure above illustrates the construction of the string of open clusters corresponding to $\psi_i$ in the case that $\psi_i$ ends at a point $z$ on $\partial B(x,4n^{-1}R^{\kappa-f})$. The red curve above represents $\phi_i$, while the full blue curve represents $\psi_i$. The dashed blue curve represents the open cluster boundaries adjacent to $\phi_i$ that are not a part of $\psi_i$. The figure above also shows the box $B_1$, which is the first box in the corridor corresponding to $\psi_i$ and it contains the point $z$. Purple curves above represent the sequence of open clusters that form an open loop with closed defects (orange touch-points) inside the annulus $2B_1\setminus B_1$. We can see above that even after truncating our sequence of clusters after the last cluster intersecting $\psi_i$, the modified sequence can still contain open clusters intersecting $\operatorname{int}\phi$ -- this can happen if these clusters touch each other at a point on the section of $\phi_i$ not covered by $\psi_i$, as on the figure above. In Lemma~\ref{lem:pivotals_are_far}, we will show that on the event $F$, there exists a closed path `shielding' the box $2B_1$ away from $\operatorname{int} \phi$ represented by the dotted red line above. This will ensure that resampling the cluster connections in $2B_1\setminus B_1$ will not influence the target distances between points $v_i,v_j$, and it will not disconnect the closed corridor from the point $v_i$.}
    \label{fig:choosing_clusters_for_psi_i}
\end{figure}

We can repeat this construction for the other three corridors. Even though the construction is identical for the open corridors, we need to take more care when dealing with clusters that could enter $\operatorname{int} \phi$. Assume without loss of generality that we are dealing with the corridor corresponding to $\psi_i$. In the case that $\phi_i$ does not intersect $\partial B(x,4n^{-1}R^{\kappa-f})$, the situation is the same as for the closed corridors -- we can show analogously as before that we can choose these clusters so that they do not intersect $\operatorname{int} \phi$. If $\phi_i$ does intersect $\partial B(x,4n^{-1}R^{\kappa-f})$, this need not be the case. Suppose that $\psi_i$ ends at some point $z\in \partial B(x,4n^{-1}R^{\kappa-f})$. Recall that in this case, we have chosen the first box in our corridor to be precisely the box $B_1$ such that $z\in \bar{B}_1$. It could indeed happen that the clusters we have chosen in $\mathscr{S}$ enter $\operatorname{int} \phi$ due to the fact that $\psi_i$ does not cover the entire length of $\phi_i$. Nonetheless, we will show in Lemma~\ref{lem:pivotals_are_far} that resampling these clusters will not affect any distances in $\operatorname{int}\phi$ between points $v_1,\ldots,v_m$ on the event $F$. Moreover, we will see that resampling all these connections will not disconnect the part of $\phi_i$ from which our $\phi_i$ corridor is starting and the point $v_i$, even though it may disconnect the path $\phi_i$ itself. See Figure \ref{fig:choosing_clusters_for_psi_i} for an illustration.

Recall that we started constructing these cluster sequences with the goal of comparing the distance between $v_i$ and $v_j$ to the distance from the a priori estimate in the box $B^+$. To do so, we had a list of requirements (see page~\pageref{list:requirements_to_compare_to_apriori}). Closing pivotals $\{v_k:k\neq i,j\}$ will address the first item on the list. Resampling around pivotals in the four cluster sequences that we have just constructed will ensure that conditions 2 and 3 are satisfied. To deal with conditions 4, 5, and 6, simply note that we can consider sequences of boxes in the grid of side length $c_1R^{-f}$: two horizontal ones corresponding to the two left-to-right crossings and four semi-circular ones corresponding to the $4$ blocking clusters around points $c_l^+$. Linking the relevant sequences of closed/open clusters in these annuli would then yields the desired paths. Importantly, note that due to the way we have constructed the $4$ corridors, the additional sequences of boxes arising from conditions 4--6 can intersect these corridors, but only if the clusters corresponding to the corridor and the sequence of boxes are both closed or both open.

Let $\mathscr{S}^{\bullet}$ denote all open touching points of the relevant sequences of adjacent closed clusters, whose existence is guaranteed by the second item in the definition of $F(A(x;n^{-1}R^{\kappa-f},n^{-1}R^{\kappa-f+1}))$. In particular, $\mathscr{S}^{\bullet}$ contains points corresponding to the corridors of $\phi_i$ and $\phi_j$, as well as touch points of clusters needed to form top and bottom closed crossings of $B^+$ and closed blocking clusters around $c_1^+$ and $c_3^+$.  Define $\mathscr{S}^{\circ}$ similarly.  
\begin{lem} Assume that the event $F$ holds. Then
\begin{enumerate}[(a)] 
\item\label{property:distance-points} No two distinct points of $\mathscr{S}^{\bullet}\cup \mathscr{S}^{\circ}\cup \{v_k:k\neq i,j\}$ are within distance $c_1c_2R^{-f}$ of each other;
\item\label{property:changing-state} Given any $z\in \mathscr{S}^{\bullet}\cup \mathscr{S}^{\circ}\cup \{v_k:k\neq i,j\}$, changing the state of sites in $B\!\left(z,\frac{1}{2}c_1c_2R^{-f}\right)$ does not change the distance between $v_i$ and $v_j$ in $\intt\phi$;
\item\label{property:resampling-connection} Upon resampling all of the open connections between clusters in $\mathscr{S}^\circ$, the clusters of $\mathscr{S}^{\bullet}$ are still connected via closed paths to both $v_i$ and $v_j$.
\end{enumerate} \label{lem:pivotals_are_far}
\end{lem}
\begin{proof} We begin by proving part~\eqref{property:distance-points}. Suppose that $z_1,z_2 \in \mathscr{S}^{\bullet}\cup \mathscr{S}^{\circ}\cup \{v_k:k\neq i,j\}$ are such that $\|z_1-z_2\|_1\leq c_1c_2R^{-f}$. By the definition of events $F_2$ and $F_3$, we know that both $z_1$ and $z_2$ have $4$ arms of alternating colors emanating from them, up to distance $a_2\cdot c_1R^{-f}$. In particular, if $w$ is a point in the grid of side length $c_1c_2R^{-f}$ such that $B(w,c_1c_2R^{-f})$ contains the point $z_1$, then there are $6$ arms of alternating colors in $A(w;2c_2c_1R^{-f},\frac{1}{2}a_2c_1R^{-f})$. Indeed, $4$ of them are the arms emanating from $z_1$. Site $z_2$ must be in one of the $4$ sector-like components of $B(w,2c_2c_1R^{-f})$ cut out by the arms of $z_1$. Two out of its $4$ arms can merge with the existing arms of $z_1$ delimiting the sector of $z_2$, but the other two must be disjoint, yielding the remaining $2$ arms in our annulus. However, this contradicts third item in the definition of the event $F$. See the left-most picture on Figure~\ref{fig:resampling_points_are_far_away}.  

\begin{figure}[ht]
    \centering
    \includegraphics[width=1\linewidth]{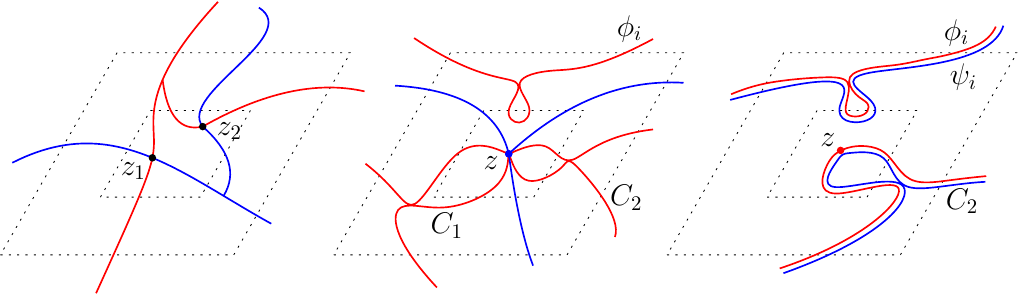}
    \caption{Left picture above illustrates the $6$ arm event from the proof of part~\eqref{property:distance-points} of the lemma. On this picture, we have chosen to illustrate the events that $z_1$ and $z_2$ are touch points of macroscopic clusters using the $4$ alternating arms emanating from both these points -- no matter if $z_\star$ is a touch-point of closed or open clusters, these clusters along with their boundaries form $4$ arms of alternating color -- these are represented by blue and red curves above. The other two pictures illustrate $6$ arms events in the proof of part~\eqref{property:changing-state} -- the central one concerns the points $z\in \mathscr{S}^{\bullet}$, and the third one depicts $z\in  \mathscr{S}^{\circ}$.}
    \label{fig:resampling_points_are_far_away}
\end{figure}

To see that part~\eqref{property:changing-state} holds, first pick $z\in \{v_k:k\neq i,j\}$. Because we are working on events $F_5(A(x;2R^{-f},4R^{-f})),F_6(A(x;2R^{-f},4R^{-f}))$ and $F_7(A(x;2R^{-f},4R^{-f}))$, we know that there exists a point $\hat{z}$ in $\intt\phi$ such that every open path from the open component in $\intt \phi$ containing $v_i$ and $v_j$ to $B(z,c_2c_1R^{-f})\cap \operatorname{int}\phi$ must go through $\hat{z}$. Thus, changing the states of sites of $B(z,c_2c_1R^{-f})$ clearly does not change the geodesic distance between $v_i$ and $v_j$ in the interior -- the shortest path will not enter $B(z,c_2c_1R^{-f})$ at all. Similarly, the resistance distance remains unchanged since any current flowing through this part of the graph has to go in and out through $\hat{z}$, making the net current flow $0$.

Assume now that $z\in \mathscr{S}^{\bullet}$ and that $B(z,c_2c_1R^{-f})$ intersects $ \operatorname{int}\phi$. The only points of $ \mathscr{S}^{\bullet}$ that could come close to $\operatorname{int}\phi$ are the ones coming from the two corridors corresponding to $\phi_i$ and $\phi_j$. Without loss of generality, assume that $z$ is in the corridor corresponding to $\phi_i$. By the way we have chosen our string of closed clusters, it must be that $\phi_i \cap B(z,c_2c_1R^{-f})\neq \emptyset$. There exists a box $B$ in our corridor such that $z$ is a touch-point of two macroscopic closed clusters in $2B\setminus B$ -- call them $C_1$ and $C_2$. By the definition of the event $\widehat{F}(B)$, we know that $B(z,a_2c_1R^{-f})\subseteq 2B\setminus B$. There are two ways in which $\operatorname{int} \phi$ can intersect $B(z,c_2c_1R^{-f})$. Either there is a segment of $\phi_i$ intersecting $B(z,c_2c_1R^{-f})$ that is in a different closed cluster to both $C_1$ and $C_2$ when restricted to the annulus $2B\setminus B$, or all segments of $\phi_i$ inside of this box intersect $C_1\cup C_2$. 

If the former is true, there are $6$ arms in the annulus $A(z;c_2c_1R^{-f},a_2c_1R^{-f})$: three closed ones corresponding to three macroscopic closed clusters that all intersect $B(z,c_2c_1R^{-f})$ ($C_1,C_2$ and the cluster of the target segment of $\phi_i$) and three open ones delimiting these clusters. See the second picture on Figure~\ref{fig:resampling_points_are_far_away}. However, this contradicts the third item in the definition of the event $F$. If the latter is true, $z$ must be a neighbor of $\phi_i$ and we know that $B(z,c_2c_1R^{-f})$ is `shielded away' by a closed path from $\operatorname{int}\phi$. This is because on the event $F$, we know that the events $F_5,F_6,F_7$ defined in Lemma \ref{lem:points_are_shielded_away} hold for the smaller annulus for which $z$ is one of the open defects in the almost-closed circuit. In particular, the only section of $\operatorname{int} \phi$ that intersects the box $B(z,c_2c_1R^{-f})$ is delimited from the rest of $\operatorname{int} \phi$ by a closed path, and hence cannot influence the geodesic/resistance distances between points $v_i$ in $\operatorname{int} \phi$. Since no other segments of $\phi$ intersect this box, changing the percolation configuration inside of it will not influence the target distance.

Finally, assume that for some $z\in \mathscr{S}^{\circ}$ we have that $B(z,c_2c_1R^{-f})$ intersects $ \operatorname{int}\phi$. Assume again without loss of generality that $z$ is in the corridor corresponding to $\psi_i$. In the case that $\psi_i$ runs along the entire length of $\phi_i$, by our choice of clusters it must be that $\psi_i$ intersects $B(z,c_2c_1R^{-f})$. Also, we know that $z$ is a touch point of two macroscopic open clusters $C_1$ and $C_2$ coming from some box $B$ in our $\psi_i$-corridor. Moreover, at most one of these clusters can intersect $\psi_i$ by construction. So suppose without loss of generality that $C_2$ does not intersect $\psi_i$.

Then there are $6$ arms inside the annulus $A(z;c_2c_1R^{-f},a_2c_1R^{-f})$: one closed one corresponding to $\phi_i$, two open ones coming from $\psi_i$, one open one corresponding to $C_2$ and two closed ones corresponding to the boundary of $C_2$. See the last picture in Figure~\ref{fig:resampling_points_are_far_away}. This cannot happen as it contradicts third item in the definition of $F$.

Assume finally that $\psi_i$ does not run along the entire length of $\psi_i$. Due to the way we have constructed these corridors, it must be that $B(z,c_2c_1R^{-f})\subseteq 2B_1$, where $B_1$ is the first box of the $\psi_i$-corridor. We claim that in this case it must be that there is a closed path inside of the annulus $\frac{b_1c_{1}^{-1}}{200}\cdot B_1\setminus 2B_1$ connecting the segment of $\phi_i$ covered by $\psi_i$, and the segment of $\phi_i$ from its last exit from the box $2B_1$ until reaching $v_{i+1}$. Indeed, if not, we would be able to find 3 disjoint arms inside $\frac{b_1c_{1}^{-1}}{200}\cdot B_1\setminus 2B_1$, and this contradicts the penultimate item in the definition of the event $F$. This closed path guarantees that changing the states of the sites in $B(z,c_2c_1R^{-f})$ does not change the distance between $v_i$ and $v_j$ in $\operatorname{int}\phi$.

Part~\eqref{property:resampling-connection} follows almost instantly from the above considerations. Simply note that if $\phi_i$ does end on $\partial B(x,4n^{-1}R^{\kappa-f})$, the closed path from the previous paragraph guarantees that the clusters of $\mathscr{S}^{\bullet}$ are still connected to $v_i$ via a closed path.
\end{proof}

\medskip \textbf{Step 3:} \textit{Resampling.} \phantomsection\label{page:resampling}

As advertised, we now have finite sequences of clusters that almost create the paths needed for the event $\mathcal{E}(B^+)$ to occur. We will use a resampling procedure to change the linking pattern around every pivotal without affecting the distances in $\operatorname{int} \phi$.

We can describe this procedure explicitly. Consider a percolation configuration $\omega$ inside $B(x,R^{-f+1})$. And suppose $c=\frac{c_1\cdot c_2}{8}>0$ is fixed as above. Consider a grid of points $w$ that are at distance $cR^{-f}$ from each other. Pick $\widetilde{M}\in \mathbf{N}_0$ randomly according to a geometric distribution with some arbitrary parameter that does not matter, independently of $\omega$. Then we repeat the following procedure $\widetilde{M}$ times:
\begin{enumerate}
\item[$\bullet$] Select a point $w$ in our grid uniformly at random and consider $A(w;cR^{-f},2cR^{-f})$. Recall that we defined $\omega_{\textrm{out}}({A(w;cR^{-f},2cR^{-f})})$ as the percolation configuration outside of $B(w,2cR^{-f})$ together with the cluster interface explorations started from $\partial B(w,2cR^{-f})$ and stopped upon hitting $\partial B(w,cR^{-f})$ or returning to $\partial B(w,2R^{-f})$. Moreover, we defined $\Lambda_{\textrm{out}}({A(w;cR^{-f},2cR^{-f})})$ to be the set of sites in $B(w,2cR^{-f})^c$ together with the explored sites in $B(w,2cR^{-f})$. Note that the conditional law of the percolation configuration in $\Lambda_{\textrm{out}}({A(w;cR^{-f},2cR^{-f})})^c$ is still i.i.d.\ Bernoulli percolation.   
\item[$\bullet$] Resample the percolation configuration $\omega$ in the unexplored region $\Lambda_{\textrm{out}}({A(w;cR^{-f},2cR^{-f})})^c$ according to i.i.d.\ Bernoulli distribution. 
\end{enumerate}

Call the resampled configuration $\omega^{\textrm{res}}$. We have determined the set of target pivotals, namely $\mathscr{S}^{\bullet}\cup \mathscr{S}^{\circ}\cup \{v_k:k\neq i,j\}$. This set can have size at most $M_0:=M_1+2R^4c_1^{-2}\cdot M_2$ on the event $F$. 
We say that \textit{the resampling is successful} if:
\begin{enumerate}
\item[$\bullet$] We have chosen $\widetilde{M}$ to be precisely the size of the set of target pivotals;
\item[$\bullet$] The points $w$ that we have chosen uniformly at random from our grid of side length $cR^{-f}$ are precisely such that the boxes $B(w,cR^{-f})$ contain the pivotal sites;
\item[$\bullet$] Within every chosen $B(w,cR^{-f})$, we can find a site $z\in \mathscr{S}^{\bullet}\cup\mathscr{S}^\circ\cup \{v_k:k\neq i,j\}$ that is a touch point of two large (of $\|\cdot\|_1$-diameter at least $a_2c_1R^{-f}$) open or closed clusters in configuration $\omega$. We require that these clusters are joined together in the resampled configuration $\omega^{\textrm{res}}$.
\end{enumerate}
We call the event that resampling is successful $E_{\textrm{res}}$. 

Recall the last item in the definition of the event $F$. This event, $E_{\textrm{sep}}=E_{\textrm{sep}}(p')$, is included in the definition of $F$ precisely because it will guarantee that the resampling is successful with positive probability.

\begin{lem} \label{lem:resampling_is_positive}
There exists a constant $g=g(p',R,M_0,c,\xi)>0$ such that for all percolation configurations $\omega$ on $B(x,n^{-1}R^{\kappa-f+1})$
\begin{equation}
    \P(E_\textrm{res} \, | \, \omega)\cdot \mathbf{1}_F(\omega)\geq g \cdot \mathbf{1}_F(\omega),
\end{equation}
whenever $\xi cR^{-f}n\geq R.$
\end{lem}
\begin{proof}
Suppose we are working on the event $F$ and assume without loss of generality that we are resampling a connection between two \textit{closed} clusters. Note first that the probability that $\widetilde{M}$ is precisely equal to the number of pivotals is bounded away from zero by a constant depending only on $M_0$, since the number of pivotals is upper bounded by $M_0$. In addition, we have a positive probability that we will pick the correct boxes to resample, again since the total number of boxes in the covering is smaller than $4c^{-2}R^4$. Moreover, point~\eqref{property:distance-points} of Lemma~\ref{lem:pivotals_are_far} tells us that by resampling the percolation configuration inside one of the boxes, we do not affect the well-separatedness of another box. Finally Lemma~\ref{le:linking_probability} tells us that for each individual box, the probability that we join the relevant clusters by resampling is bounded away from zero by a constant depending only on~$\xi$. 
\end{proof}

We are now ready to prove Proposition~\ref{lem:distancebetweencrossings}.

\begin{proof}[Proof of Proposition~\ref{lem:distancebetweencrossings}] In this proof, we will write $F$ instead of $F(A(x;n^{-1}R^{\kappa-f},n^{-1}R^{\kappa-f+1}))$ for brevity. Suppose that all of the constants are chosen so that $\P(F)\geq 1-p'$ as in Lemma~\ref{lem:bound_on_F}. All events in the definition of $F$ are measurable with respect to our original (non-resampled) percolation configuration $\omega$. Assume that we are working on $F$ and, as before, let $v_1,\ldots v_m$ be the cluster touch-points for $\widehat{F}(A(x;2n^{-1}R^{\kappa-f},4n^{-1}R^{\kappa-f}))$. Fix some $v_i$ and $v_j$ that are connected in $\operatorname{int}\phi$ and assume $i<j$. 

On the event that the resampling was successful, it must be that $\mathcal{E}(B^{+})$ holds in the resampled configuration $\omega_{\textrm{res}}$. Moreover, $v_i$ and $v_j$ are pivotal points for this event. In particular $X_{B^+}\geq D_{\omega_{\textrm{res}}}(v_i,v_j)$. Moreover, Lemma~\ref{lem:pivotals_are_far}\eqref{property:changing-state} tells us that the distance between $v_i$ and $v_j$ in the entire configuration $\omega_{\textrm{res}}$ is the same as the distance inside $\operatorname{int}\phi$ in the original configuration $\omega$. Combining these observations and using Lemma~\ref{lem:resampling_is_positive} we find that 
\begin{equation}
\begin{split}
\mathbf{P}(X_{B^+} \geq s \altq_n, \mathcal{E}(B^+))  &\geq \P\!\left(F, j\leq m, v_i \xleftrightarrow{\operatorname{int}\phi} v_j, D_{\omega}|_{\operatorname{int}\phi}(v_i,v_j)\geq s \altq_n, E_{\textrm{res}}\right) \\
&= \E\!\left[ \mathbf{1}(F, j\leq m, v_i \xleftrightarrow{\operatorname{int}\phi} v_j, D_{\omega}|_{\operatorname{int}\phi}(v_i,v_j)\geq s \altq_n) \cdot \E\!\left[ \mathbf{1}_{E_{\textrm{res}}}\,|\,\omega\right]\right] \\
&\geq \E\!\left[ \mathbf{1}(F, j\leq m, v_i \xleftrightarrow{\operatorname{int}\phi} v_j, D_{\omega}|_{\operatorname{int}\phi}(v_i,v_j)\geq s \altq_n) \cdot g\right] \\
&= g\cdot \P\!\left(F, j\leq m, v_i \xleftrightarrow{\operatorname{int}\phi} v_j, D_{\omega}|_{\operatorname{int}\phi}(v_i,v_j)\geq s \altq_n\right).
\end{split}
\end{equation}
Additionally, Lemma~\ref{lem:compXscales} gives us an upper bound on the left-hand-side of this equation. Namely, for every $\zeta>0$, there exists a constant $C=C(\alpha,\zeta)$ such that for all $R$ sufficiently large depending on $\zeta$, 
\begin{equation}
    \P(X_{B^{+}}\geq s\altq_n,\mathcal{E}(B^+))\leq \P(X_{R^{-f+1}}\geq s\altq_n) \leq Cs^{-\!\left(\frac{4}{3}+\zeta\right)\alpha}R^{-\alpha(f-1)}.
\end{equation}
In conclusion, 
\begin{equation}
\begin{split}
\P(G^{c} \, | \, F)&\leq \P(F)^{-1} \sum_{1\leq i<j\leq M_1} \P\!\left(F, j\leq m, v_i \xleftrightarrow{\operatorname{int}\phi} v_j, D_\omega |_{\operatorname{int}\phi}(v_i,v_j)\geq s\altq_n\right) \\
&\leq (1-p')^{-1} \cdot M_1^2 \cdot \frac{C}{g} s^{- \!\left( \frac{4}{3}+\zeta\right)\alpha}R^{-\alpha f}.
\end{split}
\end{equation} 
\end{proof}
 
\begin{rem} Once we establish the proof of Proposition~\ref{prop:subpolynomialconc}, we will be able to use Proposition~\ref{lem:distancebetweencrossings} with arbitrary exponent $\alpha>0$, not just $\alpha$ coming from the a priori estimate. The proof follows exactly the same -- we just replace Lemma~\ref{lem:compXscales} with Corollary~\ref{cor:subpolynomialconc}. \label{rem:distancebetweencrossings_general}
\end{rem}

\section{Superpolynomial concentration} \label{sec:superpolynomial-concentration}

In this section, we prove Proposition~\ref{prop:subpolynomialconc}. Throughout the section, we fix an $\alpha>\frac{3}{4}-\varepsilon$ with $\varepsilon>0$ small and assume that the bound~\eqref{eq:assumptionsubpoly} holds. Throughout the entire section we assume that $R$ is large enough (independent of $\alpha$) and $p$ is small enough depending on $R$ to ensure that all of the following proofs work. In particular, we want to choose $R$ sufficiently large so that Propositions~\ref{prop:apriori} and~\ref{lem:distancebetweencrossings} hold for some $\zeta>0$ to be chosen later. 

We consider percolation on $\Lambda_n$ restricted to some box $B^{(k)}\in\cB^R_{k,n}$. We first recall the definitions of $\tau_\ell$, $u,v,(x_i)_i$,  $X_{x_i,x_j}$ and $N$ from Section~\ref{subsubsec:subpolynomial-sketch}. We say a site $x$ is \emph{pivotal} if every open path from the left side of $B^{(k)}$ to the right side passes through $x$. We will also call them \emph{pinch points}. Assume the event $\mathcal{E}(B^{(k)})$ holds, so that there exist pivotal sites. Let $u,v$ be the left- and rightmost pivotal points and let $(u=x_0,x_1,\ldots,v)$ be the sequence of pivotal points from left to right. Write $X_{x_i,x_j}=D(x_i,x_j)$ to be the (geodesic or effective resistance) distance between $x_i$ and $x_j$. In particular, on the event $\mathcal{E}(B^{(k)})$, $X_{B^{(k)}}=X_{u,v}$. Define $\tau_0=0$ and for $\ell\geq 1$ define
\begin{equation} \label{eq:tau_definition2}
    \tau_{\ell}=\min\{i\geq\tau_{\ell-1}\colon X_{x_{\tau_{\ell-1},x_i}}> \altq_n\}.
\end{equation}
If the set above is empty, set $\tau_{\ell}=\infty$. Let
\begin{equation}
    N=\min\{\ell\geq 1\colon \tau_{\ell}=\infty\}.
\end{equation}
Then, because the $x_i$ are pivotal sites, we know that $D(u,v)=D(x_{\tau_0},x_{\tau_1})+D(x_{\tau_1},x_{\tau_2})+\cdots + D(x_{\tau_N-1},v)$ by the series law (property~\eqref{property:serial} of the metric). In particular,
\begin{equation}\label{eq:upperboundXdecomp2}
    X_{u,v}\leq N\altq_n+\sum_{\ell=1}^{N-1}X_{{\tau_{\ell}-1},\tau_{\ell}}.
\end{equation}

\subsection{Distance between pivotals}\label{sec:distance-between-pivotals}

In this subsection we bound the second term of~\eqref{eq:upperboundXdecomp2}. 

We bound the probability that the distance between the pinch points of a single bubble is large. {Define the \emph{bubble} in between two consecutive pivotals $x_i$ and $x_{i+1}$ to be the area enclosed by the top- and bottom-most (potentially non-simple) open paths between $x_i$ and $x_{i+1}$. Note that unless $x_i$ and $x_{i+1}$ are neighbors, there always exist two disjoint open paths between $x_i$ and $x_{i+1}$, so this is well-defined.} Let $\dprime(x_i,x_{i+1})$ be the radius $r$ of the smallest box $B(x_i,r)$ that contains the entire bubble between $x_i$ and $x_{i+1}$.

\begin{lem}\label{lem:crossingbubble}
There exists $\beta>0$ that is independent of $\alpha$ such that for all $\zeta>0$, there exists $R_0=R_0(\zeta)$ such that for all $R\geq R_0$, all $p$ sufficiently small depending on $R$, there exists a constant $C=C(R,\alpha, \zeta)$ such that 
\begin{equation}
    \P\!\left(\exists x_i\colon X_{x_i,x_{i+1}}\geq \altq_n,\,R^{-(d+1)}\leq \dprime (x_i,x_{i+1})\leq R^{-d}\right)\leq C(R,\alpha,\zeta)\cdot R^{(\frac{3}{4}+\zeta)(d-k)-(\alpha+\beta)d}
\end{equation}for all $n{,k}$ and all levels $d\geq k$ with $R^{-d}n\geq 1$.
\end{lem}
Recall the sketch below Proposition~\ref{prop:no_individual_bubble_long} in Section~\ref{subsubsec:subpolynomial-sketch}. For the proof of this lemma, we first remark that for a box of size $R^{-d}$ to contain a bubble, it must have four alternating arms. Using a union bound and the four-arm probability, we will later obtain the exponent $R^{(\frac{3}{4}+\zeta)(d-k)}$. We then proceed with an annuli covering argument, which is divided into several steps. Assume that the bubble between $x_i$ and $x_{i+1}$ is inside such a box. Our goal is to find a string of good annuli as in Proposition~\ref{lem:distancebetweencrossings} leading from $x_i$ to $x_{i+1}$, each of which can be crossed with short (geodesic or resistance) distance. To bound the probability that we fail, we distinguish between two cases. In case the two boundary segments from $x_i$ to $x_{i+1}$ are sufficiently separated, the distance has to be large ``along both boundary segments''. In case the two segments come close in between, we do not have two separated open paths from $x_i$ to $x_{i+1}$, but the pinch point forces a $5$-arm event which comes with an additional cost. In Step 1, we define two good events, $E_1$ and $E_3$ which will deal with these two cases, respectively. The event $E_1$ says that there do not exist two separate annuli inside the box that have crossings that are connected but have large (geodesic or resistance) distance from another. The event $E_3$ says that there does not exist an annulus that has connected crossings with large (geodesic or resistance) distance and also has 5 alternating arms to the boundary of the box. We show that on $E_1\cap E_3$, we have $D_n(x_i,x_{i_1})\leq \altq_n$. The event $E_1$ roughly guarantees that there can be at most one bad annulus. 

In Step 2a, we explain that if one of the bubble interfaces does not enter the bad annulus, we can cover it with good annuli and  the total distance between $x_i$ and $x_{i+1}$ is short. This is the step in which the stringing argument is first introduced. The intuition is as follows. Let $A_1$ and $A_2$ be two consecutive annuli in the covering of $\gamma$ for which $F\cap G$ holds (recall Proposition~\ref{lem:distancebetweencrossings}) and let $\phi(A_j),v_i(A_j)$ be defined as in Lemma~\ref{lem:decompositioncrossings}. Then there are points $x_+(A_1),x_-(A_1)\in V(A_1):=\{v_i(A_1)\}_i$ and $x_+(A_2),x_-(A_2)\in V(A_2)$ such that $x_{\pm}(A_j)=\gamma(t_{\pm}(A_j))$, where $t_-(A_1)\leq t_-(A_2)\leq t_+(A_1)\leq t_+(A_2)$. Now if $d_n$ were the geodesic metric, then $D_n(x_-(A_j),x_+(A_j))$ would be given by the length of some length-minimizing path between $x_-(A_j)$ and $x_+(A_j)$. If these paths for $A_1$ and $A_2$ intersected, then we could bound $D_n(x_-(A_1),x_+(A_2))$ by the triangle inequality. We continue stringing together these paths inductively and use the triangle inequality to obtain a bound on $D_n(x_i,x_{i+1})$. Of course, these paths are not guaranteed to intersect in general, and we want to give a unified proof for both the resistance and geodesic metrics. 
However, the simplified stringing argument above still gives a good idea of how the generalized stringing procedure works. Instead of using the existence of length-minimizing paths that intersect for consecutive annuli, a key ingredient will be the generalized parallel law that holds for both the geodesic and resistance metrics (recall property~\eqref{property:parallel} on page \pageref{property:serial}). Furthermore, instead of bounding $D_n(x_-(A_1),x_+(A_2))$, it will be necessary to bound $D_n(x_-(A_1),x_-(A_3))$. See Figure~\ref{fig:linking_boxes} for an illustration of the argument.

In Step 2b, we consider the case that both interfaces enter the bad annulus. In that case they have to enter a box with $5$ arms, so on $E_3$ that annulus can in fact not be the bad one. In Steps 3a and 3b, we bound the probabilities of $E_1^c$ and $E_3^c$. Using Proposition~\ref{lem:distancebetweencrossings}, we obtain that the probability of $E_1^c$ is a bit more than $R^{-2\alpha d}$, and the probability of $E_3^c$ is a bit less than $R^{-(\alpha+\beta) d}$ where $\beta>0$ comes from the $5$-arm event. In Step 4, we finally combine all bounds above together with the 4-arm bound from the beginning to obtain the desired result.

\begin{proof}[Proof of Lemma~\ref{lem:crossingbubble}]
Assume that there exist $x_i,x_{i+1}\in B^{(k)}$ such that $x_i$ and $x_{i+1}$ are consecutive pivotals and $R^{-(d+1)}\leq \dprime (x_i,x_{i+1})\leq R^{-d}$. We say that a site in $B^{(k)}\cap \Lambda_n$ is $R$-adic at level $\ell$ if it is of the form $\!\left(a\frac{R^{\kappa-\ell}}{n},b\frac{R^{\kappa-\ell}}{n}\right)$ for $a,b\in \{0,1,\ldots,R^{\ell}\}$. Note that there exists a site $x$ in $B^{(k)}\cap \Lambda_n$ at $R$-adic level $d$ 
such that $x_i\in B(x,R^{-(d-1)})$ and such that $B(x,2R^{-(d-1)})$ contains the entire bubble between $x_i$ and $x_{i+1}$. Then $B(x,2R^{-(d-1)})$ has 4 alternating arms, one to each side of $B^{(k)}$. Note that if $x$ is close to the boundary of $B^{(k)}$, it is possible that $B(x,2R^{-(d-1)})\nsubseteq B^{(k)}$, or even $B(x,2R^{-(d-1)})\nsubseteq [0,1]^2$. This does not matter to us -- if $B(x,2R^{-(d-1)})$ intersects a side of $B^{(k)}$ then the event that it has an arm up to this side is automatically satisfied (see the proof of Lemma~\ref{lem:4arms} -- we do not ask that $B(x,a)\subseteq B^{(k)}$).    
Furthermore, {we will bound the probability of the event that there exists a bubble inside $B(x,2R^{-(d-1)})$ of size $R^{-(d+1)}\leq \dprime (x_i,x_{i+1})\leq R^{-d}$ such that the (geodesic or effective resistance) distance between the two pivotals is at least $\altq_n$ by the probability of two events $E_1(x)\cup E_3(x)$ which depend only on the configuration inside $B(x,2R^{-(d-1)})\cap \Lambda_n$.}

\medskip 
\textbf{Step 1:} \textit{Definition of good events.}

Fix an $R$-adic point $x$ at level $d$. Recall the events $F,G$ from Lemma~\ref{lem:bound_on_F} and Proposition~\ref{lem:distancebetweencrossings}. We fix $\zeta=\frac{2}{3}$ 
and drop both $\zeta$ and $\alpha$ from the notation for the sake of clarity. We will also drop from the notation all of the parameters of the event $F$ except for $M_1$, as they will not be relevant for the following calculation. Fix some $\eta>0$ to be determined later.
\begin{itemize}
    \item Let $E_1(z,s,M_1)$ be the event that there exists a level $(1+\eta)d+1\leq f\leq (1+2\eta)d-2$ such that $F(M_1,f)\cap G(s,f)$ holds inside the annulus $A(z;n^{-1}R^{\kappa-f},n^{-1}R^{\kappa-f+1}) \cap \Lambda_n$. As above, note that it need not be that $B(z,n^{-1}R^{\kappa-f+1}) \subseteq [0,1]^2$. In this case the events $F$ and $G$ are not well defined. However, note that if we additionally define independent i.i.d.\ percolation on the sites of $\!\left(\bT_n \cap B(n^{-1}R^{\kappa-f+1})\right)\setminus [0,1]^2$, the distribution of distances within the relevant bubble between $x_i$ and $x_{i+1}$ remains the same. This is because these sites are fully separated from the new ones that we are adding by the bubble boundary. Because we will use events $F$ and $G$ only to control the distances between points inside the bubble, we will continue talking about $F,G$ of $A(z;n^{-1}R^{\kappa-f},n^{-1}R^{\kappa-f+1})$ for all annuli, keeping in mind that we are adding additional sites that do not influence these distances when $A(z;n^{-1}R^{\kappa-f},n^{-1}R^{\kappa-f+1}) \nsubseteq B^{(k)}$. Define
    \begin{equation}
        \begin{split}
            E_1(x):={}&\bigcap_{\|z_1-z_2\|_1\geq 4R^{-(1+\eta)d+1}}\!\left(E_1\!\left(z_1,\frac{1}{2}(M_1+1)^{-2}R^{-4\eta d},M_1\right)\right.\\
            &\left.\qquad\qquad\qquad\qquad\qquad\cup \, E_1\!\left(z_2,\frac{1}{2}(M_1+1)^{-2}R^{-4\eta d},M_1\right) \right) ,
        \end{split}
    \end{equation}
    where we take the union over all pairs of $R$-adic sites $z_1,z_2$ in $B(x,2R^{-(d-1)})$ at level $\lfloor(1+2\eta)d\rfloor-2$. 
    \item Fix some $\frac{1}{16}>\zeta'>0$ that will serve as the error in the five arm exponent and $1< C< 9/8-2\zeta'$. Let $E_2(z)$ be the event that there are 3 open and 2 closed arms (in order: open, closed, open, closed, open) from $B\!\left(z,2R^{-(1+\frac{\eta}{C})d}\right)$ to $B\!\left(z,\frac{1}{2}R^{-(d+1)}\right)$ in $B^{(k)}$ (again, we don't ask that $B\!\left(z,\frac{1}{2}R^{-(d+1)}\right)\subseteq [0,1]^2$) 
    and define
    \begin{equation}
        \begin{split}
            E_3(x)&:=\bigcap_{z}\!\left(E_2(z)^{c}\cup E'_1\!\left(z,\frac{1}{4}(M_1+1)^{-2},M_1\right)\right),
        \end{split}
    \end{equation}
    taking the union over all $R$-adic sites $z$ in $B(x,2R^{-(d-1)})$ at level $\lceil(1+\eta)d\rceil$. 
    Here $E_1'$ is defined as $E_1$ but for levels $\!\left(1+\frac{\eta}{C}\right)d+1\leq f\leq\!\left(1+\eta\right)d-2$ instead. Note that $C$ is introduced to make $E_2$ and $E_1'$ independent, as they depend on disjoint sets of vertices.\footnote{In fact, this proof directly works only when $\eta d\geq 4$ and $(1+\eta)d-2-\!\left(1+\frac{\eta}{C}\right)d-1 \geq 1$, i.e.\ when $d\geq \frac{4C}{\eta(C-1)}$ -- these relations guarantee that there is at least one level at which we can test for the existence of a good annulus in events $E_1$ and $E_1'$. We will see that we can choose $C,\eta$ deterministically, and thus our proof directly applies to $d$ larger than some deterministic constant, say $\hat{d}$. If $d\leq \hat{d}$, we can nevertheless subdivide box $B(x,2R^{-(d-1)})$ into $R$-adic boxes of level $\tilde{d}$ and carry out the proof as before. We now only get the additional constant factor of $R^{2\hat{d}}$ for the number of boxes. }
\end{itemize}

On the event $E_1(x) \cap E_3(x)$, there either does not exist a bubble inside $B(x,2R^{-(d-1)})$ of size at least $R^{-(d+1)}$ or there does exist a bubble and the distance between the left and right endpoints of the bubble is at most $\altq_n$ for the following reason. We say a box $B(z,n^{-1}R^{\kappa-f})$ is good if $F(z,M_1,f)\cap G\!\left(z,\frac{1}{2}(M_1+1)^{-2}R^{-4\eta d},f\right)$ holds. The annulus $A(z,2R^{\kappa-f},4R^{\kappa-f})$ is the corresponding good annulus. On $E_1(x)$, we can almost cover $B(x,2R^{-(d-1)})$ by good boxes centered on points $z$ at level $\lfloor(1+2\eta)d\rfloor-2$, possibly leaving a box of radius at most $4R^{-(1+\eta)d+1}$ uncovered, which we call the bad box. 

\medskip
\textbf{Step 2a:} \textit{Stringing of good boxes.}

Assume that one of the two bubble interfaces does not enter the bad box. We will show that on this event we can `string together' a short path from $x_i$ to $x_{i+1}$. Let us first introduce some notation. If $B$ is a good box, let $A$ be the good annulus around that box. Recall Lemma~\ref{lem:decompositioncrossings} -- let $V(A)=\{v_1(A),\ldots, v_{m(A)}(A)\}$ with $m(A)\leq M_1$ be the open sites that can be joined by closed paths to form a loop $\phi(A)$ in the annulus, whose existence is guaranteed by the lemma. As before, assume that these sites correspond to the outermost such loop. 
Let also $\intt \phi(A)$ denote the set of sites enclosed by and including this loop. Finally, for two sites $x,y$, $x\xleftrightarrow{\intt \phi(A)} y$ means that $x,y$ are open, belong to $\intt \phi(A)$ and there exists an open path in $\intt \phi(A)$ connecting them.

 Consider all open paths $\tilde{\gamma}$ that cross $B^{(k)}$ from left to right such that the section of $\tilde{\gamma}$ from $x_i$ to $x_{i+1}$ can be covered by interiors of good annuli corresponding to good boxes. The bubble interface that does not enter the bad box is precisely one such path. Let $\tilde{t}<\tilde{T}$ be times such that $x_i=\tilde{\gamma}(\tilde{t})$ and $x_{i+1}=\tilde{\gamma}(\tilde{T})$. For every such $\tilde{\gamma}$, consider the smallest sequence $\tilde{B}_1,\ldots \tilde{B}_{\tilde{h}}$ of (not necessarily distinct) good boxes such that their associated good annuli $\tilde{A}_1,\ldots, \tilde{A}_{\tilde{h}}$ satisfy the following:
 \begin{enumerate}
 \item[$\bullet$] For every $1\leq j\leq \tilde{h}$ there exist times $t_-(\tilde{A}_j)<t_+(\tilde{A}_j)$ such that
 \begin{equation}
     \begin{split}
         &x_-(\tilde{A}_j):=\tilde{\gamma}(t_-(\tilde{A}_j))\in V(\tilde{A}_j),\\
         &x_+(\tilde{A}_j):=\tilde{\gamma}(t_+(\tilde{A}_j))\in V(\tilde{A_j}),\\
         &\tilde{\gamma}|_{[t_-(\tilde{A}_j),t_+(\tilde{A}_j)]}\subseteq \operatorname{int} \phi(\tilde{A}_j);
     \end{split}
 \end{equation}
 \item[$\bullet$] $t_-(\tilde{A}_{j+1})\leq t_+(\tilde{A}_{j}) < t_+(\tilde{A}_{j+1})$ for all $1\leq j\leq \tilde{h}-1$ and $t_-(\tilde{A_1})\leq \tilde{t} < \tilde{T} \leq t_+(\tilde{A}_{\tilde{h}})$
 \end{enumerate}
 To see that at least one such sequence exists, one can inductively construct it starting from any annulus $\tilde{A}_1$ corresponding to a box containing $x_i$ and setting $t_\pm(\tilde{A}_1)$ to be the maximal/minimal time before/after $\tilde{t}$ such that $\tilde{\gamma}(t_\pm(\tilde{A}_1))\in V(A_1)$ and $\tilde{\gamma}|_{t_-,t_+}\subseteq \operatorname{int}\phi(A_1)$. Then find $\tilde{A}_2$ in the same way with $\tilde{\gamma}(t_+(\tilde{A}_1))$ in place of $x_i$. Because $t_+(\tilde{A}_j)$ is strictly increasing we will eventually cover $x_{i+1}=\tilde{\gamma}(\tilde{T})$.

 The reason we allow multiple occurrences of the same box $B$ is that the path $\tilde{\gamma}$ might wind around and come back into $B$ after some time. In the optimal scenario, i.e.\ for the path that minimizes the quantity $\tilde{h}$, if it does wind back into the same box, this time around it will visit sites that are not connected inside of $\operatorname{int}\phi(A)$ to the sites it has already visited. In particular, two consecutive boxes cannot be the same. 

Let $\gamma$ be such a path with the associated number of boxes $h$ the least possible. Let $A_1,\ldots A_h$ be the associated minimal string of good annuli. Note that we cannot have $x_-(A_j) \xleftrightarrow{\intt \phi(A_{j+1})}x_+(A_j)$ for any $1\leq j<h$. Suppose that were the case. Then we would be able to swap the segment of $\gamma|_{[t_-(A_j),t_+(A_j)]}$ with the open path from $x_-(A_j)$ to $x_+(A_j)$ that is inside of $\intt  \phi(A_{j+1})$. This new path is covered by the sequence $A_1,\ldots, A_{j-1},A_{j+1},\ldots, A_h$, contradicting the minimality of $\gamma$. Similarly, one sees that it cannot be that $x_-(A_j) \xleftrightarrow{\intt \phi(A_{j-1})}x_+(A_j)$ for any $j=2,\ldots,h$ either. Moreover, it must be that $t_-(A_j)\leq t_-(A_{j+1})\leq t_+(A_j)$, because otherwise we could discard $A_j$ from our list, contradicting the minimality of $h$. This means that $x_-(A_{j+1})\in \intt  \phi(A_j)$. We are now ready to string together a short path.

\begin{enumerate}
\item[$\bullet$] If $x_{i+1}$ belongs to the connected component of $\intt  \phi(A_1)$ containing $x_i$, 
then $D(x_i,x_{i+1})\leq D(x_-(A_1),x_+(A_1))\leq \frac{1}{2}(M_1+1)^{-2}R^{-4\eta d} \altq_n$ by the series law (property~\eqref{property:serial} at the beginning of Section~\ref{sec:sketch}), since $A_1$ is a good annulus and $x_i$ and $x_{i+1}$ are touch-points of our two closed clusters, so every open path from $x_-(A_1)$ to $x_+(A_1))$ must pass through these points. 

\item[$\bullet$]Assume now that $x_{i+1}$ does not belong to the connected component $\intt \phi(A_1)$ containing $x_i$. 
Consider the set 
\begin{equation}
    S(A_2)=\!\left\{v\in V(A_2): \, v\xleftrightarrow{\operatorname{int}\phi(A_2)} x_+(A_1)\right\}.
\end{equation}
It is non-empty because $x_-(A_2)\in S(A_2)$. Also, every open path from $x_-(A_1)$ to $x_{+}(A_1)$ must pass through at least one site in $S(A_2)$. Indeed, if not, we would have $x_-(A_1) \xleftrightarrow{\operatorname{int}\phi(A_2)}x_+(A_1)$, since $x_+(A_1)\in \operatorname{int}\phi(A_2)$ by construction. This contradicts the statement above. Let $v\in S(A_2)$ be such that $D(v,x_-(A_1))=\min_{u\in S(A_2)} D(u,x_-(A_1))$. By property~\eqref{property:parallel} of our metric $D$ (see the beginning of Section~\ref{sec:sketch}) we have that 
$$D(v,x_-(A_1))\leq M_1 \cdot D(x_- (A_1),x_+(A_1)).$$ 
\begin{itemize}
\item[$\star$]If $x_{i+1}$ belongs to the connected component of $\operatorname{int}\phi(A_2)$ containing $x_-(A_2)$,

because $x_i$ and $x_{i+1}$ are pinch points, we know that \begin{equation}\begin{split}D(x_i,x_{i+1})&\leq D(x_-(A_1),x_+(A_2))\leq D(x_-(A_1),v)+D(v,x_+(A_2))\\&\leq (M_1+1)\cdot \frac{1}{2}(M_1+1)^{-2}R^{-4\eta d} \altq_n. \nonumber \end{split}\end{equation} The last bound follows from the fact that $A_1$ and $A_2$ are good annuli and that $v,x_+(A_2)\in V(A_2)$ are connected in $\operatorname{int}\phi(A_2)$.
 
\item[$\star$]If $x_{i+1}$ is not in the connected component of $ \operatorname{int}\phi(A_2)$ containing $x_-(A_2)$, 
consider the set $R(A_2)=\!\left\{v\in V(A_2): \, v\xleftrightarrow{\operatorname{int}\phi(A_2)} x_-(A_3)\right\}$. 
It is again non-empty since $x_-(A_3)\in \operatorname{int}\phi(A_2)$ and so $x_+(A_2)\in R(A_2)$. If we let $w\in R(A_2)$ be such that $D(w,x_+(A_3))=\min_{u\in R(A_2)} D(u,x_+(A_3))$, then exactly as above we have that $$D(w,x_+(A_3))\leq M_1 \cdot D(x_- (A_3),x_+(A_3)).$$ Thus, since $A_2$ is good, $v,w\in V(A_2)$ and $v \xleftrightarrow{\operatorname{int}\phi(A_2)} x_+(A_1) \xleftrightarrow{\operatorname{int}\phi(A_2)} x_-(A_3) \xleftrightarrow{\operatorname{int}\phi(A_2)} w$ we have that \begin{equation} \begin{split}
D(x_-(A_1),x_-(A_3))&\leq D(x_-(A_1),v)+D(v,w)+D(w,x_+(A_3))+D(x_+(A_3),x_-(A_3))  \\&\leq (M_1+1+M_1+1) \cdot \frac{1}{2}(M_1+1)^{-2}R^{-4\eta d} \altq_n. \nonumber\end{split}\end{equation} 
See Figure~\ref{fig:linking_boxes} for an illustration of this argument.
\end{itemize}
\item[$\bullet$] Proceed recursively, starting the same procedure from $x_-(A_3)$ instead of $x_-(A_1)$.
\end{enumerate}
Using the triangle inequality, this algorithm gives that  $$D(x_i,x_{i+1})\leq \begin{cases} \!\left(\frac{h-1}{2} \cdot 2(M_1+1) + 1\right)\cdot \frac{1}{2}(M_1+1)^{-2}R^{-4\eta d} \altq_n 
 & \textrm{if }h \textrm{ is odd} \\ \!\left(\frac{h-2}{2} \cdot 2(M_1+1) + (M_1+1)\right)\cdot \frac{1}{2}(M_1+1)^{-2}R^{-4\eta d} \altq_n & \textrm{if }h \textrm{ is even.}\end{cases} $$ Either way, we obtain that $D(x_i,x_{i+1})\leq h\cdot \frac{1}{2}(M_1+1)^{-1}R^{-4\eta d} \altq_n$. Finally, $\gamma$ enters at most $M_1R^{4\eta d}$ good boxes, since the number of $R$-adic points $z$ at level $\lfloor(1+2\eta)d\rfloor-2$ inside $B(x,2R^{-(d-1)})$ is at most $R^{4\eta d}$, each good annulus associated with a good box contains at most $M_1$ entry points and we have chosen $\gamma$ to be minimal. Substituting this into the above we find that $D(x_i,x_{i+1})\leq  \frac{1}{2}\altq_n$, as wanted.

\begin{figure}[ht]
    \centering
    \includegraphics[width=1\linewidth]{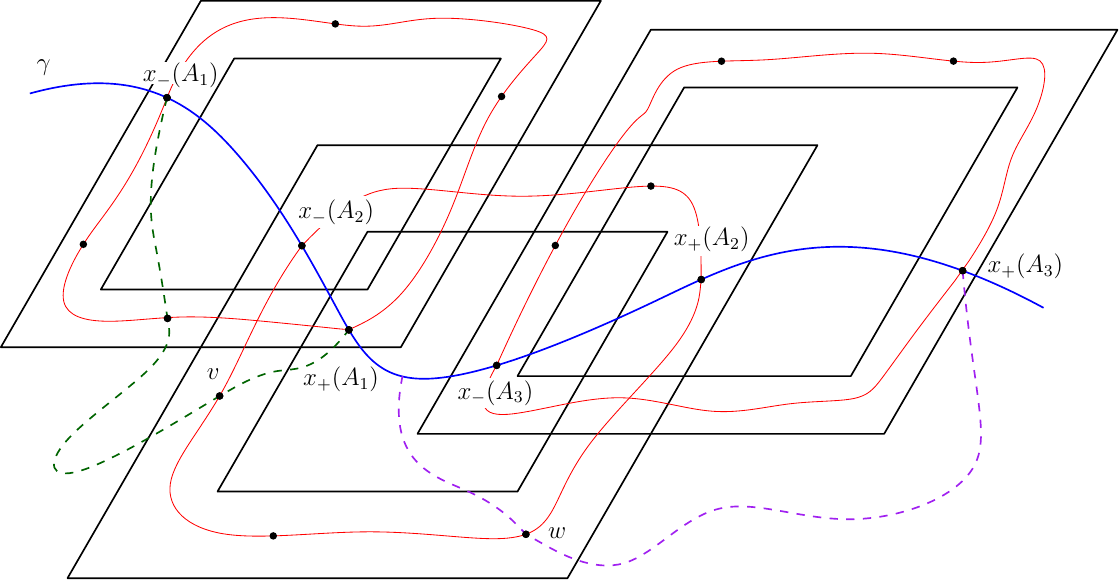}
    \caption{An illustration of the linking construction from Lemma~\ref{lem:crossingbubble}. The blue line represents the open path covered by good boxes. Black points in each of the annuli represent the sets $V(A_1), V(A_2)$ and $V(A_3)$ of our `defects' and the red loops connecting them represent the associated `almost closed loop' in each annulus.}
    \label{fig:linking_boxes}
\end{figure}

\medskip 
\textbf{Step 2b:} \textit{The bad box cannot have 5 arms.}

Now assume that both interfaces enter the bad box. Then inside this bad box, there exists an $R$-adic point $z$ of level $\lceil(1+\eta)d\rceil$ such that the bad box is contained in $B(z,R^{-(1+\eta)d+2})$. Moreover, $B(z,\frac{1}{2}R^{-(d+1)})$ cannot contain {the entire bubble between $x_i$ and $x_{i+1}$ by definition of $d$.} Hence, there are 5 arms of alternating colors in between these two boxes: two open ones coming from the inside boundaries of the bubble going towards the {part of the bubble interface that does not belong to the bigger box through or towards one of the endpoints $x_i$ or $x_{i+1}$, one open arm going through or towards the other endpoint}, and two closed arms corresponding to the outer boundaries of the bubble -- see Figure~\ref{fig:bubble_5_arms} below. In particular, the box $B(z,2R^{-(1+\frac{\eta}{C})d})$\footnote{Again, this box has $5$ arms only if its radius is larger than $R^{-(1+\eta)d+2}$. This holds if, for example $\eta d\geq \frac{3C}{C-1}$. As explained in the previous footnote, we will choose $C,\eta$ deterministically, so if $d$ is smaller than this absolute constant, we will simply further subdivide our box to this level and then carry out the proof as before with an additional constant.} has $5$ arms of alternating color to $B(z,\frac{1}{2}R^{-(d+1)})$ and so $E_2(z)$ holds. 

\begin{figure}[ht]
    \centering
    \includegraphics[scale=0.9]{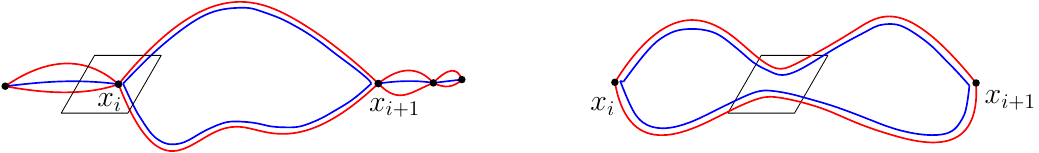}
    \caption{Two examples of both interfaces entering the bad box. In both cases, the box has 5 arms going far.}
    \label{fig:bubble_5_arms}
\end{figure}

On $E_3(x)$, there exists $\!\left(1+\frac{\eta}{C}\right)d+1\leq f\leq (1+\eta)d-2$ such that the bad box is contained in a box for which $F(M_1,f)\cap G\!\left(\frac{1}{4}(M_1+1)^{-2},f\right)$ holds, so we can alter our open path exactly as above. This time we use all the good boxes and the bad box to obtain a cover of our path $\gamma$. Moreover, in the algorithm, we will use the bad box at most $2M_1$ times. Indeed, the bad box can be entered at most $M_1$ times by the path $\gamma$. Every time the index of the box is odd in our covering, it is used in the algorithm twice and when its index is even it is used in the algorithm once. Moreover, when the box is used, it can contribute at most $(M_1+1)\cdot  \frac{(M_1+1)^{-2}}{4}\altq_n$ to the distance. Thus, 
$$D(x_i,x_{i+1})\leq \frac{1}{2}\altq_n+2M_1(M_1+1)\cdot \frac{(M_1+1)^{-2}}{4}\altq_n\leq \altq_n,$$ where the first term comes from the good boxes and the last term comes from the bad box as explained above. It follows that
\begin{equation}
    \begin{split}
        &\P\!\left(\exists x_i,x_{i+1}\colon X_{x_i,x_{i+1}}\geq \altq_n,\,R^{-(d+1)}\leq \dprime (x_{i},x_{i+1})\leq R^{-d}\right)\\
        \leq{}&\sum_{x \textrm{ of level }d} \P\!\left(\textrm{an arm from $B(x,2R^{-(d-1)})$ to each of the $4$ sides of $B^{(k)}$}\right)\\ &\cdot(\P(E_1(x)^{c})+\P(E_3(x)^{c})). \label{eq:bound_on_distance_x_i_x_i+1}
    \end{split}
\end{equation}

\medskip 
\textbf{Step 3a:} \textit{Bounding $\p(E_1^c)$.}

We first bound $\P(E_1(x)^{c})$. Let $z_1\in B(x,R^{-d})$. We will make use of Lemma~\ref{lem:bound_on_F} and Proposition~\ref{lem:distancebetweencrossings}.
Choose $0<p' \leq \frac{1}{2}R^{-\alpha/\eta}$ and suppose the parameters of $F$ are chosen so that $\P(F(M_1,f)) \geq 1-p'$ for all levels $f$ such that $\xi\frac{c_1c_2}{8} \cdot R^{\kappa-f-1}\geq R$.

Suppose first that $d$ is such that $\xi\frac{c_1c_2}{8} \cdot R^{\kappa-(1+2\eta)d-1}\leq R$. If we have $\dprime (x_i,x_{i+1})\leq R^{-d}$, then $X_{x_i,x_{i+1}}\leq 4R^{-2d}n^2$ deterministically. Indeed, if we know that sites $x_i$ and $x_{i+1}$ in a graph are connected via a path of length $l$, then both geodesic and resistance distance between $x_i$ and $x_{i+1}$ are smaller than $l$; since the graph we are considering is a box of radius $R^{-d}$, the total number of sites is $4R^{-2d}n^2$, yielding the result. In particular, this means that for $\eta<\frac{1}{8}$,
\begin{equation}
    X_{x_i,x_{i+1}}\leq 4R^{-2d}n^2\leq4\!\left(\!\left(\xi\frac{c_1c_2}{8}\right)^{-1}R^{2-\kappa}\right)^{\frac{2}{1+2\eta}}n^2\leq C(R,\alpha,\eta)n^{\frac{-2}{1+2\eta}}n^2\leq C(R,\alpha,\eta)n^{\frac{2}{5}},
\end{equation}
where the value of $C(R,\alpha,\eta)$ may change from line to line.

By Corollary~\ref{cor:lower_bound_on_m_n}, we know that there exists a constant $\tilde{C}$ depending only on $R$ such that $\altq_n\geq \tilde{C}n^{1/2}$. Thus, for all $n$ big enough (depending on $\eta,R,\alpha$) we have that $X_{x_i,x_{i+1}}< \altq_n$, and hence 
\begin{equation}
    \P\!\left(\exists x_i : X_{x_i,x_{i+1}}\geq \altq_n, R^{-(d+1)}\leq \dprime (x_i,x_{i+1})\leq R^{-d}\right)=0.
\end{equation}
For $n$ small (depending on $\eta,R,\alpha$), there are finitely many possibilities for $n,k,d$ and hence the statement of the lemma holds upon possibly making the constant $C(R,\alpha,\eta)$ bigger. Thus, the lemma holds in this case.

Assume now that $d$ satisfies $\xi\frac{c_1c_2}{8}\cdot R^{\kappa-(1+2\eta)d-1}\geq R$. Then we know that for each level $(1+\eta)d+1\leq f\leq(1+2\eta)d-2$, we have $\P(F(M_1,f))\geq 1-p'$. 
Since the event $F(z,M_1,f)$ depends only on the sites of $A(z;n^{-1}R^{\kappa-f},n^{-1}R^{\kappa-f+1})$ and these annuli are disjoint at different levels, we know that events $F(z,M_1,f)$ are independent for different $f$. Thus, the probability that there does not exist a level $(1+\eta)d+1\leq f\leq (1+2\eta)d-2$ such that $F(M_1,f)$ holds is bounded from above by $(p')^{\eta d-4}\leq (p')^{-4}\cdot R^{-\alpha d}$. Then by applying Proposition~\ref{lem:distancebetweencrossings} with $\zeta=\frac{2}{3}$, conditional on $F$ being true at level $f$, the probability of $G\!\left(\frac{1}{2}(M_1+1)^{-2}R^{-4\eta d},f\right)^{c}$ is bounded from above by
\begin{equation}
    C(R,\alpha,\eta) (M_1+1)^{4\alpha} R^{8\eta d\alpha}R^{-\alpha f}\leq C(R,\alpha,\eta)\cdot R^{-\alpha(1-7\eta)d},
\end{equation} where $C(R,\alpha,\eta)$ denotes a constant depending only on $R,\alpha,\eta$ that can increase from line to line.
We conclude that
\begin{equation}
    \begin{split}
        \P\!\left(E_1\!\left(z_1,\frac{1}{2}(M_1+1)^{-2}R^{-4\eta d},M_1\right)^c\right)\leq&\, (p')^{-4}\cdot R^{-\alpha d}+C(R,\alpha,\eta)\cdot R^{-\alpha(1-7\eta)d}\\
        \leq&\, C(R,\alpha,\eta) \cdot R^{-\alpha(1-7\eta)d}.
    \end{split}
\end{equation}
For $z_1,z_2$ such that $||z_1-z_2\|_1\geq 4R^{-(1+\eta)d+1}$, the events $E_1\!\left(z_1,\frac{1}{2}(M_1+1)^{-2}R^{-4\eta d},M_1\right)$ and $E_{1}\!\left(z_2,\frac{1}{2}(M_1+1)^{-2}R^{-4\eta d},M_1\right)$ are independent, so 
\begin{equation}
    \begin{split}
        \P(E_1(x)^{c})\leq C(R,\alpha,\eta) \cdot R^{8\eta d}R^{-2\alpha(1-7\eta)d}\leq C(R,\alpha,\eta) R^{-\frac{3}{2}\alpha d} \label{eq:bound_on_E_1}
    \end{split}
\end{equation}
for all $\eta$ small enough, say $\eta\leq 1/100$.

\medskip 
\textbf{Step 3b:} \textit{Bounding $\p(E_3^c)$.}

Next, we bound $\P(E_3(x)^{c})$. Note that the events $E_2(z)$ and $E_1'\!\left(z,\frac{(M_1+1)^{-2}}{4},M_1\right)$ are independent. By Lemma~\ref{lem:arm_exponents}, we know that for all $z$ and $1/16>\zeta'>0$, there exists $C_1=C_1(\zeta')$ such that for all $d, 1< C< 9/8-2\zeta',\eta$, 
\begin{equation}
    \P(E_2(z))\leq C_1\cdot  R^{-(2-\zeta')\frac{\eta}{C}d}.
\end{equation}
Furthermore, similarly to the above, if we choose $p'<R^{-\alpha\frac{1+\eta/C}{\eta(1-1/C)}}$, we have that
\begin{equation}
    \P\!\left(E'_1\!\left(z,\frac{(M_1+1)^{-2}}{4},M_1\right)^{c}\right)\leq C(R,\alpha,\eta)\cdot R^{-\alpha(1+\frac{\eta}{C})d},
\end{equation} for some constant $C(R,\alpha,\eta)$ depending only on these parameters. Summing over all $z$ at level $\lceil(1+\eta)d\rceil$, we get
\begin{equation}
    \P(E_3(x)^{c})\leq C(R,\alpha,\eta)\cdot R^{2\eta d}R^{-\alpha(1+\frac{\eta}{C})d}R^{-((2-\zeta')\frac{\eta}{C})d}\leq C(R,\alpha,\eta) R^{-\alpha(1+\frac{\eta}{2C})d}, \label{eq:bound_on_E_3}
\end{equation}
where the last inequality follows since we have chosen $\zeta'<1/16$ and $1<C<9/8-2\zeta'$ and because we are working in the range $\alpha \geq \frac{3}{4}-\varepsilon >1/2$.

\medskip 
\textbf{Step 4:} \textit{Final bound.}

Fixing $\eta$ to be any number in $(0,1/100)$ guarantees that all of the above steps are correct. Using equations~\eqref{eq:bound_on_E_1} and~\eqref{eq:bound_on_E_3} to get an upper bound on~\eqref{eq:bound_on_distance_x_i_x_i+1}, and noting that we are working in the range $\alpha>1/2$, we find that there exists $\beta>0$ ($\beta=\min \{\frac{1}{4},\frac{\eta}{4C}\}$ works) such that
\begin{equation}
    \begin{split}
        &\P\!\left(\exists x_i,x_{i+1}\colon X_{x_i,x_{i+1}}\geq \altq_n,\,R^{-(d+1)}\leq \dprime (x_{i},x_{i+1})\leq R^{-d}\right)\\
        &\leq\sum_{x \textrm{ at level }d}
        C(R,\alpha)R^{-(\alpha+\beta)}\cdot\P(\textrm{an arm from $B(x,2R^{-(d-1)})$ to each of the $4$ sides of $B^{(k)}$}).
    \end{split}
\end{equation}
By Lemma~\ref{lem:4arms}, for all $\zeta>0$,
\begin{equation}
    \begin{split}
        \sum_{x \textrm{ at level }d}
        \P(\textrm{an arm from $B(x,2R^{-(d-1)})$ to each of the $4$ sides of $B^{(k)}$})\leq C(R,\zeta)\cdot R^{(\frac{3}{4}+\zeta)(d-k)},
    \end{split}
\end{equation}
for some constant $C(R,\zeta)$. This concludes the proof.
\end{proof}

\begin{cor}
There exist $\beta>0$, $\varepsilon>0$ such that the following holds for all $R$ sufficiently large and $p$ sufficiently small depending on $R$. For all $\alpha>\frac{3}{4}-\varepsilon$ for which~\eqref{eq:assumptionsubpoly} holds, we also have that for all $k$ and $n>R^k$,

\begin{equation}
    \begin{split}
        \P\!\left(\exists i\colon X_{x_{i},x_{i+1}}\geq \altq_n\right)\leq C(R,\alpha)\cdot R^{-(\alpha+\beta)k},
    \end{split}
\end{equation} for some constant $C(R,\alpha)>0$. 
\label{cor:no_individual_bubble_long}
\end{cor}

\begin{proof}   
Let $\beta$ be as in the statement of Lemma~\ref{lem:crossingbubble}. Then for all $\zeta>0$, there exists $R_0=R_0(\zeta)$ such that for all $R\geq R_0$ there exists $C(R,\alpha,\zeta)$ for which
\begin{equation}
    \begin{split}
        \P(\exists i\colon X_{x_i,x_{i+1}}\geq \altq_n) &\leq\sum_{d=k}^{\infty}\P(\exists i\colon X_{x_i,x_{i+1}}\geq \altq_n,\,R^{-(d+1)}n\leq \dprime (x_i,x_{i+1})\leq R^{-d}n)\\
        &\leq C(R,\alpha,\zeta) \cdot \sum_{d=k}^\infty R^{(\frac{3}{4}+\zeta)(d-k)-(\alpha+\beta)d}.
\end{split}
\end{equation}
Choosing $\varepsilon=\beta/2$ and $\zeta=\beta/4$ makes the above series convergent, and the final expression can be bounded from above by $C(R,\alpha)\cdot R^{-(\alpha+\beta)d}$ for some constant $C(R,\alpha)$, as desired.
\end{proof}

\subsection{Bootstrapping}\label{subsec:bootstrapping}

\begin{figure}[ht]
    \centering
    \includegraphics[scale=0.85]{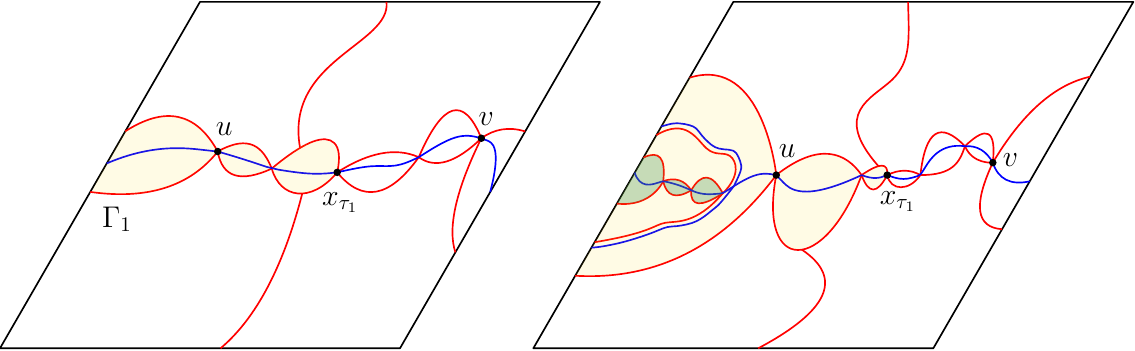}
    \caption{On the left-hand side, the area shaded in yellow represents the set $\Gamma_1$. The picture on the right illustrates why the set $\Gamma_1$ is not local. The area shaded in yellow once again represents the true set $\Gamma_1$. The area shaded in green consists of bubbles bounded in between the two closed clusters represented by the red curves that neighbor the green region. These closed clusters are fully surrounded by the open cluster of $\cE(B^{(k)})$. In this case, we do not want to consider these bubbles as they are not relevant for our $\cE(B^{(k)})$. }
    \label{fig:exploration_is_not_local}
\end{figure}

In the previous subsection we have found a bound on the second term in~\eqref{eq:upperboundXdecomp2} in terms of~$N$. Let us now turn our attention to bounding $N$. We will do so by showing that $N$ is `almost geometric'. Recall the definition~\eqref{eq:tau_definition2} of our stopping times $\tau_l$:
\begin{equation}
    \tau_0=0 \textrm{  and  }\tau_{\ell}=\min\{i\geq\tau_{\ell-1}\colon X_{x_{\tau_{\ell-1},x_i}}\geq \altq_n\} \textrm{ for }\ell\geq 1. \nonumber
\end{equation} If the set above is empty, we simply set $\tau_{\ell}=\infty$. Moreover, $N$ is given by $N=\min\{\ell : \tau_{\ell}=\infty\}$. On the event $B^{(k)}$, we will `explore the bubbles until time $\tau_\ell$' for a fixed $\ell$. Namely, let $\Gamma_\ell$ be the set of sites enclosed by and including the boundaries of clusters corresponding to our top and bottom closed crossings `until time $\tau_\ell$'. More precisely, this set consists of the first few bubbles along with the connected component containing the endpoint of the open left-to-right path on the left-hand-side of $B^{(k)}$. See Figure~\ref{fig:exploration_is_not_local} for an illustration of the set $\Gamma_1$. Because we will later condition on the shape of the set $\Gamma_\ell$, we will refer to it as the `explored domain'. This is a slight abuse of notation, since there is no exploration that would yield this set as a result. In fact, this set is not even local. {This is because for $\Gamma_1$ to be equal to some specific set $S$, we not only need to know what happens inside of $S$, but we also need to know that $S$ is connected to the top and bottom closed clusters as well as the open cluster coming from the event $\mathcal{E}(B^{(k)})$. Otherwise, it may happen that inside the set $S$ we can identify the closed/open interface on the boundary and bubble pinch-points such that the distance between the first pinch-point and the last one is larger than $\altq_n$, but that $S$ is in fact bounded in between two closed clusters that are fully surrounded by the crossing open cluster of $\cE(B^{(k)})$. In this case, the first relevant bubble pinch-point is located outside of $S$, and we must look outside the set $S$ to know `when to start looking at the distance across bubbles'. See Figure~\ref{fig:exploration_is_not_local} for an illustration.}

The general strategy of the proof will be to show that with very high probability our `unexplored domain'\footnote{Again, this terminology is a slight abuse of notation.} $C=B^{(k)}\setminus \Gamma_\ell$ is sufficiently nice, and on the event that it is nice, we will be able to string together a short path from the last pinch point of our `explored' bubbles to the right-hand side of $B^{(k)}$ with high probability.

In the rest of the proof we will fix $\ell=1$, $C=B^{(k)}\setminus \Gamma_1$ -- considering only the first time $\tau_1$ is sufficient for our purposes. However, all the results can easily be extended to hold for any time $\tau_\ell$ as above. Let $x=x_{\tau_1}$ be the last pinch point of the bubbles in this exploration. Let $\eta_t$ and $\eta_b$ be the outer boundaries of the top and bottom closed paths given by the exploration of the open cluster. And finally, let $\partial:=\partial \!\left(B^{(k)}\setminus C\right)$.

At the time $\tau_1$ of the above form, we consider the \emph{$R$-adic Whitney decomposition} of $C$. For $k+1\leq i\leq\kappa$, let $\cB_{i,n}^{R,k}$\phantomsection{}\label{def:cBRk_jn} be the set of boxes $B \in \cB_{i,n}^R$ that are contained in $B^{(k)}$. In particular, each box $B\in \cB_{i,n}^{R,k}$ contains $R^{\kappa-i}\times R^{\kappa-i}$ vertices (equivalently, its side-length is $n^{-1}R^{\kappa-i}$), and we have that $|\cB_{i,n}^{R,k}|=R^{2(i-k)}$.

Now consider the collection of boxes in $\cB^{R,k}_{i,n}$ that are fully contained in $C$ and have Euclidean distance between $\frac{1}{2}n^{-1}R^{\kappa-i-2}$ and $\frac{1}{2}n^{-1}R^{\kappa-i+2}$ from $\partial$, discarding any boxes that are fully inside another. We call this collection the $R$-adic Whitney decomposition of $C$. Note that for every site $x\in C$ there exists such a box containing it -- indeed, if $k\leq i\leq \kappa$ is such that $n^{-1}R^{\kappa-i}\leq\dist(x,\partial)<n^{-1}R^{\kappa-i+1}$, then one of the $R$-adic boxes of levels $i$ or $i+1$ containing $x$ must be a Whitney box. So the union of all Whitney boxes equals the set $C$.

We will first show that there exists a global event $H$ that happens with very high probability, and on the event $H$ it will be possible to find an open path from $x$ to the last bubble pinch-point in the box $B^{(k)}$ that can be covered by a `small' number of $R$-adic Whitney boxes of each level $i\geq k+1$. Following this, we will show that with high probability none of these boxes are individually bad. Finally, we will string together short crossings across each one of these boxes to obtain an upper bound on the remaining distance $X_{x,v}$.

\subsubsection{Controlling the explored domain}\label{subsubsec:controlling-explored-domain}
Let us first define $H$. Given a box $B \in \cB_{i,n}^{R,k}$, we say that it has $6$ semiarms of length $l\geq 0$ if there exists $0\leq \tilde{l}\leq l$ such that there are $6$ arms crossing the annulus $R^{\tilde{l}}\cdot B\setminus R^{3}\cdot B$ and $6$ arms crossing the annulus $R^l\cdot B \setminus R^{\tilde{l}+3}\cdot B$. We similarly define the concept of semiarms in the half-plane. For a box $B\in \cB_{i,n}^{R,k}$ such that $R^3B\,\cap B^{(k)}\neq \emptyset$, let $z'$ be the projection of the center of $B$ to the nearest side of $B^{(k)}$. We say that $B$ has $3$ semiarms of length $l\geq 0$ near the boundary if there exists $0\leq \tilde{l}\leq l$ such that there are $3$ arms in $B^{(k)} \cap A\!\left(z';n^{-1}R^{\kappa-i+3},n^{-1}R^{\kappa-i+\tilde{l}}\right)$ and $3$ arms in $B^{(k)} \cap A\!\left(z';n^{-1}R^{\kappa-i+\tilde{l}},n^{-1}R^{\kappa-i+l}\right)$. Note that these arms are always contained in a half-plane. 

\begin{figure}[ht]
    \centering
    \includegraphics[scale=0.7]{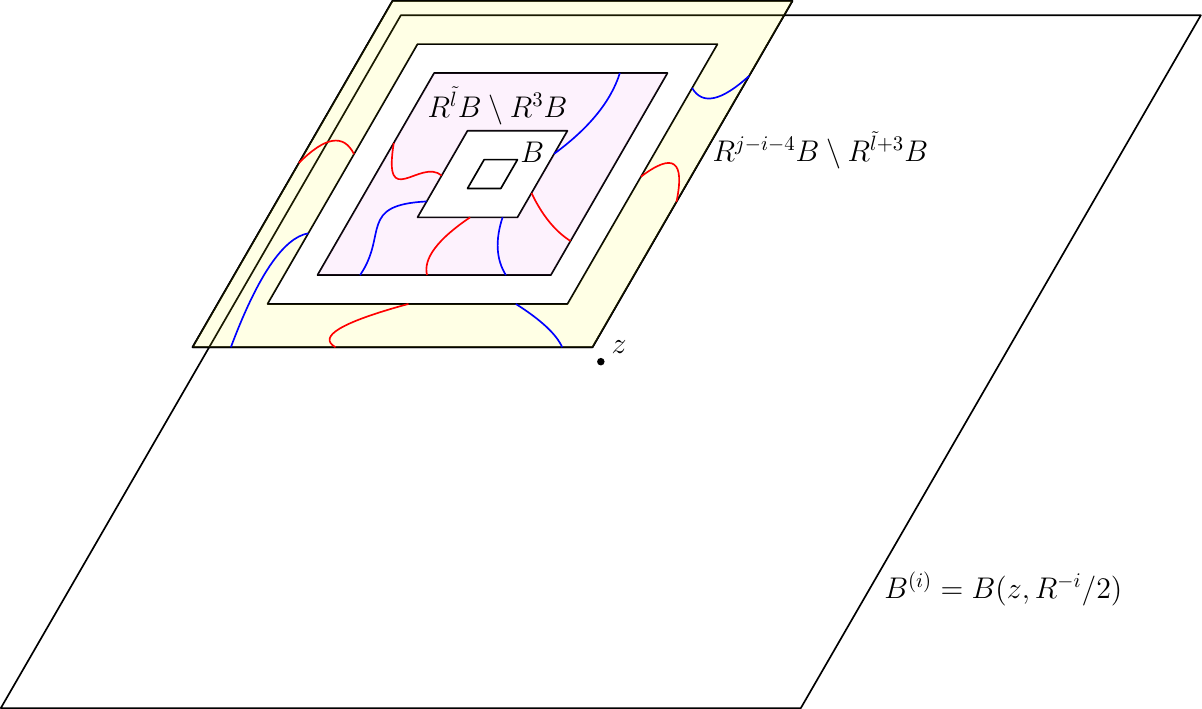} 
    \caption{On the picture above $z$ is a center of an $R$-adic box of level $i+3$ in $B^{(k)}$. One can also see a box $B$ of level $j\geq i+9$ inside $B^{(i)}=B(z,R^{-i}/2)$. The pink and yellow annuli around the box $B$ represent the two annuli in which we have $6$-arm events. Then $B$ has $6$ semiarms of length $j-i-4$. 
    }
    \label{fig:definition_of_H_i}
\end{figure}

For every $k\leq j\leq \kappa$, we define the event $H_j$ to be the following:

  \begin{enumerate}
  \item[$\bullet$] Consider some $k-2\leq i\leq j-10$ and all the points $z$ in $B^{(k)}$ that are centers of boxes in $\cB^{R,k}_{i+3,n}$. Let $B^{(i)}=B(z,\frac{1}{2}R^{-i})$ for one such $z$. Note that it need not be the case that $B^{(i)}\subseteq B^{(k)}$. For any such $i$ and $B^{(i)}$, let $N_{j,B^{(i)}}$ be the number of boxes $B\in\cB_{j,n}^{R,k}$ that are contained in $B^{(i)}$ such that $B$ has $6$ semiarms of alternating color of length $j-i-4$. 
  We require that $N_{j,B^{(i)}}\leq R^{j\nu}$ for all $k-2\leq i\leq j-10$ and all $B^{(i)}$ as above, where $\nu>0$ is a small absolute constant that we will determine later. See Figure~\ref{fig:definition_of_H_i} for an illustration. 
  \item[$\bullet$] Consider some $j\geq k+9$. Let $N_j^h$ be the number of boxes $B\in \cB_{j,n}^{R,k}$ such that $R^3B\cap \partial B^{(k)}\neq \emptyset$ that have $3$ semiarms of length $j-k-1$ near the boundary. We require that $N_j^h\leq R^{j\nu}$, where $\nu>0$ is as above.
  \end{enumerate}

\noindent Finally, define $H=\bigcap_{j= k}^{\kappa} H_j$.

\begin{lem} 
The following holds for all $R$ large enough. For any $\nu$ as in the definition of events $H_j$, there exists $C=C(\alpha,R,\nu)$ such that for all $k,n$ and $k\leq j\leq \kappa$, we have that $$\P(H_j^c)\leq C\cdot R^{-4\alpha j}.$$ In particular, we have $\mathbf{P}\!\left(H^c\right)\leq C'\cdot  R^{-4\alpha k}$ for a constant $C'=C'(\alpha,R,\nu)$.\label{lem:bound_on_global_event}
\end{lem}

The proof of Lemma~\ref{lem:bound_on_global_event} is given in Appendix~\ref{app:six-semiarms}.

\subsubsection{Path covering}
Our next goal is to show that we can find an open path from $x$ to $v$ (the last pinch point of bubbles inside of the box $B^{(k)}$ created by the clusters adjacent to the top and bottom sides of $B^{(k)}$) that can be covered by a collection of Whitney boxes whose number can be bounded by a constant plus additional Whitney boxes that either have $6$ semiarms of appropriate length or they have three semiarms near the boundary. Consider a bubble in the unexplored domain $C$, and let the pivotals corresponding to the bubble be $a$ and $b$. Let $\gamma_t$ and $\gamma_b$ be the closed paths that are boundaries of this bubble adjacent to the top and bottom sides of $B^{(k)}$ respectively.

 Consider the set $\mathcal{S}\subseteq C$ of sites $s$ in $B^{(k)}$ such that $\!\left|\textrm{dist}(s,\gamma_t) -\textrm{dist}(s,\gamma_b)\right|\leq \frac{4}{n}$, i.e.\ the sites such that their distances from the top and bottom sides of our bubble differ by at most 4 times the distance between neighboring sites in the triangular lattice $\bT_n$. We now explain that there exists a (not necessarily open) path $\gamma \subseteq \mathcal{S}$ connecting $a$ and $b$ inside the bubble. First, we know that $a,b\in \mathcal{S}$, since $a,b$ neighbor both top and bottom boundaries. If there does not exist a path from $a$ to $b$ inside of the bubble that consists only of sites in $\mathcal{S}$, then there exists a blocking path $a_1,\ldots a_j$ inside the bubble, starting at $\gamma_t$ and ending at $\gamma_b$ such that $a_{\tilde{j}}\notin \mathcal{S}$ for all $1\leq \tilde{j}\leq j$. But this is impossible, since $\textrm{dist}(a_1,\gamma_t) -\textrm{dist}(a_1,\gamma_b)<0$, $\textrm{dist}(a_j,\gamma_t) -\textrm{dist}(a_j,\gamma_b)>0$ and for two consecutive sites $a_{\tilde{j}}$ and $a_{\tilde{j}+1}$ we have that this difference of distances changes by at most
\begin{equation}
\!\left| \!\left(\textrm{dist}(a_{\tilde{j}},\gamma_t) -\textrm{dist}(a_{\tilde{j}},\gamma_b)\right)-\!\left(\textrm{dist}(a_{\tilde{j}+1},\gamma_t) -\textrm{dist}(a_{\tilde{j}+1},\gamma_b)\right)\right| \leq 2 \textrm{dist}(a_{\tilde{j}},a_{\tilde{j}+1})=\frac{2}{n}\nonumber.
\end{equation} The last inequality follows by applying the triangle inequality twice. Thus, at least one of $a_{\tilde{j}}$'s must be in $\mathcal{S}$, contradicting our assumption. So there exists a path connecting $a$ and $b$ inside the bubble using only the sites in $\mathcal{S}$. Let $\gamma$ be one such path. It is contained fully in the interior of our bubble except at the endpoints. 

By definition of a bubble, there must be an open path from $a$ to $b$ adjacent to $\gamma_t$. This path is contained in the part of our bubble bounded by $\gamma$ and $\gamma_t$. Consider the `lowest' such (not necessarily simple) open path,
i.e., the open path from $a$ to $b$ that is the closest to $\gamma$. This is precisely the outer boundary of the closed clusters attached to $\gamma$ from above. This path can contain vertices of $\gamma$ and can touch $\gamma_t$. Call it $\sigma_t$. Consider all of the $R$-adic Whitney boxes needed to cover $\sigma_t$.

Suppose that $B$ is one such box. Then $R^3 B$ intersects either $\eta_b\cup \eta_t$ or $\partial B^{(k)}$. Assume first that the former holds. Since $\eta_b\cup\eta_t$ are entirely outside the bubble bounded by $\gamma_b\cup\gamma_t$, the box $R^3B$ must also intersect at least one of $\gamma_b$ and $\gamma_t$. We will use the fact that these paths intersect $B$ to find that $B$ is one of the boxes that contribute to some $N_{j,B^{(i)}}$ in the definition of the event $H$. Let $l\geq 0$ be the largest non-negative integer such that $x\notin R^l B$. 

Let $z$ be the center of $B$, and suppose that $B\in\cB^{R,k}_{j,n}$ for some $j\geq k$. Let $k-2\leq i\leq\kappa$ be the unique integer such that $\frac{R^{-i-1}}{6}< \textrm{dist}(x,z) \leq \frac{R^{-i}}{6}$. Since $x\notin R^{l}B$, we have that $$\frac{R^{-i}}{6}\geq \dist(x,z)\geq \frac{1}{2}n^{-1}R^{\kappa-j+l}\geq \frac{1}{2}R^{-j+l},$$ and thus $j-i-1\geq l$. Similarly, since $x\in R^{l+1}\cdot B$, we know that $$\frac{R^{-i-1}}{6}\leq\textrm{dist}(x,z)\leq \frac{1}{2}n^{-1}R^{\kappa-j+l+1}<\frac{1}{2}R^{-j+l+2},$$ and thus $l\geq j-i-3$. Finally, let $x'$ be the center of the box in $\cB^{R,k}_{i+3,n}$ that contains $x$. These quantities are introduced as a parallel to the quantities in the definition of the event $H_j$.

\begin{figure}[ht]
    \centering
    \includegraphics[scale=0.7]{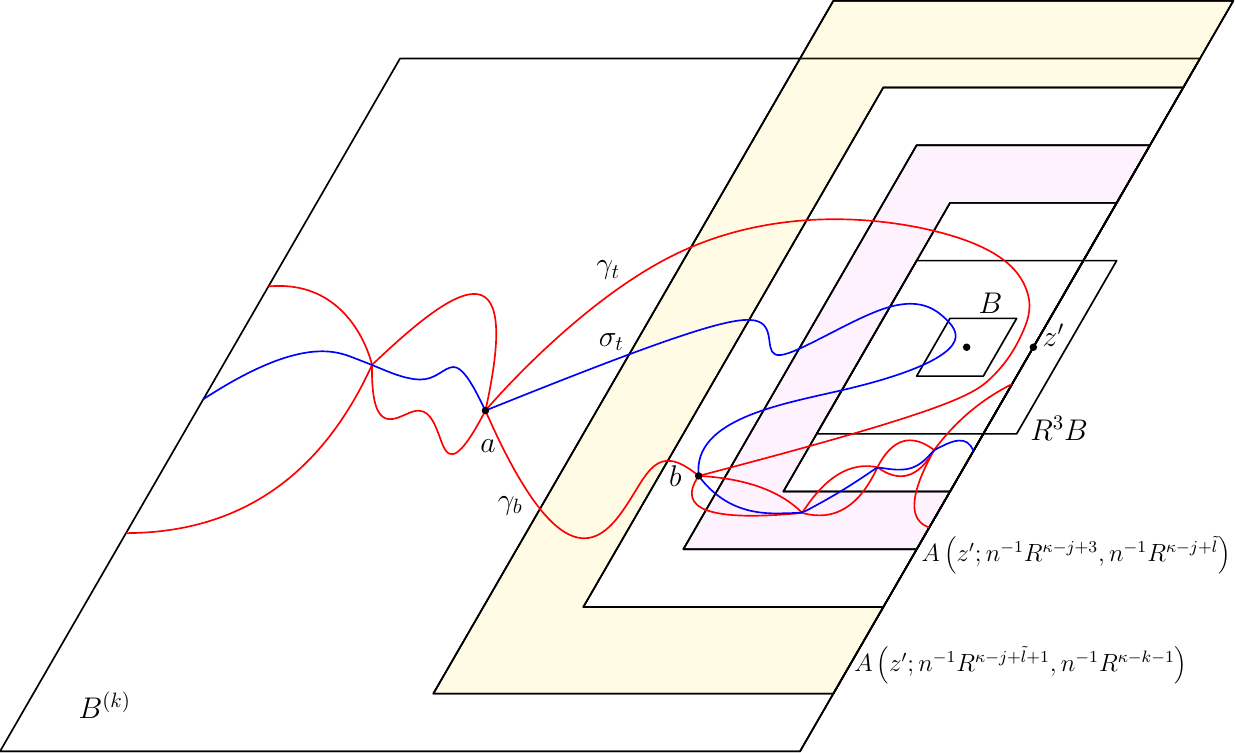}
    \caption{This figure illustrates the case when there is a Whitney box $B$ of level $j$ covering the path $\sigma_t$ that is near the boundary of the ambient box $B^{(k)}$. In the case above, the top boundary of the bubble $\gamma_t$ intersects $R^3B$. We can see that the pink semiannulus above has three arms crossing it - two closed ones corresponding to $\gamma_t$ and one open one coming from $\sigma^t$. The yellow semiannulus is the smallest one whose interior contains at least one of the two bubble endpoints, $a$ and $b$. In this case, it contains $b$. Therefore, we can also identify three arms of alternating color within the yellow semiannulus: two closed ones corresponding to $\gamma_t$ and $\gamma_b$, and one open one corresponding to $\sigma_t$.}
    \label{fig:whitney_box_near_the_boundary_has_three_semiarms}
\end{figure}

There are three cases to consider:
\begin{enumerate}[(a)]
 \item\label{property:intersect-both}  Assume $l\geq 9$ and $R^3 B$ intersects both $\gamma_b$ and $\gamma_t$. Then 
 \begin{equation}
     \dist(z,x')\leq \frac{1}{6}R^{-i}+ 2n^{-1}R^{\kappa-(i+3)}\leq \frac{1}{6}R^{-i}+ 2R^{-i-2}.
 \end{equation}
 Since $j\geq i+l\geq i+9$, this implies that $B=B(z,\frac{1}{2}n^{-1}R^{\kappa-j})\subset B(x';R^{-i}/2)$. Moreover, we know that $R^3 B$ has six arms of alternating colors that go up to distance $n^{-1}R^{\kappa}\cdot\min \{\frac{1}{2}R^{l-j},\delta R^{-k}\}$: one closed and one open arm coming from the section of the boundary of the exploration connecting the box $R^3B$ and $x$ that go at least up to distance $R^l\cdot \frac{1}{2}n^{-1}R^{\kappa-j}$, 
 and two open and two closed arms coming from the open path between the bubbles and the top/bottom closed cluster going up to distance $n^{-1}R^{\kappa}\cdot \min \{\frac{1}{2}R^{l-j},\delta R^{-k}\}$. See Figure~\ref{fig:6_semiarms_case_1} for an illustration of this event. The reason these arms are long even if the box $B$ is close to the boundary of $B^{(k)}$ is that we are working on the event $\mathcal{E}(B^{(k)})$, and on this event we know that all our arms travel at least distance $\delta n^{-1}R^{\kappa-k}$ upon leaving the bubbles. In addition, we know that the closed arm coming from the boundary of the exploration and the closed arms corresponding to the top/bottom closed clusters in $B^{(k)}$ are disjoint inside of the relevant annulus, as they are separated by the open arms corresponding to the exploration and to the open crossing of $B^{(k)}$. See Figure~\ref{fig:6_semiarms_case_1} below.  
 Recall that $l\geq j-i-3$. Moreover, $\min \{\frac{1}{2}R^{l-j},\delta R^{-k}\}\geq R^{l-j-1}$, since $B^{(k)}\subseteq B(z,n^{-1}R^{\kappa-k})$, implying that $l-j\leq -k$ (recall that $l$ is the smallest integer such that $x\in R^{l+1}\cdot B$) and since we can choose $R\geq 1/\delta$. Thus, $B$ has $6$ semiarms of alternating colors of length $l-1\geq j-i-4$.
 \begin{figure}[ht]
    \centering
    \includegraphics[scale=0.7]{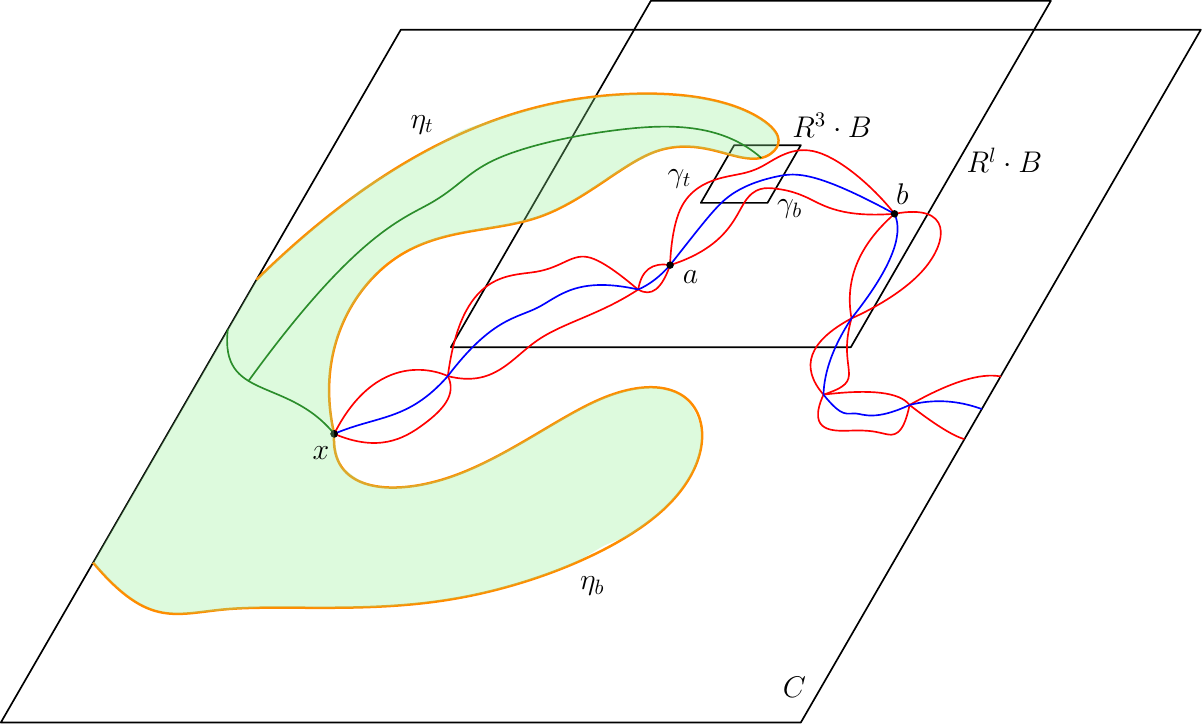}
    \caption{The shaded green area on the image represents the explored part of the open cluster; orange curve represents the closed boundary of the explored region. The red curves represent the closed boundaries of the bubbles created by closed clusters adjacent to the top and bottom of $B^{(k)}$, while the blue curve represents an open left-to-right path delimiting these clusters in the unexplored domain. This figure illustrates the case when both boundaries $\gamma_t$ and $\gamma_b$ of the target bubble intersect the box $R^3B$. The arms corresponding to the boundary of the explored region are long because $x\notin R^{l}B$. For the same reason, both open arms corresponding to the unexplored domain are disjoint from the open arm coming from the explored domain, even though they are all in the same open cluster. Finally, the arms corresponding to the unexplored crossings are long since we are working on the event $\mathcal{E}(B^{(k)})$, and on this event, we know that the open/closed paths must cross distance at least $\delta n^{-1} R^{\kappa-k}$ upon leaving the bubbles to reach the boundary.}
    \label{fig:6_semiarms_case_1}
\end{figure}
 \item\label{property:intersect-one} If $l\geq 9$, but $R^3B$ intersects only one of $\gamma_b,\gamma_t$, then again the box $B(x',R^{-i}/2)$ contains the entire box $B$. Moreover, there must be a closed path from $\sigma_t$ to $\gamma$ starting at $R^3B$, since $\sigma_t$ is the outer boundary of closed clusters adjacent to $\gamma$. This yields that there are $6$ arms in the annulus $R^{\tilde{l}-2}\cdot B\setminus R^3\cdot B$, where $\tilde{l}\geq 0$ is the smallest non-negative integer such that $R^{\tilde{l}}B$ intersects $\gamma$. Indeed, similarly as in the case~\eqref{property:intersect-both} we find that there are $6$ arms going up to distance $n^{-1}R^{\kappa}\cdot\min \{\frac{1}{2}R^{\tilde{l}-1-j},\delta R^{-k}\}$ away from $z$, the center of $B$, and again $\min \{R^{\tilde{l}-1-j},\delta R^{-k}\}\geq R^{\tilde{l}-2-j}$. As before, these arms really are long even if $\sigma_t$ is close to the boundary of $B^{(k)}$ since, again, the arms corresponding to the bubbles must cross distance at least $\delta n^{-1}R^{\kappa-k}$ upon leaving the bubble. See Figure~\ref{fig:6_semiarms_case_2} for an illustration of this event.  
 \newline Further, note that $R^{\tilde{l}}\cdot B$ contains a point in $\gamma$ and a point in $\gamma_t\cup\gamma_b$. Since $\gamma$ is `equidistant' to $\gamma_b$ and $\gamma_t$, and $R>10$, we see that $R^{\tilde{l}+1}\cdot B$ must intersect both $\gamma_b$ and $\gamma_t$. Hence, as in case~\eqref{property:intersect-both}, we find that $R^{\tilde{l}+1}\cdot B$ has 6 arms of alternating colors until the box $B(z, n^{-1}R^{\kappa}\cdot \min \{\frac{1}{2}R^{l-j},\delta R^{-k}\})\supseteq R^{l-1}\cdot B$. As before, $l\geq j-i-3$. Thus, $B$ has 6 semiarms of alternating colors of length $l-1\geq j-i-4$.
 \begin{figure}[ht]
    \centering
    \includegraphics[scale=0.7]{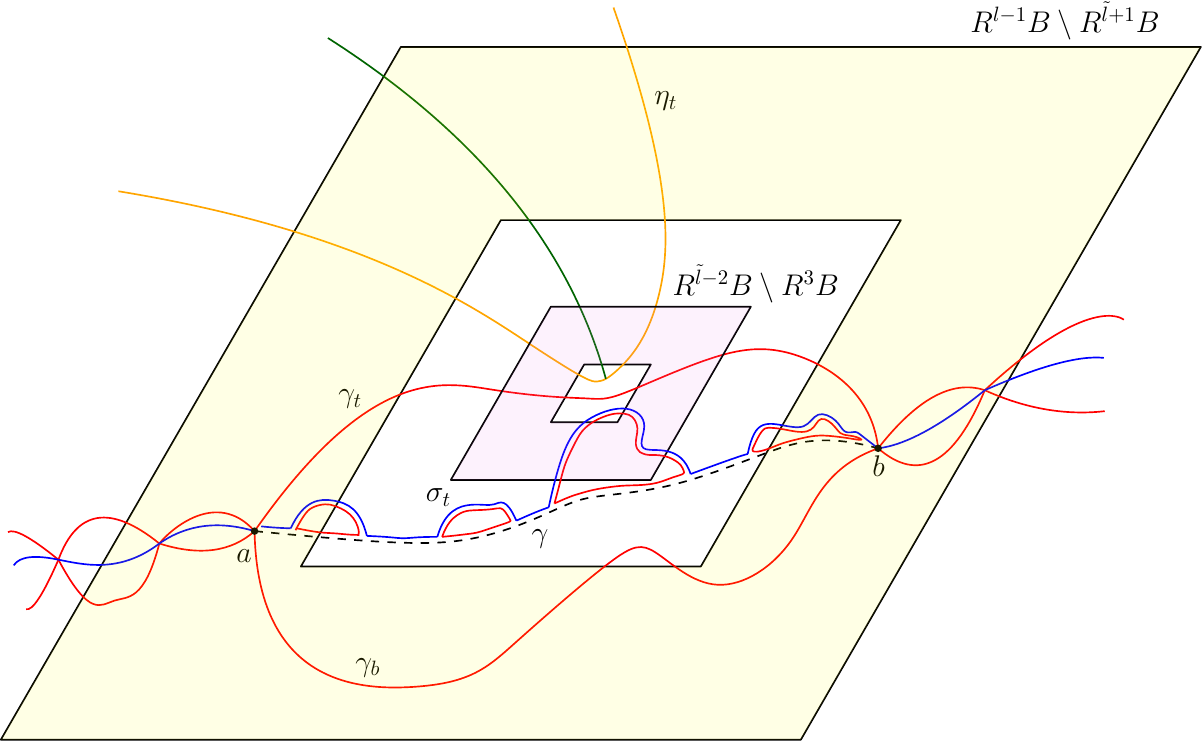}
    \caption{As before, the orange and green curves represent the closed boundary of the explored region and an open path inside that region respectively. The red curves are again the unexplored boundaries of closed clusters creating bubbles. This figure depicts the target bubble with endpoints $a$ and $b$. The black dashed line represents a sequence $\gamma$ of sites `equidistant' to the boundaries of the bubble $\gamma_t$ and $\gamma_b$ that connects $a$ and $b$. Blue curve above it is precisely $\sigma_t$, the boundary of all red clusters attached to $\gamma$ from above. The pink and yellow annuli are precisely the annuli in which we can find $6$ arms of alternating colors. Compare this to Figure~\ref{fig:definition_of_H_i}. }
    \label{fig:6_semiarms_case_2}
\end{figure}
 
\item\label{property:l<9} If $l< 9$, i.e.\ if $x\in R^9 \cdot B$, then we do not have the 6-semiarm event, but importantly, there are at most $6R^{18}$ such boxes $B$ of each level $j\geq k$. 
\end{enumerate}
On the other hand, if we assume that $B$ is such that $R^3 B$ intersects $\partial B^{(k)}$, we will show that $B$ has $3$ semiarms of length $j-k-1$ near the boundary. Since $B$ intersects the interior of our bubble, it must be that $R^3 B$ intersects $\gamma_t$ or $\gamma_b$. Without loss of generality assume it intersects $\gamma_t$. If $z'$ is the projection of the center of $B$ to the nearest side of $B^{(k)}$, then the box $B(z',n^{-1}R^{\kappa-j+3})$ also intersects $\gamma_t$. Let $\tilde{l}$ be the largest integer such that $B(z',n^{-1}R^{\kappa-j+\tilde{l}})$ does not contain the endpoints $a,b$ of our bubble. Then the semiannulus $A(z';n^{-1}R^{\kappa-j+3},n^{-1}R^{\kappa-j+\tilde{l}})\cap B^{(k)}$ has $3$ arms of alternating color -- two closed ones corresponding to the segments of $\gamma_t$ going towards $a$ and $b$ upon leaving $B(z',n^{-1}R^{\kappa-j+3})$, and one open arm corresponding to the open cluster of the bubble. Moreover, since $B(z',n^{-1}R^{\kappa-j+\tilde{l}+1})$ contains $a$ or $b$ it must be that it intersects the top and bottom closed clusters, as well as the open crossing clusters corresponding to $\mathcal{E}(B^{(k)})$. Thus, there are $3$ arms of alternating colors in $A(z';n^{-1}R^{\kappa-j+\tilde{l}},n^{-1}R^{\kappa-k-1})$. In particular, $B$ has $3$ semiarms of length $j-k-1$ near the boundary. See Figure~\ref{fig:whitney_box_near_the_boundary_has_three_semiarms} for an illustration. 

This shows that there exists an open path $\sigma_t$ from $x$ to $v$ that can be covered by Whitney boxes such that for each level $j\geq k$ all but at most $6R^{18}$ many boxes have $6$-semiarms of alternating colors of appropriate length or $3$ semiarms near the boundary. For every $j\geq k$, denote by $\mathcal{W}_j$ the set of Whitney boxes of level $j$ used to cover this path and let $\mathcal{W}_j'$ \phantomsection\label{defn:definition_of_W'} be the number of boxes in 
$\cB_{j,n}^{R,k}$ for which we know that~\eqref{property:intersect-both}, \eqref{property:intersect-one} or~\eqref{property:l<9} holds or, such that this box has $3$ semiarms of length $j-k-1$ near the boundary. One can see that on the event $H$, 
\begin{equation}\label{eq:bound-Wj}
    |\mathcal{W}_j|\leq|\mathcal{W}'_j|\leq 2R^{18}+R^{j\nu}+\sum_{k-2\leq i\leq j-9}R^{j\nu}\leq 2R^{18}+ jR^{j\nu}\leq C_1\cdot R^{\frac{3}{2}j\nu}
\end{equation}
for some $C_1=C_1(R,\nu)$ and all $j\geq k$. 

\subsubsection{Addressing the non-locality}
\label{se:nonlocality}

Our goal is ultimately to say that there exists an improved exponent $\alpha+\beta$ (where $\beta$ does not depend on $\alpha$) such that $\mathbf{P}(N\geq 3)\lesssim R^{-(\alpha+\beta)k}$. In other words, we want to show roughly that $\mathbf{P}(\tau_2<\infty,\mathcal{E}(B^{(k)}))\lesssim \mathbf{P}(\tau_1<\infty,\mathcal{E}(B^{(k)}))\cdot R^{-\beta k}$. As advertised, we will do so by conditioning on all the possible shapes of the `explored domain' $\Gamma_1$. We have already shown that on some high probability event $H$ this domain is not too `ugly', meaning that there is a small number of Whitney boxes covering the open path $\sigma_t$ from $x$ to the last pinch-point of our bubbles. Ideally, we would like to use Lemma~\ref{lem:bound_on_F} and Proposition~\ref{lem:distancebetweencrossings} to then say that the distance across each of these Whitney boxes is suitably small. We could then `string together' the short paths across all the boxes to obtain that on these events $\tau_2=\infty$, as we did in the proof of Lemma~\ref{lem:crossingbubble}. 

Recall that on the event $H$, the number of Whitney boxes of level $j$ intersecting $\sigma_t$ is bounded by $|\mathcal{W}'_j| \le R^{j\nu}$ where $\nu>0$ is a small constant. So by Proposition~\ref{lem:distancebetweencrossings} we expect that on $H$ the conditional probability given~$\Gamma_1$ that $F(B) \cap G(B)^c$ occurs for one such box $B$ is bounded by~$R^{-j\alpha+j\nu}$. However, note that the location of the boxes in $\mathcal{W}_j'$ is not determined by $\Gamma_1$, therefore the correct way of union bounding this event is by summing over \emph{all} boxes $B\in \cB_{j,n}^{R,k}$ and bounding the probability
\begin{equation}\begin{split}\label{eq:bad_box_cond_on_exploration}
&\p\!\left(\cE(B^{(k)}),\,\tau_1<\infty ,\, B \in \mathcal{W}_j' ,\, F(B) \cap G(B)^c\right) \\
&= \p\!\left(\cE(B^{(k)}),\,\tau_1<\infty ,\, B \in \mathcal{W}_j' \right) \p\!\left(F(B) \cap G(B)^c \mmiddle| \cE(B^{(k)}),\,\tau_1<\infty ,\, B \in \mathcal{W}_j' \right) . 
\end{split}\end{equation} 
We would like to use Proposition~\ref{lem:distancebetweencrossings} to bound the last term by $R^{-j\alpha}$, so that~\eqref{eq:bad_box_cond_on_exploration} is bounded by
\[
R^{-j\alpha}\, \p\!\left(\cE(B^{(k)}),\,\tau_1<\infty ,\, B \in \mathcal{W}_j'\right)
\]
which upon summing over all boxes $B\in \cB_{j,n}^{R,k}$ yields
\begin{equation}\begin{split}\label{eq:geo_exploration_sketch}
&\p\!\left(\cE(B^{(k)}),\,\tau_1<\infty ,\, F(B) \cap G(B)^c \text{ for some } B \in \mathcal{W}_j' \right)\\
&\le R^{-j\alpha}\,\E\!\left[|\mathcal{W}_j'| \, \1_{\cE(B^{(k)})} \, \1_{\tau_1<\infty} \right] \\
&\lesssim R^{-j\alpha+j\nu}\,\p\!\left(\cE(B^{(k)}) ,\, \tau_1<\infty\right) +R^{-j\alpha} \cdot R^{2j}\cdot \P(H_j^c) \qquad \text{(by~\eqref{eq:bound-Wj})}\\
&\lesssim R^{-2j\alpha+j\nu} +R^{-j(5\alpha-2)} \qquad \text{(by the a priori bound~\eqref{eq:assumptionsubpoly} and Lemma~\ref{lem:bound_on_global_event}).}
\end{split}\end{equation}
However, as we discussed at the beginning of Section~\ref{subsec:bootstrapping} (see Figure~\ref{fig:exploration_is_not_local}), for $S\subseteq B^{(k)}$, the conditional law given~$\{\Gamma_1 = S\}$ and $\cE(B^{(k)}) \cap \{\tau_1<\infty\}$ is not just an i.i.d.\ percolation in $B^{(k)}\setminus S$. We will describe the conditional law more precisely on page~\pageref{page:Ein_Eout} below. As we explained above, the event $\{\Gamma_1 = S\}$ \emph{may} depend on the sites in $B$ even if $B \subseteq B^{(k)}\setminus S$, therefore we cannot apply Lemma~\ref{lem:bound_on_F} and Proposition~\ref{lem:distancebetweencrossings} directly.

To combat this, let us fix a box $B$ as above, and let $\tilde{\omega}$ denote the percolation configuration outside~$4B$. We will define an auxiliary event $\mathscr{P}(\tilde{\omega})$ that only depends on the sites outside $B$ and simulates the event $\cE(B^{(k)}) \cap \{\Gamma_1 = S\}$. More precisely, it is the event that \emph{one can construct a configuration in $B$} that makes $\cE(B^{(k)}) \cap \{\Gamma_1 = S\}$ occur. We will then see that on another good event $\check{F}$ that depends only on the sites in the annulus $4B \setminus B$, the conditional probability of $\cE(B^{(k)}) \cap \{\Gamma_1 = S\}$ given the sites outside $4B$ is comparable to the conditional probability of~$\mathscr{P}(\tilde{\omega})$. If we now condition on just $\mathscr{P}(\tilde{\omega})$ and the sites outside $B$, we will be able to apply Proposition~\ref{lem:distancebetweencrossings} to the box $B$.

The proof proceeds in three steps.
\begin{itemize}
    \item Lemma~\ref{lem:number-crossings-check}: Here we define good events $\check{F}_1,\check{F}_2$ which state that the annulus does not have too many crossings and that the crossing clusters have enough `thickness' near the tips of their interfaces on the inner boundary. We show that these events occur with high probability.
    \item Lemma~\ref{lem:E_out_happens_with_positive_probability}: This states that on a good event $\check{F}$ defined in terms of $\check{F}_1,\check{F}_2$, the conditional probability of $\cE(B^{(k)}) \cap \{\Gamma_1 = S\}$ is comparable to the conditional probability of the auxiliary event $\mathscr{P}(\tilde{\omega})$. This will be used to argue that it suffices to bound a variant of~\eqref{eq:bad_box_cond_on_exploration} where we instead condition on the event $\mathscr{P}(\tilde{\omega})$.
    \item Lemma~\ref{lem:bound_on_stopping_time}: Here we conclude the proof of $\P(\tau_2<\infty,\cE(B^{(k)}))\lesssim\P(\tau_1<\infty,\cE(B^{(k)}))\cdot R^{-\beta k}$, following the strategy described in the paragraphs above. In particular, we use the last two steps to resolve the problem of non-locality described above. See the beginning of Section~\ref{se:geo_exploration_conclusion} for more details.
\end{itemize}

\begin{lem}\label{lem:number-crossings-check} 
Consider an annulus $A(x;r,2r)$. Let $\check{F}_1(\check{M})$ be the event that there are at most $\check{M}$ distinct clusters (open or closed) that cross this annulus. If there exist both an open crossing and a closed crossing of $A(x;r,2r)$, then for any open (resp.\ closed) cluster that crosses this annulus, one can find two (not necessarily disjoint) open (resp.\ closed) crossings of the annulus $\eta_c$ and $\eta_a$ that serve as the clockwise most and counterclockwise most boundaries of this cluster. 
Define the event $\check{F}_2(\check{b})$ as follows. If there are at least 2 crossing clusters of $A(x;r,2r)$, then for every open (resp.\ closed) crossing cluster of $A(x;r,2r)$, let $v$ be the endpoint of $\eta_a$ on the inner boundary of $A(x;r,2r)$. On $\check{F}_2$, we require that there exists an open (resp.\ closed) path from $\eta_a$ to the inner boundary of $A(x;r,2r)$ that stays at distance at least $10\check{b}r$ away from $v$. We also require the same assertion when we swap $\eta_a$ with $\eta_c$. In particular, this implies that for points $v$ that are the endpoints of the interfaces of crossing clusters at the inner boundary of $\check{A}$, no crossing clusters apart from the ones whose interface $v$ is on enter the box $B(v,10\check{b}r)$. 

The following holds for $R$ large enough. For each $\check{p}>0$, there exist $\check{M}$ large enough and $\check{b}>0$ small enough so that $$\P\!\left(\check{F}_1(\check{M})\cap \check{F}_2(\check{b})\right)\geq 1-\check{p}$$ for all annuli $A(x;r,2r)$ of radius $r\geq \frac{R}{\check{b}n}$.
\label{lem:there_exists_box_for_crossing_clusters} 
\end{lem}
\begin{proof} The proof of this lemma is very similar to the proof of Lemma~\ref{lem:number_crossings}. To start off, if there are $\check{M}>2$ different crossing clusters, then there are $\frac{\check{M}}{2}$ distinct open paths in the annulus. By the BK inequality we can find $\check{M}=\check{M}(\check{p})$ large enough so that $\P(\check{F}_1(\check{M})^c)\leq \frac{\check{p}}{2}$.

We now turn our attention to the event $\check{F}_2$. Note that on $\check{F}_2^c$, there must exist at least two crossing clusters of $A(x;r,2r)$, and there must be a crossing cluster of $A(x;r,2r)$ and a point $v$ that is an endpoint of $\eta_a$ or $\eta_c$ corresponding to this cluster such that there does not exist a path inside this cluster from $\eta_a$ (resp.\ $\eta_c$) to the inner boundary of $A(x;r,2r)$ outside of the box $B(v;10\check{b}r)$. Assume without loss of generality that we are looking at $\eta_a$, and that the relevant cluster is open. Then, if we cannot find an open path connecting $\eta_a$ and the inner boundary of $\check{A}$ outside of the box $B(v,10\check{b}r)$, then there must be a blocking path of closed sites that crosses the annulus $A(v;10\check{b}r, \frac{1}{3}r)$, say. Thus, there are three arms of alternating colors crossing this annulus -- the two coming from the interface and one blocking one. Furthermore, these arms are inside $A(x;r,2r)$, so they must belong to a sector of angle $5\pi/3$. Thus, similarly as in the proof of Lemma~\ref{lem:number_crossings}, we can use the $3$-arm exponent in the half-plane and conformal invariance\footnote{As in Lemma~\ref{lem:number_crossings}, we can also do it without using conformal invariance.} to conclude that for all $\check{b}$ small enough $\P(\check{F}_2(\check{b})^c)\leq \check{p}/2$ whenever our annulus has radius $r$ such that $r\geq \frac{R}{\check{b}n}$. This concludes the proof. 
\end{proof}

\phantomsection\label{page:bound_on_check_F}
We can now define the event $\check{F}$. Recall the definition of events $\hat{F}$ and $F$ on page~\pageref{defn:definition_of_event_F}. Fix some $\check{p}>0$. Consider an annulus $A(x;r,4r)$. Then we define $\check{F}$ as the intersection of the following events: 
\begin{enumerate}
\item[$\bullet$] $\check{F}_1(A(x;2r,4r))\cap \check{F}_2(A(x;2r,4r))$, 
where the parameters $\check{M}_1,\check{b}_1$ are chosen so that $$\P(\check{F}_1\cap \check{F}_2)\geq 1-\check{p}/4,$$ for all annuli of radius $r$ such that $\check{b}_1rn\geq R$.

\item[$\bullet$] Choose $\check{c}_1\in (0,\check{b}_1/\check{M}_1^2)$. Cover the annulus $A(x;3r/2,5r/2)$ with at most $32\check{c}_1^{-2}$ boxes of side length $\check{c}_1r$. Assume that the parameters $\check{M}_2,\check{a}_2,\check{b}_2,\check{c}_2$ (depending on $\check{p}$ and $\check{c}_1$) are such that $\P(\hat{F}(2B\setminus B))\geq 1-\frac{\check{p}\check{c}_1^2}{128}$ for any box $B$ of side length $\check{c}_1r$ satisfying $\check{c}_2\check{c_1}rn\geq R$. We assume that $\hat{F}(2B\setminus B)$ and $\hat{F}'(2B\setminus B)$ hold for all boxes $B$ in our covering. 
\item[$\bullet$] There does not exist a point $w$ in the grid of side length $\check{c}_1\check{c}_2 r$ inside $A(x;3r/2,5r/2)$ such that there are $6$ arms of alternating colors crossing $A(w;2\check{c}_2\check{c}_1r,\frac{1}{2}\check{a}_2\check{c}_1r)$. Using Lemma~\ref{lem:arm_exponents} one can see that by choosing $\check{c}_2$ small enough (depending on $\check{p},R$ and $\check{a}_2$), then the probability of this event can be made larger than $1-\check{p}/8$.
\item[$\bullet$] Consider also the grid of equally spaced points at distance $\frac{\check{c}_1\check{c}_2}{8}r$ away from each other inside $A(x;3r/2,5r/2)$. Consider all the boxes of side length $\frac{\check{c}_1\check{c}_2}{8}r$ centered at these points. Recall Lemma~\ref{lem:omega_is_well_separated} and the discussion preceding it. We will assume that for every box $B$ as above, every set of crossings of the annulus $2B\setminus B$ can be made $\check{\xi}$-well separated. We will choose $\check{\xi}>0$ depending on $\check{p},R,\check{c}_1,\check{c}_2$ so that the probability of this event is larger than $1-\check{p}/8$ for all radii $r$ such that $\check{\xi} \frac{\check{c}_1\check{c}_2}{8}rn\geq R$. Lemma~\ref{lem:omega_is_well_separated} ensures that we can make this choice.
\end{enumerate}

By a union bound, we can see that we have chosen all the parameters of $\check{F}$ so that $\P(\check{F}(A(x;r,4r)))\geq 1-\check{p}$ for all annuli of radius $r$ satisfying $\check{\xi} \frac{\check{c}_1\check{c}_2}{8}rn\geq R$. 

\phantomsection\label{page:Ein_Eout}
As advertised earlier, we will condition our percolation configuration on the shape of $\Gamma_1$. Suppose that a set $S\subseteq B^{(k)}$ is in the support of $\Gamma_1$ and $x_S$ is a designated point on the boundary of $S$ that serves as the tip of our `explored domain' $\Gamma_1$ called $x$. We will find it useful to write the event $\mathcal{E}(B^{(k)})\cap \{\tau_1<\infty\}\cap \{(\Gamma_1,x)=(S,x_S)\}$ as the intersection of two events that we define below. 

The boundary of $S$ consists of three sections -- a segment of the left-hand side of $B^{(k)}$, a section connecting the endpoint of this segment belonging to the top half of the left-hand side of $B^{(k)}$ to $x_S$, called $\eta_t^S$ and an analogous segment connecting the bottom half to $x_S$, called $\eta _b^S$. Moreover, since $S$ is in the support of $\Gamma_1$, there exists at least one site other than $x_S$ such that it is a neighbor of both $\eta_t^S$ and $\eta_b^S$.
\begin{enumerate}
\item[(1)] Define $\mathcal{E}_{\textrm{in}}(S,x_S)$ to be the event that all sites of $\eta_b^S\cup \eta_t^S$ but $x_S$ are closed, and 
have a neighboring site that is open and connected to $x_S$ via an open path inside $S$. 
Moreover, we require that there exists an open path in $S$ connecting $x_S$ to the segment of length $2\delta n^{-1}R^{\kappa-k}$ on the left-hand side of $B^{(k)}$, centered at its midpoint. Finally, we require that the geodesic distance (resp.\ effective resistance) inside of $S$ between the first touch-point of $\eta_b$ and $\eta_t$ to $x_S$ to be larger than $\altq_n$, and there isn't another touchpoint before $x_S$ that makes the distance larger than $\altq_n$. 
\item[(2)] Define the event $\mathcal{E}_{\textrm{out}}(S,x_S)$ to be the event that there exists an open path in $S^c$ connecting $x_S$ to the segment of length $2\delta n^{-1}R^{\kappa-k}$ centered at the midpoint of the right-hand side of $B^{(k)}$. The open cluster of this path in $S^c$ cannot come within distance $\delta n^{-1}R^{\kappa-k}$ from the centers of top and bottom sides of $B^{(k)}$. Further, there is a closed path in $B^{(k)}\cap S^c$ connecting $\eta_t^S$ to the segment of length $2\delta n^{-1} R^{\kappa-k}$ centered at the midpoint of the top side of $B^{(k)}$. The collection of all closed clusters in $S^c$ that are connected to the segment $\eta_t^S$ contains a left-to-right crossing of $B^{(k)}$. Moreover, this set does not come within distance $\delta n^{-1} R^{\kappa -k}$ from the midpoints of the other three sides of $B^{(k)}$. Assume further that the analogous assertion holds for $\eta_b^S$ and the bottom side of $B^{(k)}$.
\end{enumerate}
Thus, we have written $\mathcal{E}\cap \{\tau_1<\infty\}\cap \{\Gamma_1=(S,x_S)\}$ as an intersection of two events, where the first of the two events depends only on the sites of $S$ whereas the other one depends only on the sites of $S^c$. Note that the fact that the top and bottom closed clusters do not come within distance $2\delta n^{-1}R^{\kappa-k}$ of $c_4^{(k)}$ comes for free when $S$ is chosen in the support of $\Gamma_1$, as well as the property that on $\mathcal{E}(B^{(k)})$ the top and bottom closed cluster come within graph distance one of each other. 

\phantomsection\label{page:definition_of_a_nice_box}
We are now ready to state and prove the following lemma that will enable us to apply Lemma~\ref{lem:bound_on_F} and Proposition~\ref{lem:distancebetweencrossings}, as advertised. Fix the set $(S,x_S)$ as above, and consider a box $B=B(x,4r)$ of radius $4r>0$ that is contained inside of $B^{(k)}\cap S^c$. We introduce a technical condition on $B$ -- we say that $B$ is \textit{nice} in $B^{(k)}$ if the following holds. As before, denote by $c_1^{(k)},\ldots,c_4^{(k)}$ the midpoints of the sides of $B^{(k)}$ in the clockwise order, starting from the top side. If $\frac{1}{2}B$ intersects $B(c_1^{(k)},2\delta n^{-1}R^{\kappa-k})$, but is not fully contained inside this box, then there are one or two sides of $\frac{1}{2}B$ that do not intersect $B(c_1^{(k)},2\delta n^{-1}R^{\kappa-k})$. We require that these sides stay at distance at least $2\check{a}_1 r$ from $B(c_1^{(k)},2\delta n^{-1}R^{\kappa-k})$.
We assume that the same holds with $c_2^{(k)}$ and $c_3^{(k)}$ in place of $c_1^{(k)}$. 

Let $\tilde{\omega}$ denote the percolation configuration in $B^{(k)}$ outside of the box $B$. Define $\mathscr{P}=\mathscr{P}(\tilde{\omega})$ to be the event on $\check{A}:=A(x;r,4r)$ that ensures that $\mathcal{E}_{\textrm{out}}$ is possible. In other words, $\mathscr{P}(\tilde{\omega})$ is the set of all percolation configurations $\omega_0(\check{A})$ on the annulus $\check{A}$ with the property that there exists at least one percolation configuration $\omega_0(\operatorname{int}\check{A})$ in $\operatorname{int}\check{A}$ such that $\mathcal{E}_{\textrm{out}}(S,x_S)$ happens when the full percolation configuration on $B^{(k)}$ is given by $\tilde{\omega}\cup \omega_0({\check{A}})\cup \omega_0(\operatorname{int}\check{A})$. Here, $\intt\check{A}$ denotes the inner face of the annulus $\check{A}$.

\begin{lem} \label{lem:E_out_happens_with_positive_probability}
The following holds for all $R$ large enough. There exists a constant $\check{g}>0$ depending on $\check{p},R,\check{M}_1,\check{M}_2,\check{c}_1,\check{c}_2,\check{\xi}$ such that whenever $(S,x_S)$ and $B=B(x,4r)$ are as above, $$\P\!\left(\check{F}\!\left(A(x;r,4r)\right), \mathscr{P}(\tilde{\omega})\bigg|\,\tilde{\omega}\right) \leq \check{g}\cdot \P(\mathcal{E}_{\textrm{out}}(S,x_S)|\,\tilde{\omega}).$$ 
\end{lem}

The proof goes as follows. On $\mathscr{P}(\tilde{\omega})$, there exists some way to connect the clusters of $\check{A}$ inside $\intt(\check{A})$ such that the event $\cE_{\out}(S,x_S)$ occurs, i.e., such that the correct open and closed paths exist in $S^c$. This is reminiscent of Proposition~\ref{lem:distancebetweencrossings}, where clusters in the annulus needed to be connected to the sides of the box $B^{(k)}$ in the correct way in order for the event $\cE(B^{(k)})$ to hold. The proof of this lemma will use a similar resampling procedure. The proof consists of several steps. The main part will be to show that on the event $\mathscr{P}(\tilde{\omega})\cap \check{F}(A(x;r,4r)$ we can find a finite set of disjoint `corridors' in the grid of side length $\check{c}_1$, such that if we resample suitably along these corridors, the resulting percolation configuration will satisfy $\mathcal{E}_{\textrm{out}}(S,x_S)$. The first step is to examine the event $\mathscr{P}(\tilde\omega)$ and define some terminology regarding cluster connections. The second step is to construct the corridors in the case that the annulus $\check{A}$ is not near the midpoints of the sides of $B^{(k)}$. In the case when the annulus is near the sides, there is a small technical complication due to how the event $\cE(B^{(k)})$ is defined. In particular, we must not interfere with the conditions in the items~\eqref{property:E3}--\eqref{property:E5} in Section~\ref{subsubsec:apriori-sketch}. We address this in Step~3 of the proof. Finally, we show that on the event $\check{F}$, the probability that the resampling is successful is strictly positive.
\begin{proof}[Proof of Lemma~\ref{lem:E_out_happens_with_positive_probability}]

\medskip \textbf{Step 1:} \textit{Preliminary observations.}

Fix a percolation configuration $\omega_0(\check{A})\in \mathscr{P}(\tilde{\omega})$ and consider an instance of percolation $\omega_0(\operatorname{int}\check{A})$ on $\operatorname{int}\check{A}$ that makes $\mathcal{E}_{\textrm{out}}(S,x_S)$ happen.

If there exists a monochromatic circuit in $\omega_0(\check{A})$, then we claim that $\mathcal{E}_{\textrm{out}}(S,x_S)$ holds regardless of what the percolation in $\operatorname{int}\check{A}$ is. Indeed, by changing the percolation configuration in $\operatorname{int}\check{A}$, we can add/subtract certain sites of $\operatorname{int}\check{A}$ to/from the cluster of our loop and we can change the structure of the clusters strictly surrounded by this loop. All clusters outside this loop remain unchanged. Thus, if $\mathcal{E}_{\textrm{out}}(S,x_S)$ holds for one such configuration, it must hold for all others. Indeed, if a changed configuration causes $\mathcal{E}_{\textrm{out}}(S,x_S)$ to fail, it must be because there are some additional sites in the cluster of our loop that are within distance $2\delta n^{-1}R^{\kappa-k}$ from one of $c_1^{(k)},c_2^{(k)}$ and $c_3^{(k)}$. But then, the loop itself would have to intersect one of these boxes, and $\mathcal{E}_{\textrm{out}}(S,x_S)$ couldn't have held for the initial configuration either. In this case, we will say that our set of corridors is empty.

Assume now that there isn't a monochromatic circuit in $\check{A}$. Then there exists at least one open and at least one closed crossing cluster in $\check{A}$. Denote by $C_1,O_1,\ldots, C_m,O_m$ all clusters in $\check{A}$ that cross the annulus in clockwise order, starting from any crossing closed cluster. Several of these $\check{A}$-clusters can be a part of the same cluster when we consider the full percolation configuration in $B^{(k)}$ or $B$. 
In particular, we will say that two closed $\check{A}$-clusters in $\{C_1,\ldots,C_m\}$ are equivalent if and only if they are a part of the same closed cluster when considering percolation on $B$, i.e.\ if there is a closed path in $\omega_0(\operatorname{int}\check{A})$ connecting these. We analogously define the notion of equivalence of open clusters. 

Assume first that the box $B$ does not intersect any of the boxes of radii $2\delta n^{-1}R^{\kappa-k}$ centered at the points $c_1^{(k)},c_2^{(k)}$ and $c_3^{(k)}$. We claim that in this case, if we substitute $\omega_0(\intt\check{A})$ by any other percolation configuration in $\intt\check{A}$ that connects all pairs of equivalent $\check{A}$-clusters, the event $\mathcal{E}_{\textrm{out}}(S,x_S)$ still holds on $B^{(k)}$. To see that this is true, one can use an argument similar to the case when there is a monochromatic circuit in $\omega_0(\check{A})$. Changing percolation in $\intt\check{A}$ from $\omega_0(\intt\check{A})$ to any other configuration that preserves connections between equivalent $\check{A}$-clusters leaves the cluster structure outside of $B$ unchanged.\footnote{This does not mean only that the clusters in the percolation restricted to $B^{(k)}\setminus B$ remain unchanged. It really means that two sites of $B^{(k)}\setminus B$ are in the same cluster in $B^{(k)}$ (when we have $\omega_0(\textrm{int }\check{A})$ on $\textrm{int } \check{A}$) if and only if they are in the same cluster on $B^{(k)}$ with the modified percolation on $\textrm{int }\check{A}$.} The only changes that can occur in terms of the cluster structure concern the sites in $B$ that are surrounded by the loop made up of inner boundaries of $C_1,O_1,\ldots,C_m,O_m$. In particular, these sites can be joined or deleted from the $B^{(k)}$-clusters that reach outside of $B$ and they can form clusters among themselves. As before, the only way for $\mathcal{E}_{\textrm{out}}(S,x_s)$ to fail when making this change would be if we added sites of some box of radius $2\delta n^{-1} R^{\kappa-k}$ around $c_1^{(k)},c_2^{(k)},c_3^{(k)}$ to the wrong `cluster that leaves $B$'. But because in this case $B$ is disjoint from these boxes, this is impossible. Hence, $\mathcal{E}_{\textrm{out}}(S,x_S)$ holds for the modified percolation as well. 

\medskip \textbf{Step 2:} \textit{Construction of corridors.}

We will now choose a collection of `corridors' depending on percolation in $\check{A}$ and $\tilde{\omega}$ and show that when we suitably resample along these corridors, the event $\mathcal{E}_{\textrm{out}}(S,x_S)$ holds for the resulting configuration. Denote by $x_1, y_1,\ldots,x_{m}, y_m$ the endpoints on the inner boundary of $\check{A}$ of the open/closed interfaces between clusters $C_1,O_1,\ldots,C_m,O_m$ starting with the one in between $C_1$ and $O_1$. In particular, $x_{i}$ is on the interface between $C_i$ and $O_i$, while $y_{i}$ is on the interface between $O_i$ and $C_{i+1}$, where we view the indices modulo $m$. See Figure~\ref{fig:picture_of_D_cor} below. Since we are working on the event $\check{F}(A(x;r,4r))$, we know that each box $B(x_i,\check{b}_1\cdot 2r)$ is `shielded away' from all crossing clusters other than $C_i$ and $O_i$ by a closed path in $C_i$ and an open path in $O_i$. The analogous statement holds for the boxes centered at points $y_i$. Consider a grid of side length $\frac{\check{b}_1}{2} \cdot 2r$ covering the box $B$. Denote by $\mathcal{D}_{\textrm{cor}}$ the domain consisting of the following boxes: \begin{enumerate}
\item[$\bullet$] $2B$ for every box $B$ in our grid that contains one of the points $x_i$ or $y_i$;

\item[$\bullet$] Include $2B'$ for every box $B'$ in our grid such that $3B'$ is fully contained in $\textrm{int }\check{A}$, but $5B'$ is not. 
We will call these boxes the loop part of $\mathcal{D}_\textrm{cor}$; 
\item[$\bullet$] if there exists a box $B$ containing one of $x_i$ or $y_i$ that also contains one of the four corners of the inner boundary of $\check{A}$, then $2B$ only shares one corner with the loop part of $\mathcal{D}_\textrm{cor}$. In this case, we will include box $2B''$ in our set, where $2B''$ is one of the two boxes that share a side with $2B$ and with the loop part of $\mathcal{D}_\textrm{cor}$.

\end{enumerate} See Figure~\ref{fig:picture_of_D_cor} for an illustration. The reason why we consider boxes that are two times larger than the original ones is simply because when selecting a box $B$ containing some $x_i$ or $y_i$, it can happen that this $x_i$ or $y_i$ is close to the boundary of the box. In this case, it may be that a $\check{c}_1$-level box containing $x_i$ or $y_i$ (that we will later use for resampling) is not fully contained inside $B$. However, it will be fully contained within $2B$.

\begin{figure}[ht]
    \centering
    \includegraphics[scale=0.8]{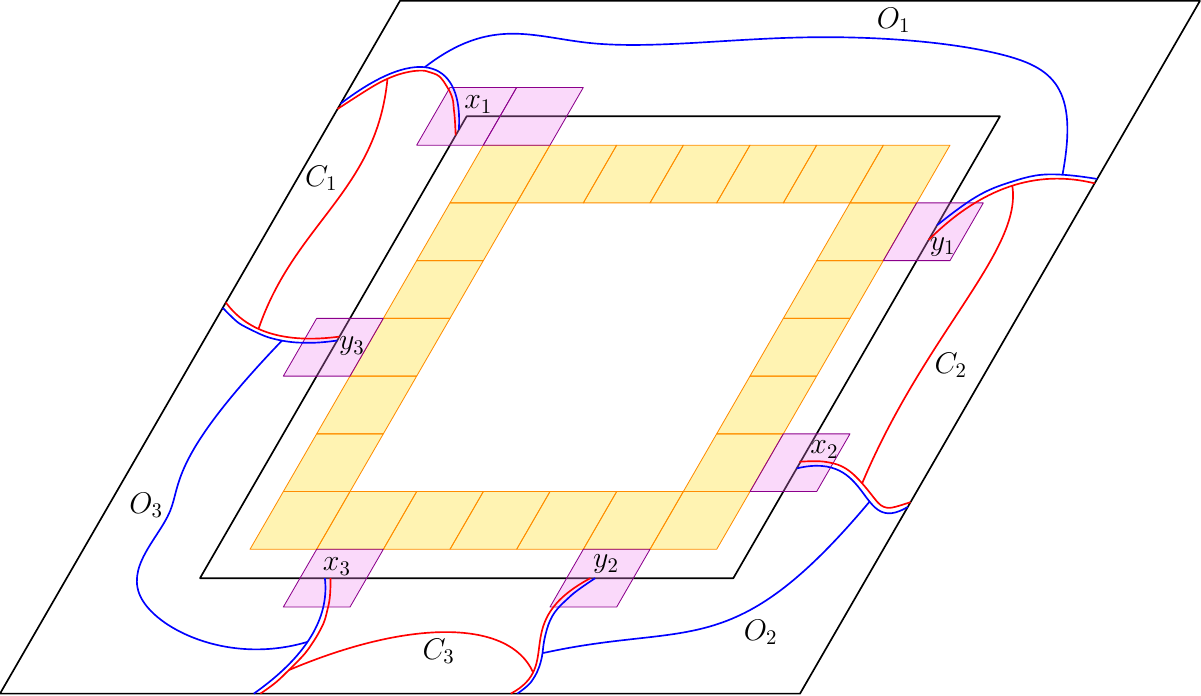}
    \caption{The picture above portrays the annulus $\check{A}$ with all of its crossing clusters. In this case there are $6$ crossing clusters - $3$ open and $3$ closed ones. The red/blue curves represent the interfaces between these clusters. The endpoints of these interfaces on the inner boundary of $\check{A}$ are denoted by $x_i$-s and $y_i$-s. The purple and yellow shaded boxes represent the domain $\mathcal{D}_{\textrm{cor}}$. Yellow boxes represent the loop part of this corridor. Note that in the example above, the box containing $x_1$ contains the corner of the inner boundary of $\check{A}$, so we have included a neighboring purple box in $\mathcal{D}_{\textrm{cor}}$. Blue and red curves connecting the consecutive interfaces depict the open (resp.\ closed) paths that `shield' the purple boxes in $\mathcal{D}_{\textrm{cor}}$ away from all the crossing clusters other than the two whose boundary endpoint belongs to that box. Their existence is guaranteed by the event $\check{F}$.}
    \label{fig:picture_of_D_cor}
\end{figure}

We will construct our collection of corridors within this domain $\mathcal{D}_{\textrm{cor}}$ using the grid of side length $\check{c}_1 \cdot 2r$. The end boxes of these corridors will be the boxes in the $\check{c}_1$-grid containing the points $x_i$ and $y_i$. Suppose that we wish to create a corridor connecting points $z,z'\in \{x_1,y_1,\ldots,x_m,y_m\}$. This corridor will consist of a sequence of boxes going `inwards' from the inner boundary of $\check{A}$, starting with the box in our grid that contains $z$ and reaching inside the loop part of $\mathcal{D}_{\textrm{cor}}$. It will then `trace along the loop part of $\mathcal{D}_{\textrm{cor}}$' and finally go straight `outwards' towards the box containing $z'$. The exception is when either of the $\check{b}_1$-boxes containing $z$ and $z'$ contain a corner of the inner boundary of $\check{A}$ -- in this case, we will first go along the boundary of $\check{A}$ to get away from the corner and then proceed to go `inwards' from there. These corridors could potentially intersect within the loop part of $\mathcal{D}_\textrm{cor}$, so we will construct different corridors at different `depths' of the loop to avoid this. Namely, since the loop part of $\mathcal{D}_{\textrm{cor}}$ has width $\check{b}_1\cdot 2r$ and $\check{b}_1 \gg \check{c}_1 \cdot \check{M}_1^2$, we will consider $\check{M}_1^2$ spaced-apart concentric sequences of boxes in the $\check{c}_1$-grid that are contained in the loop part of $\mathcal{D}_{\textrm{cor}}$ and trace the shape of a box. We will refer to these as levels, first of which will be the outermost loop. See Figure~\ref{fig:construction_of_corridors_at_different_levels} below. 

We are now ready to construct our sequence of corridors. If we consider the equivalence classes in $\{C_1,O_1,\ldots,C_m,O_m\}$, there must be a cluster (open or closed) whose equivalence class is a singleton. Indeed, consider the equivalence class of $C_1$. If it is a singleton, then we are done. If not, suppose that $i$ is the smallest index such that $C_1$ and $C_i$ are equivalent. Then, it must be that equivalence classes of elements in $\{O_1,C_2,O_2,\ldots,O_{i-1}\}$ stay within this set. We can repeat our reasoning with the cluster $O_1$, and continue recursively until we find a cluster whose equivalence class is a singleton. Note that in each step after the first one, the number of available clusters is odd and hence non-zero, thus we can always continue the procedure. We will eventually stop this recursion because the number of available clusters decreases in every step. 

To construct the first corridor pick any cluster whose equivalence class is a singleton. Without loss of generality assume that this cluster is $O_1$. Consider a corridor obtained by connecting $x_1$ and $y_1$ via the first level. There are two possible choices for the arc of the first level loop -- choose the one such that $\frac{1}{4}B$ is not surrounded by the union of this corridor and the cluster $O_1$. This ensures that all the other points $x_i,y_i$ can `reach' levels two and above without intersecting this first corridor. Declare this corridor to be closed (if we had started with a closed cluster, we would declare this corridor to be open).

If we were to make a closed (resp.\ open) connection along this corridor, 
we would be left with $m-1$ open and $m-1$ closed clusters -- two closed clusters $C_1$ and $C_2$ would merge, while $O_1$ is no longer relevant for us, as it is closed off from the rest of the clusters by the first corridor. In this new set of `clusters' there must be a cluster whose equivalence class is a singleton. Construct the second corridor analogously to the above, except that it will travel through the loop part of $\mathcal{D}_{\textrm{cor}}$ via level $2$. See Figure~\ref{fig:construction_of_corridors_at_different_levels} for an illustration of this construction.

\begin{figure}[ht]
    \centering
    \includegraphics[scale=0.8]{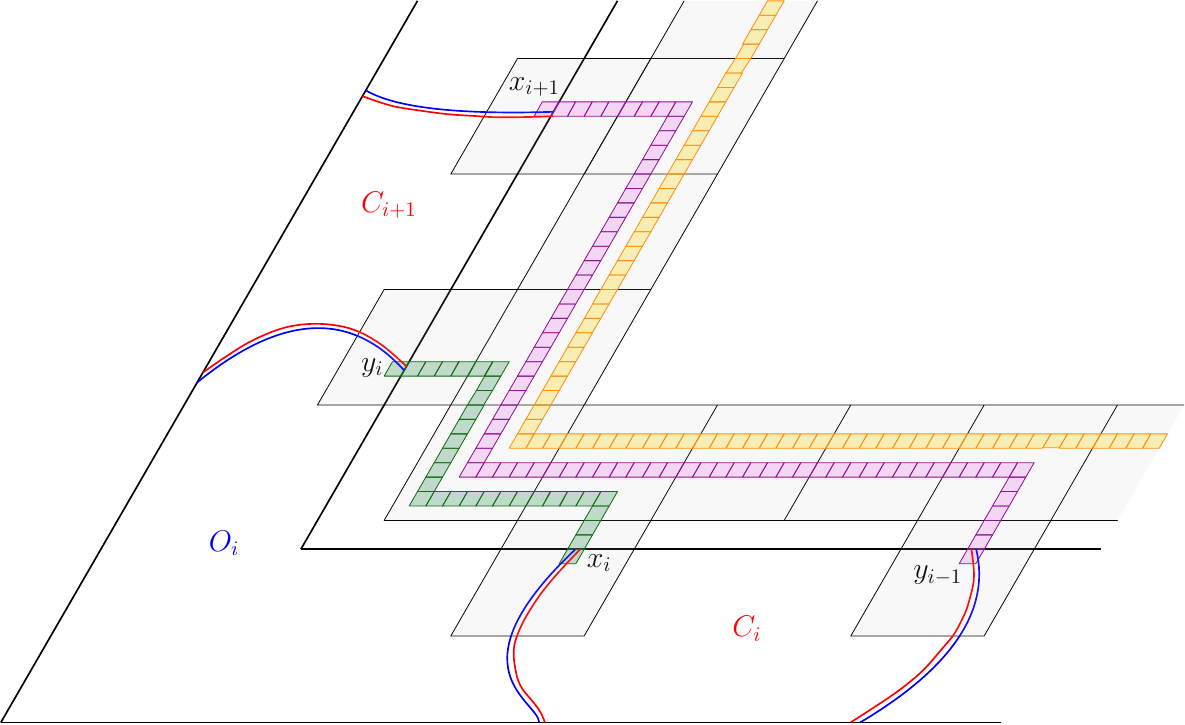}
    \caption{On the figure above, we can see a sector of the annulus $\check{A}$ together with $4$ interfaces between consecutive closed cluster represented by the red/blue curves. These interfaces delimit the clusters $C_i,O_i$ and $C_{i+1}$, as shown. The region shaded in gray represents $\mathcal{D}_{\textrm{cor}}$. In this instance, we have $3$ different levels within the loop part of $\mathcal{D}_{\textrm{cor}}$. Green, purple and yellow shaded regions represent three corridors at these three distinct levels. The green corridor connects $x_i$ to $y_i$, the purple one connects $y_{i-1}$ to $x_{i+1}$, while the endpoints of the yellow corridor are outside of the pictured sector of $\check{A}$. Note that the trajectory of each corridor is as described earlier: we first go `inwards', then we `go along the loop', and finally we return `outwards' towards the target endpoint. In the case pictured above, $O_i$ is alone in its equivalence class and the green corridor is the one constructed in the first iteration of our algorithm. It will be regarded as closed, and as such, it will serve to at the same time connect closed clusters $C_i$ and $C_{i+1}$ and prevent $O_i$ from connecting to any other open clusters in $\check{A}$.}
    \label{fig:construction_of_corridors_at_different_levels}
\end{figure}

Continue recursively. We will stop this procedure after $m-1$ steps, since in every turn the number of relevant clusters decreases by $2$. At this point, we will be left with one open and one closed `cluster', both of which must have singleton equivalence classes, and it is no longer necessary to make any connections. Note that since the equivalent class of $O_1$ is a singleton, $O_1$ must be disconnected from the open clusters $O_2$ and $O_{m-1}$ in $\omega_0(\check{A})$. This means that $C_1$ and $C_2$ must be connected in $\omega_0(\check{A})$, or in other words that $C_1$ and $C_2$ are equivalent. The same holds for all subsequent steps in the algorithm. So at each step, we construct corridors between clusters that are equivalent. Thus, if we manage to make closed connections along all closed corridors and open connections along all open corridors, we will precisely connect open/closed clusters if and only if they are equivalent. 

Note that because points $x_i$ and $y_i$ are at distance at least $\check{a}_1 \cdot 2r$ apart on the event $\check{F}$ and the different levels are suitably spaced apart, none of the resulting corridors intersect each other. In fact, they satisfy that for any two $\check{c}_1$-boxes $B$ and $B'$ in different corridors we have $2B\cap 2B'=\emptyset$. Recall that on the event $\check{F}$ we also have that all the $\check{b}_1$-boxes used in the construction of $\mathcal{D}_{\textrm{cor}}$ are suitably shielded away from all crossing clusters but the two whose interface enters the box by open/closed paths. Consider the collection of closed corridors, and consider the finite collection of closed clusters that almost form a loop in each $2B\setminus B$ for boxes $B$ in these corridors, whose existence is guaranteed by the definition of event $\check{F}$. If we were to resample the connections between these closed clusters so that we connect them, and we do the same for the open clusters, the resulting configuration satisfies $\mathcal{E}_{\textrm{out}}(S,x_S)$. 

\medskip \textbf{Step 3:} \textit{Boxes near the sides.}

Assume now that the box $B$ intersects one of the boxes of radii $2\delta n^{-1} R^{\kappa-k}$ centered at points $c_1^{(k)},c_2^{(k)},c_3^{(k)}$. Without loss of generality, assume that it intersects $B(c_1^{(k)},2\delta n^{-1} R^{\kappa-k})$. If $B$ is fully contained inside of $B(c_1^{(k)},2\delta n^{-1} R^{\kappa-k})$, 

the same construction of corridors as above works. Namely, the bottom closed cluster and the open crossing cluster of $\mathcal{E}_{\textrm{out}}(S,x_S)$ must already be realized in $\tilde{\omega}$, thus it is still true that if we change the configuration $\omega_0(\textrm{int }\check{A})$ to any other configuration that connects the equivalent crossing clusters of $\check{A}$, the event $\mathcal{E}_{\textrm{out}}(S,x_S)$ still holds for the edited percolation. 

If, on the other hand, $B$ is not contained inside of $B(c_1^{(k)},2\delta n^{-1} R^{\kappa-k})$, it is no longer sufficient to connect equivalent crossing clusters. We must additionally assume that our configuration on $\textrm{int }\check{A}$ is such that no open sites of $B\cap B(c_1^{(k)},2\delta n^{-1}R^{\kappa-k})$ are joined to the open crossing clusters that are a part of the left-to-right open crossing in $\tilde{\omega}\cup\omega_0(\check{A})\cup \omega_0(\textrm{int }\check{A})$. So, it suffices to modify the way we build our set of corridors, so that the closed corridors shield the relevant open clusters away from the interior of $B(c_1^{(k)},2\delta n^{-1} R^{\kappa-k})$. We will again construct these corridors within the set $\mathcal{D}_{\textrm{cor}}$ that is defined analogously as in the previous case. We will find it useful that $B(c_1^{(k)},2\delta n^{-1} R^{\kappa-k})$ splits the loop part of $\mathcal{D}_{\textrm{cor}}$ into two `arcs' (one of which could be empty). This means that the loop part of $\mathcal{D}_{\textrm{cor}}$ that is not in $B(c_1^{(k)},2\delta n^{-1} R^{\kappa-k})$ has enough `depth' to allow for the construction of disjoint corridors -- see Figure~\ref{fig:construction_of_corridors_when_we_are_near_the_midpoint} below. This is why we impose the technical condition that the box $B$ is nice in $B^{(k)}$.

There are possibly several equivalence classes of open crossing clusters of $\check{A}$ that are a part of the cluster corresponding to the open crossing of $\mathcal{E}_{\textrm{out}(S,x_S)}$ on $\tilde{\omega}\cup \omega_0(\check{A})\cup \omega_0(\textrm{int }\check{A})$. We will call them relevant. Without loss of generality, assume that $0\leq m_0\leq m$ is such that $O_1,\ldots,O_{m_0}$ do not intersect $B(c_1^{(k)},2\delta n^{-1} R^{\kappa-k})$ and $O_{m_0+1},\ldots,O_m$ do.  All the relevant clusters belong to $\{O_1,\ldots,O_{m_0}\}$. 

If there are no relevant clusters, the same proof as above works. Otherwise, assume that $1\leq l_1\leq k_1\leq m_0$ are such that none of the clusters in $\{O_1,\ldots,O_{l_1-1}\}$ are relevant, but $O_{l_1},O_{k_1}$ are, and moreover, the entire equivalence class of these clusters belongs to $\{O_{l_1},\ldots,O_{k_1}\}$. If there are no relevant open clusters with index larger than $k_1$, we finish our procedure. Otherwise, define $k_1<l_2\leq k_2\leq m_0$ to be such that no clusters in $\{O_{k_1+1},\ldots,O_{l_2-1}\}$ are relevant, but $O_{l_2},O_{k_2}$ are, and their entire equivalence class belongs to $\{O_{l_2},\ldots,O_{k_2}\}$. Continue until there are no relevant clusters left.

Construct a corridor connecting $x_{l_1}$ and $y_{k_1}$ at level $\check{M}_1$ that moves through the loop part of $\mathcal{D}_{\textrm{cor}}$ clockwise from $x_{l_1}$ to $y_{k_1}$ and declare it to be closed. Do the same for all pairs $(x_{l_i},y_{k_i})$. This subdivides $\textrm{int }\check{A}$ into several faces -- 1 face containing the center of $B$ and several other ones that are bounded by the corridor between $x_{l_i}$ and $y_{k_i}$ and clusters $\{O_{l_i}, C_{l_i+1},\ldots,C_{k_i},O_{k_i}\}$. For each of these faces we can construct the relevant corridors inside the face using the same algorithm as in the first case, except that for the face containing the center of $B$, we will start constructing our corridors at level $\check{M}_1+1$ instead of level $1$. Note that since we originally constructed $\check{M}_1^2$ levels, we do not run out levels to make corridors for our algorithm.

\begin{figure}[ht]
    \centering
    \includegraphics[scale=0.8]{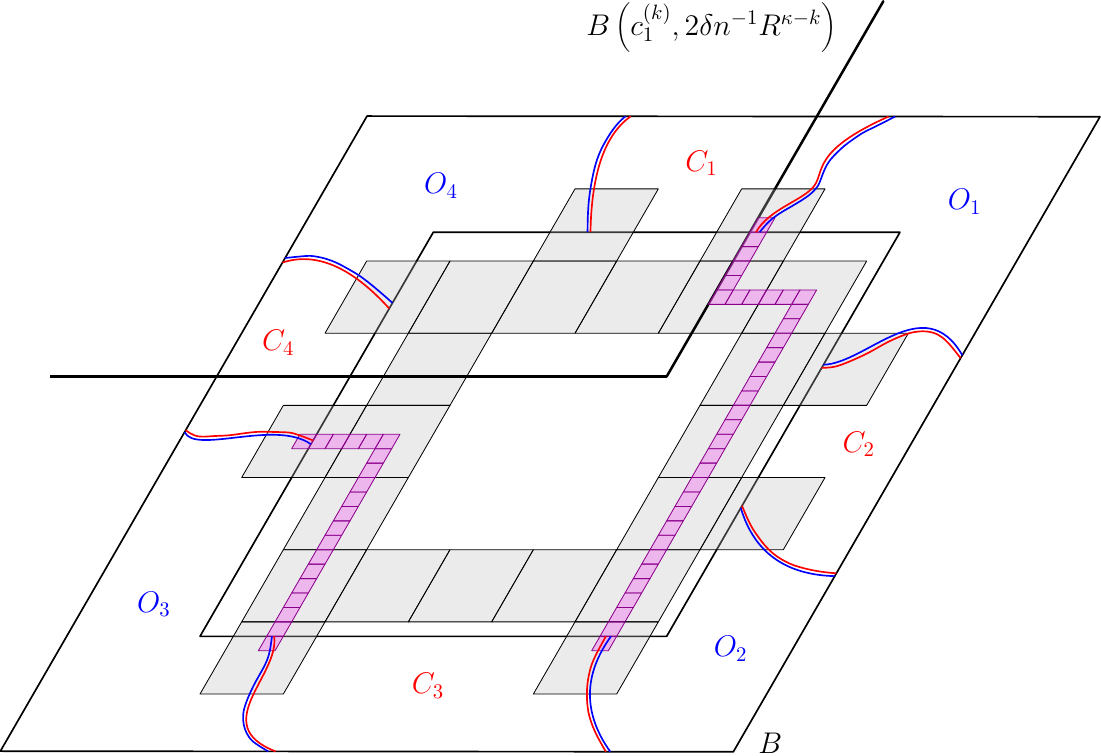}
    \caption{The figure above illustrates the case when a box $B$ that is nice in $B^{(k)}$ intersects $B\!\left(c_1^{(k)},2\delta n^{-1}R^{\kappa-k}\right)$. In this case, there are $8$ crossing clusters in total, and $C_4,O_4,C_1$ are the ones that intersect $B\!\left(c_1^{(k)},2\delta n^{-1}R^{\kappa-k}\right)$. Thus $m_0=3$. The area shaded in gray again represents the domain $\mathcal{D}_{\textrm{cor}}$. The two smaller purple corridors represent the corridors that are constructed in the first stage of our algorithm. Namely, the picture above illustrates the case that $O_1,O_2$ and $O_3$ are all relevant, and that $O_1$ is equivalent to $O_2$, but not to $O_3$. These corridors are drawn at level $\check{M}_1$, and they leave sufficiently many levels on either side to build corridors in the subsequent phase of our algorithm. In the above case, the corridors have subdivided our circuit into $3$ faces - the face containing $O_1,C_2$ and $O_2$,  the face containing $O_3$, and the face containing the center of $\check{A}$ delimited by $(C_1,C_3,C_4)$ and $O_4$. In the second phase of the algorithm we would only need to make a connection between $O_1$ and $O_2$ within the first face.
    Finally, note that by closing all of the connections between closed clusters corresponding to boxes of these two purple corridors, we would describe closed loops around the two corridors. Even though these closed paths might intersect $B\!\left(c_1^{(k)},2\delta n^{-1}R^{\kappa-k}\right)$ (as they would on the picture above for the right corridor), they will still shield off the relevant open clusters from the box $B\!\left(c_1^{(k)},2\delta n^{-1}R^{\kappa-k}\right)$ (in the case above, the cluster $O_1$ is fully separated from this box by these closed paths).}
    \label{fig:construction_of_corridors_when_we_are_near_the_midpoint}
\end{figure}

Similarly as before, if we were to resample all the connections in the corridors we would get all the necessary connections between equivalent corridors. Additionally, we get that all the relevant open clusters are `shielded' away from $B(c_1^{(k)},2\delta n^{-1}R^{\kappa-k})$ by the closed corridors. Indeed, since all the corridors and $B^{(k)}$ have sides parallel to the axes, even if some boxes of a closed corridor intersect $B(c_1^{(k)},2\delta n^{-1}R^{\kappa-k})$ the boundary of the corridor that connects $O_{l_i}$ and $O_{k_i}$ stays outside of $B(c_1^{(k)},2\delta n^{-1}R^{\kappa-k})$. See Figure~\ref{fig:construction_of_corridors_when_we_are_near_the_midpoint} for an illustration. Thus, ${\mathcal{E}_{\textrm{out}}(S,x_S)}$ holds for the resampled configuration, as wanted.

\medskip \textbf{Step 4:} \textit{Resampling.}

The resampling procedure is defined analogously as in Section~\ref{subsec:distance-crossings} (see page~\pageref{page:resampling}). Lemma~\ref{lem:pivotals_are_far}\eqref{property:distance-points} still applies to ensure that the points around which we need to resample are suitably far away. Thus, if we again denote by $E_{\textrm{res}}$ the event that the resampling is successful, by a result analogous to Lemma~\ref{lem:resampling_is_positive} we know that there exists a constant $\check{g}=\check{g}(\check{p},R,\check{M}_1,\check{M}_2,\check{c}_1,\check{c}_2,\check{\xi})>0$ such that 

$$\P (E_\textrm{res}|\omega(A(x;r,4r)),\tilde{\omega})\cdot \1_{\check{F}\cap\mathscr{P}(\tilde{\omega})}(\omega(A(x;r,4r)), \tilde{\omega})\geq \check{g}\cdot \1_{\check{F}\cap\mathscr{P}(\tilde{\omega})}(\omega(A(x;r,4r)), \tilde{\omega})$$ whenever $\check{\xi} \frac{\check{c}_1\check{c_2}}{8} \cdot 2rn\geq R$. In the above we wrote $\check{F}$ instead of $\check{F}(A(x;r,4r))$ for brevity. In particular, this means that 
\begin{equation}
    \begin{split}
        \P(\mathcal{E}_{\textrm{out}}(S,x_S) | \tilde{\omega}) &\geq \P(\check{F}, \mathscr{P}(\tilde{\omega}), E_{\textrm{res}}| \tilde{\omega}) \\
        &= \E \!\left[ \1_{\check{F}\cap\mathscr{P}(\tilde{\omega})} \cdot \E [\1_{E_{\res}}|\omega(A(x;r,4r)),\tilde{\omega}] \bigg| \tilde{\omega}\right] \\
        &\geq \check{g}\cdot \E \!\left[ \1_{\check{F}\cap\mathscr{P}(\tilde{\omega})} \bigg| \tilde{\omega}\right],
    \end{split}
\end{equation}
as wanted.
\end{proof}
\begin{rem} Note that the above proof can also be carried out without using resampling. Namely, if we explore all open/closed interfaces starting from the outer boundary of $\check{A}$ until they either hit the inner boundary, or they come back to $\textrm{ext }{\check{A}}$, we will end up exploring all the interfaces between open and closed crossing clusters of $\check{A}$. Furthermore, since we are working on the event $\check{F}_1(\check{A})\cap \check{F}_2(\check{A})\cap \check{F}_3(\check{A})$, we know that around each explored site $v$ on the inner boundary of $\check{A}$, there is a box of radius $\check{b}_1\cdot 2r$ such that any explored site belonging to this box must in fact belong to the open/closed interface ending at $v$. If we construct our corridors exactly as above, one can see that forming open (resp.\ closed) loops in the unexplored part of each $2B'\setminus B'$ for boxes $B'$ that are parts of the open (resp.\ closed) corridors. Using RSW estimates and positive correlation, one can reobtain the statement of the lemma. 
\end{rem}

\subsubsection{Concluding superpolynomial concentration}
\label{se:geo_exploration_conclusion}

We are now ready to conclude the proof of $\P(\tau_2<\infty,\cE(B^{(k)}))\lesssim\P(\tau_1<\infty,\cE(B^{(k)}))\cdot R^{-\beta k}$. Recall the definition of~$\tau_\ell$ in~\eqref{eq:tau_definition2}.
\begin{lem}
There exists $\beta>0$ such that the following holds for all $R$ sufficiently large and $p$ sufficiently small depending on~$R$. For all $\alpha$ for which~\eqref{eq:assumptionsubpoly} holds, there exists $C=C(\alpha,R)$ such that for all $k$ and $n\geq R^k$,
\[\mathbf{P}\!\left(\tau_2 <\infty, \mathcal{E}\!\left(B^{(k)}\right)\right)\leq C R^{-(\alpha+\beta)k}.\] \label{lem:bound_on_stopping_time}
\end{lem}

The proof uses a stringing argument that is similar to that of Lemma~\ref{lem:crossingbubble}. Recall that in the proof of Lemma~\ref{lem:crossingbubble}, we had two bubble interfaces and with sufficiently high probability, we could cover at least one of the interfaces with `good' annuli and string together short paths inside these annuli using Proposition~\ref{lem:distancebetweencrossings} to bound the overall distance between the bubble pinch points. In the following proof, we will follow the `good' annuli along the path $\sigma_t$.

We have already shown that on the event $H$, we can bound the number of boxes needed to cover~$\sigma_t$ (recall~\eqref{eq:bound-Wj}). In the first step of this lemma's proof, we define a `good' event $E_1$. In the second step, we give the stringing argument to show that if $H$ and the events $E_1$ hold for each covering Whitney box, then $\tau_2=\infty$. In Step~3, we finally bound the probability $\P(\tau_2<\infty,\cE(B^{(k)}),H)$. This follows the strategy outlined at the beginning of Section~\ref{se:nonlocality}. We briefly sketch how we resolve the problem of non-locality mentioned there. Recall that on page~\pageref{page:Ein_Eout} we have defined for each choice of $(S,x_S)$ the events $\mathcal{E}_{\textrm{in}}(S,x_S), \mathcal{E}_{\textrm{out}}(S,x_S)$ such that $\mathcal{E}(B^{(k)})\cap \{\tau_1<\infty\}\cap \{(\Gamma_1,x)=(S,x_S)\} = \mathcal{E}_{\textrm{in}}(S,x_S) \cap \mathcal{E}_{\textrm{out}}(S,x_S)$. In particular, the conditional law of the percolation given $\mathcal{E}(B^{(k)})\cap \{\tau_1<\infty\}\cap \{(\Gamma_1,x)=(S,x_S)\}$ is that of a percolation in~$S^c$ \emph{conditioned} on $\mathcal{E}_{\textrm{out}}(S,x_S)$. For a fixed box $B$ for which $4B \subseteq S^c$, we denoted by $\tilde{\omega}$ the percolation configuration outside of $4B$ and we have introduced the event $\mathscr{P}(\tilde\omega)$ above Lemma~\ref{lem:E_out_happens_with_positive_probability} as a proxy for~$\mathcal{E}_{\textrm{out}}(S,x_S)$, with the advantage that (unlike $\mathcal{E}_{\textrm{out}}(S,x_S)$) the event $\mathscr{P}(\tilde\omega)$ depends only on the sites outside $B$. Lemma~\ref{lem:E_out_happens_with_positive_probability} now implies that it suffices to bound a variant of~\eqref{eq:bad_box_cond_on_exploration} with $\mathcal{E}_{\textrm{out}}(S,x_S)$ replaced by $\mathscr{P}(\tilde\omega)$. Indeed, we have
\begin{equation}\begin{split}\label{eq:geo_exploration_decorrelated1}
&\p\!\left(\cE(B^{(k)}),\,\tau_1<\infty,\, (\Gamma_1,x)=(S,x_S) ,\, B \in \mathcal{W}_j' \right) \\
&= \p\!\left(\mathcal{E}_{\textrm{in}}(S,x_S) ,\, B \in \mathcal{W}_j' \right) \, \p\!\left( \mathcal{E}_{\textrm{out}}(S,x_S) \mmiddle| \mathcal{E}_{\textrm{in}}(S,x_S) ,\, B \in \mathcal{W}_j' \right) \\
&\gtrsim \p\!\left(\mathcal{E}_{\textrm{in}}(S,x_S) ,\, B \in \mathcal{W}_j' \right) \, \p\!\left( \check{F}(4B \setminus B),\, \mathscr{P}(\tilde\omega) \mmiddle| \mathcal{E}_{\textrm{in}}(S,x_S) ,\, B \in \mathcal{W}_j' \right) .
\end{split}\end{equation}
Further,
\begin{equation}\begin{split}\label{eq:geo_exploration_decorrelated2}
&\p\!\left(\cE(B^{(k)}),\,\tau_1<\infty,\, (\Gamma_1,x)=(S,x_S) ,\, B \in \mathcal{W}_j' ,\, \check{F}(4B \setminus B)  ,\, F(B) \cap G(B)^c \right) \\
&\le \p\!\left(\mathcal{E}_{\textrm{in}}(S,x_S) ,\, \mathscr{P}(\tilde\omega) ,\, B \in \mathcal{W}_j' ,\, \check{F}(4B \setminus B) ,\, F(B) \cap G(B)^c \right) \\
&= \p\!\left(\mathcal{E}_{\textrm{in}}(S,x_S) ,\, \mathscr{P}(\tilde\omega) ,\, B \in \mathcal{W}_j' ,\, \check{F}(4B \setminus B) \right) \\
&\qquad \cdot \p\!\left( F(B) \cap G(B)^c \mmiddle| \mathcal{E}_{\textrm{in}}(S,x_S) ,\, \mathscr{P}(\tilde\omega) ,\, B \in \mathcal{W}_j' ,\, \check{F}(4B \setminus B) \right) .
\end{split}\end{equation}
Now Proposition~\ref{lem:distancebetweencrossings} can be applied to the last term on the right-hand side of~\eqref{eq:geo_exploration_decorrelated2} since this conditioning depends only on the sites outside $B$. The first term on the right-hand side of~\eqref{eq:geo_exploration_decorrelated2} can be estimated by~\eqref{eq:geo_exploration_decorrelated1}, yielding an estimate similar to~\eqref{eq:geo_exploration_sketch}. This resolves the problem described at the beginning of Section~\ref{se:nonlocality}.

\begin{proof}[Proof of Lemma~\ref{lem:bound_on_stopping_time}] Note that $\mathbf{P}\!\left(\tau_2<\infty,\mathcal{E}\!\left(B^{(k)}\right)\right) \leq \mathbf{P}({H^c}) +\mathbf{P}\!\left({\tau_2<\infty,H},\mathcal{E}\!\left(B^{(k)}\right)\right)$. Since $\mathbf{P}({H^c}) \lesssim R^{-4\alpha k}$ by Lemma~\ref{lem:bound_on_global_event}, 
it suffices to bound $\mathbf{P}\!\left({\tau_2<\infty,H}, \mathcal{E}\!\left(B^{(k)}\right)\right)$. To do so, we will exploit the fact that on the event $H\cap \mathcal{E}\!\left(B^{(k)}\right)$ we can cover an open path from $x$ to the last bubble pinch point by a small number of Whitney boxes. Moreover, Lemma~\ref{lem:bound_on_F}, Proposition~\ref{lem:distancebetweencrossings} and Lemma~\ref{lem:E_out_happens_with_positive_probability} will guarantee that no individual box has long crossings (with respect to the geodesic or effective resistance metric) with high probability. In turn, we will be able to use the triangle inequality to bound the distance across the remaining bubbles by distance in each individual box, showing that on this event $\tau_2=\infty$.

\medskip \textbf{Step 1:} \textit{Path covering.}

Let $\nu>0$ be as in the definition of the event $H$. We will use a strategy similar to that in the proof of Lemma~\ref{lem:crossingbubble}. Consider a box $B\in\mathcal{W}_{j}$. Cover the box $B$ with at most $4R^{2\nu j}$ smaller boxes of side length $n^{-1}R^{\kappa-j(1+\nu)}$. Let $\mathcal{W}_{j}^\nu$ be the set of all such boxes of size $n^{-1}R^{\kappa-j(1+\nu)}$. Then again $\cup_{j\geq k} \mathcal{W}_{j}^\nu$ covers the path $\sigma_t$. Suppose $B_\nu \in \mathcal{W}_{j}^\nu$ is one such box, and let $z$ denote its center. Consider the disjoint concentric annuli $A(z;n^{-1}R^{\kappa-2f},n^{-1}R^{\kappa-2f+2})$ for $j+4\leq 2f< j(1+\nu)$. In contrast to the proof of Lemma~\ref{lem:crossingbubble}, the annuli now have radius ratio $R^2$ instead of $R$. This is because we will split these into two concentric annuli of radius ratio $R$ each: the inner one will serve the same purpose as before (we will use it to ensure that lengths across this annulus are suitably small), while the outer one will help us decorrelate $\mathcal{E}(B^{(k)})$ and the event that the lengths across the inner annulus are small. 
All of the levels $j+4\leq 2f<j(1+\nu)$ but at most 9 are such that the box $B(z,4n^{-1}R^{\kappa-2f+1})$ is nice in $B^{(k)}$ (recall this definition from page~\pageref{page:definition_of_a_nice_box}). We will only consider the `nice' levels $2f$.\footnote{The same comment applies as in Lemma~\ref{lem:crossingbubble} -- this proof only works when $j\nu > 14$, i.e.\ when there is at least one level at which we can test for the existence of a good annulus. Since we will be able to choose $\nu$ deterministically, we can simply subdivide the boxes of such levels $j$ into smaller boxes of level $\hat{j}=\lceil 14/\nu\rceil$ for which the proof works. We just get an additional constant factor of $R^{2\hat{j}}$ for the number of boxes.} Because $B_{\nu}$ was chosen as a box in the covering of a Whitney box, it must be that none of these annuli intersect the explored domain. They also stay within the box $B^{(k)}$. 

Continuing in parallel with the proof of Lemma~\ref{lem:crossingbubble}, fix some $\check{p}>0$ to be determined later. Recall the definition of the event $\check{F}$ and the discussion afterwards on page~\pageref{page:bound_on_check_F}. Choose the constants $\check{M}_1,\check{a}_1,\ldots,\check{\xi}$ depending on $\check{p}$ so that $$\mathbf{P}\!\left(\check{F}\!\left(A(z;2n^{-1}R^{\kappa-2f+1},4n^{-1}R^{\kappa-2f+1})\right)\right)\geq 1-\check{p}$$ for all levels $f$ for which $\check{\xi} \frac{\check{c}_1\check{c}_2}{8} R^{\kappa-2f+1}\geq R$. Denote this event by $\check{F}(B_\nu,2f-1)$. 

Similarly, fix some $0<p'$ to be determined later, and choose all the constants $M_1,a_1,\ldots \xi$ in Lemma~\ref{lem:bound_on_F} to be such that $$\mathbf{P}\!\left(F(A(z;n^{-1}R^{\kappa-2f},n^{-1}R^{\kappa-2f+1}))\right)\geq 1-p'$$ for all levels $2f$ for which $\xi\frac{c_1c_2}{8} R^{\kappa-2f-1}\geq R$. Fix also $\zeta=2/3$ and $s_j= \frac{(16C_1R^3)^{-1}}{(M_1+1)^2} R^{-4\nu j}$ in the definition of the event $G$ in Proposition~\ref{lem:distancebetweencrossings}, where $C_1=C_1(R,\nu)$ is as in the equation~\eqref{eq:bound-Wj}. We again drop all the parameters from the notation for the sake of brevity and denote these events by $F(B_\nu,2f)$ and $G(B_\nu,2f)$.

Define $E_1(B_\nu,s_j)$ to be the event that there exists a nice level $j+4\leq 2f< j(1+\nu)$ such that $F(B_\nu,2f)\cap G(B_\nu,2f)$ holds for the annulus $A(z;n^{-1}R^{\kappa-2f},n^{-1}R^{\kappa-2f+1})$ and $\check{F}(B_\nu,2f-1)$ holds for $A(z;2n^{-1}R^{\kappa-2f+1},4n^{-1}R^{\kappa-2f+1})$. We will ultimately use Lemma~\ref{lem:bound_on_F}, Proposition~\ref{lem:distancebetweencrossings} and Lemma~\ref{lem:E_out_happens_with_positive_probability} to bound the probabilities of such events $E_1(B_\nu,s_j)$. Recall that the bounds in these lemmas hold only when the annulus radii are large enough. In particular, these bounds hold for all annuli of level $j$ such that $\xi\frac{c_1c_2}{8}R^{\kappa-j(1+\nu)-1}\geq R$ and $\check{\xi}\frac{\check{c}_1\check{c}_2}{8}R^{\kappa-j(1+\nu)-1}\geq R$. In other words, they hold for all levels $j\leq J$, where $J=\frac{\kappa+\log_R\!\left(\min \!\left\{\xi\frac{c_1c_2}{8},\check{\xi}\frac{\check{c}_1\check{c}_2}{8}\right\}\right)-2}{1+\nu}$.

\medskip \textbf{Step 2:} \textit{Stringing.}

We claim that on the event $\mathcal{E}\!\left(B^{(k)}\right)\cap H \cap \bigcap_{k-1\leq j\leq J}\bigcap_{B_\nu\in \mathcal{W}_{j}^\nu}E_1(B_\nu,s_{j})$, the distance from $x$ to $v$ is at most $\altq_n$. To see this, we will use a stringing argument that is almost analogous to the one we have used in the proof of Lemma~\ref{lem:crossingbubble}. Recall that the open path $\sigma_t$ connects $x$ to $v$. It is covered by $\bigcup_{k-1\leq j\leq J} \mathcal{W}_{j}^{\nu} \cup \bigcup_{j\geq J} \mathcal{W}_{j}$.

To start off, we need to make sure that we can deal with boxes in $\mathcal{W}_{j}$ for $j\geq J$, i.e.\ the ones for which the discussion above does not apply. Namely, for such $B\in \mathcal{W}_{j}$, note that we know deterministically that the distance between any two sites in $2B$ that are connected by an open path inside of this box 
is smaller than the total number of sites in $2B$ (regardless of whether $2B$ intersects the explored domain). In particular, their distance is bounded from above by $$\!\left(2R^{\kappa-j}\right)^2\leq 4 R^{2\kappa}\!\left(\frac{R^{2-\kappa}}{\min \!\left\{\xi\frac{c_1c_2}{8},\xi''\frac{c_1''c_2''}{8}\right\}}\right)^{\frac{2}{1+\nu}}\leq C(\nu,R,\alpha) \cdot n^{\frac{2\nu}{1+\nu}}.$$ 
At the same time, this also provides an upper bound on the number of times our path $\sigma_t$ enters the box $2B$ -- if it enters through the same site twice we can just discard the portion of $\sigma_t$ in between different visits to this site.

We can now proceed with our linking argument. The argument follows in almost exactly the same way as the linking argument of Lemma~\ref{lem:crossingbubble}. Consider all open paths from $x$ to $v$ that can be covered with boxes in $\bigcup_{k-1\leq j\leq J} \mathcal{W}_{j}^{\nu} \cup \bigcup_{j\geq J} \mathcal{W}_{j}$, and let $B_1,\ldots,B_a$ be a minimal sequence of consecutive boxes needed to cover an open path from $x$ to $v$. Our sequence will consist of several strings of boxes alternating between a string of boxes in $\bigcup_{j\geq J} \mathcal{W}_{j}$ and a string of boxes in $\bigcup_{k-1\leq j\leq J} \mathcal{W}_{j}^{\nu}$. For strings of the second type, we can bound the distance between the entry and exit points of our path into these strings in exactly the same way as we did in Lemma~\ref{lem:crossingbubble} -- every box is entered at most $M_1$ times, and each time it is entered it can contribute to the overall length by at most $2(M_1+1)s_{j}\cdot \altq_n$ if the box belongs to $\mathcal{W}_j^{\nu}$. The segments of our open path in between these `good' strings are covered by boxes in $\bigcup_{j\geq J} \mathcal{W}_{j}$. Thus, we can bound the distance from the exit point of one `good' string to the entry point of the following `good' string by the deterministic bound in terms of the number of sites in these boxes that we established in the previous paragraph. Hence, on the event $\mathcal{E}\!\left(B^{(k)}\right) \cap H \cap \bigcap_{k-1\leq j\leq J}\bigcap_{B_\nu\in \mathcal{W}_{j}^\nu}E_1(B_\nu,s_{j})$, using the triangle inequality, we can bound $D(x,v)$ from above by
\begin{equation}
\begin{split}\label{eq:upper_bound_on_distance_x_v}
D(x,v) &\leq \sum_{j=k-1}^{J} |\mathcal{W}_{j}^{\nu}|\cdot M_1\cdot 2(M_1+1)\cdot \frac{(16C_1R^3)^{-1}}{(M_1+1)^2} R^{-4\nu j}\cdot \altq_n + \sum_{j=J}^{\kappa} |\mathcal{W}_{j}| \cdot C(\nu, R,\alpha)^2 n^{\frac{4\nu} {1+\nu}} \\
&\leq \sum_{j=k-1}^{J} C_1 R^{\frac{3}{2}\nu j} \cdot 4R^{2\nu j}\cdot (8C_1R^3)^{-1}R^{-4\nu j} \altq_n + \sum_{j=J}^{\kappa} C_1 R^{\frac{3}{2}\nu j} \cdot C(\nu, R,\alpha)^2 n^{\frac{4\nu} {1+\nu}} \\
&\leq \frac{R^{-3}}{2} \cdot \!\left(\sum_{j\geq k-1}R^{-\frac{1}{2}\nu j}\right)\cdot \altq_n+ C(\nu, R,\alpha) \cdot \log_R n\cdot n^{\frac{3}{2}\nu+\frac{4\nu}{1+\nu}},
\end{split}
\end{equation}

where $C(\nu,R,\alpha)$ is a constant that depends only on $\nu,R,\alpha$ and can increase from line to line. The constant $C_1=C_1(R,\nu)$ above comes from equation~\eqref{eq:bound-Wj} that we have used to bound the numbers of boxes in $\cW_j,\cW^\eta_j$. 

Fix $\nu=1/100$. The series in the first term is summable and it sums up to $\frac{R^{-\nu(k-1)/2}}{1-R^{-\nu/2}}$, which is smaller than $R^3$ for $R$ large enough (in terms of $\nu=1/100$). Moreover, $\log_R n\cdot n^{\frac{3}{2}\nu+\frac{4\nu}{1+\nu}} \lesssim n^{6\nu}$. Corollary~\ref{cor:lower_bound_on_m_n} tells us that there exists a constant $\tilde{C}(R)$ such that $\altq_n\geq \tilde{C}(R)\cdot n^{1/2}$. In particular, this means that the second term in the final inequality above can be bounded from above by $C'(\nu,R,\alpha)\cdot n^{6\nu-1/2} \cdot \altq_n$. Combining all these observations and substituting them in the equation above we find that $$D(x,v)\leq  \frac{1}{2}\altq_n + C'(\nu,R,\alpha) n^{6\nu-1/2} \altq_n \leq \frac{1}{2}\altq_n+\frac{1}{2}\altq_n=\altq_n ,$$ for all $n$ large enough (in terms of $R,\alpha$ and $\nu=1/100$). Thus, in this case $\tau_2=\infty$. From now on, we assume that $n$ is large enough so that this bound holds. It suffices to prove the lemma for such $n$, since by enlarging the constant $C=C(\alpha,R)$ in the statement, the claim also follows for finitely many small values of $n$.

\medskip \textbf{Step 3:} \textit{Concluding the probability bound.}

We will now use this to bound $\mathbf{P}\!\left({\tau_2<\infty,H},\mathcal{E}\!\left(B^{(k)}\right)\right)$ from above. In what follows, we will write $\mathcal{E}$ in place of $\mathcal{E}\!\left(B^{(k)}\right)$ to simplify the notation. Also, whenever we are taking the union or sum over levels $j+4\leq 2f<j(1+\nu)$, it is understood that we are taking the union/sum only over the nice levels (so the index set will depend on the choice of the box $B_\nu$). We omit this to alleviate the notation. Recall that $\Gamma_1$ denote the sites in the `explored region' before $\tau_1$, which exists and is a subset of $B^{(k)}$ on the event $\{\tau_1<\infty\}\cap\cE$. Then, since the event $\mathcal{E}\cap H\cap \bigcap_{k-1\leq j \leq J} \bigcap_{B_{\nu}\in \mathcal{W}_{j}^{\nu}} E_1(B_\nu,s_{j})$ implies $\mathcal{E}\cap H \cap \{\tau_2=\infty\}$, we have
\begin{align}
&\mathbf{P}\!\left(\tau_2<\infty, \mathcal{E}, H\right) \leq \P \!\left(\mathcal{E},\tau_1<\infty, H, \bigcup_{j= k-1}^{J}\bigcup_{B_\nu\in \mathcal{W}_{j}^\nu}E_1(B_\nu,s_{j})^c \right) \notag\\
& \leq \P \!\left(\bigcup_{j= k-1}^{J}\bigcup_{B\in \cB_{j,n}^{R,k}}\bigcup_{B_\nu\textrm{ covers }B} \!\left(\bigcap_{j+4\leq 2f<j(1+\nu)} \check{F}(B_\nu,2f-1)^c \cup F(B_\nu,2f)^c\right) \right) + \notag\\
&\quad \P \!\left(\mathcal{E},\tau_1<\infty, \bigcup_{j= k-1}^{J}\bigcup_{B_\nu\in \mathcal{W}_{j}^\nu} \!\left(\bigcup_{j+4\leq 2f<j(1+\nu)} \check{F}(B_\nu,2f-1) \cap F(B_\nu,2f)\cap G(B_\nu,2f)^{c}\right) \right) \notag\\
& \leq \sum_{j=k-1}^J R^{2(j+1)} \cdot R^{2\nu +2} \cdot (p''+p')^{j\nu/2-14} + \sum_{(S,x_S)}\P \!\left[\mathcal{E},\tau_1<\infty, \Gamma_1=(S,x_S), \phantom{\!\left(\bigcup_{B_\nu}^J\right)} \right. \notag\\
&\quad \hspace{2cm} \left.  \bigcup_{j= k-1}^{J}\bigcup_{B_\nu\in \mathcal{W}_{j}^\nu} \!\left(\bigcup_{j+4\leq 2f<j(1+\nu)}\check{F}(B_\nu,2f-1) \cap F(B_\nu,2f)\cap G(B_\nu,2f)^{c}\right) \right],
\label{eq:first_bound_on_tau_1<infty}
\end{align}
where we are summing over all appropriate pairs $(S,x_S)$, where $S\subseteq B^{(k)}$ is in the support of $\Gamma_1$ and $x_S$ is a designated point on the boundary of $S$ that serves as the tip of our `explored domain' $\Gamma_1$ called $x$. In the last inequality we have used that events $\check{F}(B_\nu,2f-1)\cap F(B_\nu,2f)$ are independent for different levels $2f$, since they depend on disjoint sets of sites. The exponent $j\nu/2-14$ is the lower bound on the number of nice levels $j+4\leq 2f<j(1+\nu)$ that holds no matter where $B_\nu$ is. Choose $\check{p}=\check{p}(R,\alpha,\nu)$ and $p'=p'(R,\alpha,\nu)$ to be such that $\check{p},p'\leq \frac{1}{2}R^{-2(3+\alpha)/\nu}$. With these constraints, the first term in the expression above can be bounded from above by $C(R,\alpha,\nu)\cdot R^{-(1+\alpha)k}$.

Denote the second term by $\mathit{II}$. 
To bound $\mathit{II}$, recall that we can write the event $\mathcal{E}\cap \{\tau_1<\infty\}\cap \{\Gamma_1=(S,x_S)\}$ as $\mathcal{E}_{\textrm{in}}(S,x_S) \cap \mathcal{E}_{\textrm{out}}(S,x_S)$ where the first of the two events depends only on the sites of $S$ whereas the other one depends only on the sites of $S^c$. Recall also the definition of $\mathcal{W}_j'$ from page~\pageref{defn:definition_of_W'} -- it consists of all boxes in $\cB^{R,k}_{j,n}$ that are either close to the tip of the exploration, $x$ ($x=x_S$ on the event $\Gamma_1=(S,x_S)$) or that are contained in some annulus `around $x$' and have $6$ semiarms of appropriate length or that have 3 semiarms near the boundary. Define $\!\left(\mathcal{W}_j'\right)^\nu$ analogously to $\mathcal{W}_j^\nu$ with $\mathcal{W}_j'$ in place of $\mathcal{W}_j$. Then, because $\mathcal{W}_j\subseteq \mathcal{W}_j'$, we can take the union over $\!\left(\mathcal{W}_j'\right)^\nu$ in the equation~\eqref{eq:first_bound_on_tau_1<infty} to obtain a further upper bound. Substituting all this into equation~\eqref{eq:first_bound_on_tau_1<infty} we find that we can bound $\mathit{II}$ from above by
\begin{equation}
    \begin{split}
        \mathit{II}\leq&\, \sum_{(S,x_S)}\P \Biggl(\cE_{\textrm{in}}(S,x_S),\cE_{\textrm{out}}(S,x_S),\\
        &\left.\hspace{1.5cm}\bigcup_{j= k-1}^{J} \bigcup_{B_\nu\in \!\left(\mathcal{W}_{j}'\right)^\nu} \!\left(\bigcup_{j+4\leq 2f<j(1+\nu)}\check{F}(B_\nu,2f-1) \cap F(B_\nu,2f)\cap G(B_\nu,2f)^{c}\right) \right) \\
     \leq&\, \sum_{j=k-1}^{J} \sum_{(S,x_S)} \sum_{B\in \cB_{j,n}^{R,k}}\sum_{B_\nu \textrm{ covers }B}
     \sum_{j+4\leq 2f<j(1+\nu)}\P \bigg(\cE_{\textrm{in}}(S,x_S),\cE_{\textrm{out}}(S,x_S), B_\nu\in \!\left(\mathcal{W}_{j}'\right)^\nu, \\ &\phantom{=} \hspace{8cm} \check{F}(B_\nu,2f-1) \cap F(B_\nu,2f)\cap G(B_\nu,2f)^{c}  \bigg). \label{eq:first_upper_bound_on_I}
    \end{split}
\end{equation}
We now bound the above probability for a fixed box $B_\nu$ as in the sum above and a fixed nice level $j+4\leq 2f<j(1+\nu)$. For the sake of brevity, we will drop $B_\nu$ and $f$ from the notation in the following paragraph. Denote by $B_\nu^f$ the box centered at the center of $B_\nu$ of radius $4n^{-1}R^{\kappa-2f+1}$. Because $B_\nu$ is always chosen to cover some $S$-Whitney box, we know that $S\subseteq \!\left(B_\nu^f\right)^c$. In particular, conditioning on the percolation configuration $\tilde{\omega}:=\omega\!\left(\!\left(B_\nu^f\right)^c\right)$ outside of $B_\nu^f$ we get that
\begin{equation}
    \begin{split}
        \P \!\left(\cE_{\textrm{in}}(S,x_S),\cE_{\textrm{out}}(S,x_S), B_\nu\in \!\left(\mathcal{W}_{j}'\right)^\nu, \check{F}(B_\nu,2f-1)\cap F(B_\nu,2f)\cap G(B_\nu,2f)^{c}  \right) \\
        = \mathbf{E} \!\left[ \mathbf{1}_{\mathcal{E}_{\textrm{in}}(S,x_S)\cap B_\nu\in \!\left(\mathcal{W}_{j}'\right)^\nu} \cdot \mathbf{E}\!\left[ \mathbf{1}_{\check{F}\cap F\cap G^{c}} \cdot \mathbf{1}_{\cE_{\textrm{out}}(S,x_S)} \big| \tilde{\omega}\right]\right],
    \end{split}
    \label{eq:upper_nound_on_probability_conditional}
\end{equation}
where we have used that $\cE_{\textrm{in}}(S,x_S)$ and $B_\nu\in (\mathcal{W}_j' )^\nu$ are measurable with respect to $\tilde{\omega}$. Further, using Lemma~\ref{lem:E_out_happens_with_positive_probability} on the (nice) box $B_\nu^f$ that is fully contained in $S^c$ yields that
\begin{equation}\label{eq:decorrelating-Eout}
    \begin{split}
    \P\!\left(\check{F},F,G^{c},\cE_{\textrm{out}}(S,x_S)|\tilde{\omega}\right)&= \P\!\left(\check{F},F,G^{c},\cE_{\textrm{out}}(S,x_S),\mathscr{P}(\tilde{\omega})|\,\tilde{\omega}\right) \\
       &= \P\!\left(F,G^{c},\cE_{\textrm{out}}(S,x_S) |\check{F},\mathscr{P}(\tilde{\omega}),\tilde{\omega}\right)\cdot \P\!\left(\check{F},\mathscr{P}(\tilde{\omega})|\,\tilde{\omega}\right) \\
       &\leq \P\!\left(F,G^{c}|\check{F},\mathscr{P}(\tilde{\omega}),\tilde{\omega}\right)\cdot \P\!\left(\check{F},\mathscr{P}(\tilde{\omega})|\,\tilde{\omega}\right) \\
       &\leq C(R,\nu,\alpha) \cdot \P(F,G^{c})\cdot \P\!\left(\mathcal{E}_{\textrm{out}}(S,x_S)|\,\tilde{\omega}\right).
    \end{split}
\end{equation} In the last line, we have also used that $\check{F},\tilde{\omega}, \mathscr{P}$ and $F,G$ depend on disjoint sets of sites. We can now apply Proposition~\ref{lem:distancebetweencrossings} with $\zeta=\frac{2}{3}$ to find that 
    \begin{equation}
    \begin{split}
    \P\!\left(\check{F}\cap F\cap G^{c} \cap \cE_{\textrm{out}}(S,x_S) \bigg| \tilde{\omega} \right) &\leq C(R,\nu,\alpha) \cdot\P\!\left(\cE_{\textrm{out}}(S,x_S)\big|\,\tilde{\omega}\right) \cdot  s_j^{-2\alpha} R^{-2\alpha f} \\  &\leq C(R,\nu,\alpha) R^{8\nu \alpha j} R^{-\alpha j} \cdot \P\!\left(\cE_{\textrm{out}}(S,x_S)\big|\tilde{\omega}\right),
    \end{split}
    \end{equation}
    where the constant $C(\alpha,R,\nu)$ can increase from line to line. The last inequality follows by substituting the value $s_j=\frac{(16C_1R^3)^{-1}}{(M_1+1)^2} R^{-4\nu j}$ and noting that $2f\geq j+4$.

Plugging this back into equation~\eqref{eq:upper_nound_on_probability_conditional}, we find that 
\begin{align*}
        &\P \!\left(\cE_{\textrm{in}}(S,x_S),\cE_{\textrm{out}}(S,x_S), B_\nu\in \!\left(\mathcal{W}_{j}'\right)^\nu, \check{F}(B_\nu,2f-1)\cap F(B_\nu,2f)\cap G(B_\nu,2f)^{c}  \right) \\
        &\leq C(R,\nu,\alpha)\cdot s_j^{-2\alpha} R^{-\alpha f} \cdot \mathbf{E} \!\left[ \mathbf{1}_{\mathcal{E}_{\textrm{in}}(S,x_S) \cap B_\nu\in \!\left(\mathcal{W}_{j}'\right)^\nu} \cdot \mathbf{E}\!\left[ \mathbf{1}_{\cE_{\textrm{out}}(S,x_S)} \big| \omega\!\left(\!\left(B_\nu^f\right)^c\right)\right]\right] \\
        &= C(R,\nu,\alpha)\cdot s_j^{-2\alpha} R^{-\alpha f} \cdot \mathbf{P} \!\left(\mathcal{E}_{\textrm{in}}(S,x_S),B_\nu\in \!\left(\mathcal{W}_{j}'\right)^\nu, \cE_{\textrm{out}}(S,x_S) \right) \\
        &\leq C(R,\nu,\alpha)\cdot R^{8\nu\alpha j} R^{-\alpha j} \cdot \mathbf{P} \!\left(\mathcal{E}\!\left(B^{(k)}\right),\tau_1<\infty, \Gamma_1=(S,x_S),B_\nu\in \!\left(\mathcal{W}_{j}'\right)^\nu \right),
\end{align*}
where the constant $C(R,\nu,\alpha)$ can increase from line to line.
Substituting everything back into equation~\eqref{eq:first_upper_bound_on_I}, we get that 
\begin{align*}
     \mathit{II} \leq{}&   C(R,\nu,\alpha) \cdot \sum_{j=k-1}^{J} R^{-\alpha j(1-8\nu)} \sum_{(S,x_S)} \sum_{B\in \cB_{j,n}^{R,k}}\sum_{B_\nu \textrm{ covers }B}\sum_{j+4\leq 2f<j(1+\nu)} \\
     & \hspace{4cm} \mathbf{P} \!\left(\mathcal{E}\!\left(B^{(k)}\right),\tau_1<\infty, \Gamma_1=(S,x_S),B_\nu\in \!\left(\mathcal{W}_{j}'\right)^\nu \right) \\
     \leq{}&   C(R,\nu,\alpha) \cdot \sum_{j=k-1}^{J} R^{-\alpha j(1-8\nu)} \sum_{(S,x_S)} \sum_{B\in \cB_{j,n}^{R,k}}\nu jR^{2j\nu} \cdot\mathbf{P} \!\left(\mathcal{E}\!\left(B^{(k)}\right),\tau_1<\infty, \Gamma_1=(S,x_S),B\in \mathcal{W}_{j}' \right) \\
     \leq{}& C(R,\nu,\alpha)\cdot \sum_{j=k-1}^J jR^{-j(\alpha-8\nu \alpha-2\nu)} \sum_{(S,x_S)} \bigg\{ \mathbf{P}\!\left(\mathcal{E}\!\left(B^{(k)}\right),\tau_1<\infty, \Gamma_1=(S,x_S)\right) \\ 
     & \hspace{4cm}\cdot \E\!\left[ |\mathcal{W}_j'| \cdot (\mathbf{1}_{H_j} +\mathbf{1}_{H_j^c}) \, \big| \mathcal{E}\!\left(B^{(k)}\right),\tau_1<\infty, \Gamma_1=(S,x_S)\right] \bigg\} \\
     \leq{}& C(R,\nu,\alpha) \cdot \sum_{j=k-1}^J jR^{-j(\alpha-8\nu \alpha-2\nu)}\cdot \bigg[\sum_{(S,x_S)} C_1 R^{\frac{3}{2}j\nu} \cdot \mathbf{P}\!\left(\mathcal{E}\!\left(B^{(k)}\right),\tau_1<\infty, \Gamma_1=(S,x_S)\right) \\
     &+ R^{2j} \cdot \mathbf{P}\!\left(H_j^c,\mathcal{E}\!\left(B^{(k)}\right),\tau_1<\infty, \Gamma_1=(S,x_S)\right)\bigg] \\
     \leq{}& C(R,\nu,\alpha) \cdot\sum_{j=k-1}^J R^{-j(\alpha-8\nu \alpha-2\nu)}\!\left[ R^{\frac{3}{2}j\nu} \cdot \mathbf{P}\!\left(\mathcal{E}\!\left(B^{(k)}\right),\tau_1<\infty\right) + R^{2j} \cdot \P(H_j^c)\right].
\end{align*}
In the fourth inequality above, we have again used equation~\eqref{eq:bound-Wj} saying that on $H_j$ there exists a constant $C_1=C_1(R,\nu)$ such that $|\mathcal{W}_{j}'|\leq C_1R^{\frac{3}{2}\nu j}$. Using the a priori estimate, and Lemma~\ref{lem:bound_on_global_event} we find that 
\begin{equation}
    \begin{split}
        \mathit{II} &\leq C(R,\nu,\alpha) \cdot \sum_{j=k-1}^J R^{-j(\alpha-8\nu \alpha-2\nu)}\!\left( R^{\frac{3}{2}}j\nu \cdot R^{-\alpha k} + R^{2j} \cdot R^{-4\alpha j} \right).         
    \end{split}
\end{equation}
Recall that we are working in the range where $\alpha>3/4-\varepsilon$, where $\varepsilon$ can be chosen to be smaller than $1/8$, for example. With this in mind, choosing $\nu=1/100$ ensures that both series above are summable and yield an improved exponent $\alpha+\beta$ for some $\beta>0$ that is independent of $\alpha$. Indeed, the first series above is summable with the exponent $\alpha-8\nu\alpha-\frac{7}{2}\nu>\alpha/2>1/4$ for all $\alpha>1/2$. The second series is summable with the exponent $5\alpha-8\nu\alpha-2\nu-2>\frac{6}{5}\alpha$. This concludes the proof.
\end{proof}

We now finish the proof of Proposition~\ref{prop:subpolynomialconc}. 
\begin{proof}[Proof of Proposition~\ref{prop:subpolynomialconc}]
Let $\varepsilon>0$ and assume that~\eqref{eq:assumptionsubpoly} holds for some $\alpha>\frac{3}{4}-\varepsilon$. By Corollary~\ref{cor:no_individual_bubble_long} and Lemma~\ref{lem:bound_on_stopping_time}, there exists some $\varepsilon_0,\beta>0$ and $C=C(R,\alpha)$ such that if $\varepsilon<\varepsilon_0$, then 
\begin{align*}
        \P(X_{R^{\kappa-k}}\geq 6\altq_n)&\leq \P\!\left(\cE\!\left(B^{(k)}\right)\right)^{-1}\cdot \P\!\left(N \altq_n+\sum_{\ell=1}^N X_{\tau_{\ell}-1,\tau_\ell}\geq 6\altq_n\right) \\
        &\leq \P\!\left(\cE\!\left(B^{(k)}\right)\right)^{-1}\cdot\!\left(\P(2N\altq_n\geq 6\altq_n)+\P\!\left(\exists i\colon X_{x_{i},x_{i+1}}\geq \altq_n\right)\right)
        \\
        &\leq \P\!\left(\cE\!\left(B^{(k)}\right)\right)^{-1}\cdot\!\left(\P(\tau_2<\infty)+\P\!\left(\exists i\colon X_{x_{i},x_{i+1}}\geq \altq_n\right)\right)\\ 
        &\leq C(R,\alpha) R^{-(\alpha+\beta)k}.
\end{align*}
By applying Lemma~\ref{lem:compXscales} (specifically Eq.~\eqref{eq:comp_mn_scales}) with, say, $\zeta=1/4$, there exists $k_0=k_0(R)$ such that for all $n$ sufficiently large,
\begin{equation}
    \altq_n\geq 6\altq_{R^{\kappa-k_0}}.
\end{equation}
Thus,
\begin{align*}
        \P(X_{R^{\kappa-k}}\geq \altq_n)\leq&\,\P(X_{R^{\kappa-k}}\geq 6\altq_{R^{\kappa-k_0}})\\
        =&\,\P(X_{R^{\kappa-k_0-(k-k_0)}}\geq 6\altq_{R^{\kappa-k_0}})\\
        \leq&\,C(R,\alpha)R^{-(\alpha+\beta)(k-k_0)}\\
        \leq&\,C(R,\alpha)R^{-(\alpha+\beta)k},
\end{align*}
where the value of $C(R,\alpha)$ may change from line to line.
\end{proof}

Combining Proposition~\ref{prop:subpolynomialconc} and Lemma~\ref{lem:compXscales}, we obtain the following result.
\begin{cor}\label{cor:subpolynomialconc}
For all $\alpha,\zeta>0$, $R$ sufficiently large depending on $\zeta$, and all $p$ sufficiently small depending on $R,\zeta$, there exists $C=C(\alpha,\zeta,R,p)>0$ such that for all $k\geq 1$, $n>R^k$  
and $0\leq s\leq 1$,
\begin{equation}
\P(X_{R^{\kappa-k}}\geq s\altq_n)\leq Cs^{-(\frac{4}{3}+\zeta)\alpha}R^{-\alpha k}.
\end{equation}
\end{cor}

\section{Tightness}\label{sec:tightness}

The goal of this section is to prove Propositions~\ref{prop:GHf-tightness-dU} and~\ref{prop:main-conv-dU-GH}.  We show Proposition~\ref{prop:holderct} in Section~\ref{sec:holder}, using the superpolynomial concentration proved in Sections~\ref{sec:polynomial_tails} and~\ref{sec:superpolynomial-concentration}. Finally, we complete the proofs of Propositions~\ref{prop:GHf-tightness-dU} and~\ref{prop:main-conv-dU-GH} in Section~\ref{sec:conv-CLE}.

\subsection{H\"older continuity}\label{sec:holder}
The goal this section is to prove Proposition~\ref{prop:holderct}. In fact, we prove the following stronger statement.

\begin{prop}\label{prop:holderct-general}
 For all $\eta<\frac{5}{8}$ and $\varepsilon>0$, there exist $p=p(\eta)$ and $C=C(\eta,\varepsilon)>0$ such that for all $n$,
\begin{equation}\label{eq:holder-bound-path-distance}
    \P\!\left(\forall x,y\in \Lambda_n\colon\,\altq_n^{-1}D_n(x,y)\leq C\dpt_n(x,y)^\eta\right)>1-\varepsilon.
\end{equation}
Furthermore, let $\gamma$ denote any open path in $\Lambda_n$. Then for any $\nu>0$ and $\alpha>0$, there exist $p=p(\nu)$, $C_1=C_1(\nu,\alpha),C_2=C_2(\nu,\alpha)$ such that for all $n$ and $M$ sufficiently large,
\begin{equation}\label{eq:holder-bound-away-from-edge}
    \P\!\left(\exists \gamma: \altq_n^{-1}D_n(\gamma(0),\gamma(1))> C_1(M+\dist(\gamma,\partial\Lambda_n)^{-\frac{1}{4}-\nu})\right)\leq C_2M^{-\alpha} .
\end{equation}
\end{prop}

For the purposes of proving that the scaling limit of $d_n$ exists, we only need the first bound~\eqref{eq:holder-bound-path-distance} in terms of $\dpt_n$. The second bound~\eqref{eq:holder-bound-away-from-edge} is used to prove Corollary~\ref{cor:scaling-exponent-distance}. We use it in Section~\ref{sec:scaling-exponent-proof} to bound the expected length of a left-to-right crossing. Most of this section will be aimed at proving~\eqref{eq:holder-bound-path-distance}. The proof of~\eqref{eq:holder-bound-away-from-edge} then follows by a slight modification.

Let $\cB_{k,n}^{R,+}$ be a set of boxes such that $\bigcup_{B\in\cB_{k,n}^{R,+}}B\supset[0,1]^2$, every box $B\in\cB_{k,n}^{R,+}$ contains $R^{\kappa-k}\times R^{\kappa-k}$ vertices, or equivalently, the side length of $B$ is $n^{-1}R^{\kappa-k}$, and such that $|\cB_{k,n}^{R,+}|\leq\lceil nR^{k-
\kappa}\rceil^2$. Note that the boxes in $\cB_{k,n}^{R,+}$ may overlap. If $B\in\cB^{R,+}_{k,n}$ is such that $\dist(B,\partial \Lambda_n)\geq n^{-1}R^{\kappa-k-2}$, we say that $B$ is an \textit{inner} box.  

The strategy of the proof is similar to that in the proof of Lemma~\ref{lem:bound_on_stopping_time}. We will construct a global event for percolation on $\Lambda_n$ that will ensure two things: \begin{enumerate}
    \item [(1)] We can use suitably small number of \textit{inner} $R$-adic boxes of each level $k$ to cover an open path from any two points $x,y$ in $\Lambda_n$.
    \item[(2)] No \textit{inner} box is bad on its own -- meaning that no inner box can have a crossing that is too long.
\end{enumerate}
The reason why we are considering only inner boxes is that our only tool to control lengths of crossings are Lemma~\ref{lem:bound_on_F} and Proposition~\ref{lem:distancebetweencrossings}, which require existence of an annulus around our box. This means that we can only apply these lemmas for boxes that are away from the boundary. On the above events we will be able to use our linking procedure to show that ${d_n}$ is comparable to some power of ${\dpt_n}$. 

We first construct an event $J$ that addresses the first item above.
Let $\nu\in (0,1)$ to be determined later. Denote by $S(a)$ the annulus (or \emph{strip}) inside of $\Lambda_n$ of width $a$ around $\partial \Lambda_n$. Given any $\frac{3}{\nu}\leq k\leq \kappa$, let $J_k^1$ be the following event: if there exists an open path of diameter $n^{-1}R^{\kappa-k}$, then this path must have a point outside of the strip $S\!\left(n^{-1}R^{\kappa-(1+\nu)k+3}\right)$. Let $J_k^2$ be the event that there are at most $R^{\!\left(\frac{1}{3} +\nu\right)k}$ many $R$-adic boxes $B$ of level $k$ such that $R^3\cdot B \cap \partial \Lambda_n\neq \emptyset$ and such that there are $3$ arms of alternating colors crossing $ R^{3+k\!\left(\frac{1}{3} +\nu\right)}\cdot B\setminus R^3\cdot B$.

\begin{lem} Given any $\frac{1}{2}>\nu>0$, there exists 
$C=C(\nu,R)$ large enough such that for all $K\geq\frac{4}{\nu}$, $$\P\!\left(\bigcup_{k= K}^{\kappa}\!\left((J^1_k)^c\cup (J_k^2)^c\right)\right)\leq CR^{-\nu K}.$$ 
\label{lem:bound_on_first_event_holder}
\end{lem}
Note that the events $J_k^1$ and $J_k^2$ are only relevant for the proof of~\eqref{eq:holder-bound-path-distance}. The bound $R^{-\nu K}$ is not strong enough to obtain the necessary concentration for~\eqref{eq:holder-bound-away-from-edge}. This is precisely why we require the paths in~\eqref{eq:holder-bound-away-from-edge} to be away from the boundary of the box.
\begin{proof}[Proof of Lemma~\ref{lem:bound_on_first_event_holder}]
Let us first bound $\P(J_k^1)$. Consider at most $4\left\lceil nR^{{\!\left(1+\nu\right)}k-3-\kappa}\right\rceil$ (not necessarily disjoint) boxes of side length $n^{-1} R^{\kappa-\!\left({1+\nu}\right)k+3}$ that cover the strip $S\!\left(n^{-1} R^{\kappa-\!\left({1+\nu}\right)k+3}\right)$. On the event $(J_k^1)^c$, there is an open path fully contained in this strip that moves distance at least $\frac{1}{2}n^{-1}R^{{\kappa-k}}$. This means that there exists a box $B'$ in our covering such that if we consider the next (in clockwise order) $\left\lfloor\frac{1}{2} R^{\nu k-3}\right\rfloor$ \textit{pairwise disjoint} boxes in our covering (excluding the 4 corner boxes), each of them has either a top-to-bottom or a left-to-right open crossing -- see Figure~\ref{fig:event_J_1} for an illustration. Since the probability of having such a crossing is $\frac{1}{2}$, we can bound $\P\!\left((J_k^1)^c\right)$ as follows:
\begin{equation}
    \begin{split}
        \P\!\left((J_k^1)^c\right) \leq 4\left\lceil n R^{\!\left({1+\nu}\right)k-\kappa}\right\rceil \cdot \!\left(\frac{1}{2}\right)^{\left\lfloor\frac{1}{2} R^{\nu k-3}\right\rfloor} 
    \end{split}
\end{equation}
Note that $(\frac{1}{2})^{\frac{1}{2}R^{\nu k-3}}\leq C(R)R^{-2 k}$ for all $k\geq 4 /\nu$ and $R$. Thus, 
\begin{equation}
    \P((J_k^1)^c)\leq C(R)R^{-\nu k}
\end{equation}
for all $k\geq 4/\nu$ and $\nu<\frac{1}{2}$.

\begin{figure}[ht]
    \centering
    \includegraphics[scale=0.7]{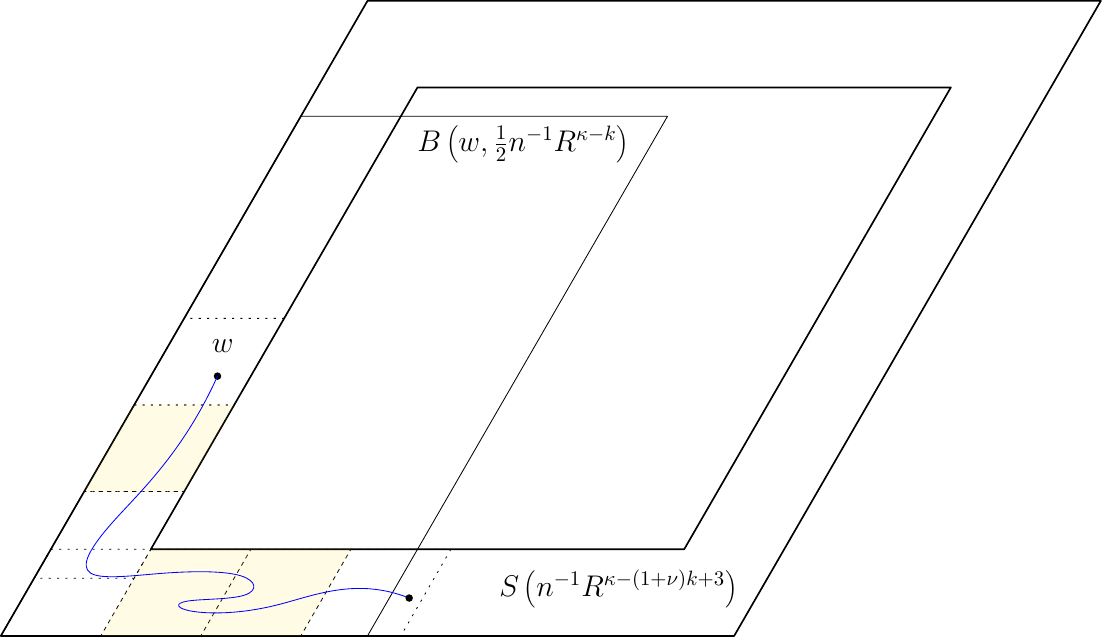}
    \caption{This figure illustrates an instance of the event $(J_k^1)^c$. The blue curve represents an open path of diameter at least $n^{-1}R^{\kappa-k}$. This path must leave the box $B\!\left(w,\frac{1}{2}n^{-1}R^{\kappa-k}\right)$, where $w$ is one of its endpoints. Moreover, it stays within the strip $S\!\left({n^{-1}\cdot R^{\kappa-\!\left(1+\nu\right)k+3}}\right)$. The dashed boxes above represent the not necessarily disjoint boxes we used to cover this strip. The boxes shaded in yellow are the ones that must have either a top-to-bottom or a left-to-right open crossing. Note that the boxes in the corner need not have such crossings, so we skip them in our analysis.}
    
    \label{fig:event_J_1}
\end{figure}

On the other hand, we can bound $\P((J_k^2)^c)$ using Markov's inequality and Lemma~\ref{lem:arm_exponents}. We find that for any error in the exponent $\zeta>0$, there exists a constant $C(\zeta)$ such that $$\P\!\left((J_k^2)^c\right) \leq C(\zeta)\frac{4R^4\cdot R^{k}\cdot \!\left(R^{k\!\left( \frac{1}{3}+\nu \right)}\right)^{-2+\zeta}}{R^{\!\left(\frac{1}{3}+\nu\right)k}} = C(\zeta)4R^4\cdot R^{(\frac{1}{3}+\nu)\zeta k} R^{-3\nu k}.$$
The desired result now follows by taking $\zeta=\frac{2\nu}{\frac{1}{3}+\nu}$.
\end{proof}
Define the event $J=J(\nu,K)$ to be $\bigcap_{k=K}^{\kappa} (J_k^1\cap J_k^2)$. 

To construct an event that will deal with the second requirement above, we resort to Remark~\ref{rem:distancebetweencrossings_general}, and Lemma~\ref{lem:bound_on_F} and Proposition~\ref{lem:distancebetweencrossings}.
Recall the definitions of events $F$ and $G$. Fix $\zeta>0$ very small, to be determined later. We will choose all the parameters of $F,G$ later depending on $R,\alpha$. Let also $\nu$ be a parameter in $(0,1)$ (the same one as in the previous lemma). 

Let $B$ be an inner $R$-adic box of level $k$. Consider all $R$-adic boxes $B^\nu\in\cB_{\left\lceil (1+\nu)k\right\rceil+4,n}^{R,+}$ needed to cover $B$. Let $E_1(B^\nu,s_k)$ be the event that there exists a level $k+4\leq f<(1+\nu)k+5$ such that $F(M_1,f)\cap G(s_k,f)$ holds for the annulus ${A\!\left(z;n^{-1}R^{\kappa-f},n^{-1}R^{\kappa-f+1}\right)}$, where $z$ is the center of $B^\nu$.

Like in the proof of Lemma~\ref{lem:crossingbubble}, Lemma~\ref{lem:bound_on_F} and Proposition~\ref{lem:distancebetweencrossings} tell us that for any $p'>0$, we can choose the parameters of $F,G$ in such a way that for $R$ sufficiently large depending on $\zeta$, $$\P(E_1(B^{\nu},s_k)^c)\leq (p')^{\nu k}+C(\alpha,M_1,R) s_k^{-\!\left(\frac{4}{3}+\zeta\right)\alpha}R^{-\alpha k}$$ for all levels $k$ such that $\xi\frac{c_1c_2}{8}\cdot R^{\kappa - (1+\nu)k-{6}}\geq R$.
Choosing $p'\leq R^{-\alpha/2}$ ensures that the first term can be bounded from above by $R^{-\nu \alpha k/2}$. 
Further, setting $s_k=\max \{C(\alpha,M_1,R)^{\frac{3}{4\alpha}},1\}\cdot R^{-\!\left(\frac{3}{4}-\nu\right)k}$ and $\zeta=\nu$, we find that the second term in the above sum is smaller than $R^{-\nu\alpha k/2}$. Thus,
\begin{equation}\label{eq:bound-holder-prob-one}
    \P\!\left({\bigcup_{B^\nu\in\cB_{{\lceil(1+\nu)k \rceil+4},n}^{R,+}}}E_1(B^\nu,s_k)^c\right) \leq C(R) \cdot R^{2\!\left(1+\nu\right)k}\cdot R^{-\nu\alpha k/2}
\end{equation}

for all levels $k\leq \widehat{K}$, where  $\widehat{K}=\frac{1}{1+\nu}\!\left(\kappa -\log_R \!\left( \frac{R^{{7}}}{\xi\frac{c_1c_2}{8}}\right)\right)$.

Choosing $\alpha>8/\nu$ makes the above summable (since we have chosen $\nu<1$). 
Since we can choose~$\alpha$ to be any large number by Corollary \ref{cor:subpolynomialconc-general}, the following lemma follows readily from the above.

\begin{lem}
Given any $\nu\in (0,1)$ and $\alpha'$, for all $R$ large enough depending on $\nu,$ and $p$ small enough depending on $R,\nu$ there exist parameters of events $F,G$ (depending only on $R,\alpha',\nu$) and a constant $C_1=C_1(R)$ such that for all $K$,
\begin{equation}\label{eq:bound-holder-prob-two}
    \P \!\left(\bigcup_{k=K}^{\widehat{K}}{\bigcup_{B^\nu\in\cB_{{\lceil(1+\nu)k \rceil+4},n}^{R,+}}} E_1(B^\nu,s_k)^c\right)\leq C_1R^{-\alpha' K}
\end{equation}

where $\widehat{K}=\frac{1}{1+\nu}\!\left(\kappa -\log_R \frac{R^{{7}}}{\xi\frac{c_1c_2}{8}}\right)$ and $s_k=C_2(R,\alpha',\nu) R^{-\!\left(\frac{3}{4}-\nu\right)k}$, for some constant $C_2(R,\alpha',\nu)$ depending only on $R,\alpha',\nu$. 
\label{lem:bound_on_second_event_holder}
\end{lem} 

Note that in order to obtain the exponent $\alpha'$ in~\eqref{eq:bound-holder-prob-two} from~\eqref{eq:bound-holder-prob-one}, we need to set the parameter $\alpha$ for $s_k$ and $p'$ to be sufficiently large depending on $\alpha'$ ($\alpha>8\alpha'/\nu$ would work). Further note that the value for $s_k$ is chosen to optimize the exponent in~\eqref{eq:holder-bound-path-distance}. The same $s_k$ works to show~\eqref{eq:holder-bound-away-from-edge}, and we do not attempt to optimize it for the sake of simplicity. 

We have now constructed the appropriate global events, so all that remains to prove~\eqref{eq:holder-bound-path-distance} is to show that on these events we really do have H\"older continuity of the desired distances. This is the content of the following lemma.

\begin{lem}\label{lem:holder_bound_on_event}
    Using the same parameters as before, for all $\nu<\frac{5}{8\cdot12}$, 
    on the event $$J(\varepsilon,\nu) \, \cap \, \bigcap_{k=K}^{\widehat{K}}{\bigcap_{B^\nu\in\cB_{{\lceil(1+\nu)k \rceil+4},n}^{R,+}}} E_1(B^\nu,s_k),$$ we have 
    \begin{equation}\label{eq:holder-path-on-event}
        {D_n}(x,y)\leq C(R,\alpha,\nu,K)\cdot {\dpt_n}(x,y)^{\frac{5}{8}-12\nu}\altq_n,
    \end{equation}
    for some constant $C(R,\alpha,\nu,K)$.
\end{lem}

The proof of this Lemma uses a similar path covering and stringing argument to that in the proofs of Lemmas~\ref{lem:crossingbubble} and~\ref{lem:bound_on_stopping_time}. The proof consists of several steps. We let $\gamma$ be a diameter-minimizing path from $x$ to $y$. In Step 1, we define the Whitney decomposition of $\Lambda_n\setminus\{x,y\}$. In order for the covering of $\gamma$ not to have too many boxes, we want $\gamma$ to stay away from the boundary of $\Lambda_n$. Therefore, in Step 2, we use the events $J^1_k$ to define a modification of $\gamma$ that mostly stays away from the boundary. In the parts of the modified $\gamma$ that do come near the boundary, there must be a 3-arm event. In Step 3, we bound the number of Whitney boxes needed to cover the modified path, using the event $J^2$. We then give the final bound on $D_n$ in Step 4. Recall the explanation in Section \ref{subsubsec:holder-sketch} on how we obtain the exponent $\frac{5}{8}$.

\begin{proof}[Proof of Lemma~\ref{lem:holder_bound_on_event}]
\medskip \textbf{Step 1:} \textit{Whitney box decomposition.}

Fix any two distinct points $x$ and $y$ inside of $\Lambda_n$, and let $\gamma$ be an open path between $x$ and $y$ in $\Lambda_n$ with minimal ({Euclidean}) diameter. Consider a box $B^{\textrm{ref}}\subseteq \Lambda_n$ of side length ${\dpt_n}(x,y)$ that contains the entire path $\gamma$.

We will consider the Whitney box decomposition of the set $\Lambda_n\setminus \{x,y\}$, \emph{started at level $K$}. Recall that the boundary of $\Lambda_n$ is the set of sites in $\bT_n$ that neighbor $\Lambda_n$ but are not contained in it. This ensures that the vertices on the boundary of $[0,1)^2$ are contained in a box in the Whitney decomposition. 
In particular, consider all $R$-adic boxes $B{\in\cB_{k,n}^{R,+}}$ of levels $K\leq k \leq \kappa$ such that $n^{-1}R^{\kappa-k-2}\leq \textrm{dist}(B, \partial \Lambda_n\cup \{x,y\})\leq n^{-1}R^{\kappa-k+2}$, as well as all boxes $B\in\cB^{R,+}_{K,n}$ of level $K$ such that $\dist(B,\partial\Lambda_n\cup\{x,y\})\geq n^{-1}R^{\kappa-K-2}$. As before, to see that every point $z\in \Lambda_n \setminus \{x,y\}$ is covered by at least one such box, note that there exists $0\leq i\leq \kappa$ such that $n^{-1}R^{\kappa-i}\leq\textrm{dist}(z,\partial \Lambda_n\cup \{x,y\})<n^{-1}R^{\kappa-i+1}$. If $i\geq K$, then at least one of the two $R$-adic boxes of levels $i$ and $i+1$ containing $z$ must be a Whitney box. If $i<K$, then one of the $R$-adic boxes of level $K$ containing $z$ must be a Whitney box.

As in the previous sections, our goal is to first show that we can cover (a modification of) $\gamma$ by a small number of Whitney boxes. Then, we can `string' together short crossings across these boxes, allowing us to compare the (geodesic or resistance) distance between $x$ and $y$ to the radius of the reference box $B^{\textrm{ref}}$.

To do so, we will first alter the path $\gamma$ within the box $B^{\textrm{ref}}$ so that the new path from $x$ to $y$ consists of boundaries of some well-chosen closed clusters. On the event $J( K,\nu)$, we will be able to cover the modified path by a `small' number of Whitney boxes, as advertised.

Assume that $0\leq k_0\leq \kappa$ is such that $n^{-1}R^{\kappa-k_0}\leq {\dpt_n}(x,y) < n^{-1}R^{\kappa-k_0+1}$. Suppose that we are working on the event $J(K,\nu)$. Let ${k_1}=\max\!\left\{K,\left\lceil\frac{k_0}{\frac{2}{3}-2\nu}\right\rceil\right\}$. Note that the diameter of $\gamma$ is at least $n^{-1}R^{\kappa-(\frac{2}{3}-2\nu)k_1}$.

\medskip \textbf{Step 2:} \textit{Path modification.}

Since we are working on the event $J(K,\nu)$, this means that there exists a point of $\gamma$ outside of the strip $S\!\left(n^{-1}R^{\kappa-(1+\nu)\!\left(\frac{2}{3}-2\nu\right){k_1}+3}\right)$. Let $x_{{k_1}}$ be the first point of $\gamma$ starting from $x$ that leaves this strip, and define $y_{{k_1}}$ similarly. Note that if $x$ is not inside the strip, then $x_{{k_1}}=x$ and similarly for $y$.

The box $B^{\textrm{ref}}$ contains the entire straight line segment $[x_{{k_1}},y_{{k_1}}]$ between $x_{{k_1}}$ and $y_{{k_1}}$. By a slight abuse of notation, we denote the part of $\gamma$ between $x_{k_1}$ and $y_{k_1}$ by $\gamma|_{[x_{{k_1}},y_{{k_1}}]}$. This segment partitions the path $\gamma|_{[x_{{k_1}},y_{{k_1}}]}$ into chunks $\gamma^{{1}},\ldots,\gamma^{{t}}$, so that the endpoint of each $\gamma^{{\tilde{t}}}$ and the starting point of $\gamma^{{\tilde{t}+1}}$ are neighbors, and such that no $\gamma^{{\tilde{t}}}$ `crosses' $[x_{{k_1}},y_{{k_1}}]$ {--} $\gamma^{{\tilde{t}}}$ can contain sites that belong to $[x_{{k_1}},y_{{k_1}}]$, but it cannot contain {neighboring} sites {that are} strictly on other sides of $[x_{{k_1}},y_{{k_1}}]$. 
For each $1\leq \tilde{t}\leq t$, $\gamma^{{\tilde{t}}}$ together with $[x_{{k_1}},y_{{k_1}}]$ bounds a domain in $\Lambda_n$. Note that it is possible that $\gamma^{\tilde t}$ starts off at one side of $[x_k,y_k]$, then wraps around behind $x_k$ or $y_k$ and hits $[x_k,y_k]$ on the other side. This is not a problem, since $\gamma^{\tilde t}$ and $[x_k,y_k]$ still bound a domain. This domain must be fully contained inside of $B^{\textrm{ref}}$ since this box is convex. Consider the open {path} connecting the beginning to the endpoint of $\gamma^{\tilde{t}}$ inside of this domain that is closest to the segment $[x_{{k_1}},y_{{k_1}}]$ -- this is the boundary of all the closed clusters attached to this section of $[x_{{k_1}},y_{{k_1}}]$. Update this section of $\gamma$ to be precisely this boundary. Upon doing this for all $1\leq \tilde{t}\leq t$, we find that every point on $\gamma|_{[x_{{k_1}},y_{{k_1}}]}$ is connected by a closed path to the segment $[x_{{k_1}},y_{{k_1}}]$.  See Figure~\ref{fig:modification_of_gamma_x_y.pdf}.

\begin{figure}[ht]
    \centering
    \includegraphics[scale=0.8]{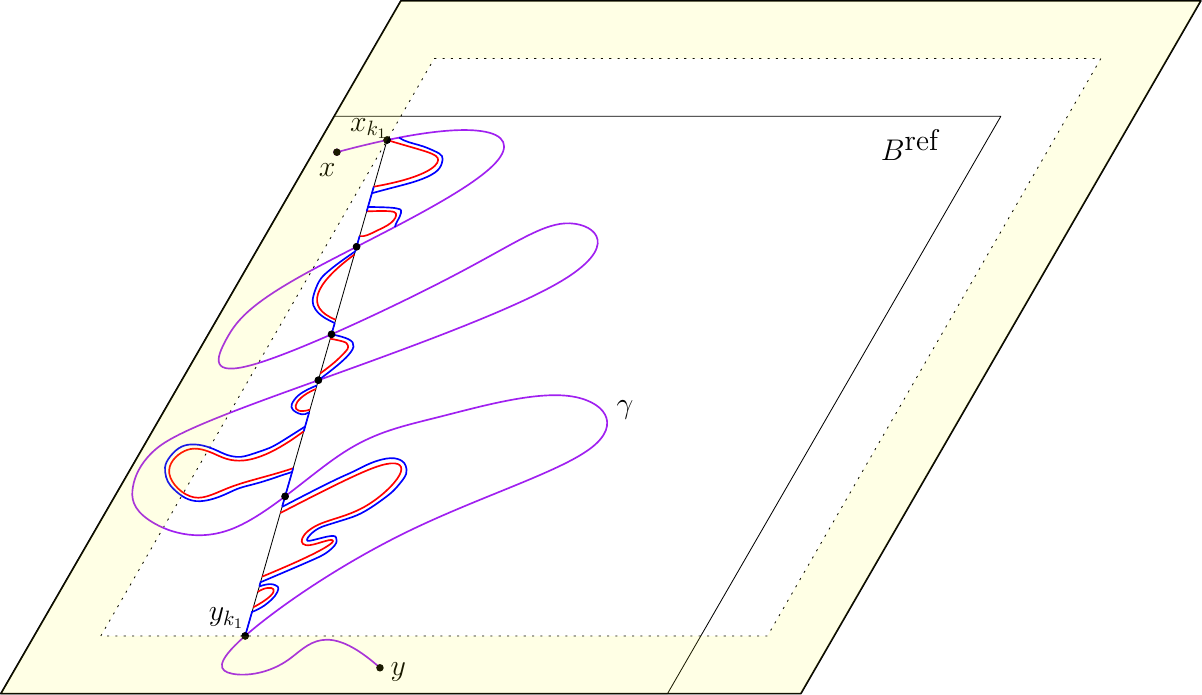}
    \caption{This figure depicts the first step of our modification of the path $\gamma$ connecting $x$ and $y$. The original path $\gamma$ is represented by the purple curve. This path is contained inside the reference box $B^{\textrm{ref}}$ marked above. Points $x_{{k_1}}$ and $y_{{k_1}}$ are the first exit points of $\gamma$ from the strip $S\!\left(n^{-1}R^{\kappa-(1+\nu)\!\left(\frac{2}{3}-2\nu\right){{k_1}}+3}\right)$ (shaded in yellow above) starting from $x$ and $y$ respectively. One can see that the segment $[x_{{k_1}},y_{{k_1}}]$ partitions $\gamma$ into several chunks and that these, together with the curve $\gamma$, enclose domains inside $\Lambda_n$. The blue-red curves inside of these domains represent precisely the cluster boundaries of all closed clusters inside of these domains attached to the segment $[x_{k_1},y_{k_1}]$. In particular, the blue curve is our modification of $\gamma$. Note that this modification still belongs to $B^{\textrm{ref}}$.}
    \label{fig:modification_of_gamma_x_y.pdf}
\end{figure}

Then either $\gamma|_{[x,x_{{k_1}}]}$ leaves the box $B\!\left(x, n^{-1}R^{\kappa-\!\left(\frac{2}{3}-2\nu\right)({{k_1}}+1)}\right)$ 
or not. If it does, then $\gamma|_{[x,x_{{k_1}}]}$ has diameter at least $n^{-1}R^{\kappa-\!\left(\frac{2}{3}-2\nu\right)({{k_1}}+1)}$, and thus it cannot be fully contained in the strip $S\!\left(n^{-1}R^{\kappa-(1+\nu)(\frac{2}{3}-2\nu)(k_1+1)+3}\right)$ on the event $J$. Let $x_{{{k_1}}+1}$ be the first point of $\gamma$ after $x$ that leaves this strip. Then edit $\gamma_{[x_{{{k_1}}+1},x_{{k_1}}]}$ so that every point on it has a closed path connecting it to the segment $[x_{{{k_1}}+1},x_{{k_1}}]$. If $\gamma$ does not leave $B\!\left(x,n^{-1} R^{\kappa-\!\left(\frac{2}{3}-2\nu\right)({k_1}+1)}\right)$, simply set $x_{{k_1}+1}=x_{k_1}$. Do the same for $\gamma|_{[y_{k_1},y]}$. Inductively define $x_k,y_k$ as above for all $k_1\leq k\leq\kappa$ and alter $\gamma_{[x_{k+1},x_k]}$ and $\gamma_{[y_{k+1},y_k]}$ similarly. 

\medskip \textbf{Step 3:} \textit{Number of Whitney boxes.}

Having constructed the altered path $\gamma$, we can now bound the number of Whitney boxes needed to cover $\gamma$. Since $n^{-1}R^{\kappa-k_0}\leq \dpt_n(x,y) < n^{-1}R^{\kappa-k_0+1}$, all Whitney boxes needed to cover $B^{\textrm{ref}}\subseteq \Lambda_n$ have level at least $k_0-4$. 
First assume that $\max\{K,k_0-4\}\leq k< k_1$. If $k=K$ (in case $K>k_0-4$), we can bound the number of Whitney boxes needed to cover $\gamma$ by $4R^2\cdot R^{2(K-k_0)}$, which is the total number of $R$-adic boxes of level $K$ constrained in a box of size $n^{-1}R^{\kappa-k_0+1}$. Otherwise, if $k>K$, there are at most $R^{k+9-k_0}$ 
Whitney boxes of level $k$ in the covering of $\gamma$, since $\gamma$ is contained in a box of size $n^{-1}R^{\kappa-k_0+1}$. 
Now assume that $k\geq k_1$. Then any Whitney box $B$ in the covering of $\gamma$ is either at distance at most $R^2\cdot n^{-1}R^{\kappa-\!\left(\frac{2}{3}-2\nu\right)k}$ from $x$ or $y$, or it covers a section of $\gamma$ in between $x_k$ and $y_k$. Assume that the latter is true. Then since $x,y\not\in R^3B$, the box $R^3B$ must intersect $\partial\Lambda_n$ since $B$ is a Whitney box. 
Note that $R^3 B$ has two disjoint open arms to $x_k$ and $y_k$ and a closed arm to $[x_k,y_k]$. Since $x_k,y_k\not\in S(n^{-1}R^{\kappa-(\frac{2}{3}-2\nu)k})$, it follows that there are 3 arms of alternating colors crossing $R^{3+k(\frac{1}{3}+\nu)}\cdot B\setminus R^3\cdot B$. On the event $J$, the number of such boxes $B$ is at most $R^{(\frac{1}{3}+\nu)k}$. Furthermore, the number of boxes of the first type of level $k$ an be at most $32R^4+16R^{\!\left(\frac{1}{3}+2\nu\right)k+4}$  -- there are at most $32R^4$ boxes $B$ such that $\dist(B,\{x,y\})\leq n^{-1}R^{\kappa+2-k}$ and at most $2\cdot 4\cdot 2R^2\cdot R^{\!\left(\frac{1}{3}+2\nu\right)k+2}$ boxes $B$ such that $\dist(B,\partial \Lambda_n)\leq n^{-1}R^{\kappa+2-k}$ and $\dist(B,\{x,y\}) \le R^2\cdot n^{-1}R^{\kappa-\!\left(\frac{2}{3}-2\nu\right)k}$.

Thus, if we denote by $\mathcal{W}_k$ the set of all Whitney boxes of level $k$ needed to cover the altered path $\gamma$, the above shows that on the event $J(K,\nu)$, we have $|\mathcal{W}_k|\leq 32R^{4}+ 32 R^{\!\left(\frac{1}{3}+2\nu\right)k+4}$ for all $k\geq k_1$. For $\max\{K,k_0-4\} < k<k_1$, we have $|\cW_k|\leq R^{k+9-k_0}$. For $k=K$ (in case $K > k_0 -4$) we have that $|\mathcal{W}_K|\leq 4R^2\cdot R^{2(K-k_0)}$, which is an upper bound on the total number of all $R$-adic boxes of level $K$ contained in a box of size $n^{-1}R^{\kappa-k_0+1}$.

Let $\mathcal{W}_k^*$ denote the set of all $R$-adic boxes of level $\lceil(1+\nu) k\rceil+4$ needed to cover boxes of $\mathcal{W}_k$. 
On the event $\bigcap_{k=K}^{\widehat{K}}\bigcap_{B^\nu\in\cB_{\lceil(1+\nu)k\rceil+4,n}^{R,+}} E_1(B^\nu,s_k)$, we can control the lengths of all open crossings of the relevant annuli of levels $K\leq k\leq \widehat{K}$. For boxes $B\in \mathcal{W}_k$ of level $k\geq \hat{K}$, we can bound both geodesic and resistance distances between any two points inside of $2B$ by the number of sites inside the box $2B$. In particular we can bound this distance from above by $$4R^{2(\kappa-k)}\leq 4R^{2\kappa} \cdot \!\left(\frac{R^6}{\min \{\frac{c_1c_2}{8},\xi\}}\cdot R^{-\kappa}\right)^{2\cdot \frac{1}{1+\nu}}\leq C(R,\alpha,\nu)\cdot R^{\frac{2\nu}{1+\nu}\kappa},$$ where $C(R,\alpha,\nu)$ is some constant depending only on $R,\alpha,\nu$ (through $\xi\frac{c_1c_2}{8}$). At the same time, this quantity provides an upper bound on how many times the same box can appear in the covering of the path $\gamma$, since $\gamma$ will not enter the same site of the box twice.

\medskip \textbf{Step 4:} \textit{Bounding the distance.}

We are now ready to establish the bound ~\eqref{eq:holder-path-on-event} on $D_n$. Assume first that $(3/2+5\nu)k_0\leq\kappa$. Using an analogous linking argument to the one used in Lemmas~\ref{lem:crossingbubble} and~\ref{lem:bound_on_stopping_time} (cf.\ \eqref{eq:upper_bound_on_distance_x_v}), we find that the distance between $x$ and $y$ is at most 
\begin{align}\label{eq:final-bound-D_n}
D_n(x,y)\leq{}& \sum_{k\geq \max\{K,k_0-4\}}^{\widehat{K}} |\mathcal{W}_k^*|\cdot M_1\cdot 2(M_1+1)\cdot s_k\altq_n+ \sum_{k=\widehat{K}}^{\kappa} |\mathcal{W}_k|\cdot \!\left(C(R,\alpha,\nu) R^{\frac{2\nu}{1+\nu}\kappa}\right)^2 \notag\\
\leq{}& C(R,\alpha,\nu)\cdot R^{2\nu K}\cdot R^{2(K-k_0)}\cdot R^{-(\frac{3}{4}-\nu)K}\altq_n\1_{K> k_0-4}\notag\\
&+ \sum_{k\geq k_0-4}^{k_1} C(R,\alpha,\nu)R^{2\nu k}R^{k-k_0}R^{-(\frac{3}{4}-\nu)k}\altq_n\notag\\
&+\sum_{k=k_1}^{\widehat{K}}C(R,\alpha,\nu)\cdot R^{2\nu k}\cdot R^{\!\left(\frac{1}{3}+2\nu\right)k}\cdot R^{-\!\left(\frac{3}{4}-\nu\right)k}\altq_n+\sum_{k=\widehat{K}}^{\kappa} C(R,\alpha,\nu)\cdot R^{\!\left(\frac{1}{3}+2\nu\right)k} R^{\frac{4\nu}{1+\nu}\kappa}\notag\\
\leq{}& C(R,\alpha,\nu,K)R^{-2k_0}\altq_n+ C(R,\alpha,\nu)R^{(\frac{1}{4}+3\nu)k_1-k_0}\altq_n \notag\\
&+ C(R,\alpha,\nu)\cdot\!\left[\!\left(\sum_{k\geq k_1}  R^{-\!\left(\frac{5}{12}-5\nu\right)k}\altq_n\right) + \kappa R^{\!\left(\frac{1}{3}+2\nu+4\nu\right)\kappa}\right]\notag\\
\leq{}& C(R,\alpha,\nu,K) \!\left(R^{-(\frac{5}{8}-6\nu)k_0}\altq_n + R^{-(\frac{5}{12}-5\nu)k_1}\altq_n+R^{-(\frac{5}{12}-8\nu)k_{1}}\cdot R^{\!\left(\frac{3}{4}-
\nu\right)\kappa}\right)\notag\\
\leq{}& C(R,\alpha,\nu,K)\!\left(R^{-(\frac{5}{8}-8\nu)k_0}\altq_n+R^{-(\frac{5}{8}-12\nu){k_0}}\cdot R^{\!\left(\frac{3}{4}-\nu\right)\kappa}\right)
\end{align}
where the constants can increase from row to row, but depend only on the highlighted values. The penultimate inequality is obtained by multiplying the final term in the previous row by $R^{-\!\left(\frac{5}{12}-8\nu\right)k_{1}}\cdot R^{\!\left(\frac{5}{12}-8\nu\right)\kappa}\geq C(R,\nu,K)$ and rearranging. Here and in the final line, we use that $\frac{3}{2}k_0\leq k_1\leq (\frac{3}{2}+5\nu)k_0+K\leq\kappa+K$  for $\nu\leq\frac{5}{12\cdot8}$.

Recall from Corollary~\ref{cor:lower_bound_on_m_n} that there exists a constant $C=C(\nu,R)$ such that $\altq_n\geq CR^{(\frac{3}{4}-\nu)\kappa}$. Finally, because $R^{-k_0}\leq n^{-1}\cdot R^{\kappa-k_0}\leq {\dpt_n}(x,y)$, we obtain that $${D_n}(x,y)\leq C(R,\alpha,\nu,K) \!\left({\dpt_n}(x,y)\right)^{\frac{5}{8}-12\nu}\altq_n.$$ 

Now assume that $(3/2+5\nu)k_0\geq\kappa$. Then for $k\geq \widehat K$, we bound the number of Whitney boxes at level $k$ by $R^{k+9-k_0}$, like in the case $K<k_0\leq k_1$. Then using the same bounds as above, we obtain a bound of
\begin{equation}
    \begin{split}
        D_n(x,y)\leq{}& C(R,\alpha,\nu,K)\left(R^{-(\frac{5}{8}-6\nu)k_0}\altq_n+R^{-(\frac{5}{12}-5\nu)k_1}\altq_n+\sum_{k=\widehat{K}}^\kappa R^{k-k_0}R^{\frac{4\nu}{1+\nu}\kappa}\right)\\
        \leq{}&C(R,\alpha,\nu,K)\left(R^{-(\frac{5}{8}-8\nu)k_0}\altq_n+R^{\kappa-k_0}R^{\frac{4\nu}{1+\nu}\kappa}\right)\\
        \leq{}&C(R,\alpha,\nu,K)\left(R^{-(\frac{5}{8}-8\nu)k_0}\altq_n+R^{(\frac{3}{4}-\nu)\kappa}R^{(\frac{1}{4}+5\nu)\kappa-k_0}\right)\\
        \leq{}&C(R,\alpha,\nu,K)\left(R^{-(\frac{5}{8}-8\nu)k_0}\altq_n+R^{(\frac{1}{4}+5\nu)(\frac{3}{2}+5\nu)k_0-k_0}\altq_n\right)\\
        \leq{}&C(R,\alpha,\nu,K)R^{-(\frac{5}{8}-12\nu)k_0}\altq_n\\
        \leq{}&C(R,\alpha,\nu,K)\!\left({\dpt_n}(x,y)\right)^{\frac{5}{8}-12\nu}\altq_n
    \end{split}
\end{equation}
as desired.
\end{proof}

We now prove~\eqref{eq:holder-bound-away-from-edge}.
\begin{lem}\label{lem:holder-bound-away-from-edge}
Using the same parameters as before, on the event
\begin{equation}
    \bigcap_{k=K}^{\widehat{K}}{\bigcap_{B^\nu\in\cB_{{\lceil(1+\nu)k \rceil+4},n}^{R,+}}}E_1(B^\nu,s_k),
\end{equation}
the following holds. Let $x,y\in\Lambda_n$ and let $\gamma$ be a path from $x$ to $y$ with minimal diameter. Assume that $\dist(\gamma,\partial\Lambda_n)\geq R^{-L}$. Then 
\begin{equation}
    D_n(x,y)\leq C(R,\alpha,\nu) \!\left(R^{(\frac{5}{4}+3\nu)K}+R^{(\frac{1}{4}+3\nu)L}\right) \altq_n.
\end{equation}
\end{lem}
\begin{proof}
The proof is analogous to that of the previous lemma, except that we do not need to alter the path $\gamma$ since it does not come near the boundary. Let $k_0$ be such that $n^{-1}R^{\kappa-k_0}\leq \diam(\gamma) < n^{-1}R^{\kappa-k_0+1}$. As before, the number of Whitney boxes of level $K$ needed to cover $\gamma$ is at most $C(R)R^{2K}$. For $K < k \le L+3$, the number of Whitney boxes of level $k$ needed to cover $\gamma$ is at most $C(R)R^{k}$. Finally, since $\dist(\gamma,\partial\Lambda_n)\geq R^{-L}$, the only Whitney boxes of level $k>L+3$ in the covering of $\gamma$ are at distance at most $n^{-1}R^{\kappa+2-k}$ from $\{x,y\}$. As stated before, there are at most $32R^4$ such boxes. Thus, (recall~\eqref{eq:final-bound-D_n}),
\begin{equation}
    \begin{split}
        D_n(x,y)\leq
        {}& C(R,\alpha,\nu)\cdot R^{2\nu K}\cdot R^{2K}\cdot R^{-(\frac{3}{4}-\nu)K}\altq_n +C(R,\alpha,\nu)\sum_{k=\max \{K,k_0-4\}}^{L+3} R^{2\nu k}\cdot R^{k}\cdot R^{-(\frac{3}{4}-\nu)k}\altq_n \\
        &+\sum_{k= L+4}^{\widehat{K}} C(R,{\alpha},\nu)\cdot R^{2\nu k}\cdot32R^4\cdot R^{-\!\left(\frac{3}{4}-\nu\right)k}\altq_n + \sum_{k=\widehat{K}+1}^{\kappa} C(R,\alpha,\nu)\cdot 32R^4\cdot R^{\frac{4\nu}{1+\nu}\kappa} \\
        \leq{}& C(R,\alpha,\nu)\!\left(R^{(\frac{5}{4}+3\nu)K}+R^{(\frac{1}{4}+3\nu)L}\right)\altq_n,
    \end{split}
\end{equation}

which completes the proof. For the last inequality, we use that $\altq_n\gtrsim R^{(\frac{3}{4}-\nu)\kappa}$ to bound the last sum as explained in the proof of Lemma~\ref{lem:holder_bound_on_event}. 
\end{proof}
\begin{proof}[Proof of Proposition~\ref{prop:holderct-general}] 
Let $\nu>0$, $\eta=\frac{5}{8}-12\nu$, and $\varepsilon>0$. Set $\alpha'=1$ (say), choose the parameters for the events $F,G$ as in Lemma~\ref{lem:bound_on_second_event_holder}, and let $R$ be sufficiently large depending on $\nu$, and let $p$ be sufficiently small depending on $\nu,R$ as in Lemma~\ref{lem:bound_on_second_event_holder}. Let $K$ be large enough that $CR^{-\nu K}<\varepsilon/2$ and $C_1R^{-K}<\varepsilon/2$. Then the proof of~\eqref{eq:holder-bound-path-distance} follows immediately from Lemmas~\ref{lem:bound_on_first_event_holder},~\ref{lem:bound_on_second_event_holder} and~\ref{lem:holder_bound_on_event}.

For the proof of~\eqref{eq:holder-bound-away-from-edge}, let $\nu,\alpha'>0$. Let $R$ be large enough depending on $\nu,$ and $p$ sufficiently small depending on $\nu,R$ as in Lemma~\ref{lem:bound_on_second_event_holder}. Let $M>0$ and let $K$ be the largest integer such that $R^{(5/4+3\nu)K}\leq M$. Then by Lemmas~\ref{lem:bound_on_second_event_holder} and~\ref{lem:holder-bound-away-from-edge}, with probability at least $1-C_1R^{-\alpha'K}\geq1-C_1M^{-\alpha'(5/4+3\nu)^{-1}}$, every path $\gamma$ satisfies $d_n(\gamma(0),\gamma(1))\leq C(R,\alpha',\nu)(M+\dist(\gamma,\partial\Lambda_n)^{-(1/4+3\nu)})$. Since $\alpha',\nu$ were arbitrary, the proof is complete.
\end{proof}
\begin{rem}\label{rem:generalization-holder}
Note that using the above proofs we can obtain H\"older continuity of these two metrics in any domain with sufficiently nice boundary, not just $[0,1)^2$. {In particular, the proof works to show the analogue of Proposition~\ref{prop:holderct-general} for the metrics $d_{Q,n}$ and $\dpt_{Q,n}$ for $Q\in\cQ$.}
\end{rem}

\subsection{Proof of tightness}\label{sec:conv-CLE}
In this section, we prove Propositions~\ref{prop:GHf-tightness-dU} and~\ref{prop:main-conv-dU-GH}.

\begin{proof}[Proof of Proposition~\ref{prop:GHf-tightness-dU}]
Let $N'\colon[0,\infty)\to\N$ be some function. Let $C,\eta>0$. Let $\mathcal{K}_0=\mathcal{K}_0(N')$ be the set of all compact metric spaces $(X,d)$ such that for all $\varepsilon>0$, $N_{(X,d)}(\varepsilon)\leq N'(\varepsilon)$, with $N_{(X,d)}(\varepsilon)$ the number of $\varepsilon$-balls required to cover $K$. Define $\mathcal{K}=\mathcal{K}(N',C,\eta)$ to be the set of all metric-function spaces $(X,d,f)$ such that $(X,d) \in \mathcal{K}_0$, $f\colon X\to Y$, with $Y$ as defined at the beginning of Section~\ref{sec:tightness-sketch}, and for all $(x,y)\in X$, we have $d(x,y)\leq Cd_Y(f(x),f(y))^{\eta}$. Then since $\mathcal{K}_0$ is relatively compact in the GH topology, the space $\mathcal{K}$ is relatively compact in the GHf topology by \cite[Proposition A.4]{amy2025tightness}.

First note that by Proposition~\ref{prop:convergence-ghf-path}, there exists a coupling and a random function $N'$ such that $((\cX_{Q,n}^{(m)})^2,(\dpt_{Q,n})^2)\in\mathcal{K}_0(N')$ for all $m,n$. Furthermore, the map $\pi_n$ is trivially H\"older continuous with respect to $\dpt_{Q,n}$ with exponent $1$. Recall that $\median_n$ was chosen in Section~\ref{subsubsec:holder-sketch} such that $d_{Q,n}$ is uniformly H\"older continuous with a fixed exponent (e.g.~$\eta=1/4$) with high probability. 
It then follows from Proposition~\ref{prop:holderct} and Remark~\ref{rem:generalization-holder} that for all $\varepsilon>0$, 
there exists $C=C(\varepsilon,\eta)$ such that $\P(((\cX_{Q,n}^{(m)})^2,(\dpt_{Q,n})^2,d_{Q,n},\pi_n)\in\mathcal{K}(N',C,\eta))>1-\varepsilon$. This implies that $(\Lambda_{Q,n}^2,(\dpt_{Q,n})^2,d_{Q,n},\pi_n)$ is tight. 
\end{proof}

\begin{proof}[Proof of Proposition~\ref{prop:main-conv-dU-GH}]
Let $Q\in\cQ$, let $m\in\N$ and let $I\subset\N$ be a subsequence such that the sequence $((\cX^{(m)}_{Q,n})^2,(\dpt_{Q,n})^2,f_{Q,n})_{n\in I}$ converges almost surely in GHf under a coupling with limit $((\cX^{(m)}_{Q})^2,(\dpt_{Q})^2,f_{Q})$, where $f_Q(x,y)=(\pi(x),\pi(y),d_Q(x,y))$, with $\pi$ the embedding of $\cX^{(m)}_Q$ into $[0,1]^2$ and $d_Q$ some function. We show that $d_Q$ is a pseudo-metric and that $(\cX^{(m)}_{Q,n},d_{Q,n},\pi_n)_{n\in I}$ converges to $(\cX^{(m)}_{Q},d_{Q},\pi)$ in GHf distance.

We claim that by the GHf convergence of $((\cX^{(m)}_{Q,n})^2,(\dpt_{Q,n})^2,f_{Q,n})_{n\in I}$ and the $\GH$ convergence of $(\cXmun,\dpt_{Q,n})_{n\in\N}$, upon switching to a further subsequence, there exist isometric embeddings $\psi\colon \cXmu \to W$, $\psi_n\colon \cXmun \to W$ such that for each $x,y\in \cX^{(m)}_Q$ there are sequences $(x_n)_{n\in I},(y_n)_{n\in I}$ such that $x_n,y_n\in\cXmun$, $d_W(\psi_n(x_n),\psi(x)) \to 0$, $d_W(\psi_n(y_n),\psi(y)) \to 0$, and $\lim_{n\in I}d_{Q,n}(x_n,y_n)=d_Q(x,y)$. To see this, let $(W,d_W)$ be a compact space and $\psi_n,\psi$ isometric embeddings such that $\Delta_{\Haus}(\psi_n(\cXmun),\psi(\cXmu)) \to 0$ (cf.\ \cite[Lemma A.2]{amy2025tightness}). Then the sequence\footnote{Here the notation $(f,\psi(u),\psi(v))$ is shorthand for the function $(u,v) \mapsto (f(u,v),\psi(u),\psi(v))$.} \[ ((\cX^{(m)}_{Q,n})^2,(\dpt_{Q,n})^2,f_{Q,n},\psi_n(u),\psi_n(v))_{n\in I} \] is relatively compact, and converges in GHf along a further subsequence to $$((\cX^{(m)}_{Q})^2,(\dpt_{Q})^2,f_{Q},\psi(u),\psi(v)).$$ This implies the claim. Further, note that the continuity of $d_Q$ implies that $$d_Q(x,y) = \lim_{n\in I} d_{Q,n}(x_n,y_n)$$ for \emph{any} sequence $x_n \to x$, $y_n \to y$ in $d_W$. It follows that $d_Q$ satisfies the triangle inequality. So $d_Q$ is indeed a pseudo-metric.

We now show the GHf convergence $(\cXmun,d_{Q,n},\pi_n) \to (\cXmu,d_Q,\pi)$ along $I$. Since the sequence $((\cX^{(m)}_{Q,n})^2,(\dpt_{Q,n})^2,f_{Q,n})_{n\in I}$ converges to $((\cX^{(m)}_{Q})^2,(\dpt_{Q})^2,f_{Q})$, it implies in particular that the sequence $(\cXmun,d_{Q,n},\pi_n)_{n \in I}$ is relatively compact in GHf. 
We argue that any subsequential limit is isometric to $(\cXmu,d_Q,\pi)$. By the argument in the previous paragraph, we have that $$((\cX^{(m)}_{Q,n})^2,(\dpt_{Q,n})^2,f_{Q,n},\psi_n(u),\psi_n(v)) \to ((\cX^{(m)}_{Q})^2,(\dpt_{Q})^2,f_{Q},\psi(u),\psi(v))$$ along a further subsequence. Let $n$ be sufficiently large that the GHf distance is at most $\varepsilon$. Then for each $(x_n,y_n)\in(\cXmun)^2$, there exists $(x_{(x_n,y_n)},y_{(x_n,y_n)})\in(\cXmu)^2$ such that we have $d_W(\psi_n(x_n),\psi(x_{(x_n,y_n)}))<\varepsilon$, $d_W(\psi_n(y_n),\psi(y_{(x_n,y_n)}))<\varepsilon$, $|\pi_n(x_n)-\pi(x_{(x_n,y_n)})| < \varepsilon$, $|\pi_n(y_n)-\pi(y_{(x_n,y_n)})| < \varepsilon$, and $|d_{Q,n}(x_n,y_n)-d_Q(x_{(x_n,y_n)},y_{(x_n,y_n)})|<\varepsilon$.

Let $\cR \subseteq \cXmun \times \cXmu$ be the correspondence given by
\[ \cR = \{(x_n,x) : d_W(\psi_n(x_n),\psi(x)) \le \varepsilon ,\, |\pi_n(x_n)-\pi(x)|\le\varepsilon\} .\]
Suppose that $(x_n,x),(y_n,y)\in\cR$. Let $(x_{(x_n,y_n)},y_{(x_n,y_n)}) \in (\cXmu)^2$ be as in the paragraph above, and note that in particular $(x_n, x_{(x_n,y_n)}) \in \cR$, $(y_n, y_{(x_n,y_n)}) \in \cR$, so that $\dpt_{Q}(x,x_{(x_n,y_n)}) \le 2\varepsilon$, $\dpt_{Q}(y,y_{(x_n,y_n)}) \le 2\varepsilon$. Then
\begin{equation}
    \begin{split}
        &|d_{Q,n}(x_n,y_n)-d_Q(x,y)|\\
        &\quad \leq |d_{Q,n}(x_n,y_n)-d_Q(x_{(x_n,y_n)},y_{(x_n,y_n)})|+|d_Q(x_{(x_n,y_n)},y_{(x_n,y_n)})-d_Q(x,y)|\\
        &\quad \leq \varepsilon+d_Q(x_{(x_n,y_n)},x)+d_Q(y_{(x_n,y_n)},y)\\
        &\quad \leq o(1) \quad\text{as } \varepsilon \searrow 0
    \end{split}
\end{equation}
where in the last line we used that $d_Q$ is continuous on $(\cXmu,\dpt_Q)$ (in fact, we have $d_Q(u,v)\leq C(\omega)\dpt_Q(u,v)^{\eta}$ by Proposition~\ref{prop:holderct} and Remark~\ref{rem:generalization-holder}). This shows that $(\cXmun,d_{Q,n},\pi_n)$ converges to $(\cXmu,d_Q,\pi)$ in GHf along the subsequence and uniquely identifies any subsequential limit.
\end{proof}

\addtocontents{toc}{
  \protect\contentsline{part}{\hyperref[sec:uniqueness-proof]{Part 3}}{}{}
}

\section{Identification of the limit}\label{sec:uniqueness-proof}

In this section, we prove Theorem~\ref{th:dU-is-CLE-metric} and Corollary~\ref{cor:final-convergence-metrics}. Let $(\fd_V)_{V\in\fC}$ be a collection of subsequential limiting metrics as defined in Section~\ref{sec:uniqueness-sketch}. We show that this is a $\CLE_6$ metric as defined in Section~\ref{subsec:cle-metrics}. We show compatibility in Section~\ref{sec:compatibility}. In Section~\ref{sec:series-parallel}, we show the series and generalized parallel laws. Finally, Markovianity and translation invariance are proved in Sections~\ref{sec:markov} and~\ref{sec:translation-invariance} respectively. This completes the proof that $\dg$ and $\dr$ are $\CLE_6$ metrics. We then show that $\dg$ is a geodesic $\CLE_6$ metric in Section~\ref{sec:geodesic-cle} and show that $\dr$ is a resistance $\CLE_6$ metric in Section~\ref{sec:edge-weight-compatibility} by showing that it satisfies the edge weight locality property. This completes the proof of Theorem~\ref{th:dU-is-CLE-metric}. We give the proof of Corollary~\ref{cor:final-convergence-metrics} in Section~\ref{sec:final-convergence-metrics}, thus completing the proof of GH convergence.

Throughout this section, we let $J \subseteq \N$ be a subsequence so that for each $Q \in \cQ$, the sequence $((\Lambda_{Q,n},\dpt_{Q,n},d_{Q,n},\pi_n))_{n \in J}$ converges weakly to $(\Lambda_{Q},\dpt_Q,d_Q,\pi)$ in the GHf topology as defined in Section~\ref{sec:tightness-sketch}. In the proofs we will work with some coupling so that $(\Lambda_{Q,n},\dpt_{Q,n},d_{Q,n},\pi_n) \to (\Lambda_{Q},\dpt_Q,d_Q,\pi)$ almost surely along $J$ for each $Q \in \cQ$. Recall from Section~\ref{subsec:topologies} and from the proof of Proposition~\ref{prop:main-conv-dU-GH} that this means that there are isometric embeddings $\psi_n\colon \cX_n \to W$, $\psi\colon \cX \to W$ in some ambient space $(W,d_W)$ so that they converge in $\Delta_f$ with respect to $d_W$. In the following, we will say `$x_n$ converges to $x$' if they converge in $d_W$ under an embedding. Finally, we will often assume that the coupling satisfies the property in Lemma~\ref{lem:touching-loops}.

\subsection{Compatibility}\label{sec:compatibility}
\begin{prop}\label{prop:compatibility}
Let $V,V'\in\fC$, $V\subset V'$, and $x,y\in\cX\cap\overline{V}$ such that for every $u\in V'\setminus V$ there is a point $z\in\overline{V}$ that separates $u$ from $x,y$ in $\cX\cap\overline{V'}$. Then
\begin{equation}
    \fd_{V'}(x,y)=\fd_V(x,y).
\end{equation}
\end{prop}

\begin{proof}
Let $Q,Q'\in\cQ$ such that $Q\Supset V$, $Q'\Supset V'$, and $Q \Subset Q'$. By the definition of~$\fC$, the boundary of $V'$ is part of finitely many loops of $\Gamma$. Therefore, if $Q'$ is sufficiently close to~$V'$, then every $u \in \cX \cap Q' \setminus \overline{V'}$ is separated from $V'$ in $\cX \cap Q'$ by a single point $z \in \partial V'$. In particular, every $u \in \cX \cap Q' \setminus Q$ is separated from $x,y$ in $\cX \cap Q'$ by a single point $z \in \overline{V}$. Indeed, if $u \in \overline{V'}$, then by the assumption there is a point $z \in \overline{V}$ separating $u$ from $x,y$ in $\cX \cap \overline{V'}$. Since $\cX \cap  Q' \setminus \overline{V'}$ consists only of `dead ends', the point $z$ also separates $u$ from $x,y$ in $\cX \cap Q'$. On the other hand, if $u \in Q' \setminus \overline{V'}$, then it lies in one of the `dead ends', so the same property holds again. See Figure \ref{fig:compatibility_property} below.

Let $x_n,y_n \in \cX_n \cap Q$ be converging to $x,y \in \cX \cap \overline{V}$, respectively, so that $d_Q(x,y) = \lim_n d_{Q,n}(x_n,y_n)$ and $d_{Q'}(x,y) = \lim_n d_{Q',n}(x_n,y_n)$. (Note that by the GHf convergence and since $d_Q,d_{Q'}$ are continuous, the limits do not depend on the choice of the approximating sequence of $x_n,y_n$.) Ultimately, we want to use the compatibility property on the discrete metrics $d_{Q,n}$ and $d_{Q',n}$ to say that the appropriate distances between $x_n$ and $y_n$ match. This will allow us to deduce compatibility for the limit metrics $d_Q$ and $d_{Q'}$, and in turn for $\fd_V$ and $\fd_{V'}$. To do so, we need that for every point $u\in \cX_n\cap Q'\setminus Q$, there exists a point $z \in \cX_n\cap Q$ so that every open path from $x_n$ or $y_n$ to $u$ passes through $z$. 

By the previous paragraph, this property holds in the continuum. We would like to use Lemma~\ref{lem:touching-loops} to argue that it holds also in the discrete for $n$ large enough. Since $V \Subset Q$, it is possible to find a \emph{finite} set of cut points $z_1,z_2,\ldots \in \cX \cap \overline{V}$ such that every $u \in \cX \cap Q' \setminus Q$ is separated from $x,y$ in $\cX \cap Q'$ by some $z_i$ (otherwise there would be infinitely many crossings of interface loops between $\partial V$ and $\partial Q$ which is impossible). 

Therefore, by Lemma~\ref{lem:touching-loops}, under a suitable coupling, there are points $z_{1,n},z_{2,n},\ldots \in \cX_n \cap Q$ such that the same property holds for $\cX_n$ for all $n$ sufficiently large. Therefore, by the analogous compatibility property for the geodesic resp.\ resistance metric on percolation, we have $d_{Q',n}(x_n,y_n) = d_{Q,n}(x_n,y_n)$ for sufficiently large $n$, hence $d_{Q'}(x,y) = d_Q(x,y)$. Since this holds for all $Q,Q'$ sufficiently close to $V,V'$, respectively, the result follows by~\eqref{eq:def-dU}.
\begin{figure}[ht]
    \centering
    \includegraphics[scale=0.8]{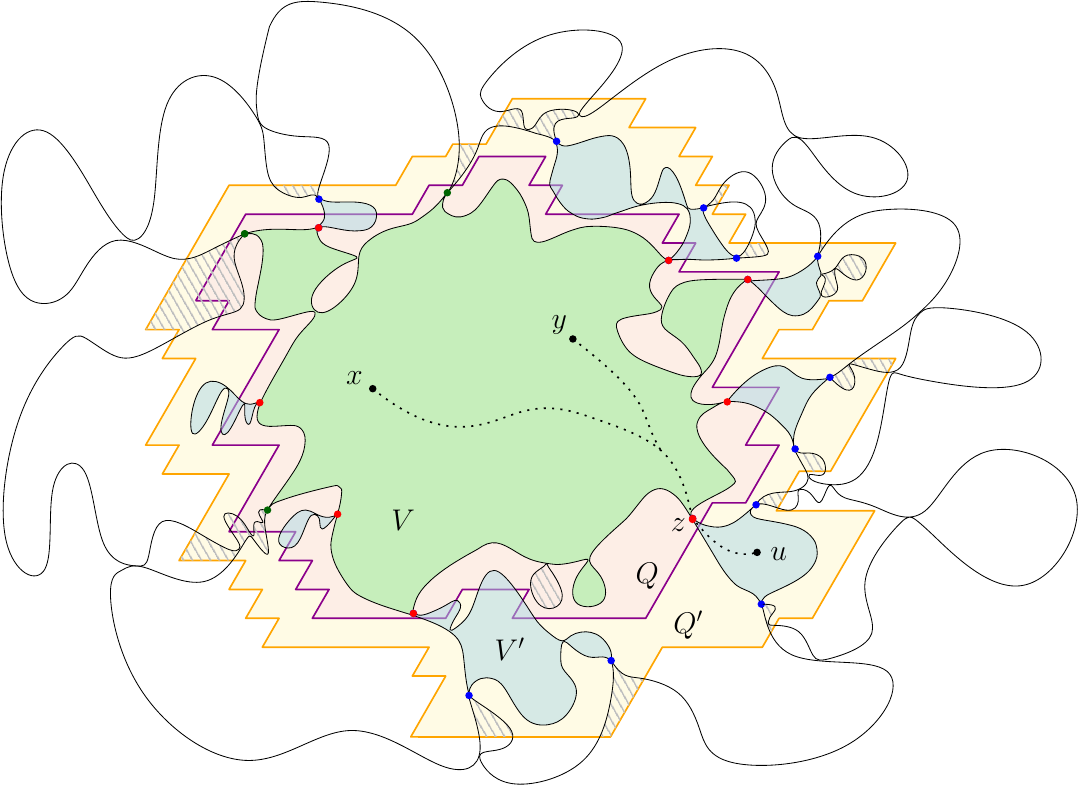}
    \caption{This figure illustrates the compatibility property for regions $V\subseteq V'$ in~$\fC$. The black loops in the picture above represent the loops of $\Gamma$. The set $V$ is represented by the green area; its boundary consists of segments of finitely many loops $\cG \subseteq \Gamma$. The set $V'\supseteq V$ consists of the blue and green areas above; we can see that any point $u\in V' \setminus V$ is separated from $x$ and $y$ in $V'$ by exactly one of the red points pictured above. The purple and yellow shaded areas represent two sets $Q,Q'$ in $\cQ$ such that $V\Subset Q$, $V'\Subset Q'$ and $Q\Subset Q'$. 
    Because the boundaries of $V$ and $V'$ consist of finitely many loops $\cG$ resp.\ $\cG'$, all boundary points except for the finitely many intersection points of distinct loops in $\cG$ resp.\ $\cG'$ (indicated in green resp.\ blue) lead to `dead ends' in~$\cX$. Further, since $\overline{V}$ and $\overline{V'}$ are required to be simply connected with respect to $\dpt_{\overline{V},\cG}$ (recall Section~\ref{subsec:cle-metrics}), these intersection points have positive Euclidean distance from another. This ensures that if the Euclidean Hausdorff distance of $Q$ resp.\ $Q'$ to $V$ resp.\ $V'$ is sufficiently small, the regions that lie behind these intersection points (indicated in gray) are also `dead ends' in $\cX \cap Q$ resp.\ $\cX \cap Q'$. Note that the total set of `dead ends' need not be finite since the loops comprising the boundary of $V$ resp.\ $V'$ have infinitely many self-intersection points. However, since $V \Subset Q$, only finitely many of them will exit~$Q$. Only these are relevant for comparing $d_{Q,n}$ and $d_{Q',n}$ -- the dead ends fully contained in $Q$ are also fully contained in $Q'$, so their discrete approximations do not incur a change between $d_{Q,n}(x_n,y_n)$ and $d_{Q',n}(x_n,y_n)$.}
    \label{fig:compatibility_property}
\end{figure}
\end{proof}

\subsection{Series and generalized parallel laws}\label{sec:series-parallel}
\begin{prop}\label{prop:series}
Let $V\in\fC$ and let $x,y,z\in \cX\cap\overline{V}$ such that $z$ separates $x$ from $y$ in $\cX\cap\overline{V}$. Then
\begin{equation}
    \fd_V(x,y)=\fd_V(x,z)+\fd_V(z,y).
\end{equation}
\end{prop}
\begin{proof}
Let $Q\in\cQ$. We first prove the series law for $d_Q$. Let $x_n,y_n\in\Lambda_n\cap Q$ such that $x_n\rightarrow x$ and $y_n\rightarrow y$. In the proof of Proposition~\ref{prop:compatibility}, it was shown that if $x,y\in Q$ are separated by $z\in Q$, then almost surely under a suitable coupling, there exists $n_0$ such that for all $n\geq n_0$, there exists $z_n\in\Lambda_n\cap Q$ such that $z_n$ separates $x_n$ and $y_n$ and $z_n\rightarrow z$. It then immediately follows from the series law for $d_n$ that
\begin{equation}
    \begin{split}
        d_Q(x,y)&=\lim_{n\rightarrow\infty}d_{Q,n}(x_n,y_n)\\
        &=\lim_{n\rightarrow\infty}d_{Q,n}(x_n,z_n)+d_{Q,n}(z_n,y_n)\\
        &=d_Q(x,z)+d_Q(z,y).
    \end{split}
\end{equation}
The statement for $V\in\fC$ follows from~\eqref{eq:def-dU}.
\end{proof}

\begin{prop}\label{prop:parallel}
Let $V\in\fC$ and let $x,y,z_1,\ldots,z_M\in V$ such that $x$ and $y$ are separated in $\cX\cap\overline{V}\setminus\{z_1,\ldots,z_M\}$. Let $K_x$ be the connected component of $\cX\cap\overline{V}\setminus\{z_1,\ldots,z_M\}$ containing $x$, and let $V_x \in \fC$ be the smallest region 
such that $\overline{V_x} \supseteq K_x$.\footnote{In other words, $\overline{V_x}$ is the union of $\overline{K_x}$ with the components bounded by the loops within $\overline{K_x}$.} Then 
\begin{equation}
    \fd_V(x,y)\geq\frac{1}{M}\min_i\fd_{V_x}(x,z_i).
\end{equation}
\end{prop}
\begin{proof}
Let $Q,Q'\in\cQ$ such that $Q\Supset V_x$ and $Q'\Supset V$. By the argument in the first paragraph of the proof of Proposition~\ref{prop:compatibility}, for all $Q'$ small enough, $x$ and $y$ are also separated in $\cX\cap Q'\setminus\{z_1,\ldots,z_M\}$. Then under a suitable coupling, there exist sequences $x_n,y_n,z_{1,n},\ldots,z_{M,n}\in\cX_n$ that converge to $x,y,z_1,\ldots,z_M$ respectively. Furthermore, by Lemma~\ref{lem:touching-loops}, the $z_i$ can be chosen such that $x_n$ and $y_n$ are separated by $z_{1,n},\ldots,z_{M,n} \in Q$ for all $n$ sufficiently large. Assume also that $Q$ contains the connected component of $\cX_n \cap Q' \setminus \{z_{1,n},\ldots,z_{M,n}\}$ containing $x_n$ (by making $Q'$ smaller if necessary, and letting $n$ be sufficiently large). 
Then by the parallel law for the percolation geodesic and resistance metrics, and the monotonicity,
\begin{equation}
    d_{Q'}(x,y)=\lim_{n\rightarrow\infty}d_{Q',n}(x_n,y_n)\geq\lim_{n\rightarrow\infty}\frac{1}{M}\min_i d_{Q,n}(x,z_i)=\frac{1}{M}\min_id_Q(x,z_i).
\end{equation}

Since this holds for all sufficiently small $Q,Q'$, the proof is complete.
\end{proof}

\subsection{Markovian Property}\label{sec:markov} 

We prove the Markovian property. Let $U\subset[0,1]^2$ be open and simply connected and recall the definitions of $\cL$, $\Gamma_{\cL}$, $\fC_U$, $U^*$ and $\Gamma_{U^*}$ from Section~\ref{subsec:cle-metrics}.
\begin{prop}\label{prop:markov}
Almost surely, the conditional law of the collection $(\cX\cap \overline{V},\fd_{V})_{V\in\fC_U}$ given $\{\cL\}\cup\Gamma_{\cL}\setminus\Gamma_{U^*}$ and $(\cX\cap \overline{V'},\fd_{V'})_{V'\in\fC_{[0,1]^2\setminus U}}$ is measurable with respect to $U^*$.
\end{prop}

Note that the characterization of the CLE$_6$ metric in Section~\ref{subsec:cle-metrics} is given in terms of the non-nested CLE$_6$ corresponding to the cluster $\cX$, but in the proof of Proposition~\ref{prop:markov} we will make use of the $\dg$ resp.\ $\dr$ metric defined on all clusters. This has the advantage that their local laws in disjoint regions are exactly independent (see Lemma~\ref{lem:indep-metrics-cle} below). The main technical work is to deduce from this the stated Markovian property for $(\fd_V)$.

We will need a few preparations. Suppose that $U \in \cQ$. For each $\ell\in\N$ we let
\begin{align*}
    W_1^*(\ell)&:=\bigcup\{Q\in\cQ_{\ell}: Q\subset U^*\} ,\\
    W_2^*(\ell)&:=W_1^*(\ell)\bigcup\{Q \cap U : Q\in\cB_{\ell},\, Q\cap \partial U^*\neq\varnothing\} .
\end{align*}
Note that for $\ell_1<\ell_2$, we have $W_1^*(\ell_{1})\subset W_1^*(\ell_{2}) \subset W_2^*(\ell_2) \subset W_2^*(\ell_1)$. Furthermore, as $\ell\rightarrow\infty$, the Hausdorff distance between $W_1^*(\ell)$ and $W_2^*(\ell)$ tends to $0$. We define the following sigma-algebras:
\[
    \cG_\ell = \sigma(\Gamma_{W_2^*(\ell)\setminus W_1^*(\ell)}, U^*) ,\qquad
    \cG = \bigcap_\ell \cG_\ell .
\]

\begin{lem}\label{le:ustar_infinitesimal_sigma_algebra}
    We have $\cG = \sigma(U^*)$ up to a.s.\ equivalence.
\end{lem}
\begin{proof}
A slight technical point in the proof is that because we selected $\cL$ to be the outermost loop surrounding the point $(1/2,1/2)$, this induces a non-trivial conditioning in the unbounded complementary connected component of $\cX$. To deal with this, we introduce a conditioning that leaves the law of $(\Gamma_U,U^*)$ absolutely continuous. 
Let us first consider the restriction of the law to the event that $U^*$ touches $\partial U$. Let $\widetilde{\p}$ be the law of $\p$ conditional on the event that $\cL$ surrounds~$U$. By resampling the connections of the loops outside $U$ (see~\cite[Section~5.2]{amy2025resampling}), the $\p$-law of~$(\Gamma_U,U^*)$ restricted to event that $U^*$ touches $\partial U$ is absolutely continuous with respect to its $\widetilde{\p}$-law. Therefore it suffices to prove the lemma for $\widetilde{\p}$. Under $\widetilde{\p}$, each complementary connected component of $\{\cL\}\cup\Gamma_{\cL}$ intersecting $U$ is bounded, and the conditional law of $\Gamma$ in each of them is just a conditionally independent CLE$_6$. Therefore (see Lemma~\ref{lem:measurability-cle}) the conditional law of $\Gamma_U$ in the connected components of $U \setminus \partial U^*$ is a conditionally independent collection of CLE$_6$, in particular, it depends only on $U^*$. 
    
Since $\cG$ is measurable with respect to $U^*$ and $\Gamma_{U}$, it suffices to show that $\Gamma_{U}$ is conditionally independent of $\cG$ given $U^*$. The measurability of CLE$_6$ with respect to the quad-crossing sigma algebra (see Section~\ref{subsubsec:convergence-boundary-loops} and \cite[Theorem 1.13(2)]{schramm2011scaling}) implies that $\sigma(\Gamma_{U},U^*)$ is generated by $\sigma(U^*)$ and the collections $(\Gamma_{Q_i}\1_{Q_i \Subset U \setminus \partial U^*})_{i=1,\ldots,n}$\footnote{Here, $\Gamma_{Q}\1_{Q \Subset U \setminus \partial U^*}$ denotes the collection $\Gamma_Q$ if $Q\Subset U \setminus \partial U^*$, and the empty collection otherwise.} where $Q_1,\ldots,Q_n \in \cQ$ and $Q_i \subseteq U$. By Lemma~\ref{lem:measurability-cle} and the previous paragraph (since the relative boundaries of $W_1^*(\ell),W_2^*(\ell)$ in $U$ are piecewise smooth), we have $\sigma(\Gamma_{U},U^*) = \cG_\ell \vee \sigma(\Gamma_{W_1^*(\ell)}) \vee \sigma(\Gamma_{U \setminus W_2^*(\ell)})$, and the latter three are conditionally independent given~$\sigma(U^*)$. Therefore for each choice of $Q_1,\ldots,Q_n \in \cQ$ with $Q_i \subseteq U$, we have that $(\Gamma_{Q_i}\1_{\dist(Q_i, \partial U^*) > 2^{-\ell}})_{i=1,\ldots,n}$ is conditionally independent of~$\cG$. Since this holds for each $\ell \in \N$, the claim follows.
    
Finally, let us address the case when $\cL \subseteq U$. Let $D_1$ be the collection of bounded Euclidean connected components of $\C\setminus\cL$, and let $D_2$ be the unbounded component. For $j=1,2$, let $\Gamma_{D_j}$ denote the collection of loops of $\Gamma$ in $D_j$. It suffices to argue that for each $E \in \cG$ and each $j=1,2$ we have $\p[E \mid \cL,\, \Gamma_{D_{3-j}}] \in \{0,1\}$ almost surely. For $j=1$, since the conditional law of $\Gamma_{D_1}$ is just a conditionally independent collection of CLE$_6$, the claim follows already from the previous paragraph. For $j=2$, we use another change of measure where we sample a ball $B(z,r)$ according to some probability measure (independently of $\Gamma$) and consider the innermost loop~$\widetilde{\cL}$ surrounding~$B(z,r)$. Then, given $\cL \subseteq U$, the probability of selecting $B(z,r)$ so that $\widetilde{\cL} = \cL$ is positive. Now, the conditional law of $\Gamma$ outside $\widetilde{\cL}$ is exactly a CLE$_6$, hence a similar argument as before completes the proof of this case.
\end{proof}

\begin{lem}\label{lem:measurability-outside}
For each $\ell\in\N$, the collection $(\cX\cap\overline{V},\fd_V)_{V\in\fC_U}$ is measurable with respect to the random variable $(U^*,\Gamma_{W_2^*(\ell)},(d_Q)_{Q\Subset W_2^*(\ell),Q\in\cQ})$.
\end{lem}
\begin{proof}
Note that $\cX \cap U$ consists of the points in $\overline{U^*}$ that are connected to $\partial U^*$ by an open path. Also, the loops of $\Gamma_{U^*}$ in the next level inside $\cL$ are exactly the ones that are connected to $\partial U^*$ by an open path in $\overline{U^*}$. 
These loops together with $\partial U^*$ determine the collection $\fC_U$. Therefore the collection $(\cX\cap \overline{V})_{V\in\fC_U}$ is measurable with respect to $U^*$ and $\Gamma_{U^*}$, which is measurable with respect to $U^*$ and $\Gamma_{W_2^*(\ell)}$. Furthermore, by~\eqref{eq:def-dU} the random variable $\fd_V$ is measurable with respect to $(d_Q)_{Q\Subset W_2^*(\ell),Q\in\cQ}$, which completes the proof.
\end{proof}

We begin by noting that the metrics $d_Q$ on the full CLE$_6$ configurations $\Gamma_Q$ are exactly spatially independent (unlike their restrictions to the cluster $\cX$ which introduces a bias).

\begin{lem}\label{lem:indep-metrics-cle}
Let $Q_1,\ldots,Q_k\in\cQ$ and define $W:=\bigcup_{i=1}^kQ_i$. Let $Q'_1,\ldots,Q'_{k'} \Subset W^c$. Then the random variables $d_{Q_1},\ldots,d_{Q_k}$ and $\Gamma_{W}$ are independent of $\Gamma_{W^c}$ and $d_{Q'_1},\ldots,d_{Q'_{k'}}$.
\end{lem}
\begin{proof}
The percolation configurations in $W$ and $W^c$ are independent, so the joint law of the random variables $(\Gamma_{W,n},d_{Q_1,n},\ldots,d_{Q_k,n})$ and $(\Gamma_{W^c,n},d_{Q'_1,n},\ldots,d_{Q'_k,n})$ is given by some product measure $\mu_n\otimes\nu_n$. Without loss of generality and using Proposition~\ref{prop:GHf-tightness-dU} and Theorem~\ref{thm:convergence-cle-percolation-full}, 
we can choose a subsequence $J$ such that along $J$ we have that $\mu_n$ and $\nu_n$ converge weakly to the laws of $(\Gamma_{W},d_{Q_1},\ldots,d_{Q_k})$ and $(\Gamma_{W^c},d_{Q'_1},\ldots,d_{Q'_k})$ respectively in the product topology of the topologies induced by the quad-crossing topology and $\Delta_{\GHf}$ (recall Sections~\ref{subsubsec:convergence-boundary-loops} and~\ref{sec:tightness-sketch} for the precise definition of the topologies). 
Since these topologies are separable, it follows that $\mu_n\otimes\nu_n$ converges to a product measure (see, e.g.\ \cite[Theorem 2.8]{billingsley1999convergence}), which completes the proof.
\end{proof}

\begin{lem}\label{lem:cond-independence-metrics-cle}
Let $W_1\subset W_2$ be finite unions of sets in $\cQ$. Let $Q_1,\ldots,Q_k,Q'_1,\ldots,Q'_{k'}\in\cQ$ such that $Q_i\subset W_2$ and $Q'_i\subset W_2^c$. Then $(\Gamma_{W_2},(d_{Q_i})_i)$ and $(\Gamma_{W_1^c},(d_{Q'_i})_i)$ are conditionally independent given~$\Gamma_{W_2\setminus W_1}$.
\end{lem}
\begin{proof}
By Lemma~\ref{lem:measurability-cle}, $\Gamma_{W_1^c}$ is measurable with respect to $(\Gamma_{W_2^c},\Gamma_{W_2\setminus W_1})$. Moreover, by Lemma~\ref{lem:indep-metrics-cle}, $(\Gamma_{W_2^c},(d_{Q'_i})_i)$ is independent of $(\Gamma_{W_2},(d_{Q_i})_i)$, and since $\Gamma_{W_2\setminus W_1}$ is determined by $\Gamma_{W_2}$, they are also conditionally independent given~$\Gamma_{W_2\setminus W_1}$. Combining these facts gives the result.
\end{proof}

We are now ready to complete the proof of Proposition~\ref{prop:markov}. 

\begin{proof}[Proof of Proposition~\ref{prop:markov}]

We can and will assume that $U \in \cQ$ since the collection $(\cX\cap \overline{V},\fd_{V})_{V\in\fC_{U}}$ is determined by $\{(\cX\cap \overline{V},\fd_{V})_{V\in\fC_{U_1}}\}_{U_1 \in \cQ,\,U_1 \Subset U}$. Indeed, upon showing this, for general $U$ we have that $(U_1^*)_{U_1 \in \cQ,\, U_1 \Subset U}$ is a function of $U^*$ and $\Gamma_{U^*}$, and therefore conditionally independent of $\Gamma \setminus \Gamma_{U^*}$. This implies also that the metrics outside $U$ are independent of $\Gamma_{U^*}$ given $U^*$. Therefore $(U_1^*)_{U_1 \in \cQ,\, U_1 \Subset U}$ and hence also the collection of metrics in all the $U_1$ are conditionally independent of the metrics outside $U$ given $U^*$, concluding the proof in the general case.

Note that the collection $(\cX\cap \overline{V'},\fd_{V'})_{V'\in\fC_{U^c}}$ is determined by the countable collection of random variables $(\cX\cap Q',d_{Q'}|_{\cX})_{Q'\in\cQ,Q'\cap \overline{U}=\varnothing}$. This is generated by the finite collections $(\cX\cap Q_i',d_{Q_i'}|_{\cX})_{i=1}^k$ for $Q'_i\cap \overline{U}=\varnothing$ and is intersection stable. Thus, it suffices to show that $(\cX\cap \overline{V},\fd_V)_{V\in\fC_U}$ and $(\cX\cap Q_i',d_{Q_i'}|_{\cX})_{i=1}^k$ are conditionally independent given $U^*$ (note that the collection of loops $\{\cL\}\cup\Gamma_{\cL}\setminus\Gamma_{U^*}$ is already determined by $(U^*,\cX \setminus \overline{U})$).

Let $\ell$ be large enough such that $\dist(Q'_i,U)>2^{-\ell}$ for all $1\leq i\leq k$. We first show that $\cX\cap Q'_i$ is measurable with respect to $\Gamma_{W_1^*(\ell)^c}$. Note that $x\in Q'_i$ is in $\cX$ if and only if $x$ is inside $\cL$ and $x$ is not surrounded by any loop in $\Gamma_{\cL}$. Note that the random variable $\cL$ is measurable with respect to $\Gamma_{W_1^*(\ell)^c}$. Indeed, since $U^*$ is surrounded by $\cL$, we must have $\cL\subset W_1^*(\ell)^c$. Then $\cL$ is simply the outermost loop in $\Gamma_{W_1^*(\ell)^c}$ that surrounds $(\frac{1}{2},\frac{1}{2})$. Furthermore, if $\gamma$ surrounds $x$, then $\gamma$ is a loop that is not fully contained in $U$, so $\gamma\subset W_1^*(\ell)^c$. We conclude that the set of loops surrounding $x$ is in $\Gamma_{W_1^*(\ell)^c}$, which suffices to show that $\cX\cap Q'_i$ is measurable with respect to $\Gamma_{W_1^*(\ell)^c}$.

By Lemma~\ref{lem:measurability-outside}, $(\cX\cap \overline{V},\fd_V)_{V\in\fC_U}$ is measurable with respect to $(U^*,\Gamma_{W_2^*(\ell)},d_Q\1_{Q\Subset W_2^*(\ell)})_{Q\in\cQ}$. This is generated by the finite collections $(d_{Q_i}\1_{Q_i\Subset W_2^*(\ell)})_{i=1}^k$ for $Q_i\in\cQ$ such that $Q_i \Subset U$ and is intersection stable.

Let $I^*$ be the set of $i\in\{1,\ldots,k\}$ for which $Q_i\Subset W_2^*(\ell)$. We claim that $(U^*,\Gamma_{W_2^*(\ell)},(d_{Q_i})_{i\in I^*})$ is conditionally independent of $(\Gamma_{W_1^*(\ell)^c},(d_{Q'_i})_i)$ given~$(\Gamma_{W_2^*(\ell)\setminus W_1^*(\ell)},U^*)$. Note that there are finitely many possible choices of $W_1^*(\ell),W_2^*(\ell)$. Let $W_1,W_2 \in \cQ_\ell$ be a possible choice. By Lemma~\ref{lem:cond-independence-metrics-cle}, the random variables $(\Gamma_{W_2},(d_{Q})_{Q \subset W_2})$ and $(\Gamma_{W_1^c},(d_{Q'_i})_i)$ are conditionally independent given~$\Gamma_{W_2\setminus W_1}$. Note that the event $\{W_2^*(\ell)\setminus W_1^*(\ell) = W_2\setminus W_1\}$ is measurable with respect to $\Gamma_{W_1^c}$, and on this event $U^*$ is measurable with respect to $\Gamma_{W_1^*(\ell)^c} =\Gamma_{W_1^c}$ and $I^*$ is measurable with respect to $U^*$. It follows that on the event $\{W_2^*(\ell)\setminus W_1^*(\ell) = W_2\setminus W_1\}$, the random variables $(\cX\cap \overline{V},\fd_V)_{V\in\fC_U}$ and $(\Gamma_{W_1^*(\ell)^c},(d_{Q'_i})_i)$ are conditionally independent given~$(\Gamma_{W_2^*(\ell)\setminus W_1^*(\ell)},U^*)$. This holds for each possible choice of $W_1,W_2 \in \cQ_\ell$, implying the claim.

Now let $X=(\cX\cap\overline{V},\fd_V)_{V\in\fC_U}$, $Y=(\cX\cap Q'_i,d_{Q'_i}|_{\cX})_i)$. By the above, $X$ and $Y$ are conditionally independent given $\cG_\ell$. Let $f,g$ be bounded measurable functions of $X$ and $Y$ respectively. By conditional independence of $X$ and $Y$ given $\cG_\ell$, we have
\[
    \E[f(X)g(Y)\mid\cG_\ell] = \E[f(X)\mid\cG_\ell]\,\E[g(Y)\mid\cG_\ell] .
\]

We have that $(\cG_\ell)_{\ell\in\N}$ is a reverse filtration since $\Gamma_A$ is measurable with respect to $\Gamma_B$ for $A\subset B$. Further, by Lemma~\ref{le:ustar_infinitesimal_sigma_algebra}, we have $\bigcap_{\ell\in\N}\cG_\ell=\sigma(U^*)$. Then by L\'evy's reverse martingale convergence theorem, we get
\[\begin{split}
\E[f(X)g(Y)\mid U^*] 
&= \lim_{\ell\to\infty}\E[f(X)g(Y)\mid\cG_\ell] \\
&= \lim_{\ell\to\infty}\E[f(X)\mid\cG_\ell]\,\E[g(Y)\mid\cG_\ell] \\
&= \E[f(X)\mid U^*]\,\E[g(Y)\mid U^*]
\end{split}\]
which was what we wanted to show.
\end{proof}

\subsection{Translation invariance}\label{sec:translation-invariance}
\begin{prop}\label{prop:translation-invariance}
Let $U\subset[0,1]^2$ be open and simply connected. Let $D_1$ be the support of $(\cX\cap \overline{V},\fd_V)_{V\in\fC_U}$ and $D_2$ be the support of $U^*$. Let $\mu_U(\cdot\mid\cdot)$ be a probability kernel from $D_2$ to $D_1$ such that $\mu_U(\cdot\mid U^*)$ is a.s.\ equal to the conditional law of $(\cX\cap \overline{V},\fd_V)_{V\in\fC_U}$ given $U^*$. Then for each $z \in \C$ such that $U+z\subset[0,1]^2$ one can choose versions of $\mu_U$ and $\mu_{U+z}$ such that for all events $A$ in $D_1$ and $u \in D_2$ we have
\begin{equation*}
    \mu_U(A\mid u)=\mu_{U+z}(A+z\mid u+z).
\end{equation*}
where $A+z$ denotes the event that the translated metrics $(\cX\cap \overline{V}-z,\fd_{V}(\cdot+z,\cdot+z))_{V\in\fC_{U+z}}$ are in the event $A$.
\end{prop}

As in Section~\ref{sec:markov}, we will make use of the $\dg$ resp.\ $\dr$ metric defined on all clusters. This has the advantage that their local laws are exactly translation invariant which we argue in Lemma~\ref{lem:translation-invariance-unconditional} below.

Recall from Section~\ref{subsec:cle-metrics} that $\fC_U$ contains only regions whose boundary loops are in $\{\cL\} \cup \Gamma_\cL$ and that the metric $\fd_V$ is only defined on the set $\cX\cap\overline{V}$. For the following, it will be useful to consider the set $\fC'_U$ of regions that is defined like $\fC_U$, except that the boundary loops $\{\gamma_1,\ldots,\gamma_m\}$ of $V\in\fC'_U$ can be any set in $\Gamma$. That is, the regions $V \in \fC'_U$ can belong to any of the clusters of the nested $\CLE_6$ $\Gamma$ and not just $\cX$. 
Furthermore, we extend the definition of $\fd_V$ to each $V \in \fC'_U$ by
\[
    \fd_V(x,y)=\sup_{\substack{Q\Supset V\\Q\in\cQ}}d_Q(x,y), \quad x,y \in \overline{V} .
\]

Below, we will write $\fd-z$ to denote the metric $\fd(\cdot+z,\cdot+z)$.

\begin{lem}\label{lem:translation-invariance-unconditional}
Let $U' \in \cQ$ and $z \in \R^2$ be such that $U'+z \subseteq [0,1]^2.$ Then
\begin{equation}
    (\Gamma_{U'},(\fd_{V})_{V\in\fC'_{U'}})\overset{d}{=}(\Gamma_{U'+z},(\fd_{V})_{V\in\fC'_{U'+z}})-z.
\end{equation}
\end{lem}
This result is intuitively obvious, since the percolation configuration itself is translation invariant. However, there are two subtleties. First of all, $\Lambda_{U',n}$ and $\Lambda_{U'+z,n}$ may actually not be the same graphs, but instead instead differ by one row or column. This extra row or column does not matter for the $\CLE$ convergence: both configurations can be coupled to converge to the same limit because the boundary $3$-arm exponent is larger than $1$. Furthermore, for a $Q\in\cQ$ with $Q\Subset U'$, it is actually not true that $Q+z\in\cQ$ if $z$ is not a dyadic number. However, we can still approximate $V\in\fC_{U'}$ and $V+z\in\fC_{U'+z}$ by $Q_1\in\cQ$, $Q_1\Subset U'$ and $Q_2\in\cQ$, $Q_2\Subset U'+z$ in such a way that the limiting metrics $\fd_V$ and $\fd_{V+z}$ satisfy $\fd_{V}=\fd_{V+z}-z$.
 
\begin{proof}[Proof of Lemma~\ref{lem:translation-invariance-unconditional}]
Let $z\in[0,1]^2$. Note that $\Lambda_{U',n}$ and $\Lambda_{U'+z,n}$ are polygon-shaped grids whose side lengths differ by at most 1 vertex. Couple them such that the percolation configurations agree except for possibly the rightmost or topmost strip, where percolation is sampled i.i.d.\ Then $\Gamma_{U',n}$ and $\Gamma_{U'+z,n}$ can be chosen such that they converge almost surely to $\Gamma_{U'}$ and $\Gamma_{U'+z}$, with $\Gamma_{U'}=\Gamma_{U'+z}-z$ because there are almost surely no boundary pivotals since the half-plane $3$-arm exponent is larger than $1$ (cf.\ Lemma~\ref{lem:arm_exponents} and the proof of Lemma~\ref{lem:touching-loops}). 

Now let $V\in\fC'_{U'}$. Then $V+z\in\fC'_{U'+z}$. Let $Q\in\cQ$ with $U'\Supset Q\Supset V$, and let $Q_{\pm}\in\cQ$ such that $U'+z\Supset Q_{+}\supset Q+z\supset Q_{-}\Supset V+z$. Note that if $Q+z\in\cQ$, we can simply choose $Q_{\pm}=Q+z$, but this is not guaranteed to be the case if $z$ is not dyadic. Now, let $x,y\in V$. Then there exist $x_n,y_n\in\Lambda_{U',n}$ and $(x+z)_n,(y+z)_n\in\Lambda_{U'+z,n}$ such that $x_n,y_n,(x+z)_n,(y+z)_n$ converge to $x,y\in\Lambda_{U'}$ and $x+z,y+z\in\Lambda_{U'+z}$ respectively. Note that $x_n+z$ is not necessarily in the grid $\Lambda_{U'+z,n}$, so we choose $(x+z)_n$ to be the vertex in $\Lambda_{U'+z,n}$ that is coupled with $x_n$, and similarly for $y_n$ (recall that every vertex in the interior of $\Lambda_{U'+z,n}$ (resp.\ $\Lambda_{U',n}$) is coupled with some vertex in $\Lambda_{U',n}$ (resp.\ $\Lambda_{U'+z,n}$)). 
Then by monotonicity we get that
\begin{equation*}
    d_{Q_+,n}((x+z)_n,(y+z)_n)\leq d_{Q,n}(x_n,y_n)\leq d_{Q_-,n}((x+z)_n,(y+z)_n).
\end{equation*}
Since $Q_{\pm}\in\cQ$, we can take a subsequential limit and obtain
\begin{equation*}
    d_{Q_+}(x+z,y+z)\leq d_{Q}(x,y)\leq d_{Q_-}(x+z,y+z).
\end{equation*}
Taking $Q\downarrow V$, we can also let $Q_{\pm}\downarrow V+z$, which shows that
\begin{equation*}
    \fd_V(x,y)=\fd_{V+z}(x+z,y+z),
\end{equation*}
which completes the proof.
\end{proof}

We are now ready to prove that the conditional law of $(\cX\cap\overline{V},\fd_V)_{V\in\fC_U}$ is translation invariant. For the proof, we first note that $(\cX\cap \overline{V},\fd_V)_{V\in\fC_U}$ is a translation invariant measurable function of $(\Gamma_{U'},(\fd_{V})_{V\in\fC'_U})$ and $U^*$. The former is translation invariant by Lemma~\ref{lem:translation-invariance-unconditional}. To deal with the conditioning on $U^*$, we first define a selection procedure of a local cluster 
that depends only on $\Gamma_{U'}$ and is translation-invariant, and then use that the conditional law given $U^*$ is the same as the conditional given $\widetilde{U}^*$ on the positive probability event that $\widetilde{U}^*=U^*$.

\begin{proof}[Proof of Proposition~\ref{prop:translation-invariance}]
First, we note that it suffices to assume that $\dist(U,\partial[0,1]^2) \ge \varepsilon$ and $\dist(U+z,\partial[0,1]^2) \ge \varepsilon$ for some $\varepsilon > 0$. Indeed, this is because the collection $(\cX\cap \overline{V},\fd_{V})_{V\in\fC_{U}}$ is determined by $\{(\cX\cap \overline{V},\fd_{V})_{V\in\fC_{U_1}}\}_{U_1 \in \cQ,\,U_1 \Subset U}$, and the conditional law of $(U_1^*)_{U_1 \in \cQ,\,U_1 \Subset U}$ given~$U^*$ is translation-invariant.

Let $U' \in \cQ$ be such that $U \subseteq U'$ and $U'+z \subseteq [0,1]^2$, and let $A' = \{u \in [0,1]^2 \mid 0<\dist(u,U') < \varepsilon\}$. Let $\Gamma_{U'},\Gamma_{A'}$ be as defined in Section~\ref{subsubsec:convergence-boundary-loops}. By Lemma~\ref{lem:translation-invariance-unconditional}, we have that the law of $(\Gamma_{A'},\Gamma_{U'},(\fd_{V})_{V\in\fC'_U})$ agrees with the law of the translation of $(\Gamma_{A'+z},\Gamma_{U'+z},(\fd_{V})_{V\in\fC'_{U+z}})$. 
Since CLE is locally finite, $\Gamma_{A'}$ a.s.\ contains only finitely many strands that cross $A'$. These strands divide the annulus $A'$ into a finite number of components which are alternatingly open and closed (in the percolation sense). Let us select a random collection $\cI$ of open components in a translation-invariant manner (e.g.\ by flipping a coin for each open component). We let $\widetilde{U}^*$ be defined analogously to $U^*$ but for the clusters that are connected to $\cI$ inside $U'$, i.e.\ $\widetilde{U}^*$ is the region containing the loops of $\Gamma$ that are entirely contained in $U$ and are connected to $\cI$ by an open path within $U'$. (In the case $\cI=\varnothing$, we can e.g.\ take instead the open cluster inside $U'$ whose outer boundary has the smallest Hausdorff distance to $\partial U'$.) Further, we let $\widetilde{\fC}_U \subseteq \fC'_U$ be the regions contained in $\widetilde{U}^*$ where $\fC'_U$ is defined above Lemma~\ref{lem:translation-invariance-unconditional}. Also let $\widetilde{\cX}$ be the set of points inside $\overline{\widetilde{U}^*}$ that are connected to $\cI$ by an open path within $U'$. Since the function that associates to $(\Gamma_{A'},\Gamma_{U'},\cI)$ the region $\widetilde{U}^*$ is translation-invariant, we have that the law of $(\Gamma_{A'},\Gamma_{U'},\widetilde{U}^*,(\widetilde{\cX}\cap \overline{V},\fd_{V})_{V\in\widetilde{\fC}_{U}})$ is translation-invariant. Hence, its disintegration along $\widetilde{U}^*$ is also translation-invariant. See Figure \ref{fig:translation_invariance} for an illustration. 

\begin{figure}[ht]
    \centering
    \includegraphics[scale=0.8]{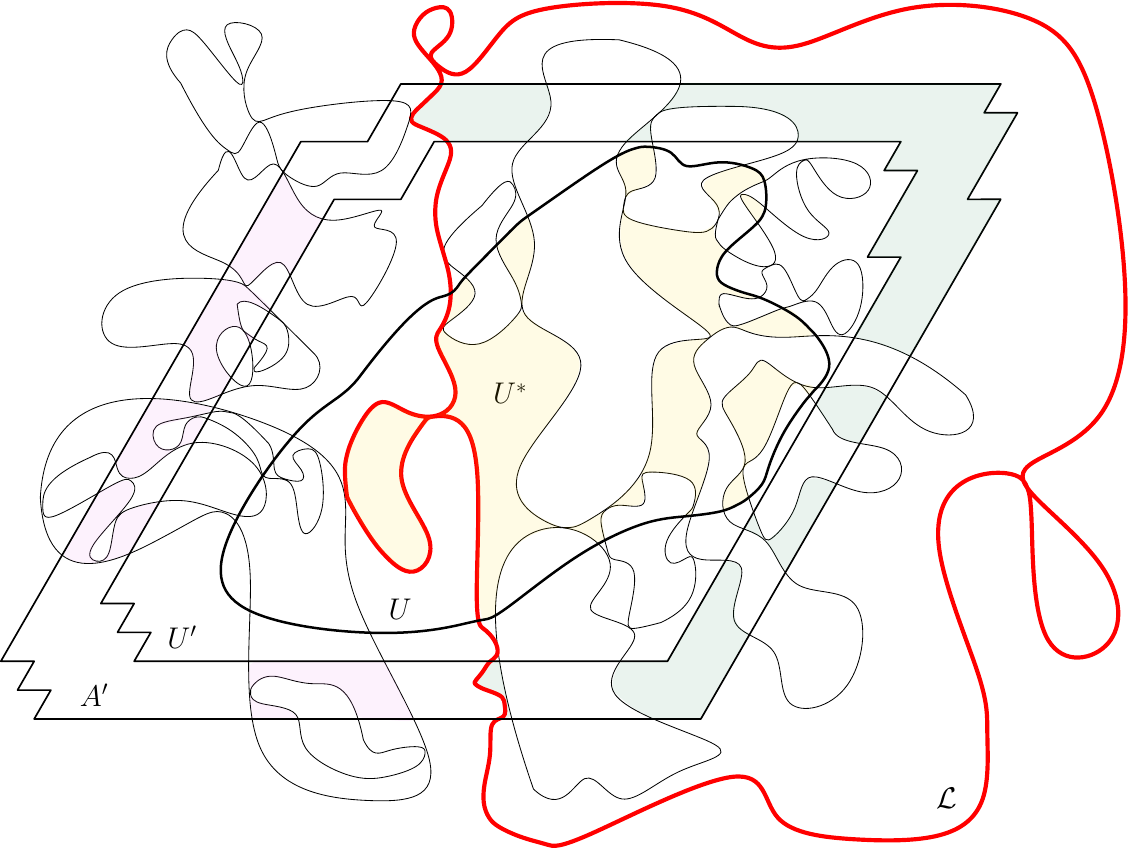}
    \caption{The figure above illustrates the construction of the translation-invariant set $\widetilde{U}^*$ that equals $U^*$ with positive probability. The figure portrays a set $U$ whose boundary is represented by the bold black curve above. The set $U'\in \cQ$ is such that $U\subseteq U'$, and $A'$ is the set of all points that are not in $U'$ but are within `triangular' $L^1$-distance $\varepsilon$ from $U'$. The outermost loop~$\cL$ surrounding the origin is represented by the red curve above. The yellow area inside the set $U$ represents the set $U^*$. The areas shaded in green and pink inside the annulus $A'$ represent the open (in a percolation sense) clusters that cross $A'$. In the case above, for the set $\widetilde{U}^*$ to be equal to $U^*$, we need to select $\cI$ to be precisely the collection of the green crossing clusters above. Note that in contrast to $U^*$, the law of~$\widetilde{U}^*$ is translation-invariant since it does not `see' which of the loops is~$\cL$.}
    \label{fig:translation_invariance}
\end{figure}

Let $F(U')$ be the event that {$\cX \cap U' \neq \varnothing$ and} $\cX$ is either not entirely contained in $U'$ or that $\cX \subseteq U'$ but is the cluster whose outer boundary has the smallest Hausdorff distance to $\partial U'$ (note that the case when $\cX \cap U' = \varnothing$ is trivial). On $F(U')$, given $U^*$ and $\Gamma \setminus \Gamma_{U^*}$, the conditional probability of the event $E = \{\widetilde{U}^* = U^*\}$ is positive. Therefore, on~$E$, the conditional law of $(\cX\cap \overline{V},\fd_{V})_{V\in\fC_{U}}$ given $\Gamma \setminus \Gamma_{U^*}$ agrees a.s.\ with the conditional law of $(\widetilde{\cX}\cap \overline{V},\fd_{V})_{V\in\widetilde{\fC}_{U}}$ given $\Gamma \setminus \Gamma_{\widetilde{U}^*}$ which is the same as the conditional law given just $\Gamma_{A' \cup U'} \setminus \Gamma_{\widetilde{U}^*}$. But since $\cI$ is sampled conditionally independently of $\Gamma \setminus \Gamma_{U^*}$, the conditional law of $(\cX\cap \overline{V},\fd_{V})_{V\in\fC_{U}}$ given $\Gamma \setminus \Gamma_{U^*}$ does not depend on $E$. By Proposition~\ref{prop:markov}, this is the same as the conditional law of $(\cX\cap \overline{V},\fd_{V})_{V\in\fC_{U}}$ given just $U^*$. Therefore, its translation invariance follows from the translation invariance of the conditional law of $(\widetilde{\cX}\cap \overline{V},\fd_{V})_{V\in\widetilde{\fC}_{U}}$ given $\Gamma_{A' \cup U'} \setminus \Gamma_{\widetilde{U}^*}$. This shows the result on the event $F(U')$. The general case follows since the collection $(\cX\cap \overline{V},\fd_{V})_{V\in\fC_{U}}$ is determined by $\{(\cX\cap \overline{V},\fd_{V})_{V\in\fC_{U_1}} \1_{F(U_1)} \}_{U_1 \in \cQ,\,U_1 \Subset U}$, and the conditional law of $(U_1^*)_{U_1 \in \cQ,\,U_1 \Subset U}$ given~$U^*$ is translation-invariant.
\end{proof}

\subsection{Non-degeneracy}

We have now verified that $(\cX \cap \overline{V}, \fd_V)_{V \in \fC}$ is a CLE$_6$ metric in the sense defined in Section~\ref{subsec:cle-metrics}. Note that the definition of a CLE$_6$ metric does not rule out the case that $\fd_V=0$ for each $V$. We now show that this is not the case, and in fact each $\fd_V$ is a true metric defined on $\cX \cap \overline{V}$.

\begin{lem}\label{le:nondegeneracy}
    Let $(\Lambda_{Q},\dpt_Q,d_Q,\pi)_{Q\in\cQ}$ be a subsequential limit of the $(\Lambda_{Q,n},\dpt_{Q,n},d_{Q,n},\pi_n)_{Q\in\cQ}$, and let $(\cX \cap \overline{V}, \fd_V)_{V \in \fC}$ be defined as in~\eqref{eq:def-dU}. We almost surely have $\fd_V(x,y) > 0$ for each $V \in \fC$ and $x,y \in \cX \cap \overline{V}$ with $x \neq y$.
\end{lem}

\begin{proof}
    By Proposition~\ref{prop:positivity}, we only need to rule out that $\fd_V=0$ for each $V \in \fC$. Suppose that this is the case. By a local absolute continuity argument (see \cite[Section~7.2]{my2025metric}), this implies also that the global metric $d$ satisfies $d=0$ a.s. So we just need to rule out this case.
    
    Let $\cE(\Lambda)$ be the analogous event to $\cE(\Lambda_n)$, but for the $\CLE_6$ configuration rather than percolation (see Section~\ref{subsubsec:apriori-sketch}). It is straightforward to show that this has positive probability. 
    Assume this event holds. Then by Lemma~\ref{lem:touching-loops}, under a suitable coupling, almost surely there exists $n_0$ such that for all $n\geq n_0$, the event $\cE(\Lambda_n)$ holds. Let $u_n$ and $v_n$ be the left- and rightmost touching points as described in Section~\ref{subsubsec:apriori-sketch}. Then, again by Lemma~\ref{lem:touching-loops}, $u_n$ and $v_n$ converge in $d$ to the left- and rightmost touching points of the closed touching $\CLE_6$ clusters that are connected to the top and bottom sides, which we call $u$ and $v$. By definition of the constants $\median_n$, there exists a constant $p\in(0,1)$ such that $1$ is the conditional $p$-th quantile of $d_n(u_n,v_n)$. Since $d(u,v)=\lim_{n\in J}d_n(u_n,v_n)$, we have that $1$ is also the conditional $p$-th quantile of $d(u,v)$. In particular, the latter is not identically $0$.
\end{proof}

\subsection{Geodesic metric}\label{sec:geodesic-cle}

We have shown that $(\fd_V)_{V\in\fC}$ is a non-degenerate $\CLE$ metric. We now show that if $d=\dg$, then it is a geodesic $\CLE$ metric. We first show in Lemma~\ref{lem:geo} that $\fd_V$ is a geodesic metric in the usual sense of the word, i.e., the distance between two points equals the $\fd_V$-length of the shortest path. We then show in Proposition~\ref{prop:geo-cle} that the $\fd_V$-length in fact equals the $d$-length. Recall that $d=d_{[0,1]^2}$ simply denotes the subsequential limit of the global metrics $d_n$. Also recall the notation $P(x,y;\overline{V}) = P(x,y;\overline{V};\Gamma)$ for the set of admissible paths in $\overline{V}$ from $x$ to $y$ (equivalently, continuous paths in $(\overline{V},\dpt)$) introduced above~\eqref{eq:def-dptU}.

\begin{lem}\label{lem:geo}
Let $d=\dg$. For each $V\in\fC$, the metric $(\cX\cap\overline{V},\fd_{V})$ is a geodesic metric space, i.e., for all $x,y\in\cX\cap\overline{V}$,
\begin{equation}\label{eq:d_v-length-metric}
    \fd_V(x,y)=\inf_{\gamma\in P(x,y;\overline{V})}L_{\fd_V}(\gamma),
\end{equation}
and the infimum above is attained by some path $\gamma\in  P(x,y;\overline{V})$.
\end{lem}
\begin{proof}
Throughout this proof, we write $d=\dg$. Let $Q\in\cQ$. We say a metric is a length metric if~\eqref{eq:d_v-length-metric} holds but if the infimum is not necessarily attained for all pairs of points in the space. Then, since the GH limit of length metrics is a length metric, it follows that $(\cX\cap Q,d_Q)$ is a length metric space. 
By Proposition~\ref{prop:GHf-tightness-dU}, Lemma~\ref{le:nondegeneracy}, and the compactness of $(\cX \cap \overline{V}, \dpt)$ by \cite[Lemma~1.10]{amy2025tightness}, we have that $\dg$ is topologically equivalent to $\dpt$, 
hence the continuous paths in $\dg$ are the same as the continuous paths in~$\dpt$.
Additionally, we have argued in the proof of Proposition~\ref{prop:compatibility} that for $Q\in\cQ$ with $Q\Supset V$ in a small enough neighborhood of $V$, all simple admissible paths in $Q$ from $x\in\overline{V}$ to $y\in\overline{V}$ must lie entirely inside $\overline{V}$. Thus, by the exact same proof as in Proposition~\ref{prop:compatibility} we have that $d_Q(x,y)=\fd_V(x,y)$. It follows that $(\cX\cap\overline{V},\fd_V)$ is also a compact space and
\begin{equation}\label{eq:d_Q-length-metric}
    \fd_V(x,y)=d_Q(x,y)=\inf_{\gamma\in P(x,y;Q)}L_{d_Q}(\gamma)=\inf_{\gamma\in P(x,y;\overline{V})}L_{\fd_V}(\gamma).
\end{equation}
The Hopf-Rinow theorem then implies that it is in fact a geodesic metric space.
\end{proof}

\begin{prop}\label{prop:geo-cle}
Let $d=\dg$. 
For every $V\in\fC$ and $x,y\in \cX\cap\overline{V}$, we have
\begin{equation}
    \fd_V(x,y)=\inf_{\gamma\in P(x,y;\overline{V})}L_{d}(\gamma).
\end{equation}
\end{prop}
\begin{proof}
Let $x,y\in \cX\cap\overline{V}$. By Lemma~\ref{lem:geo}, it suffices to show that for every $\gamma\in P(x,y;\overline{V})$, $L_{\fd_V}(\gamma)=L_{d_{[0,1]^2}}(\gamma)$. In particular, by~\eqref{eq:d_Q-length-metric}, it suffices to show that for every $Q\Supset V$, $Q\in\cQ$, we have $L_{d_Q}(\gamma)=L_{d_{[0,1]^2}}(\gamma)$.

We have shown in Lemma~\ref{le:nondegeneracy} that $d_{[0,1]^2}\big|_{\cX}=\fd_D$ is a true metric on $\cX$. Therefore $\delta:=d_{[0,1]^2}(\gamma,\partial Q)>0$. 
Let $t_0,\ldots,t_M$ such that $d_{[0,1]^2}(\gamma(t_{i-1}),\gamma(t_{i}))<\delta/2$ for all $1\leq i<M$. Let $x^0_n,\ldots,x^M_n\in\Lambda_{n}$ that such that $(x^i_n)_{n\in\N}$ converges to $\gamma(t_i)$. Then for $n$ sufficiently large, we have $d_{Q,n}(x^{i-1}_n,x^i_n)<\delta/2$ and $d_{Q,n}(x_n^i,\partial\Lambda_{Q,n})>3\delta/4$. This means the graph-length minimizing path between $x^{i-1}_n$ and $x_n^{i}$ lies entirely inside $Q$ and in particular $d_{Q,n}(x_n^{i-1},x_n^{i})=d_n(x_n^i,x_n^{i-1})$. Thus,
\begin{equation}
    \sum_{i=1}^Md_{Q}(\gamma(t_{i-1}),\gamma(t_{i}))=\lim_{n\rightarrow\infty}\sum_{i=1}^Md_{Q,n}(x^{i-1}_n,x^i_n)=\lim_{n\rightarrow\infty}\sum_{i=1}^Md_{n}(x^{i-1}_n,x^i_n)=\sum_{i=1}^Md_{[0,1]^2}(\gamma(t_{i-1}),\gamma(t_{i})).
\end{equation}
This implies that $L_{d_Q}(\gamma)=L_{d_{[0,1]^2}}(\gamma)$, which completes the proof. 
\end{proof}

\subsection{Resistance metric}\label{sec:edge-weight-compatibility}
We show that if $d=\dr$, then $(\fd_V)_{V\in\fC}$ is a $\CLE_6$ resistance metric. We first show in Proposition~\ref{prop:res} that $(\cX\cap\overline V,\fd_V)$ is a resistance metric space. We then prove the edge weight locality property in Proposition~\ref{prop:edge-weight-locality}.
\begin{prop}\label{prop:res}
Let $d=\dr$. For each $V\in\fC$, the metric $(\cX\cap\overline{V},\fd_V)$ is a resistance metric space.
\end{prop}
\begin{proof}
Throughout this proof, we write $d=\dr$. Let $V\in\fC$ and $Q\in\cQ$ with $Q\Supset V$ be small enough such that for all $x,y\in\overline{V}$, each simple admissible path in $Q$ from $x$ to $y$ must lie entirely inside $\overline{V}$, as in the previous proof. Then $(Q,d_Q)$ is a resistance metric space, since the GH limit of resistance metric spaces is a resistance metric space \cite{kasue2010convergence}. Furthermore, the exact same proof as in Proposition~\ref{prop:compatibility} shows that $\fd_V(x,y)=d_Q(x,y)$ on $\cX\cap\overline{V}$, and by Lemma~\ref{le:nondegeneracy} we have that $\fd_V$ is indeed a metric on $\cX\cap\overline{V}$.
So $(\cX\cap\overline{V},\fd_V)$ is a resistance metric space.
\end{proof}

\begin{prop}\label{prop:edge-weight-locality}
Let $d=\dr$. Suppose that $V,V_1,V_2\in\fC$ are such that $\cX\cap\overline{V} = (\cX\cap\overline{V_1}) \cup (\cX\cap\overline{V_2})$ and $(\cX\cap\overline{V_1}) \cap (\cX\cap\overline{V_2}) = \{a,b\}$. There exists a countable set $\cZ_V=\{z_1,z_2,\ldots\} \subseteq \cX\cap\overline{V}$ that contains $a,b$, is dense with respect to $\dpt_{\overline{V}}$, and such that the following holds for sufficiently large~$m$. Let $A = \{z_1,\ldots,z_m\}$, $A_1 = A \cap \overline{V_1}$, $A_2 = A \cap \overline{V_2}$, and let $w$ (resp.\ $w_1$, $w_2$) be the weight function associated with $\fd_V(\cdot,\cdot)|_{A}$ (resp.\ $\fd_{V_1}(\cdot,\cdot)|_{A_1}$, $\fd_{V_2}(\cdot,\cdot)|_{A_2}$). Then
\begin{itemize}
    \item $w(x,y)=w_1(x,y)$ for each $x,y\in A_1$;
    \item $w(x,y)=w_2(x,y)$ for each $x,y\in A_2$;
    \item $w(x,y)=0$ for each $x\in A_1 \setminus \{a,b\}$, $y\in A_2 \setminus \{a,b\}$;
    \item $w(a,b)=0$.
\end{itemize}
\end{prop}

\begin{proof}
Recall that $(\overline{V},\dpt_{\overline{V}})$ is compact and therefore separable. Let $\gamma,\gamma'$ be loops inside $\overline{V}$. Then the sets of touching points $\gamma\cap \gamma'$ and $\gamma \cap \partial V$ are separable. Further, for each $\gamma$, the set of self-intersection points is separable. Let $\cZ_V$ be any countable set containing a dense subset of $\gamma\cap\gamma'$ and of $\gamma \cap \partial V$ for each such $\gamma\neq\gamma'$ and a dense subset of the set of self-intersection points of each~$\gamma$. We show that this set satisfies the properties stated in the lemma.

Let $Q \in \cQ$ with $Q \Supset V$ be small enough so that $\fd=\dr_Q$ on $\cX\cap\overline{V}$, as in the previous proof. Recall that if $A \subseteq \cX \cap Q$ is a finite set and $w$ is the weight function associated with~$\dr_Q|_{A}$, then $w(x,y) = \lim_n w_n(x_n,y_n)$ where $x_n,y_n$ converge respectively to $x,y$ in the coupling described at the beginning of the section and $w_n$ is the weight function associated with $\dr_{Q,n}|_{A_n}$ where $A_n = \{x_n \mid x \in A\} \subseteq \cX_n \cap Q$. (See the proof of \cite[Lemma~3.29]{my2025diffusion}.)

Let $V_1,V_2$ and $\{a,b\}$ be as in the lemma and let $m$ be large enough such that $\{a,b\}\subset A$ and such that $a$ and $b$ are separated by points in $A\setminus\{a,b\}$ in both $\cX \cap \overline{V_1}$ and $\cX \cap \overline{V_2}$. 
Note that $w(x,y)=0$ if $x$ and $y$ are separated by $A \setminus \{x,y\}$ (by the previous paragraph and the analogous property on graphs together with Lemma~\ref{lem:touching-loops}). 
Then the fourth point is automatically true. Furthermore, for each $x\in A_1\setminus\{a,b\}$ and $y\in A_2\setminus\{a,b\}$, the points $x$ and $y$ are separated inside $\cX\cap\overline{V}$ by $a$ and $b$. So the third point is also true. For the first and second points, note that $\cX\cap \overline{V_1}$ and $\cX\cap \overline{V_2}$ are separated in $\cX \cap \overline{V}$ by $a$ and $b$, and again apply Lemma~\ref{lem:touching-loops}. 
\end{proof}

\subsection{Proof of uniqueness}\label{sec:final-convergence-metrics}
\begin{proof}[Proof of Corollary~\ref{cor:final-convergence-metrics}]
By Proposition~\ref{prop:main-conv-dU-GH}, for each subsequence $I\subset\N$, there exists a further subsequence $J\subset I\subset\N$ such that on a suitable coupling $(\Lambda_{Q,n},d_{Q,n},\pi_n)_{n\in J}$ converges for each $Q\in\cQ$. Define the subsequential limits $(\cX\cap\overline{V},\fd_{V})_{V\in\fC}$ as in~\eqref{eq:def-dU}. Then by Theorem~\ref{th:dU-is-CLE-metric}, this collection is a $\CLE_6$ metric. Now let $(\Lambda,d)$ and $(\Lambda,\widetilde{d})$ be two subsequential limits of $(\cX_n,d_n)$ with corresponding subsequences $J$ and $\widetilde{J}$, coupled with the same $\CLE_6$. By Theorems~\ref{th:uniqueness-geodesic-cle} and~\ref{th:uniqueness-resistance-cle}, there exists a constant $c$ such that $\fd_V=c\widetilde{\fd}_V$ for each $V \in \fC$. By a local absolute continuity argument (see \cite[Section~7.2]{my2025metric}), this implies also that $d=c\widetilde{d}$. 
We have seen in the proof of Lemma~\ref{le:nondegeneracy} that $1$ is the conditional $p$-th quantile of both $d(u,v)$ and $\widetilde{d}(u,v)$ on the event $\cE(\Lambda)$. Therefore $c=1$. 
\end{proof}

\section{Measure convergence}\label{sec:ghpf-proof}

In this section, we prove the convergence of the vertex counting measures $\mu^{(m)}_n$ defined in~\eqref{eq:counting_measure_cluster} to the CLE$_6$ gasket measure $\mu^{(m)}$ on $\cX^{(m)}$, thus completing the last step of the proof of Theorem~\ref{thm:main-convergence}. Throughout this entire section, we write $\cX$ (resp.\ $\cX_n$) to denote any of the clusters $\cX^{(m)}$ (resp.\ $\cX^{(m)}_n$) in order to ease the notation, and do the same thing for $\mu$\phantomsection{}\label{def:mu} (resp.\ $\mu_n$). We first show that $(\cX_n,d_n,\mu_n)_{n\in\N}$ converges in the GHP distance to some limiting metric-measure space $(\cX,d,\widetilde{\mu})$ in Section~\ref{subsec:ghp-proof}. We then prove that $\widetilde{\mu}$ is a multiple of the CLE$_6$ gasket measure $\mu$ by showing that it satisfies~\eqref{eq:expression-mu}, and conclude the proof of Theorem~\ref{thm:main-convergence}.

\subsection{GHP convergence}\label{subsec:ghp-proof}
Recall from Section~\ref{subsubsec:measure-conv} that the measure convergence results from \cite{gps} were only stated for measures of the form $\mu^A$, where $A$ is some annulus. Therefore, we first express the desired measures on $\cX_n$ and $\cX$ in terms of annulus measures. In order to do this, we introduce the following tiling of $[0,1]^2$ by annuli.

Let $\cH_\ell$\phantomsection{}\label{def:cHl} be the set of annuli of the form $2^{-\ell}\!\left([i-1,i+2]\times[j-1,j+2]\right)\setminus([i,i+1]\times[j,j+1])$ for $0\leq i,j< 2^{\ell}$. Then $\cH_\ell$ is a tiling of $[0,1]^2$, i.e., each point in $[0,1]^2$ is contained in the inner face of exactly one annulus in $\cH_\ell$. 
Define $\mu_n^{\cH_\ell}:=\sum_{A\in\cH_\ell}\mu_n^A$\phantomsection{}\label{def:mu_ncH}, where $\mu_n^A$ is defined as in Section~\ref{subsubsec:measure-conv}. It follows from the convergence in Theorem~\ref{th:gps} that there is a coupling such that $\mu_n^{\cH_\ell}$ converges almost surely in Prokhorov distance on $[0,1]^2$ to some limiting measure $\mu^{\cH_\ell}$\phantomsection{}\label{def:mucH} that is a function of the $\CLE_6$ configuration $\Gamma$. Define the measure $\widetilde{\mu}^\ell:=\mu^{\cH_\ell}\circ\pi$\phantomsection{}\label{def:tilde-muell} on $\cX$ where $\pi\colon\cX\to[0,1]^2$ is the embedding map. (Note that this measure is well-defined since the set of double points of $\pi(\cX)$ has measure $0$ by Lemma~\ref{lem:measure_2arm}.)

Recall the definition of the measure $\mu_n=\frac{1}{n^2A(1,n)}\sum_{x\in\cX_n}\delta_x$ on $\cX_n$. Suppose that $\ell$ is large enough so that $2^{-\ell+1}<\diam(\cX_n)$. Then $\mu_n=\mu_n^{\cH_\ell}\circ\pi_n$ where $\pi_n\colon\cX_n\to[0,1]^2$ is the embedding map. Indeed, note that since the map $\pi_n$ is defined only on the cluster $\cX_n$, the measure $\mu_n^{\cH_\ell}\circ\pi_n$ counts only the points in $\cX_n$). Further, a point~$x$ in $\cX_n$ is counted by $\mu_n^{\cH_\ell}\circ\pi_n$ if the cluster $\cX_n$ crosses the annulus in $\cH_\ell$ whose inner face contains~$x$. If $\diam(\cX_n)>2^{-\ell+1}$, then every $x\in\cX_n$ is counted, so it follows that $\mu_n=\mu_n^{\cH_\ell}\circ\pi_n$.

Suppose that we are on the event that $\diam(\cX) > 2^{-\ell+1}$. Since $\cX_n$ converges to $\cX$ in Hausdorff distance, it follows that $\diam(\cX_n)>2^{-\ell+1}$ for $n$ sufficiently large almost surely. The goal of this section is to deduce the convergence of $\mu_n$ to $\widetilde{\mu}^\ell$ in the GHP distance, and in Section~\ref{sec:char-limit-measure-proof} we identify~$\widetilde{\mu}^\ell$ with the CLE$_6$ gasket measure~$\mu$ defined in \cite{miller2024existence}. Note that this implies in particular that $\widetilde{\mu}^\ell$ does not depend on $\ell$ as long as we are on the event that $\diam(\cX) > 2^{-\ell+1}$.

\begin{prop}
Let $d \in \{\dpt,\dg,\dr\}$. There exists a coupling under which, on the event $\{\diam(\cX) > 2^{-\ell+1}\}$, the sequence $(\cX_n,d_n,\mu_n)_{n\in I}$ converges to $(\cX,d,\widetilde{\mu}^\ell)$ in the GHP distance.
\end{prop}

By the weak convergence of $(\cX_n,d_n,\pi_n) \to (\cX,d,\pi)$ in the GHf topology (Corollary~\ref{cor:final-convergence-metrics}), we can further couple the $\cX_n$ and $\cX$ so that we have both $\mu_n^{\cH_\ell} \to \mu^{\cH_\ell}$ and $(\cX_n,d_n,\pi_n) \to (\cX,d,\pi)$ almost surely. In particular, there exists a metric space $(Y,d_Y)$ and isometries $\psi_n\colon\cX_n\to Y$ and $\psi\colon\cX\to Y$ such that $(\psi_n(\cX_n))_{n\in I}$ converges to $\psi(\cX)$ in Hausdorff distance. We show that the pushforward measure $(\mu_n\circ\psi^{-1}_n)_{n}$ converges to $\widetilde{\mu}^\ell\circ\psi^{-1}$ in the L\'evy-Prokhorov distance. This then implies convergence of $(\cX_n,d_n,\mu_n)_{n}$ to $(\cX,d,\widetilde{\mu}^\ell)$ in the GHP distance. The proposition follows immediately from the following.

For a set $A$ and some metric $d$, we define $A^{\varepsilon}:=\{x:d(x,A)<\varepsilon\}$ to be the \emph{$\varepsilon$-fattening} of $A$.

\begin{lem}\label{lem:prokhorov}
The following holds almost surely on the event $\{\diam(\cX)>2^{-\ell+1}\}$ under the coupling described above. Let $\varepsilon>0$. Then there exists $n_0$ such that for all $n\geq n_0$ and all measurable $A\subset Y$,
\[ 
\mu_n(\psi_n^{-1}(A))\leq\widetilde{\mu}^\ell(\psi^{-1}(A^\varepsilon))+\varepsilon
\quad\text{and}\quad
\widetilde{\mu}^\ell(\psi^{-1}(A))\leq\mu_n(\psi_n^{-1}(A^\varepsilon))+\varepsilon .
\]
\end{lem}
We first prove the following result. We point out that in the statements below, the notation $A^\delta$ refers to the $d_Y$-metric whereas $\pi(A)^\delta$ refers to the Euclidean metric.
\begin{lem}\label{lem:closeness-projections}
For every $\delta>0$, there exists $n_0=n_0(\delta)$ such that for all $n\geq n_0$ and measurable $A\subset Y$,
\begin{equation*}
    \pi_n(\psi_n^{-1}(A))\subset\pi(\psi^{-1}(A^{\delta}))^{\delta} ,\qquad \pi(\psi^{-1}(A))\subset\pi_n(\psi_n^{-1}(A^\delta))^\delta.
\end{equation*}
\end{lem}
\begin{proof}

It follows from the GHf convergence $(\cX_n,d_n,\pi_n) \to (\cX,d,\pi)$ that there exists $n_0=n_0(\delta)$ such that for all $x\in\cX_n$, there exists $x'\in\cX$ such that $d_Y(\psi_n(x),\psi(x'))<\delta$ and furthermore, $\|\pi_n(x)-\pi(x')\|_2\leq\delta$. So for $\psi_n(x)\in A$, we have $\psi(x')\in A^{\delta}$ and $\pi_n(x)\in\pi(\psi^{-1}(A^\delta))^{\delta}$, which shows the first inclusion. Furthermore, for all $x'\in\cX$, there exists $x\in\cX_n$ such that $d_Y(\psi_n(x),\psi(x'))<\delta$ and furthermore, $\|\pi_n(x)-\pi(x')\|_2\leq\delta$. So for $\psi(x')\in A$, we have $\psi_n(x)\in A^\delta$ and $\pi(x')\in\pi_n(\psi^{-1}(A^\delta))^\delta$, which shows the second inclusion.
\end{proof}

Before we prove Lemma~\ref{lem:prokhorov}, we need another auxiliary lemma. This is because we have assumed $\mu_n^{\cH_\ell} \to \mu^{\cH_\ell}$ in the Prokhorov distance with respect to the Euclidean metric, but for the statement of Lemma~\ref{lem:prokhorov} we need the upper bound in terms of the mass of $\pi(\psi^{-1}(A^\varepsilon))$ where the fattening is in the $d_Y$-metric, so its mass can be smaller than that of $\pi(\psi^{-1}(A))^\delta$. The next lemma ensures that for $\delta>0$ sufficiently small it is not much smaller.

\begin{lem}\label{lem:two-arm-vs-one-arm}
For every $\varepsilon>0$, there exists $\delta=\delta(\varepsilon)>0$ and $n_0=n_0(\varepsilon)$ such that for all measurable $A\subset Y$ and $n\geq n_0$, we have
\begin{equation*}
    \mu^{\cH_\ell}\!\left(\pi(\psi^{-1}(A))^{\delta}\setminus\pi(\psi^{-1}(A^\varepsilon))\right)<\varepsilon ,\qquad \mu_n^{\cH_\ell}\!\left(\pi_n(\psi_n^{-1}(A))^\delta\setminus\pi_n(\psi_n^{-1}(A^\varepsilon))\right)<\varepsilon
\end{equation*}
\end{lem}
\begin{proof}
We first prove the first inequality. The idea is as follows. Suppose that $x \in \pi(\psi^{-1}(A))^{\delta}$ is in the support of $\mu^{\cH_\ell}$, and let $x' \in \pi(\psi^{-1}(A))$ be the closest point in Euclidean distance. In case (the preimage of) $x$ is also close to (the preimage of) $x'$ in the $d$-distance, then it is contained in $A^\varepsilon$ when $\delta$ is small enough. 
If not, there are at least two disjoint open arms from $x,x'$, respectively, that are separated by two closed arms. However, the measure $\mu^{\cH_\ell}$ assigns zero mass to $4$-arm points.

Assume without loss of generality that $\psi^{-1}(A)$ is non-empty. Let $x\in[0,1]^2$ and let $B\in\cH_\ell$ be the annulus with inner face containing $x$. Using Proposition~\ref{prop:char-gps} it is straightforward to check that $\mu^{\cH_\ell}$ is indeed only supported on points that have an open arm to the boundary of $B$. Assume that $x$ has an open arm to the outer boundary of $B$, and $\|x-\pi(\psi^{-1}(A))\|_2<\delta$, and $x\not\in\pi(\psi^{-1}(A^\varepsilon))$.

Consider the Euclidean box $B(x;\delta)$ of radius $\delta$ around $x$. Then $B(x;\delta)$ intersects $\pi(\psi^{-1}(A))$, so $B(x;\delta)$ has at least one open arm to distance 
$\diam(\pi(\cX))/2-\delta\geq2^{-\ell-1}-\delta$ via $\cX$. Recall that $x$ was assumed to have an open arm to the boundary of $B$. Note that $x$ can be on the same cluster $\cX$ or on a different cluster with Euclidean diameter at least $2^{-\ell}$. Assume first that $x\not\in\pi(\cX)$. Then $B(x;\delta)$ has another disjoint open arm to distance $2^{-\ell}\geq 2^{-\ell-1}-\delta\geq\varepsilon-\delta$ for $\varepsilon$ sufficiently small. Assume now that $x\in\pi(\cX)$. Then since $x\not\in\pi(\psi^{-1}(A^\varepsilon))$, $x$ has another disjoint open arm that goes at least distance $\varepsilon/2$ with respect to the metric $d$. By H\"older continuity of $d$ with respect to $\dpt$, there exists some constant $c\geq1$ and exponent $\alpha>0$ such that almost surely for all $x$, this arm goes at least Euclidean distance $c\varepsilon^{\alpha}$. We conclude that for $\varepsilon,\delta$ sufficiently small, the box $B(x,\delta)$ has two disjoint open arms that go at least Euclidean distance $c\varepsilon^\alpha-\delta$ (which are separated by two closed arms).

Since the 2-arm exponent is strictly larger than the 1-arm exponent, it follows from Lemma~\ref{lem:measure_2arm} that
\begin{equation*}
    \begin{split}
        \mu^{\cH_\ell}\!\left(\pi(\psi^{-1}(A))^{\delta}\setminus\pi(\psi^{-1}(A^\varepsilon))\right)\leq\mu^{\cH_\ell}\!\left({\{x:B(x,\delta)\text{ has 2 arms to distance $c\varepsilon^\alpha-\delta$}\}}\right)\rightarrow0
    \end{split}
\end{equation*}
as $\delta\downarrow0$, which completes the proof.

The proof of the the second inequality is similar. Let $n$ be large enough so that the GHf distance between $(\cX_n,d_n,\pi_n)$ and $(\cX,d,\pi)$ is smaller than $\delta$. Let $x_n$ be such that $x_n$ has an open arm in $\Lambda_n$ to the outer boundary of $B$, and $\|x_n-\pi_n(\psi_n^{-1}(A))\|_2<\delta$, and $x_n\not\in\pi_n(\psi_n^{-1}(A^\varepsilon))$. Then there are points $x,x' \in B(x_n;2\delta)$ with $d(x,x') > \varepsilon-\delta$, and we have two 
arms from $B(x;2\delta)$ to distance $c\varepsilon^\alpha-2\delta$.

Therefore
\begin{equation*}
    \mu_n^{\cH_\ell}\!\left(\pi_n(\psi_n^{-1}(A))^{\delta}\setminus\pi_n(\psi_n^{-1}(A^\varepsilon))\right)\leq\mu_n^{\cH_\ell}\!\left(\{x:B(x;2\delta)\text{ has 2 arms to distance $c\varepsilon^\alpha-\delta$}\}\right).
\end{equation*}
for some constants $c,\alpha>0$. Again, using Lemma~\ref{lem:measure_2arm} and the convergence $\mu_n^{\cH_\ell} \to \mu^{\cH_\ell}$, we can make the last part arbitrarily small by choosing $\delta$ small enough and $n$ large enough. 
\end{proof}

We are now ready to prove Lemma~\ref{lem:prokhorov}.
\begin{proof}[Proof of Lemma~\ref{lem:prokhorov}]
Let $\varepsilon>0$ and let $A\subset Y$ be measurable. Then by Prokhorov convergence of $\mu_{n}^{\cH_\ell}$ to $\mu^{\cH_\ell}$, on the event that $\diam(\cX)>2^{-\ell+1}$, we have that for all $\delta>0$, there exists $n_0$ such that for all $n\geq n_0(\delta)$, 
\begin{equation*}
    \begin{split}
        \mu_n(\psi_n^{-1}(A))=\mu_{n}^{\cH_\ell}(\pi_n(\psi_n^{-1}(A)))\leq\mu^{\cH_\ell}(\pi_n(\psi_n^{-1}(A))^\delta)+\delta
    \end{split}
\end{equation*}
and
\begin{equation*}
    \widetilde{\mu}^\ell(\psi^{-1}(A))=\mu^{\cH_\ell}(\pi(\psi^{-1}(A)))\leq\mu_n^{\cH_\ell}(\pi(\psi^{-1}(A))^\delta)+\delta.
\end{equation*}
Now by Lemma~\ref{lem:closeness-projections}, there exists $n_1=n_1(\delta)$ such that for all $n\geq n_1(\delta)$,
\begin{equation*}
    \mu^{\cH_\ell}(\pi_n(\psi_n^{-1}(A))^\delta)\leq\mu^{\cH_\ell}(\pi(\psi^{-1}(A^\delta))^{2\delta})
\end{equation*}
and
\begin{equation*}
    \mu_n^{\cH_\ell}(\pi(\psi^{-1}(A))^\delta)\leq\mu^{\cH_\ell}_n(\pi_n(\psi_n^{-1}(A^\delta))^{2\delta}).
\end{equation*}
Choosing $\delta=\delta(\varepsilon/2) \wedge \varepsilon/2$ and $n_2=n_2(\varepsilon/2)$ from Lemma~\ref{lem:two-arm-vs-one-arm}, we obtain that for all $n\geq n_2$
\begin{equation*}
    \begin{split}
        \mu^{\cH_\ell}(\pi(\psi^{-1}(A^\delta))^{2\delta})\leq\mu^{\cH_\ell}(\pi(\psi^{-1}(A^{\varepsilon/2+\delta})))+\varepsilon/2\leq\widetilde{\mu}^\ell(\psi^{-1}(A^\varepsilon))+\varepsilon
    \end{split}
\end{equation*}
and 
\begin{equation*}
    \mu_n^{\cH_\ell}(\pi_n(\psi_n^{-1}(A^\delta))^{2\delta})\leq\mu_n^{\cH_\ell}(\pi_n(\psi_n^{-1}(A^{\varepsilon/2+\delta})))+\varepsilon/2 \le \mu_n(\psi_n^{-1}(A^\varepsilon))+\varepsilon
\end{equation*}
which completes the proof.
\end{proof}

\subsection{Characterization of the limit}\label{sec:char-limit-measure-proof}

Let $B\subset[0,1]^2$ be a ball or a box, let $k\in\N$ and recall the definitions of $Y_k^A(B)$ and $Y^{(m)}_k(B)$ from Section~\ref{subsubsec:measure-conv}. For the latter, we again drop the superscript and just write $Y_k(B)$\phantomsection{}\label{def:Y_k} for the sake of brevity. 
Recall that we have defined the measure $\widetilde{\mu}^\ell$ in Section~\ref{subsec:ghp-proof} from the annulus measures $\mu^{\cH_\ell}$, and shown that it is the weak limit of the counting measure $\mu_n$ on the event $\{\diam(\cX) > 2^{-\ell+1}\}$. In particular, we can consistently define the measure $\widetilde{\mu} := \widetilde{\mu}^\ell$ where $\ell$ can be chosen to be any random number with $\ell > -\log_2(\diam(\cX))+1$. 

Our goal now is to identify $\widetilde{\mu}$ with the conformally covariant CLE$_6$ gasket measure defined in \cite{miller2024existence}. For this, we prove that $\widetilde{\mu}$ satisfies the same box-counting characterization as the annulus measure stated in Proposition~\ref{prop:char-gps} and the CLE$_6$ gasket measure stated in Proposition~\ref{prop:char-measure-my}.

\begin{prop}
\label{prop:char-mu-tilde} 
There exists a constant $c>0$ such that for all balls or boxes $B\subset[0,1]^2$,
\begin{equation*}
    \widetilde{\mu}(\pi^{-1}(B)) = \lim_{k\rightarrow\infty}\frac{cY_k(B)}{2^{2k}\AC_1(2^{-k},1)} \qquad\text{in probability}.
\end{equation*}
\end{prop}

Recall that we have a similar characterization for the annulus measures $\mu^A$ from Proposition~\ref{prop:char-gps} for domains $B$ that have rectifiable boundary. However, since $\pi(\cX)$ does not have rectifiable boundary, it is not immediately clear that the same holds for the measure $\widetilde{\mu}^\ell$, which is supported on $\pi(\cX)$. The proof of this proposition follows from two key ingredients. We first show that we can approximate $\pi(\cX)$ by rectifiable domains in an appropriate way. We then finish the proof by using the fact that the measure does not give mass to the set of $2$-arm points (Lemma~\ref{lem:measure_2arm}).

\begin{proof}[Proof of Proposition~\ref{prop:char-mu-tilde}]
For the ease of notation, we let $\ell \in \N$ be arbitrary and implicitly restrict to the event $\{\diam(\cX) \ge 2^{-\ell+1}\}$ throughout the proof. Since the probability of this event goes to~$1$ as $\ell \to \infty$, the desired convergence in probability then follows.

Define $Y_{k}^{\cH_\ell}(B):=\sum_{A\in\cH_\ell}Y^A_k(B)$. It follows from Proposition~\ref{prop:char-gps} that for all balls or boxes $B$,
\begin{equation}\label{eq:char-mu-Bl-B}
    \mu^{\cH_\ell}(B)=\lim_{k\rightarrow\infty}\frac{cY^{\cH_\ell}_k(B)}{2^{2k}\AC_1(2^{-k},1)}\qquad\text{in $L^2$}.
\end{equation}
Note that the formula above can be defined for any set that is a finite union of disjoint balls or boxes. Recall also that for each $B$ we have $\mu^{\cH_\ell}(\partial B) = 0$ a.s.\ by \cite[Corollary~4.16]{gps}. We would like to conclude that~\eqref{eq:char-mu-Bl-B} remains true when $B$ is replaced by $B \cap \pi(\cX)$. Since the boundary of $\cX$ is not rectifiable, we argue as follows.

Fix $m > \ell$. Let $B_{m}^{\cX}$ be the union of the set of boxes $B'$ in $\cB_{m}$ such that $B'$ has non-empty intersection with $\pi(\cX)$. Then it follows from~\eqref{eq:char-mu-Bl-B} that
\[
    \mu^{\cH_\ell}(B^\cX_m \cap B)=\lim_{k\rightarrow\infty}\frac{cY^{\cH_\ell}_k(B^\cX_m \cap B)}{2^{2k}\AC_1(2^{-k},1)}\qquad\text{in $L^2$}.
\]
Let $B_1\in\cB_k$ be such that $B_1\subset B$ and $2B_1$ intersects $\pi(\cX)$. Then for $k>m$, $B_1\subset B_m^{\cX}$ and $2B_1$ has an open arm that goes distance at least $2^{-\ell}$. So $Y_k(B)\leq Y_k^{\cH_\ell}(B_m^{\cX}\cap B)$. Now
suppose that a box $B_1\in\cB_k$ inside $B_m^{\cX}\cap B$ is counted in $Y^{\cH_\ell}_k(B^\cX_m)$ but $2B_1$ does not intersect $\pi(\cX)$. Let $B' \in \cB_m$ be the box containing $B_1$ and having non-empty intersection with $\pi(\cX)$. Then $B'$ has at least two open arms that go to distance $2^{-\ell}$ and are separated by two closed arms. Thus, for all $\ell$ and all $k>m$ 
\begin{align*}
        &Y_k^{\cH_\ell}(B_m^\cX\cap B) \\
        & \geq Y_k(B)\\
        &\geq Y_k^{\cH_\ell}(B_m^\cX\cap B)
        -Y_k^{\cH_\ell}\!\left(\bigcup\{B'\in\cB_m:\text{$B'$ has four arms that go distance $2^{-\ell}$}\}\right).
\end{align*}
For each $\varepsilon > 0$, by Lemma~\ref{lem:measure_2arm} we can choose $m$ sufficiently large so that
\[
\limsup_{k\to\infty}\E\!\left[ \left\lvert \frac{cY_k(B)}{2^{2k}\AC_1(2^{-k},1)} - \frac{cY^{\cH_\ell}_k(B_m^\cX\cap B)}{2^{2k}\AC_1(2^{-k},1)} \right\rvert\1_{\{\diam(\cX)\geq 2^{-\ell}\}} \right] \le \varepsilon
\]
Similarly, again by Lemma~\ref{lem:measure_2arm}, we have
\[\begin{split}
&\E[ \lvert \widetilde{\mu}^\ell(\pi^{-1}(B))-\mu^{\cH_\ell}(B^\cX_m \cap B) \rvert\1_{\{\diam(\cX)\geq 2^{-\ell}\}} ]\\ 
&\le \E\!\left[ \mu^{\cH_\ell}\!\left(\bigcup\{B'\in\cB_m:\text{$B'$ has four arms that go distance $2^{-\ell}$}\}\right) \right] \\
&\le \varepsilon .
\end{split}\]
This concludes the proof.
\end{proof}

\begin{proof}[Proof of Theorem~\ref{thm:main-convergence} and Corollary~\ref{co:cle_one_arm}]
Let $\widetilde{\mu}=\mu^{\cH_\ell}\circ\pi$ be the measure defined above, and let $\mu$ be the conformally covariant gasket measure on $\cX$. 
Let $(B_i)_{i\in\N}$ be a countable basis of boxes in $[0,1]^2$. By Lemma~\ref{lem:measure_2arm} and \cite[Proposition~2.6]{my2025diffusion}, both measures give $0$ mass to the set of double points of $\cX$, hence they are characterized by their values on $(\widetilde{\mu}(\pi^{-1}(B_i))_{i\in\N}$. By Proposition~\ref{prop:char-mu-tilde} and Proposition~\ref{prop:char-measure-my}, both $\widetilde{\mu}(\pi^{-1}(B_i))$ and $\mu(\pi^{-1}(B_i))$ are given by the following limits in probability 
\begin{align*}
    \widetilde{\mu}(\pi^{-1}(B_i))&=\lim_{k\to\infty}\frac{cY_k(B_i)}{2^{2k}\AC_1(2^{-k},1)} ,\\
    \mu(\pi^{-1}(B_i))&=\lim_{k\to\infty}\frac{\widetilde{c}Y_k(B_i)}{2^{(91/48)k}} .
\end{align*}
It follows that $\frac{2^{2k}\AC_1(2^{-k},1)}{2^{(91/48)k}} \to c$ for some constant $c>0$ (since both measures are non-trivial). In particular, $c\widetilde{\mu}(\pi^{-1}(B_i)) = \mu(\pi^{-1}(B_i))$ almost surely for each $i$, hence $c\widetilde{\mu}=\mu$. It is then straightforward to check that the GHP convergence of $(\cX_n,d_n,\mu_n)_{n\in\N}$ and the GHf convergence of $(\cX_n,\dpt_n,d_n,\pi_n)_{n\in\N}$ imply the GHPf convergence of $(\cX_n,d_n,\mu_n,\pi_n)_{n\in\N}$.
\end{proof}

\section{Distance, resistance, and displacement exponents}\label{sec:scaling-exponent-proof}

In this section, we are going to prove Corollaries~\ref{cor:scaling-exponent-distance}--\ref{co:spectral-dimension}.

\begin{proof}[Proof of Corollary~\ref{cor:scaling-exponent-distance}]
We first show show that $\E[S_n\mid \partial_L\Lambda_n\leftrightarrow\partial_R\Lambda_n]\lesssim \median_n$ as $n\rightarrow\infty$. Let $R\in\N$ be large, $\nu>0$ be small, and consider percolation on the enlarged box $B^+=[-1,2]\times[0,1]$. Note that both the chemical distance as well as the effective resistance between $\partial_L\Lambda_n$ and $\partial_R\Lambda_n$ is the same on $[0,1)^2$ as on $B^+$, since adding extra sites to the left and right does not influence the distances. Let $A_{K}$ be the event that there exists a left-to-right crossing of $\Lambda_n$ that does not come within distance $R^{-2K+3}$ of the top and bottom of $[0,1]^2$. If $\partial_L\Lambda_n\leftrightarrow\partial_R\Lambda_n$ but $A_K$ does not hold, then there must be a closed path that comes within distance $R^{-2K+3}$ of the top and bottom side of $[0,1]^2$, and an open path from left to right that fits above or below this closed path. Thus, there exists a box of size $R^{-2K+3}$ that shares one side with the top (resp.\ bottom) side of $[0,1]^2$, and that has two open arms and one closed arm coming out of it that go distance at least~$\frac{1}{2}$ in the half-plane. The number of such boxes is at most $\lesssim R^{2K}$, and the probability of such an arm event occurring for one box is at most $\lesssim R^{-(4-\nu)K}$ for every $\nu>0$ 
by Lemma~\ref{lem:arm_exponents}, since the half-plane $3$-arm exponent is $2$. On the other hand, if there exists a left-to-right crossing $\gamma$ of $\Lambda_n$ that does not come within distance $R^{-2K+3}$ of the top and bottom, then $\gamma$ is a path between two points $x,y$ in $B^+$ such that $\gamma$ does not come within distance $R^{-2K+3}$ of $\partial B^+$. Furthermore, $\median_n^{-1}S_n$ is bounded from above by $d_n(x,y)$. Let $\alpha>1$ and let $C_1=C_1(\alpha,\nu)$ from Proposition~\ref{prop:holderct-general}. Then there exists $C_3$ such that $R^{K-1}\geq C_1(R^{K-C_3}+\dist(\gamma,\partial B^+)^{-\frac{1}{4}-\nu})$ for all $K$ sufficiently large, since $R^{\frac{1}{4}\cdot 2K}<R^{K-1}$. Thus, by Proposition~\ref{prop:holderct-general}, we have $\P(\median_n^{-1}S_n>R^{K-1},A_{K})\lesssim R^{-\alpha K}$. Here we use that $\median_n=\altq_n(\fp)$ for some sufficiently small $\fp$. 
Therefore,
\begin{align*}
        &\,\E[\median_n^{-1}S_n\mid \partial_L\Lambda_n\leftrightarrow\partial_R\Lambda_n]\\
        \leq&\,\sum_{K\geq 1}R^K\P(\median_n^{-1}S_n>R^{K-1}\mid \partial_L\Lambda_n\leftrightarrow\partial_R\Lambda_n)\\
        \leq&\,\sum_{K\geq1}R^{K}\!\left(\P(A^c_{K}\mid \partial_L\Lambda_n\leftrightarrow\partial_R\Lambda_n)+\P(\median_n^{-1}S_n>R^{K-1},A_{K}\mid \partial_L\Lambda_n\leftrightarrow\partial_R\Lambda_n)\right)\\
        \leq&\,\sum_{K\geq1}C(R,\alpha,\nu)R^K(R^{2K}R^{-(4-\nu)K}+R^{-\alpha K})\leq C(R,\alpha,\nu)        
\end{align*}
for $\alpha$ large enough, for some constant $C(R)$ whose value may change from line to line. Here we also use that $\P(\partial_L\Lambda_n\leftrightarrow\partial_R\Lambda_n)$ is bounded from below uniformly in $n$.

We now show that $\operatorname{quant}_q[S_n\mid \partial_L\Lambda_n\leftrightarrow\partial_R\Lambda_n]\gtrsim \median_n$. Let $A_{M,n}$ be the event that there is a set of at most $M$ points $u_{1,n},\ldots,u_{M,n}$ that separate $\partial_L\Lambda_n$ from the line $\{\frac{1}{4}\}\times[0,1]$ and a set of at most $M$ points $v_{1,n},\ldots,v_{M,n}$ that separate $\partial_R\Lambda_n$ from the line $\{\frac{3}{4}\}\times[0,1]$. Define $A_M$, $u_1,\ldots,u_M$ and $v_1,\ldots,v_M$ similarly for $\Lambda$. By Lemma~\ref{lem:touching-loops}, we can find a coupling such that almost surely on the event $A_M \cap \{\partial_L[0,1]^2 \leftrightarrow \partial_R[0,1]^2\}$ the events $A_{M,n} \cap \{\partial_L\Lambda_n \leftrightarrow \partial_R\Lambda_n\}$ hold for $n$ large enough, with the $u_{i,n},v_{i,n}$ converging to $u_i,v_i$, respectively, for all $i$. Note that by the BK inequality, there exists $p_0<1$ such that $\P(A_{M,n}^c)\leq p_0^M$ for all $n$. Furthermore, by the parallel law, we have $\median_n^{-1}S_n\geq\frac{1}{M}\min_{i,j}d_n(u_{i,n},v_{j,n})$. Under the coupling above, the latter converges almost surely to $\frac{1}{M}\min_{i,j}d(u_i,v_j)$, which is a.s.\ positive since the Euclidean distance between the $u_i$ and $v_j$ is at least $\frac{1}{2}$. Thus,
\begin{align*}
        &\lim_{n\rightarrow\infty}\P(\median_n^{-1}S_n\leq x\mid \partial_L\Lambda_n\leftrightarrow\partial_R\Lambda_n)\\
        \leq&\, \lim_{n\rightarrow\infty}\!\left[\P\!\left(\frac{1}{M}\min_{i,j}d_n(u_{i,n},v_{j,n})\leq x,\,A_{M,n}\mid \partial_L\Lambda_n\leftrightarrow\partial_R\Lambda_n\right) + \P(A_{M,n}^c\mid \partial_L\Lambda_n\leftrightarrow\partial_R\Lambda_n)\right]\\
        \leq&\,\P\!\left(\frac{1}{M}\min_{i,j}d(u_i,v_j)\leq x,\,A_M\mid\partial_L\Lambda\leftrightarrow\partial_R\Lambda\right)+Cp_0^M.
\end{align*}
for some $C>0$. The last line can be made smaller than $q$ upon choosing $M$ sufficiently large depending on $q$ and $x$ sufficiently small depending on $q$ and $M$. We conclude that $\operatorname{quant}_q[S_n\mid\partial_L\Lambda_n\leftrightarrow\partial_R\Lambda_n]\geq x\median_n$.

It remains to show that $\median_n=n^{\beta+o(1)}$. Let $R\in\N$ and let $Q\subset[0,1]^2$ 
be a parallelogram that is small enough such that $RQ\subset[0,1]^2$. Recall that $X_{Q,n}$ is defined analogously to $X_n$, but on the parallelogram $Q$ and conditional on the event $\cE(\Lambda_{Q,n})$. Then 
\begin{equation*}
    \lim_{n\rightarrow\infty}\median_n^{-1}\med(X_{Q,n})=\lim_{n\rightarrow\infty}\median_{Rn}^{-1}\med(X_{Q,Rn})=\med(X_Q).
\end{equation*}
If we let $\Gamma$ denote a $\CLE_6$ in $Q$ and $d_Q(\cdot,\cdot;\Gamma)$ the metric $d_Q$ with respect to $\Gamma$, then by \cite[Theorem 1.2]{my2025diffusion} and \cite[Theorem 1.5]{my2025metric}, we have that $d_{RQ}(R\cdot,R\cdot;R\Gamma) = R^\beta d_Q(\cdot,\cdot;\Gamma)$. Thus,
\begin{align*}
        \lim_{n\to\infty} \median_{Rn}^{-1} \operatorname{median}(X_{Q,Rn}) 
        &= \lim_{n\to\infty} \frac{\median_n}{\median_{Rn}} \median_{n}^{-1} \operatorname{median}(X_{RQ,n}) \\
        &= \!\left( \lim_{n\to\infty} \frac{\median_n}{\median_{Rn}} \right) \operatorname{median}(X_{RQ}) \\
        &= \!\left( \lim_{n\to\infty} \frac{\median_n}{\median_{Rn}} \right) R^\beta \operatorname{median}(X_{Q}).
\end{align*}
This shows that for each $\varepsilon > 0$ and $R>1$, there exists $n_0=n_0(\varepsilon,R)$ such that 
\begin{equation}\label{eq:polynomial-growth-qn}
    R^{\beta-\varepsilon}\median_n \le \median_{Rn} \le R^{\beta+\varepsilon}\median_n
\end{equation}
for all $n\geq n_0$. Note that $\median_n$ is not obviously increasing, so we need to do a bit more work in order to conclude that $\median_n = n^{\beta+o(1)}$ as $n \to \infty$. Recall that $\kappa=\kappa(n)$ was defined such that $R^{\kappa-1}\leq n\leq R^{\kappa}$. We show that $\median_n\asymp \median_{R^{\kappa-1}}\asymp \median_{R^\kappa}$, which will imply the desired result.

Write $Q=[0,q]^2$ and let $Q'=[0,\frac{q}{2R}]\times[0,2Rq]$. Then for every $n$, we can embed $\Lambda_{Q',n}$ and $\Lambda_{Q,R^{\kappa}}$ into $\Lambda_n$ in such a way that $\Lambda_{Q',n}$ is taller than $\Lambda_{Q,R^{\kappa}}$ but less wide. Furthermore, let $S^1_{Q,n}=\inf\{d_n(x_n,y_n) \colon x\in\partial_L\Lambda_{Q,n},\,y\in\partial_R\Lambda_{Q,n}\}$, and define $S^1_{Q}$ analogously. Since $\p[\partial_L\Lambda_{Q,n}\leftrightarrow\partial_R\Lambda_{Q,n}] = \frac{1}{2}$, the $1/4$-quantiles of $S^1_{Q,n}$ and $S^1_{Q',n}$ are finite.

Then we have the stochastic domination $S^1_{Q',n}\leq S^1_{Q,R^{\kappa}}$ for all $n$. 
By Theorem~\ref{thm:main-convergence} and Lemma~\ref{lem:touching-loops} we have
\begin{align*}
    \lim_{n\to\infty} \median_n^{-1}\operatorname{quant}_{1/4}(S^1_{Q,n}) &= \operatorname{quant}_{1/4}(S^1_{Q}) ,\\
    \lim_{n\to\infty} \median_n^{-1}\operatorname{quant}_{1/4}(S^1_{Q',n}) &= \operatorname{quant}_{1/4}(S^1_{Q'}) .
\end{align*}

Thus,
\begin{align*}
        \operatorname{quant}_{1/4}(S^1_Q)=&\,\lim_{n\to\infty}\median_{R^\kappa}^{-1}\operatorname{quant}_{1/4}(S^1_{Q,R^\kappa})\\=&\,\lim_{n\rightarrow\infty}\frac{\median_n}{\median_{R^\kappa}}\median_n^{-1}\operatorname{quant}_{1/4}(S^1_{Q,R^\kappa})\\
        \geq&\,\limsup_{n\to\infty}\frac{\median_n}{\median_{R^\kappa}}\median_n^{-1}\operatorname{quant}_{1/4}(S^1_{Q',n})=\limsup_{n\to\infty}\!\left(\frac{\median_n}{\median_{R^{\kappa}}}\right)\operatorname{quant}_{1/4}(S^1_{Q'}).
\end{align*}
So there exists a constant $C=C(R)$ such that for all $n$, $\median_n\leq C\median_{R^\kappa}$. Similarly, in Lemma~\ref{lem:compXscales}, it was shown that there exists a constant $c=c(R)$ such that for all $n$, $\median_n\geq c\median_{R^{\kappa-1}}$. Let $\kappa_0$ be such that $R^{\kappa_0-1}<n_0\leq R^{\kappa_0}$ It follows from~\eqref{eq:polynomial-growth-qn} and the above that
\begin{equation*}
        \median_n\geq c\median_{R^{\kappa-1}}\geq c\median_{n_0R^{\kappa-1-\kappa_0}}\geq cR^{(\kappa-1-\kappa_0)(\beta-\varepsilon)}\geq n^{\beta-\varepsilon}
\end{equation*}
for $n$ sufficiently large, where the value of the constant $c$ may change from line to line. Similarly,
\begin{equation*}
    \median_n\leq C\median_{R^\kappa}\leq C\median_{n_0R^{\kappa-k_0}}\leq CR^{\kappa-\kappa_0(\beta+\varepsilon)}\leq n^{\beta+\varepsilon}
\end{equation*}
for $n$ sufficiently large, where the value of the constant $C$ may change from line to line. We conclude that $\median_n=n^{\beta+o(1)}$ as $n\rightarrow\infty$.
\end{proof}

To transfer the result to the IIC, we use the following result proved in \cite{kesten1986incipient,gps}. We denote by $\cA_1(r,R)$ the event that there exists an open crossing of the annulus $B(0,R) \setminus B(0,r)$. Recall from~\cite{lawler2002onearm} that $A_1(r,R) = \p[\cA_1(r,R)] = (r/R)^{5/48+o(1)}$.

\begin{lem}[{\cite[Proposition~5.2]{gps}}]\label{le:iic}
    There exists $c>0$ such that the following holds. For each $10 < r < R/100$, total variation between the restriction to $B(0,r)$ of $\piic$ and $\p[\cdot \mid \cA_1(1,R)]$ is at most $(r/R)^c$.
\end{lem}

\begin{proof}[Proof of Corollary~\ref{co:exponent_iic}]
We first show the upper bound. Let $\varepsilon > 0$, and let $\alpha = \alpha(\varepsilon) > 1$ be large. By applying Lemma~\ref{le:iic} with $r=n$ and $R=n^\alpha$, and Proposition~\ref{prop:holderct-general} on the box $B(0,2n)$ with exponent $\alpha^2$, we have 
\begin{equation}\label{eq:iic_estimate}
    \begin{split}
        \piic[ S_{n} > n^{\varepsilon}\median_{2n}] 
        &\le n^{-c(\alpha-1)} + \frac{\p[S_{n} > n^{\varepsilon }\median_{2n} ,\, \cA_1(1,n^{\alpha})]}{\p[\cA_1(1,n^{\alpha})]} \\
        &\lesssim n^{-c(\alpha-1)} + n^{-\alpha^2\varepsilon} n^{(5/48+\varepsilon)\alpha} 
    \end{split}
\end{equation}
for some absolute $c>0$. Since $\median_{2n} \le n^{\beta+o(1)}$ by Corollary~\ref{cor:scaling-exponent-distance}, this shows that $\piic$-almost surely we have
\[
S_{n} \le n^{\beta+o(1)}
\quad\text{as } n \to \infty .
\]

We turn to the lower bound. Let $\alpha > 1$ be arbitrary and let $M = M(\varepsilon,\alpha) \in \N$. For $k \in \N$ let $A_{M,k}$ be the event that there are at most $M$ disjoint crossings of the annulus $A(0;2^{k-1},2^{k})$. By RSW and the BK inequality we have $\p[(A_{M,k})^c] \le p_1^M$ for some $p_1 \in (0,1)$. On the event $A_{M,k-2} \cap A_{M,k}$, let $u_1,\ldots,u_M$ (resp.\ $v_1,\ldots,v_M$) be a collection of points in $A(0;2^{k-3},2^{k-2})$ (resp.\ $A(0;2^{k-1},2^{k})$) that separate its inner from its outer boundary in the percolation. Let
\[ Z_k = \min_{i,j} D_{A(0;2^{k-3},2^{k})}(u_i,v_j) \,\1_{A_{M,k-2} \cap A_{M,k}} \]
where $D_A$ denotes the unnormalized geodesic or resistance metric on percolation on $\bT \cap A$. We then have $S_{2^{k}} \ge M^{-2}Z_k$ by the parallel law on the event $\cA_1(1,2^k)$.

By Theorem~\ref{thm:main-convergence} and Lemma~\ref{lem:touching-loops}, the laws of $\median_{2^{k}}^{-1}Z_k$ converge to their CLE$_6$ analogue, and since the CLE$_6$ metric is a.s.\ positive, there exists for each $p_2 = p_2(\varepsilon,\alpha) > 0$ some $c_1 = c_1(\varepsilon,\alpha,p_2) > 0$ such that $\limsup_{k \to \infty} \p[ Z_k \le c_1 \median_{2^{k}} \mid A_{M,k-2} \cap A_{M,k} ] \le p_2$. Therefore
\[
\p[ Z_k \le c_1 \median_{2^{k}} ] \le p_1^M+p_2 =: p_0 ,
\]
and by independence 
\[
\p[ S_{2^{k+\varepsilon k}} \le M^{-2}c_1 \median_{2^{k}},\, \cA_1(1,2^{k+\varepsilon k}) ] \le p_0^{\varepsilon k} 
\quad\text{for $k$ sufficiently large.}
\]

The same argument as in~\eqref{eq:iic_estimate} then implies that for all $\alpha > 1$,
\[
\piic[ S_{2^{k+\varepsilon k}} \le M^{-2}c_1 \median_{2^{k}} ]\leq 2^{-c(1-\varepsilon)\alpha k} + 2^{(5/48+\varepsilon)\alpha k}p_0^{\varepsilon k}
\quad\text{for $k$ sufficiently large.}
\]
We can choose $p_0 = p_0(\varepsilon,\alpha)$ sufficiently small by choosing $M(\varepsilon,\alpha)$ sufficiently large and $p_2(\varepsilon,\alpha)$ sufficiently small, thus obtaining for all $\alpha > 1$,
\[
\piic[ S_{2^{k+\varepsilon k}} \le M^{-2}c_1 \median_{2^{k}} ] \le 2^{-c\alpha k} 
\quad\text{for $k$ sufficiently large.}
\]
Let $k$ be the largest integer such that $2^{k+\varepsilon k}\leq n$. Then $S_n\geq S_{2^{k+\varepsilon k}}$, and so for all $\varepsilon'$ small enough and $n$ large enough depending on $\varepsilon'$, we have
\[
    \piic[S_n\leq n^{-\varepsilon'}\median_n]\leq\piic[S_{2^{k+\varepsilon k}}\leq M^{-2}c_1\median_{2^k}]\lesssim n^{-c\alpha},
\]
as required.
\end{proof}

\begin{proof}[Proof of Corollary~\ref{co:spectral-dimension}]
Let $\mu$ be the vertex counting measure of the percolation cluster containing the origin (resp.\ the IIC). We first show that $\piic$-a.s.\ we have $\mu(B(0,n)) = n^{91/48+o(1)}$ as $n \to \infty$.

The lower bound for $\mu(B(0,n))$ is similar to the lower bound for the metric in the proof of Corollary~\ref{co:exponent_iic} above. Let $\varepsilon>0$, let $\alpha>1$ be arbitrary and let $M=M(\varepsilon,\alpha)\in\N$. Let $A_{M,k}$ be the event that there are $m\leq M$ disjoint open clusters $\cC^k_1,\ldots,\cC^k_m$ crossing the annulus $A(0;2^{k-1},2^k)$. Let
\begin{equation*}
    Z_k=\min_i|\cC^k_i| \,\1_{A_{M,k}} .
\end{equation*}
Then on the event $\cA_1(1,2^k)$, we have $\mu(B(0,2^k))\geq Z_k$. By convergence of the (renormalized) percolation measure to the $\CLE_6$ gasket measure, $2^{-2k}A(1,2^k)^{-1}Z_k$ converges to the analogous quantity for $\CLE_6$. Thus, for any $p_0>0$ there is $c_1 > 0$ such that 
\begin{equation*}
     \P[Z_k\leq c_12^{2k}A(1,2^k)]\leq p_0
\end{equation*}
for $k$ sufficiently large. By independence, we have
\begin{equation*}
    \P[\mu(B(0,2^{k+\varepsilon k}))\leq c_1 2^{2k}A(1,2^{k}), \,\cA_1(1,2^{k+\varepsilon k})]\leq p_0^{\varepsilon k} \quad \text{for $k$ sufficiently large.}
\end{equation*}
Now in exactly the same way as in the proof of Corollary~\ref{co:exponent_iic}, using the fact that $A(1,n)=n^{-5/48+o(1)}$, we obtain that for all $\nu>0$ and $\alpha>1$,
\[
\piic[\mu(B(0,n))\leq n^{91/48-\nu}] = O(n^{-\alpha}).
\]

For the upper bound, note that on the event $\cA_1(1,n)$, every site in the cluster containing the origin has an arm going Euclidean distance at least $n/2$. Thus, Lemma~\ref{lem:one-arm-measure-concentration} immediately implies that $$\p[\mu(B(0,n)) \geq n^{91/48+\nu} ,\, \cA_1(1,n) ] = O(n^{-\alpha})$$ as $n\to\infty$ for any $\alpha>1$. Again, it follows that
\[
\piic[\mu(B(0,n))\geq n^{91/48+\nu}] = O(n^{-\alpha})
\]
as $n\rightarrow\infty$, which concludes the proof that $\mu(B(0,n)) = n^{91/48+o(1)}$ holds $\piic$-a.s.

The proof of Corollary~\ref{co:spectral-dimension} then follows from \cite[Propositions~3.1 and~3.2]{kumagai2008heat}, applied with $v(n)=n^{91/48}$, $r(n)=n^\beta$, $\mathcal{I}(n)=n^{1/(91/48+\beta)}$, $\lambda=\lambda(n)=n^{o(1)}$, $m=m(n)=n^{o(1)}$, and $\varepsilon=\varepsilon(n)=n^{-o(1)}$.
\end{proof}

\addtocontents{toc}{
  \protect\contentsline{part}{\hyperref[app:arms]{Back matter}}{}{}
}

\appendix

\section{Arm calculations}
\label{app:arms}

We prove some arm bounds that are needed in Sections~\ref{sec:superpolynomial-concentration} and~\ref{sec:scaling-exponent-proof}.

\subsection{A four-arm bound}

We prove the following lemma, which is used in Section~\ref{sec:distance-between-pivotals}. Let $B^{(k)}$ be a box of side length $\frac{1}{n}R^{\kappa-k} \asymp R^{-k}$, and consider the percolation restricted to the sites in $B^{(k)}$.
\begin{lem}\label{lem:4arms}
For all $\zeta>0$, there exists $R=R(\zeta)$ large enough and a constant $C=C(R,\zeta)$ such that for all $n,k$, $x\in B^{(k)}$, and $a>0$ such that $R^{\kappa-k}\geq an\geq R$, we have
\begin{equation}\label{eq:4arm_box}
    \P\!\left(\textrm{an arm from $B(x,a)$ to each of the $4$ sides of $B^{(k)}$}\right)\leq C\cdot\!\left(\frac{a}{R^{-k}}\right)^{\frac{5}{4}-\zeta}.
\end{equation}
\end{lem}
The event in~\eqref{eq:4arm_box} denotes the event that there exist four disjoint arms $\eta_1,\ldots,\eta_4$ inside $B^{(k)}\setminus B(x,a)$ such that $\eta_i$ connects $B(x,a)$ to the $i$th side of the box $B^{(k)}$, where the sides of $B^{(k)}$ are numbered clockwise starting from the top. Note that we do not require that $B(x,a)\subset B^{(k)}$. The bound given on the right-hand side of~\eqref{eq:4arm_box} is just the $4$-arm probability from radius $a$ to $R^{-k}$, so the result is immediate when $B(x,a)$ is away from the boundary of $B^{(k)}$. The proof explains that the same (and in fact a stronger) bound is valid when $B(x,a)$ is near the boundary of $B^{(k)}$.

\begin{proof}[Proof of Lemma~\ref{lem:4arms}] 
Recall that we are working in a coordinate system determined by the triangular lattice. In this proof, when we consider distances of points from various lines parallel to the coordinate axes, we refer to the 1-Minkowski distance with respect to this coordinate system and not their Euclidean distance.

Let $d_1>0$ be the distance from $x$ to the closest side of $B^{(k)}$ and let $x_1$ be the lattice-parallel projection of $x$ onto this side. Similarly, let $d_2$ be the distance from $x$ to the second closest side and let $x_2$ be the vertex of $B^{(k)}$ connecting the two sides of $B^{(k)}$ closest to $x$. Note that if $B(x,a)$ has an arm to each side of $B^{(k)}$, then $B(x,a)$ has 4 arms to $B(x,d_1)$, $B(x_1,2d_1)\cap B^{(k)}$ has 3 arms in the half-plane to $B(x_1,d_2)\cap B^{(k)}$, and $B(x_2,2d_2)\cap B^{(k)}$ has 2 arms in the wedge domain of angle $\frac{\pi}{3}$ or $\frac{2\pi}{3}$ to $B(x_2,R^{-k})\cap B^{(k)}$. We can now conclude the proof by using arm estimates like in Lemma~\ref{lem:number_crossings}. That is, given $\zeta>0$, there exists $R=R(\zeta)$ large enough so that all the necessary arm exponent estimates apply. Thus, we find that for all $n,k,x,a$ such that $R^{\kappa-k}\geq an\geq R$ we have
\begin{equation}
    \begin{split}
        \P(\textrm{an arm from $B(x,a)$ to each of the $4$ sides of $B^{(k)}$})&\leq \!\left(\frac{a}{d_1}\right)^{\frac{5}{4}-\zeta}\!\left(\frac{2d_1}{d_2}\right)^{2-\zeta}\!\left(\frac{2d_2}{n^{-1}R^{\kappa-k}}\right)^{\frac{3}{2}-\zeta}\\ 
        &\lesssim \!\left(\frac{a}{R^{-k}}\right)^{\frac{5}{4}-\zeta} \cdot \!\left(\frac{d_1}{d_2}\right)^{3/4} \cdot \!\left(\frac{d_2}{R^{-k}}\right)^{1/4} 
        \\&\leq\!\left(\frac{a}{R^{-k}}\right)^{\frac{5}{4}-\zeta}, \nonumber 
    \end{split}
\end{equation} 
as required.
\end{proof}

\subsection{Six semi-arm concentration}\label{app:six-semiarms}
We prove the following result, which is used in Section~\ref{subsubsec:controlling-explored-domain}. Recall the definitions of the events $H_j$ and $H$ from Section~\ref{subsubsec:controlling-explored-domain}.

\begin{lem} 
The following holds for all $R$ large enough. For any $\alpha>0$ 
and $\nu$ as in the definition of events $H_j$, there exists $C=C(\alpha,R,\nu)$ such that for all $k,n$ and $k\leq j\leq \kappa$, we have that $$\P(H_j^c)\leq C\cdot R^{-4\alpha j}.$$ In particular, we have $\mathbf{P}\!\left(H^c\right)\leq C'\cdot  R^{-4\alpha k}$, for a constant $C'=C'(\alpha,R,\nu)$. 
\label{lem:bound_on_global_event_appendix}
\end{lem}

In order to prove this lemma, we want to bound the probability that $m$ different boxes of level $j$ have six semiarms of length $i-j-4$. To achieve this, we want to decorrelate the six semiarm events for different boxes. To do so, we use the following lemma, whose proof we omit since it is analogous to \cite[Lemma A.2]{kavvadias2023conformal}. The lemma states that around each box, we can find a sequence of annuli that are all pairwise disjoint and such that the overall ratio between the inner and outer boundary distances is upper bounded by a quantity involving the distance between these boxes. Then we upper bound the six semiarm events for a set of $m$ boxes with the probability that each annulus in the sequence has 6 arms crossing it. If the distance between the boxes is large, this probability will be small by Lemma~\ref{lem:arm_exponents}.

\begin{lem}  \label{lem:annuli_lemma}
Fix $m\in \mathbf{N}$, and let $i,j\in \mathbf{N}$ with $i\leq j-10$. Let $B(0,R^{-i}/2)$ be the box centered at the origin with side length $R^{-i}$, and tile it with boxes of side length $n^{-1}R^{\kappa-j}$. Let $K$ be the set of centers of such boxes that are fully contained inside $B(0,R^{-i}/2)$. Suppose that $z_1,\ldots, z_m\in K$ are distinct. Then there exist $n_{u}\in \{1,\ldots,m\}$ 
for each $1\leq u\leq m$ and $s_{u,l},r_{u,l} \in \mathbf{R}_{\geq 0}$ such that $0<s_{u,l}<r_{u,l} \leq n^{-1}R^{\kappa-i}$ for $1\leq l\leq n_{u}-1$ and $0=s_{u,n_{u}}< \frac{1}{2}n^{-1}R^{\kappa-j}\leq r_{u,n_{u}}$ so that all the annuli $A(z_u;s_{u,l},r_{u,l})$ are pairwise disjoint and
\[ \prod_{u=1}^{m} \frac{\prod_{l=1}^{n_{u}-1} s_{u,l}}{\prod_{l=1}^{n_{u}} r_{u,l} } \leq C(R,m)\prod_{u=1}^{m} \frac{1}{R^{-i} \wedge \min_{1\leq v<u} |z_{u}-z_v|} , \]
where $C(R,m)$ is a constant depending only on $R,m$. We use the convention that $\prod_{l=1}^{n_{u}-1} s_{u,l}=1$ if $n_{u}=1$, and the minimum of an empty set is infinite. 
\end{lem}

\begin{proof}[Proof of Lemma~\ref{lem:bound_on_global_event_appendix}] Fix $k\leq j\leq \kappa$. Using the union bound, we can bound the probability $\P(H_j^c)$ from above by
\begin{equation} \begin{split} \mathbf{P}\!\left(H_j^c\right) &\leq P\!\left(N_j^h\geq R^{j\nu}\right)+  \sum_{i=k-2}^{j-10}  \sum_{ \substack{B^{(i)}}} \mathbf{P}\!\left(N_{j,B^{(i)}}\geq R^{j\nu}\right)\\ &\leq P\!\left(N_j^h\geq R^{j\nu}\right)+\sum_{i=k-2}^{j-10}R^{2i-2k+6} \cdot \mathbf{P}(N_{j,B^{(i)}}\geq R^{j\nu}), \label{eq:union_bound} 
\end{split}
\end{equation}
where the sum in the first line is over $R^{2i-2k+6}$ boxes $B^{(i)}$ that have side length $R^{-i}$ and are centered at centers of boxes in $\cB^{R,k}_{i+3,n}$. In the second line, we are regarding $B^{(i)}$ as an abstract box of side length $R^{-i}$ rather than a subbox of $B^{(k)}$. This quantity provides an upper bound on $\P(N_{j,B^{(i)}}\geq R^{i\nu})$ for all boxes $B^{(i)}$ in the above sum, even if these boxes are not fully contained inside $B^{(k)}$ -- sites of $B^{(i)}\setminus B^{(k)}$ can only make $N_{j,B^{(i)}}$ larger. 
  
We start by examining the behavior of $N_{j,B^{(i)}}$. Fix any $k-2\leq i\leq j-10$. Fix a box $B^{(i)}=B(0,R^{-i}/2)$ of side length $R^{-i}$ as described above. For brevity of notation, in the remaining section of the proof we will write $N_j$ instead of $N_{j,B^{(i)}}$. We will find a bound on the $m$'th moment of $N_{j}$ for all $m\in \mathbf{N}$, and then use $m$ large enough to obtain the desired upper bound on~\eqref{eq:union_bound}. To do so, we will need to analyze the probability that we have the $6$-semiarm event for any $m$ fixed level $j$ boxes. Given any $m$ such boxes in $B^{(i)}$, Lemma~\ref{lem:annuli_lemma} gives us a way to find a collection of pairwise disjoint annuli centered around these boxes. Suppose $z_1,\ldots z_m$ are centers of these boxes, and denote the boxes by $B_u=B\!\left(z_u,\frac{1}{2}n^{-1} R^{\kappa-j}\right)$ for $u=1,\ldots,m$. Moreover, let $\tilde{l}_u$ denote the level of the gap in the semiarms of box $B_u$. Let $s_{u,l},r_{u,l}$ be as in Lemma~\ref{lem:annuli_lemma}. Then if there are six semiarms in the $u$'th box, each of the annuli $$A(z_u;s_{u,l},r_{u,l})\cap \!\left( A(z_u;\frac{1}{2}n^{-1}R^{\kappa-j+3},\frac{1}{2}n^{-1}R^{\kappa-j+\tilde{l}_u}) \cup A(z_u;\frac{1}{2}n^{-1}R^{\kappa-j+\tilde{l}_u+3},\frac{1}{2}n^{-1}R^{\kappa-i-4}) \right)$$ for $1\leq u\leq m$ and $1\leq l\leq n_u$ 
has six alternating arms. To account for the fact that our semiarms have gaps, we have intersected our annuli $A(z_u;s_{u,l},r_{u,l})$ above by the set on which we know we have an alternating $6$ arm event. 

Furthermore, all of these annuli are pairwise disjoint. So we can upper bound the event that we have $6$ semiarms of length $j-i-4$ around each of the $m$ boxes by the \textit{independent} events that tell us that there are 6 arms of alternating colors in each of the annuli above. 

Recall from Lemma~\ref{lem:arm_exponents} that for any $\zeta>0$ we can choose $R=R(\zeta)$ large enough so that for any\footnote{This bound ensures that the inner graph radius of the annulus $A(z_u;s',r')$ is at least $R$.} $n^{-1}\cdot R\leq s'\leq r'$ the 6-arm probability for any annulus $A(z_u;s',r')$ is bounded from above by $(\frac{s'}{r'})^{35/12-\zeta}$. Note that the smallest annulus radius relevant for us is $n^{-1}R^{\kappa-j+3}/2\geq n^{-1}R$ for all $j\leq \kappa$, thus this lemma applies. Finally, note that the inverse of the product of the aspect ratios of the new annuli centered at $z_u$ can be bounded from above by $$2R^4 \cdot R^3 \cdot \frac{\frac{1}{2}n^{-1}R^{\kappa-j+3}}{r_{u,n_u}}\prod_{l=1}^{n_u-1} \frac{s_{u,l}}{r_{u,l}}.$$
To see this, note that the $\frac{1}{2}n^{-1}R^{\kappa-j+3}$ term is precisely the smallest radius of our annuli, the $R^3$ term comes from the gap in our arms that has aspect ratio $R^3$ and the $2R^4$ term stems from the fact that our arms have length at most $\frac{1}{2}n^{-1}R^{\kappa-i-4}$, while $r_{u,1}$ is chosen to be smaller than $n^{-1}R^{\kappa-i}$.

Combining the above observations, we have
\begin{equation*}
\begin{split}
&\mathbf{P}\!\left( B_u \textrm{ has 6 semiarms of length } j-i-4  \textrm{ for each } 1\leq u\leq m\right) \\ 
&\leq C(m,R)\sum_{0\leq  \tilde{l}_1< j-i} \cdots \sum_{0\leq \tilde{l}_m< j-i}\!\left(\prod_{u=1}^m R^{11}\frac{R^{-j}}{r_{u,n_u}}\prod_{l=1}^{n_u-1}\frac{ s_{u,l} }{ r_{u,l}}\right)^{35/12-\zeta},
\end{split}
\end{equation*}
where $C(m,R)$ denotes a constant depending only on $m,R$ and can increase from line to line, $35/12$ is the 6-arm exponent and $\zeta>0$ is a constant that will be chosen later. Now by Lemma~\ref{lem:annuli_lemma}, we have
\begin{equation*}
\begin{split}
&\mathbf{P}\!\left( B_u \textrm{ has 6 semiarms of length } j-i-4  \textrm{ for each } 1\leq u\leq m\right) \\ 
&\leq C(m,R) \sum_{0\leq  \tilde{l}_1< j-i} \cdots \sum_{0\leq \tilde{l}_m< j-i}\!\left( \prod_{u=1}^m \frac{R^{11} \cdot R^{-j}}{R^{-i} \wedge \min_{v<u} |z_v-z_u|}\right)^{35/12-\zeta}  \\
&\leq C(m,R)(j-i)^m \!\left( \prod_{u=1}^m \frac{R^{11} \cdot R^{-j}}{R^{-i} \wedge \min_{1\leq v<u} |z_v-z_u|}\right)^{35/12-\zeta} \\ 
&\leq C(m,R) (j-i)^m \cdot R^{-jm(35/12-\zeta)} \!\left(\prod_{u=1}^m \frac{1}{R^{-i} \wedge \min_{1\leq v<u} |z_v-z_u|}\right)^{35/12-\zeta}.
\end{split}
\end{equation*}

Let $\beta=35/12-\zeta$. We can use the above to bound $\mathbf{E}[N_j^m]$ as follows:
 \begin{equation*}
 \begin{split}
 \mathbf{E}[N_j^m] &= \mathbf{E} \!\left[ \!\left(\sum_{\substack{B\subseteq B^{(i)}\\ B \textrm{ has level }j}} \mathbf{1}_{\!\left\{B \textrm{ has 6 semiarms of length } j-i-4 \right\}}\right)^m\right]\\
 &= \sum_{1\leq \tilde{m} \leq m} c_{\tilde{m}} \sum_{\substack{B_{1},\ldots, B_{\tilde{m}}\\ \textrm{ distinct}}} \mathbf{P} \!\left( B_{1},\ldots, B_{{\tilde{m}}} \textrm{ all have 6 semiarms of length } j-i-4 \right) \\
 &\leq \sum_{1\leq \tilde{m} \leq m} c_{\tilde{m}} \sum_{\substack{B_{1},\ldots, B_{\tilde{m}}\\ \textrm{ distinct}}} C(\tilde{m},R) (j-i)^{\tilde{m}} R^{-j\tilde{m}\beta} \!\left(\prod_{u=1}^{\tilde{m}} \frac{1}{R^{-i} \wedge \min_{1\leq v<u} |z_v-z_u|}\right)^{\beta}.
 \end{split}
 \end{equation*}
The constants $c_{\tilde{m}}$ in the above come from the multinomial expansion. Noting that $c_{\tilde{m}}$ and $C(\tilde{m},R)$ are bounded by some constant depending only on $m,R$ and splitting the sum over $\tilde{m}$ boxes into $\tilde{m}$ separate sums, we can further bound $\E[ N_j^m]$ as
 \begin{equation*}
 \begin{split}
 \mathbf{E}[N_j^m]&\leq C(m,R) \sum_{1\leq \tilde{m} \leq m} (j-i)^{\tilde{m}}  R^{-j\tilde{m}\beta} \frac{1}{\tilde{m}!}\sum_{B_1} \bigg[ \!\left(\frac{1}{R^{-i}}\right)^{\beta} \\ 
 &\cdot  \sum_{B_2} \!\left(\frac{1}{|z_1-z_2|}\right)^{\beta} \cdots \sum_{B_{\tilde{m}}} \!\left(\frac{1}{\min_{v<\tilde{m}} \{|z_{\tilde{m}}-z_v|\}}\right)^{\beta} \bigg] \\
 &\leq C(m,R) \sum_{1\leq \tilde{m} \leq m} (j-i)^{\tilde{m}} R^{-j\tilde{m}\beta} \frac{1}{\tilde{m}!}\sum_{B_1} \bigg[ R^{i\beta} \sum_{a_2=i}^{j} \sum_{B_2 : |z_1-z_2| \leq R^{-a_2}} R^{a_2\beta} \\
 &\cdot \cdots \sum_{a_{\tilde{m}}=i}^{j} \sum_{B_{\tilde{m}}: \min_{v<\tilde{m}}|z_{\tilde{m}}-z_v|\leq R^{-a_{\tilde{m}}}} R^{a_{\tilde{m}}\beta} \bigg].
 \end{split}
\end{equation*}
Now note that the number of boxes $B_u$ such that $\min_{v<u} |z_u-z_v| \leq R^{-a_u}$ is at most $u R^{2(j-a_u)+2}$, since such a level $j$ box must intersect at least one of the $u$ boxes centered around $z_0,\ldots z_{u-1}$ of radius $R^{-a_u}$. This implies that 
 \begin{equation*}
 \begin{split}
 \E[N_j^m]& \leq C(m,R) \sum_{1\leq \tilde{m} \leq m} (j-i)^{\tilde{m}}  R^{-j\tilde{m}\beta} \frac{1}{\tilde{m}!}\sum_{B_1} R^{i\beta} \sum_{a_2=i}^{j} 2R^{2(j-a_2)} R^{a_2\beta} \cdots \sum_{a_{\tilde{m}}=i}^{j} \tilde{m} R^{2(j-a_{\tilde{m}})} R^{a_{\tilde{m}}\beta} \\
 &\leq C(m,R) \sum_{1\leq \tilde{m}\leq m} (j-i)^{\tilde{m}} R^{-j\tilde{m}\beta} R^{2(j-i) + i\beta} ((j-i) R^{2j+j(\beta-2)})^{\tilde{m}-1} \\
 &\leq  C(m,R) (j-i)^{2m} R^{-(\beta-2)(j-i)} \leq C(m,R) R^{-(11/12-2\zeta)(j-i)},
 \end{split}
 \end{equation*}
In the above, the constant $C(m,R)$ increases from line to line, but only depends on $m$ and $R$ throughout.

We can analogously obtain that for all $m\geq 1$ there exists a constant $C(m,R)$ such that $$\E\!\left[(N_j^h)^m\right]\leq C(m,R) \cdot R^{-(1-2\zeta)(j-k)}.$$ Indeed, we can redo the above calculation with $k$ in place of $i$ and $\beta=2-\zeta$ in place of $\beta=35/12-\zeta$. This corresponds to the $3$-arm exponent in the half-plane (see Lemma~\ref{lem:arm_exponents}). To complete the calculation, one just needs to note that now the number of boxes $B_u$ that are suitably close to the boundary and such that $\min_{v<u}|z_u-z_v|\leq R^{-a_u}$ is at most $4R^{9}\cdot u R^{j-a_u+1}$.

Finally, going back to equation~\eqref{eq:union_bound} and applying Markov's inequality, we get that for any $m\in \mathbf{N}$
 \begin{equation}
 \begin{split}
 \mathbf{P}(H_j^c) &\leq  \frac{\E\!\left[ \!\left(N_j^h\right)^m\right]}{R^{mj\nu}}+ \sum_{i=k-2}^{j-10} R^{2i-2k+6} \cdot \frac{\mathbf{E}\!\left[N_{j,B^{(i)}}^m\right]}{R^{mj\nu}} \\
&\leq C(m,R,\nu) \cdot \bigg[ \frac{R^{-(1-2\nu)(j-k)}}{R^{mj\nu}}+ \sum_{i=k-2}^{j-10} \frac{R^{2i-2k}}{R^{mj\nu}} R^{-(11/12-2\zeta)(j-i)} \bigg] \\
&\leq C(m,R,\nu) \cdot j R^{(2-m\nu)j},\nonumber 
 \end{split}
 \end{equation}
 where we have chosen $\zeta=1/6$ for example to make our bounds valid. Setting $m=m(\alpha,\nu)= \lceil \frac{4+4\alpha}{\nu}\rceil +1$ in the above yields $$\P(H_j^c)\leq C(\alpha,\nu, R)\cdot R^{-4\alpha j}.$$
Summing this formula over all $j\geq k$, we find that $\mathbf{P}(H^c)\leq C(\alpha,R, \nu) \cdot R^{-4\alpha k}$ for a possibly larger choice of the constant $C(\alpha,R,\nu)$, as wanted. 
 \end{proof}

\subsection{One-arm concentration}
We prove the following lemma which we need for the volume upper bound in the proof of Corollary~\ref{co:spectral-dimension}. The lemma is a one-arm analogue of Lemma~\ref{lem:bound_on_global_event_appendix} for the case $j=\kappa$.

\begin{lem}\label{lem:one-arm-measure-concentration}
Consider critical percolation on $\Lambda_n$. Let $N_n$ be the number of open sites $z\in \Lambda_n$ such that $z$ has an open arm that goes Euclidean distance at least $1/2$. Then for any $\alpha>0$ and $\nu>0$, there exists $C=C(\alpha,\nu)$ such that for all $n$, we have
\begin{equation}
    \P\left(N_n\geq n^{\frac{91}{48}+\nu}\right)\leq Cn^{-\alpha}.
\end{equation}
\end{lem}
The proof of Lemma~\ref{lem:one-arm-measure-concentration} is an easier version of Lemma~\ref{lem:bound_on_global_event_appendix}.
\begin{proof}[Proof of Lemma~\ref{lem:one-arm-measure-concentration}]
Like before, we need to show that $n^{-91/48-\nu}N_n$ has bounded moments. Let $m\in\N$ and let $z_1,\ldots,z_m\in\Lambda_n$ be distinct such that each $z_u$ has an open arm that goes Euclidean distance at least $1/2$. We apply Lemma~\ref{lem:annuli_lemma} with $i=0$ and $j=\kappa$. Let $s_{u,l},r_{u,l}$ be as in Lemma~\ref{lem:annuli_lemma}. Note that the $r_{u,l}$ can be chosen such that $s_{u,l}<r_{u,l}\leq1/4$. Then each annulus $A(z_u;s_{u,l},r_{u,l})$ has an open arm crossing it and the annuli are all disjoint. Therefore, by independence and because the 1-arm exponent is $5/48$ (recall \cite{lawler2002onearm}), for each $\zeta>0$, there exists $C(\zeta,m)$ such that for all $n$,
\begin{equation*}
    \begin{split}
        {}&\P\!\left(z_1,\ldots,z_m\text{ each have an open arm that goes Euclidean distance at least $\frac{1}{2}$}\right)\\
        \leq{}& C(\zeta,m)\left(\prod_{u=1}^m\frac{n^{-1}}{r_{u,n_u}}\cdot\frac{\prod^{n_u-1}_{l=1}s_{u,l}}{\prod_{l=1}^{n_u-1}r_{u,l}}\right)^{\frac{5}{48}-\zeta}\\
        \leq{}& C(\zeta,m)\left(\prod_{u=1}^m\frac{n^{-1}}{\min_{1\leq v<u}|z_u-z_v|}\right)^{\frac{5}{48}-\zeta},
    \end{split}
\end{equation*}
where the value of $C(\zeta,m)$ may change from line to line. The term $\frac{n^{-1}}{r_{u,n_u}}$ corresponds to the one-arm bound for the inner annulus of radii $s_{u,n_0}=0$ and $r_{u,n_u}$. The final inequality is due to Lemma~\ref{lem:annuli_lemma}. It follows that for all $n$,
\begin{equation*}
    \begin{split}
        &\E[N_n^m]\\
        \leq{}&\sum_{1\leq\widetilde{m}\leq m}c_{\widetilde{m}} \sum_{\substack{z_1,\ldots,z_{\widetilde{m}}\\ \text{distinct}}} \P\!\left(z_1,\ldots,z_{\widetilde{m}}\text{ each have an open arm that goes Euclidean distance at least $\frac{1}{2}$}\right)\\
        \leq{}&\sum_{1\leq\widetilde{m}\leq m}c_{\widetilde{m}} \sum_{\substack{z_1,\ldots,z_{\widetilde{m}}\\ \text{distinct}}}C(\zeta,\widetilde{m})\left(\prod_{u=1}^{\widetilde{m}}\frac{n^{-1}}{\min_{1\leq v<u}|z_u-z_v|}\right)^{\frac{5}{48}-\zeta}\\
        \leq{}&\sum_{1\leq\widetilde m\leq m}C(\zeta,\widetilde{m})n^{(2-\frac{5}{48}+\zeta)\widetilde{m}}\\
        \leq{}&C(\zeta,m)n^{(\frac{91}{48}+\zeta)m},
    \end{split}
\end{equation*}
where $c_{\widetilde m}$ is as in the proof of Lemma~\ref{lem:bound_on_global_event_appendix} and the values of the constants $C(\zeta,m),C(\zeta,\widetilde m)$ may change from line to line. The third inequality is a standard summation inequality and follows as in the proof of Lemma~\ref{lem:bound_on_global_event_appendix}.

Now by Markov's inequality,
\begin{equation*}
    \begin{split}
        \P\!\left(N_n\geq n^{\frac{91}{48}+\nu}\right)\leq{}&C(\zeta,m)n^{-(\frac{91}{48}+\nu)m+(\frac{91}{48}+\zeta)m}=C(\zeta,m)n^{(\zeta-\nu)m}.
    \end{split}
\end{equation*}
The desired result follows from choosing $\zeta=\nu/2$ and $m=\lceil2\alpha/\nu\rceil$.
\end{proof}

\section*{Index of notation}

\renewcommand{\arraystretch}{1.15}
\setlength{\LTleft}{0pt}
\setlength{\LTright}{0pt}

\begin{longtable}{@{}>{\raggedright\arraybackslash}p{0.20\textwidth}
                  >{\raggedright\arraybackslash}p{0.55\textwidth}
                  >{\raggedright\arraybackslash}p{0.21\textwidth}@{}}
\toprule
\textbf{Notation} & \textbf{Meaning} & \textbf{Defined in} \\
\midrule
\endfirsthead

\toprule
\textbf{Notation} & \textbf{Meaning} & \textbf{Defined in} \\
\midrule
\endhead

\midrule
\multicolumn{3}{r}{\emph{Continued on next page}} \\
\endfoot

\bottomrule
\endlastfoot
\multicolumn{3}{@{}l}{\textbf{Lattices, domains, and $\CLE_6$}}\\[0.3em]
$\bT$ & The triangular lattice centered at $0$ & \hyperref[def:bT]{Sec.~\ref*{def:bT}}\\
$\bT_n := n^{-1}\bT$ & The rescaled triangular lattice & \hyperref[def:bT_n]{Sec.~\ref*{def:bT_n}} \\
$[0,1]^2$ & The triangular unit parallelogram, i.e., the region enclosed by the Euclidean vertices $(0,0)$, $(1,0)$, $(1/2,\sqrt{3}/2)$, $(3/2,\sqrt{3}/2)$ & \hyperref[def:unit-square]{Sec.~\ref*{def:unit-square}}\\
$\Lambda_n$ & The subgraph of $\bT_n$ induced by vertices in $\bT_n \cap [0,1)^2$; by abuse of notation also its vertex set & \hyperref[def:Lambda_n]{Sec.~\ref*{def:Lambda_n}} \\
$\pi_n$ & The embedding of $\Lambda_n$ into $[0,1]^2$ & \hyperref[def:pin]{Sec.~\ref*{def:pin}} \\

$\cX_n^{(1)},\cX_n^{(2)},\ldots$ & The open and closed clusters of $\Lambda_n$, ordered by decreasing Euclidean diameter & \hyperref[def:cXnm]{Sec.~\ref*{def:cXnm}} \\
$\cX_n$ & The outermost open percolation cluster surrounding $(1/2,1/2)$ & \hyperref[def:cXn-first]{Sec.~\ref{def:cXn-first}}, \hyperref[def:cXn-second]{Sec.~\ref{def:cXn-second}}\\[0.4em]
$\Gamma_n$ & The collection of interface loops of a percolation configuration on $\Lambda_n$ & \hyperref[def:Gamma_n]{Sec.~\ref*{def:Gamma_n}} \\
$\Gamma$ & A nested $\CLE_6$ in $(0,1)^2$ & \hyperref[def:Gamma]{Sec.~\ref*{def:Gamma}} \\
$\Lambda=\Lambda(\Gamma)$ & The space of prime ends of $\Gamma$ & \hyperref[def:Lambda]{Sec.~\ref*{def:Lambda}} \\
$\pi$ & The embedding of $\Lambda$ into $[0,1]^2$ & \hyperref[def:pi]{Sec.~\ref*{def:pi}} \\[0.4em]
$\cX^{(1)},\cX^{(2)},\ldots$ & The connected components of $\Lambda$, equivalently the interior gaskets of $\Gamma$, ordered by decreasing Euclidean diameter & \hyperref[def:cXm]{Sec.~\ref*{def:cXm}} \\
$\cL$ & The outermost loop in $\Gamma$ surrounding $(1/2,1/2)$ & \hyperref[def:cL]{Sec.~\ref*{def:cL}}\\
$\cX$ & The outermost open cluster of $\Lambda$ surrounding $(1/2,1/2)$, equivalently the gasket of $\Gamma$ whose exterior boundary is $\cL$ & \hyperref[def:cX-first]{Sec.~\ref*{def:cX-first}}, \hyperref[def:cX-second]{Sec.~\ref*{def:cX-second}} \\[0.4em]
$\Gamma_{\cL}$ & The collection of loops in $\Gamma$ that are surrounded by $\cL$ but no other loops & \hyperref[def:Gamma_cL]{Sec.~\ref*{def:Gamma_cL}}\\
$U$ & An open and simply connected subset of $[0,1]^2$ & \hyperref[def:U]{Sec.~\ref*{def:U}}\\
$\Gamma_{U^*}$ & The collection of loops in $\Gamma$ that are inside $\cL$ and entirely contained in $U$ & \hyperref[def:GammaU*]{Sec.~\ref*{def:GammaU*}}\\
$U^*$ & The set of points in $U$ that are inside $\cL$ and not on or inside any other loop of $\Gamma_{\cL}\setminus\Gamma_{U^*}$ & \hyperref[def:U*]{Sec.~\ref*{def:U*}}\\

\multicolumn{3}{@{}l}{\textbf{Metrics}}\\[0.3em]
$\dist$ & The Euclidean distance & \hyperref[def:dist]{Sec.~\ref*{def:dist}} \\
$\diam$ & The Euclidean diameter & \hyperref[def:diam]{Sec.~\ref*{def:diam}}, \hyperref[def:diam-second]{Sec.~\ref*{def:diam-second}} \\
$P_n(x,y)$ & The set of open (resp.\ closed) paths from $x$ to $y$ in the relevant discrete cluster & \hyperref[def:Pn]{Sec.~\ref*{def:Pn}} \\
$P(x,y)$ & The set of admissible paths from $x$ to $y$ with respect to $\Gamma$ & \hyperref[def:P(x,y)]{Sec.~\ref*{def:P(x,y)}} \\
$\dpt_n$ & The discrete path metric & \hyperref[eq:dpt_n]{Eq.~\eqref{eq:dpt_n}} \\
$\dpt$ & The continuum path metric on $\Lambda$ & \hyperref[eq:dpt-first]{Eq.~\eqref{eq:dpt-first}} \\
$\Dg_n$, $\Dr_n$ & The unnormalized geodesic and resistance percolation metrics on $\Lambda_n$ & \hyperref[def:Dgeon]{Sec.~\ref*{def:Dgeon}} \\

$D_n$ or $D$ & Refers to either $\Dg_n$ or $\Dr_n$. The notation $D$ is used for the metric on certain subgraphs of $\Lambda_n$; sometimes the subscript $n$ is also dropped when the underlying graph is clear from the context. & \hyperref[def:D_n]{Sec.~\ref*{def:D_n}} \\

$\median_n^{\geo}$, $\median_n^{\res}$ & Normalizing constant for the geodesic and resistance percolation  metrics & \hyperref[def:mediangeo]{Sec.~\ref*{def:mediangeo}} \\

$\median_n$ & Refers to either $\median_n^{\geo}$ or $\median_n^{\res}$. Defined as $\median_n=\altq_n(\fp)$, see \hyperlink{proofs-index-notation}{Proofs}. & \hyperref[def:median_n]{Sec.~\ref*{def:median_n}}, \hyperref[eq:def-median_n-second]{Eq.~\eqref{eq:def-median_n-second}}, \hyperref[def:median_n-hoelder]{Sec.~\ref*{def:median_n-hoelder}} \\
$\dg_n:=\median_n^{-1}D^{\geo}_n$ & The normalized geodesic percolation metric & \hyperref[eq:def_dgeon_dresn]{Eq.~\eqref{eq:def_dgeon_dresn}} \\
$\dr_n:=\median_n^{-1}D^{\res}_n$ & The normalized resistance percolation metric & \hyperref[eq:def_dgeon_dresn]{Eq.~\eqref{eq:def_dgeon_dresn}}\\

$d_n$ & Refers to either $\dg_n$ or $\dr_n$ & \hyperref[def:d_n]{Sec.~\ref*{def:d_n}} \\
$\dg$, $\dr$ & The geodesic and resistance $\CLE_6$ metrics on $\Lambda$; also the (subsequential) scaling limits of $\dg_n$ and $\dr_n$ & \hyperref[def:dg]{Sec.~\ref*{def:dg}}, \hyperref[def:dg_scaling_limit]{Sec.~\ref*{def:dg_scaling_limit}} \\

$d$ & Refers to either $\dg$ or $\dr$ & \hyperref[def:def-d]{Sec.~\ref*{def:def-d}} \\
$\fC$ & The class of admissible subregions $V$ used in the characterization of $\CLE_6$ metrics & \hyperref[def:fC]{Sec.~\ref*{def:fC}} \\
$\fd_{V}$ & The geodesic or resistance $\CLE_6$ metrics on $\cX\cap \overline{V}$ for $V\in\fC$ & \hyperref[def:fd-first]{Sec.~\ref*{def:fd-first}}, \hyperref[eq:def-dU]{Eq.~\eqref{eq:def-dU}} \\[0.4em]

\multicolumn{3}{@{}l}{\textbf{Measures}}\\[0.3em]
$\mu_n^{(m)}$, $\mu_n$ & The normalized counting measure on $\cX_n^{(m)}$. The superscript~$(m)$ is dropped in Sec.~\ref{sec:ghpf-proof}. & \hyperref[eq:counting_measure_cluster]{Eq.~\eqref{eq:counting_measure_cluster}}, \hyperref[def:mu]{Sec.~\ref*{def:mu}} \\
$\mu_n^A$ & The normalized counting measure of all vertices in $\Lambda_n$ in the inner face of $A$ that are connected by an open path to the outer boundary of $A$ & \hyperref[eq:annulus_measure]{Eq.~\eqref{eq:annulus_measure}}\\
$\mu^{(m)},\mu$ & The $\CLE_6$ gasket measure on $\cX^{(m)}$. The superscript~$(m)$ is dropped in Sec.~\ref{sec:ghpf-proof}. & \hyperref[def:mum]{Sec.~\ref*{def:mum}}, \hyperref[def:mum-second]{Sec.~\ref*{def:mum-second}}, \hyperref[def:mu]{Sec.~\ref*{def:mu}}\\
$\mu^A$ & The $\CLE_6$ measure on all points in $\Lambda$ in the inner face of $A$ that are connected by an admissible path to the outer boundary of $A$; the scaling limit of $\mu_n^A$ & \hyperref[th:gps]{Th.~\ref{th:gps}}, \hyperref[eq:muA(B)]{Eq.~\eqref{eq:muA(B)}}\\
$Y^{(m)}_{k}$, $Y_{k}$ & The number of dyadic boxes $B$ at level $k$ such that $2B$ intersects $\cX^{(m)}$. The superscript~$(m)$ is dropped in Sec.~\ref{sec:ghpf-proof}. & \hyperref[def:Ymk]{Sec.~\ref*{def:Ymk}}, \hyperref[def:Y_k]{Sec.~\ref*{def:Y_k}}\\
$Y^{A}_{k,n}$, $Y^A_{k}$ & The number of dyadic boxes $B$ at level $k$ in the inner face of $A$ such that $2B$ is connected to the outer face of $A$ in the percolation or $\CLE_6$, respectively. & \hyperref[def:YA_k,n]{Sec.~\ref*{def:YA_k,n}}, \hyperref[def:YA_k]{Sec.~\ref*{def:YA_k}}\\
$\cH_\ell$ & A tiling of $[0,1]^2$ by dyadic annuli & \hyperref[def:cHl]{Sec.~\ref*{def:cHl}}\\
$\mu_n^{\cH_\ell}:=\sum_{A\in\cH_\ell}\mu_n^A$ & The percolation annulus measure extended to a tiling~$\cH_\ell$ & \hyperref[def:mu_ncH]{Sec.~\ref*{def:mu_ncH}}\\
$\mu^{\cH_\ell}$ & The scaling limit of $\mu_n^{\cH_\ell}$ & \hyperref[def:mucH]{Sec.~\ref*{def:mucH}}\\
$\widetilde{\mu}^{\ell}:=\mu^{\cH_\ell}\circ\pi$ & The pullback measure of $\mu^{\cH_\ell}$ on $\cX$ under the projection $\pi$ & \hyperref[def:tilde-muell]{Sec.~\ref*{def:tilde-muell}}\\[0.4em]

\multicolumn{3}{@{}l}{\textbf{$R$-adic geometry}}\\[0.3em]
$\kappa$ & The unique integer such that $R^{\kappa-1}<n<R^{\kappa}$ & \hyperref[def:kappa]{Sec.~\ref*{def:kappa}} \\
$\cD_k^R$ & The set of $R$-adic grid points at level $k$ & \hyperref[def:D_kR]{Sec.~\ref*{def:D_kR}} \\
$\cB_k^R$ & The collection of $R$-adic boxes at level $k$ & \hyperref[def:B_kR]{Sec.~\ref*{def:B_kR}} \\
$\cQ_k^R$ & The simply connected unions of boxes from $\cB_k^R$ & \hyperref[def:cQ_kR]{Sec.~\ref*{def:cQ_kR}} \\
$\cD^R := \bigcup_k \cD_k^R$ & The set of all $R$-adic grid points & \hyperref[def:cDR]{Sec.~\ref*{def:cDR}} \\
$\cB^R := \bigcup_k \cB_k^R$ & The set of all $R$-adic boxes & \hyperref[def:cBR]{Sec.~\ref*{def:cBR}} \\
$\cQ^R := \bigcup_k \cQ_k^R$ & The set of all simply connected $R$-adic domains\\
$\cB_{k,n}^R, \cD_{k,n}^R$ & The $n$-dependent $R$-adic boxes of side length $n^{-1}R^{\kappa-k}$ and corresponding grid points & \hyperref[def:cBR_k,n]{Sec.~\ref*{def:cBR_k,n}}\\
$B^{(k)}$ & In Section~\ref{sec:apriori}: the $R$-adic superbox at level $k$ of a box~$B$. In Section~\ref{sec:superpolynomial-concentration}: some box in $\cB^R_{k,n}$ on which we apply the bootstrapping procedure to obtain superpolynomial concentration for $X_{R^{\kappa-k}}$. & \hyperref[def:Bk-first]{Sec.~\ref{def:Bk-first}}, \hyperref[def:Bk-second]{Sec.~\ref{def:Bk-second}}\\
$\cB_{i,n}^{R,k}$ & The set of boxes in $\cB^{R}_{i,n}$ that are contained in $B^{(k)}$ &  \hyperref[def:cBRk_jn]{Sec.~\ref*{def:cBRk_jn}}\\[0.4em]
& When $R=2$, the superscript $R$ is suppressed in all of the above. & \hyperref[def:cQR]{Sec.~\ref*{def:cQR}}\\[0.4em]

\multicolumn{3}{@{}l}{\textbf{Local versions on subdomains}}\\[0.3em]
$Q\in\cQ$ & A simply connected domain built from dyadic boxes & \hyperref[def:Q-in-cQ]{Sec.~\ref*{def:Q-in-cQ}} \\
$\Lambda_{Q,n}$ & The discrete set of points in $\Lambda_{n}$ inside $Q$ & \hyperref[def:Lambda_Q,n]{Sec.~\ref*{def:Lambda_Q,n}}\\
$\Gamma_{Q,n}$ & The collection of interface loops of the percolation configuration on $\Lambda_{Q,n}$ & \hyperref[def:Gamma_Q,n]{Sec.~\ref*{def:Gamma_Q,n}}\\

$\Gamma_Q$ & The collection of $\CLE_6$ loops in $Q$ induced by $\Gamma$ & \hyperref[def:Gamma_Q]{Sec.~\ref*{def:Gamma_Q}} \\
$\Lambda_Q=\Lambda_Q(\Gamma_Q)$ & The prime-end space associated with $\Gamma_Q$ & \hyperref[def:Lambda_Q]{Sec.~\ref*{def:Lambda_Q}} \\
$\cX_{Q,n}^{(1)},\cX_{Q,n}^{(2)},\ldots$ & The open and closed clusters of $\Lambda_{Q,n}$, ordered by decreasing Euclidean diameter & \hyperref[def:cXQmn]{Sec.~\ref*{def:cXQmn}} \\
$\cX_Q^{(1)},\cX_Q^{(2)},\ldots$ & The connected components of $\Lambda_Q$, equivalently the interior gaskets of $\Gamma_Q$, ordered by decreasing Euclidean diameter & \hyperref[def:cXQm]{Sec.~\ref*{def:cXQm}} \\
$\dpt_{Q,n}, \dg_{Q,n}, \dr_{Q,n}$ & The discrete percolation path, geodesic, and resistance metrics on $\Lambda_{Q,n}$ & \hyperref[def:dpt_{Q,n}]{Eq.~\eqref{def:dpt_{Q,n}}}, \hyperref[eq:def_dQn]{Eq.~\eqref{eq:def_dQn}}\\
$\dpt_Q, \dg_Q, \dr_Q$ & The $\CLE_6$ path, geodesic, and resistance metrics metrics on $\Lambda_Q$ & \hyperref[eq:def-dptU]{Eq.~\eqref{eq:def-dptU}}, \hyperref[def:dg_Q]{Sec.~\ref*{def:dg_Q}} \\[0.4em]

\multicolumn{3}{@{}l}{\textbf{Boxes, annuli, and arm events}}\\[0.3em]
$B(x,r)$ & The parallelogram (also called box) $x+[-r,r)^2$ & \hyperref[def:box]{Sec.~\ref*{def:box}} \\
$A(x;r_1,r_2):=B(x,r_2)\setminus B(x,r_1)$ & The parallelogram-shaped annulus between radii $r_1$ and $r_2$ & \hyperref[def:annulus]{Sec.~\ref*{def:annulus}} \\
$A_{j,\sigma}(r,R)$ & The probability of a $j$-arm event for percolation of color pattern $\sigma$ in the annulus $A(0;r,R)$ & \hyperref[def:arm-percolation]{Sec.~\ref*{def:arm-percolation}} \\
$A^{h}_{j,\sigma}(r,R)$ & The corresponding half-plane arm probability & \hyperref[def:arm-percolation-half]{Sec.~\ref*{def:arm-percolation-half}} \\
$\AC_{j,\sigma}(r,R)$ & The corresponding $\CLE_6$ arm-event probability & \hyperref[def:arm-cle]{Sec.~\ref*{def:arm-cle}} \\
$\alpha_j, \alpha_j^h, \alpha_j'$ & The full-plane polychromatic, half-plane polychromatic, and full-plane monochromatic arm exponents & \hyperref[lem:arm_exponents]{Lem.~\ref{lem:arm_exponents}}, \hyperref[lem:arm_exponents_monochromatic]{Lem.~\ref{lem:arm_exponents_monochromatic}} \\[0.4em]

\multicolumn{3}{@{}l}{\textbf{Gromov--Hausdorff-type distances and loop topologies}}\\[0.3em]
$\Delta_{\Haus}$ & The Hausdorff distance & \hyperref[eq:hausdorff]{Eq.~\eqref{eq:hausdorff}} \\
$\Delta_{\GH}$ & The Gromov-Hausdorff distance & \hyperref[eq:def-gh-embeddings]{Eq.~\eqref{eq:def-gh-embeddings}} \\
$\Delta_{\GHf}$ & The Gromov-Hausdorff-function distance & \hyperref[eq:ghf]{Eq.~\eqref{eq:ghf}} \\

$\Delta_{\Pro}$ & The Prokhorov distance & \hyperref[eq:ghp]{Eq.~\eqref{eq:ghp}} \\
$\Delta_{\GHP}$ & The Gromov-Hausdorff-Prokhorov distance & \hyperref[eq:ghp]{Eq.~\eqref{eq:ghp}} \\
$\Delta_{\GHPf}$ & The Gromov-Hausdorff-Prokhorov-function distance & \hyperref[eq:ghpf]{Eq.~\eqref{eq:ghpf}} \\
$\Delta_{\rGHPf}$ & The rooted Gromov-Hausdorff-Prokhorov-function distance & \hyperref[eq:rghpf]{Eq.~\eqref{eq:rghpf}} \\
$\Delta_{\cur}$ & The uniform distance modulo reparametrisation between curves or loops & \hyperref[eq:delta-cur-paths]{Eq.~\eqref{eq:delta-cur-paths}}, \hyperref[eq:delta-cur-loops]{Eq.~\eqref{eq:delta-cur-loops}} \\
$\Delta_{\Hc}$ & The Hausdorff distance between collections of curves or loops with respect to $\Delta_{\cur}$ & \hyperref[def:delta-Hc]{Sec.~\ref*{def:delta-Hc}} \\

\multicolumn{3}{@{}l}{\hypertarget{proofs-index-notation}{\textbf{Proofs}}}\\[0.3em]
$\cE$ & The auxiliary event on which the a priori bound holds & \hyperref[def:cE]{Sec.~\ref*{def:cE}}\\
$\delta=\frac{1}{10000}$ & The distance that arms need to away from or near the centers of $[0,1]^2$ in the event $\cE$ & \hyperref[def:delta]{Sec.~\ref*{def:delta}}\\
$X_n$ & The $d_n$-distance between the left- and rightmost pivotals on the event $\cE(\Lambda_n)$ & \hyperref[def:Xn]{Sec.~\ref*{def:Xn}}\\
$\altq_n=\altq_n(p)$ & The auxiliary normalizing constant used to define $\median_n$. It is the $(1-p)$th quantile of $X_n$. & \hyperref[eq:def-median_n]{Eq.~\eqref{eq:def-median_n}}\\
$p$ & The parameter for $\altq_n=\altq_n(p)$ & \hyperref[eq:def-median_n]{Eq.~\eqref{eq:def-median_n}}\\
$\fp$ & The value of $p$ such that $\median_n=\altq_n(\fp)$ & \hyperref[eq:def-median_n-second]{Eq.~\eqref{eq:def-median_n-second}}, \hyperref[def:median_n-hoelder]{Sec.~\ref*{def:median_n-hoelder}}\\
$\alpha$ & The exponent for which the a priori bound holds in the bootstrapping procedure & \hyperref[def:alpha-first]{Sec.~\ref*{def:alpha-first}}, \hyperref[prop:subpolynomialconc]{Prop.~\ref{prop:subpolynomialconc}}, \hyperref[cor:subpolynomialconc-general]{Cor.~\ref{cor:subpolynomialconc-general}}\\
\end{longtable}

\bibliographystyle{hmralphaabbrv}
\bibliography{references}

\end{document}